\documentclass[a4paper,10pt,twoside]{book}

\usepackage[utf8]{inputenc}
\usepackage[T1]{fontenc}
\usepackage[francais]{babel}
\usepackage{amsmath, amssymb, amsthm}
\usepackage{amstext, amsfonts, a4}
\usepackage{latexsym}
\usepackage{mathrsfs}
\usepackage[all]{xy}
\usepackage{graphicx}
\usepackage{enumerate}

\usepackage[pdftex,a4paper,left=3cm,right=3cm,top=3cm,bottom=3cm]{geometry}

\def\no{\noindent}

\def\ra{\rightarrow}

\def\rra{\rightrightarrows}

\def\lra{\longrightarrow}

\def\hra{\hookrightarrow}
\def\lmapsto{\longmapsto}

\def\bbD{\mathbb{D}}
\def\bbF{\mathbb{F}}\def\bbG{\mathbb{G}}
\def\bbL{\mathbb{L}}
\def\bbN{\mathbb{N}}\def\bbP{\mathbb{P}}
\def\bbQ{\mathbb{Q}}
\def\bbW{\mathbb{W}}
\def\bbZ{\mathbb{Z}}

\def\cA{\mathcal{A}}\def\cC{\mathcal{C}}
\def\cG{\mathcal{G}}\def\cH{\mathcal{H}}
\def\cI{\mathcal{I}}
\def\cN{\mathcal{N}}\def\cO{\mathcal{O}}\def\cP{\mathcal{P}}
\def\cQ{\mathcal{Q}}\def\cS{\mathcal{S}}\def\cT{\mathcal{T}}
\def\cU{\mathcal{U}}\def\cV{\mathcal{V}}\def\cW{\mathcal{W}}

\def\fgg{\mathfrak{g}}

\def\fm{\mathfrak{m}}

\def\om{\overline{m}}

\def\ov{\overline{v}}

\DeclareMathOperator{\ad}{ad}

\DeclareMathOperator{\can}{can}
\DeclareMathOperator{\Coker}{Coker}
\DeclareMathOperator{\cont}{cont}

\DeclareMathOperator{\deriv}{deriv}
\DeclareMathOperator{\disc}{disc}
\DeclareMathOperator{\dlog}{dlog}

\DeclareMathOperator{\Ext}{Ext}
\DeclareMathOperator{\Et}{\textrm{\'E}t}

\DeclareMathOperator{\et}{\textrm{é}t}
\DeclareMathOperator{\ev}{ev}

\DeclareMathOperator{\FPPF}{FPPF}
\DeclareMathOperator{\Fil}{Fil}

\DeclareMathOperator{\Frac}{Frac}

\DeclareMathOperator{\Gal}{Gal}

\DeclareMathOperator{\Gr}{Gr}
\DeclareMathOperator{\Groupes}{Groupes}

\DeclareMathOperator{\Hom}{Hom}

\DeclareMathOperator{\Id}{Id}
\DeclareMathOperator{\Ind}{Ind}
\DeclareMathOperator{\indpro}{indpro}
\DeclareMathOperator{\Ima}{Im}

\DeclareMathOperator{\Ker}{Ker}

\DeclareMathOperator{\klin}{k-lin.}

\DeclareMathOperator{\kschgr}{k-sch.gr.}

\DeclareMathOperator{\li}{\varinjlim}
\DeclareMathOperator{\Lie}{Lie}

\DeclareMathOperator{\lo}{long}

\DeclareMathOperator{\lp}{\varprojlim}

\DeclareMathOperator{\modd}{mod}

\DeclareMathOperator{\N}{N}

\DeclareMathOperator{\parf}{parf}
\DeclareMathOperator{\pf}{pf}

\DeclareMathOperator{\Pro}{Pro}
\DeclareMathOperator{\proind}{proind}

\DeclareMathOperator{\pizf}{\pi_0-f}

\DeclareMathOperator{\pr}{pr}

\DeclareMathOperator{\fQf}{\overrightarrow{\underleftarrow{\cQ}}}
\DeclareMathOperator{\fQfzz}{\overrightarrow{\underleftarrow{\cQ}}^{00}}

\DeclareMathOperator{\rat}{rat}
\DeclareMathOperator{\red}{red}
\DeclareMathOperator{\Res}{Res}

\DeclareMathOperator{\Rrat}{R-rat}

\DeclareMathOperator{\Sch}{Sch}
\DeclareMathOperator{\SchR}{(Sch/R)}
\DeclareMathOperator{\SchgrR}{(Sch. gr./R)}
\DeclareMathOperator{\SchK}{(Sch/K)}
\DeclareMathOperator{\SchgrK}{(Sch. gr./K)}

\DeclareMathOperator{\Spec}{Spec}

\DeclareMathOperator{\sym}{sym}

\DeclareMathOperator{\tf}{tf}

\DeclareMathOperator{\tors}{tors}

\DeclareMathOperator{\uExt}{\underline{Ext}}
\DeclareMathOperator{\uHom}{\underline{Hom}}

\DeclareMathOperator{\fUf}{\overrightarrow{\underleftarrow{\cU}}}
\DeclareMathOperator{\fUnf}{\overrightarrow{\underleftarrow{\cU(n)}}} 
\DeclareMathOperator{\fUuf}{\overrightarrow{\underleftarrow{\cU(1)}}}
\DeclareMathOperator{\fUfzz}{\overrightarrow{\underleftarrow{\cU}}^{00}}
\DeclareMathOperator{\fVf}{\overrightarrow{\underleftarrow{\cV}}}

\DeclareMathOperator{\VIB}{VI_B}

\DeclareMathOperator{\ZAR}{ZAR}

\newtheorem{counter}[subsection]{$\!\!$}
\newenvironment{Def}{\begin{counter} {\bf Définition.}}{\end{counter}}
\newenvironment{Not}{\begin{counter} \rm {\bf Notation.}}{\end{counter}}
\newenvironment{Nots}{\begin{counter} \rm {\bf Notations.}}{\end{counter}}

\newenvironment{Notss}{\begin{counter} \rm {\bf Notations.}}{\end{counter}}
\newenvironment{Prop}{\begin{counter} {\bf Proposition.}}{\end{counter}}
\newenvironment{Lem}{\begin{counter} {\bf Lemme.}}{\end{counter}}
\newenvironment{Cor}{\begin{counter} {\bf Corollaire.}}{\end{counter}}
\newenvironment{Th}{\begin{counter} {\bf Théorème.}}{\end{counter}}
\newenvironment{Rem}{\begin{counter} \rm {\bf Remarque.}}{\end{counter}}
\newenvironment{Ex}{\begin{counter} \rm {\bf Exemple.}}{\end{counter}}
\newenvironment{Exs}{\begin{counter} \rm {\bf Exemples.}}{\end{counter}}
\newenvironment{Pt}{\begin{counter} \rm}{\end{counter}}

\newtheorem{counter*}[subsubsection]{$\!\!$}
\newenvironment{Def*}{\begin{counter*} {\bf Définition.}}{\end{counter*}}
\newenvironment{Not*}{\begin{counter*} \rm {\bf Notation.}}{\end{counter*}}
\newenvironment{Notss*}{\begin{counter*} \rm {\bf Notations.}}{\end{counter*}}
\newenvironment{DefNot*}{\begin{counter*} \rm {\bf Définition-Notation.}}{\end{counter*}}
\newenvironment{Nots*}{\begin{counter*} \rm {\bf Notations.}}{\end{counter*}}
\newenvironment{Prop*}{\begin{counter*} {\bf Proposition.}}{\end{counter*}}
\newenvironment{Lem*}{\begin{counter*} {\bf Lemme.}}{\end{counter*}}
\newenvironment{Cor*}{\begin{counter*} {\bf Corollaire.}}{\end{counter*}}
\newenvironment{Th*}{\begin{counter*} {\bf Théorème.}}{\end{counter*}}
\newenvironment{Rem*}{\begin{counter*} \rm {\bf Remarque.}}{\end{counter*}}
\newenvironment{Ex*}{\begin{counter*} \rm {\bf Exemple.}}{\end{counter*}}
\newenvironment{Exs*}{\begin{counter*} \rm {\bf Exemples.}}{\end{counter*}}
\newenvironment{Pt*}{\begin{counter*} \rm}{\end{counter*}}
\newenvironment{Q*}{\begin{counter*} \rm {\bf Question.}}{\end{counter*}}

\title{\huge{\textbf{Dualité sur un corps local de caractéristique positive à corps résiduel algébriquement clos}}}
\author{Cédric Pépin}
\date{\today}

\begin{document}

\maketitle

\clearpage
\pagestyle{empty}
\phantom{placeholder}
\clearpage
\pagestyle{plain}

\chapter*{Introduction}
\addcontentsline{toc}{chapter}{Introduction}

Soit $K$ un corps discrètement valué complet à corps résiduel $k$ de caractéristique $p>0$. On dispose d'une théorie de dualité pour la cohomologie des $K$-schémas en groupes (commutatifs) finis dans les trois cas suivants : $K$ est de caractéristique $0$ et $k$ fini (Tate \cite{Ta}, 1962), $K$ est de caractéristique $p$ et $k$ fini (Shatz \cite{Sh}, 1964), $K$ est de caractéristique $0$ et $k$ algébriquement clos (Bégueri \cite{Beg}, 1980). On traite ici le cas où $K$ est de caractéristique $p$ et $k$ algébriquement clos. Une approche indépendante a été donnée récemment par Suzuki \cite{Su} 2.7.6 (1) (a).

La théorie cohomologique que l'on considère est à défaut la cohomologie fppf. La catégorie des coefficients est celle des $K$-schémas en groupes finis, qui est munie de la dualité de Cartier. Lorsque $K$ est de caractéristique $0$ et $k$ fini, la cohomologie est alors à valeurs dans la catégorie des groupes finis, qui est munie de la dualité de Pontryagin, et le théorème de dualité énonce que la dualité de Cartier sur les coefficients induit la dualité de Pontryagin sur la cohomologie. Lorsque $K$ est de caractéristique $p$ et $k$ fini, la cohomologie est à valeurs dans la catégorie des groupes topologiques localement compacts, à laquelle la dualité de Pontryagin des groupes finis s'étend, et le théorème de dualité s'énonce alors comme en caractéristique $0$. Lorsque $K$ est de caractéristique $0$ et $k$ algébriquement clos, la cohomologie prend ses valeurs dans la catégorie des groupes quasi-algébriques sur $k$ \cite{S2}, à composante neutre unipotente, qui satisfont à une dualité de type Pontryagin-Serre, et la cohomologie transforme la dualité de Cartier en la dualité de Pontryagin-Serre.

Lorsque $K$ est de caractéristique $p$ et $k$ algébriquement clos, on définit une catégorie $\fQf$ de \emph{groupes $\cQ$-admissibles}, qui sont une version quasi-algébrique sur $k$ des groupes topologiques localement compacts. La dualité de Serre des groupes quasi-algébriques unipotents connexes s'étend à la catégorie $\fUfzz$ des objets \emph{unipotents strictement connexes} de $\fQf$, et le théorème de dualité est le suivant :

\medskip

\textbf{Théorème.} \emph{Soit $K$ un corps discrètement valué complet de caractéristique $p>0$ à corps résiduel algébriquement clos. Soit $G_K$ un $K$-schéma en groupes commutatif fini.}
\begin{enumerate}

\medskip

\item \emph{On a $H^i(K,G_K)=0$ pour tout $i>1$.}

\medskip

\item \emph{Le groupe $H^1(K,G_K)$ possède une structure naturelle d'objet de $\fQf$, et sa composante neutre $H^1(K,G_K)^0$ est un objet de $\fUfzz$.}

\medskip

\item \emph{L'accouplement de dualité de Cartier $G_K\times_K G_{K}^*\ra\bbG_{m,K}$ induit des isomorphismes}
$$
\xymatrix{
H^1(K,G_{K}^*)/H^1(K,G_{K}^*)^0\ar[r]^<<<<<{\sim} & H^0(K,G_K)^{\vee_P},& H^1(K,G_{K}^*)^0\ar[r]^>>>>>{\sim} &(H^1(K,G_K)^0)^{\vee_S}, 
}
$$
\emph{où $(\cdot)^{\vee_P}$ désigne la dualité de Pontryagin et $(\cdot)^{\vee_S}$ celle de Serre.}
\end{enumerate}
\medskip

Les groupes $\cQ$-admissibles admettent deux descriptions : l'une \emph{indproalgébrique}, l'autre \emph{proindalgébrique}. Les groupes proalgébriques ont été introduits par Serre dans \cite{S2} : il s'agit de groupes ordinaires, qui sont limites projectives de groupes quasi-algébriques sur $k$ ; ils forment une catégorie abélienne $\cP$. Les groupes indproalgébriques sont alors des groupes ordinaires, qui sont limites inductives de groupes proalgébriques ; ils forment une catégorie $\cI\cP$, qui n'est plus abélienne, mais qui admet un plongement épais dans la ind-catégorie de $\cP$, en faisant une catégorie exacte. Les groupes indalgébriques sont des groupes ordinaires limites inductives de groupes quasi-algébriques, formant une catégorie exacte $\cI$, et les groupes proindalgébriques en sont les limites projectives, formant une catégorie exacte $\cP\cI$. La catégorie $\fQf$ est alors une sous-catégorie commune aux catégories $\cI\cP$ et $\cP\cI$, dont les objets satisfont à une certaine condition d'admissibilité. La dualité de Serre des groupes quasi-algébriques unipotents connexes s'étend aux sous-catégories des objets unipotents strictement connexes de $\cI\cP$ et $\cP\cI$, donnant lieu à des foncteurs $\cI(\cP^{00})(\cU)^{\circ}\ra\cP^{00}(\cI)(\cU)$ et $\cP^{00}(\cI)(\cU)^{\circ}\ra\cI(\cP^{00})(\cU)$ inverses l'un de l'autre, qui coïncident sur $\fUfzz$, et définissent donc une bidualité sur celle-ci.

Pour munir le $H^1$ des $K$-schémas en groupes finis d'une structure de groupe $\cQ$-admissible, on passe par la réalisation indproalgébrique de $\fQf$. Soit $R$ l'anneau des entiers de $K$ et $\fm$ son idéal maximal. Le groupe $R$, limite projective des tronqués $R/\fm^n$, $n\in\bbN$, est naturellement un groupe proalgébrique sur $k$ (c'est en particulier un espace topologique, mais dont la topologie n'est pas la topologie $\fm$-adique). Plus généralement, tout $R$-schéma en groupes de type fini se projette par le foncteur de Greenberg sur un pro-schéma en groupes de type fini sur $k$, induisant en particulier un groupe proalgébrique sur $k$. Partons maintenant d'un $K$-schéma en groupes fini $G_K$. Celui-ci possède des résolutions de longueur 1 par des $K$-schémas en groupes de type fini (et lisses), et $k$ étant algébriquement clos, certaines de ces résolutions sont à composantes $\Gamma(K,\cdot)$-acycliques ; soit $f_K:L_{K,0}\ra L_{K,1}$ une telle résolution. Supposons que $f_K$ s'étende en un morphisme de $R$-schémas en groupes de type fini $f:L_0\ra L_1$ et que $\Gamma(R,L_i)=\Gamma(K,L_{K,i})$ pour $i=0,1$. Alors le conoyau de la projection de $f$ par le foncteur de Greenberg fournit une structure de groupe proalgébrique sur le groupe ordinaire $H^1(K,G_K)$. Mais un tel morphisme $f$ n'existe pas en général : plus précisément, il n'existe pas de foncteur
$$
j_*:\textrm{($K$-schémas en groupes de type fini lisses)} \ra \textrm{($R$-schémas en groupes localement de type fini)}
$$
tel que $(j_*L_K)_K=L_K$ et $(j_*L_K)(R)=L_K(K)$ pour tout $L_K$, et $(j_*f_K)(R)=f_K(K)$ pour tout $f_K:L_{K,0}\ra L_{K,1}$. L'obstruction à l'existence d'un tel foncteur provient du groupe additif $\bbG_{a,K}$. Cependant, cette obstruction disparaît si l'on remplace la catégorie but par sa \emph{ind}-catégorie :

\medskip  

\textbf{Théorème.} \emph{Soit $K$ un corps discrètement valué complet à corps résiduel algébriquement clos. Soit $L_K$ un $K$-schéma en groupes lisse de type fini. Alors il existe un ind-objet $(L_i)_{i\in I}$ de la catégorie des $R$-schémas en groupes de fibre générique $L_K$ ayant les propriétés suivantes :}
\begin{enumerate}

\medskip

\item \emph{Pour tout $i\in I$, $L_i$ est lisse et séparé sur $R$.}

\medskip

\item \emph{Pour tout $i\in I$, le groupe ordinaire des composantes connexes de la fibre spéciale de $L_i$ est de type fini.}

\medskip

\item \emph{Pour tout $R$-schéma lisse de type fini $X$, tout morphisme de $K$-schémas $X\otimes_R K\ra L_K$ se prolonge de manière unique en un morphisme de ind-objets de la catégorie des $R$-schémas $X\ra (L_i)_{i\in I}$.}

\medskip

\item \emph{Pour tout $K$-schéma en groupes lisse de type fini $L_{K}'$, et pour tout ind-objet $(L_{i'}')_{i'\in I'}$ de la catégorie des $R$-schémas en groupes de fibre générique $L_{K}'$ ayant les propriétés 1. et 2., tout morphisme de $K$-schéma en groupes $L_{K}'\ra L_K$ se prolonge de manière unique en un morphisme de ind-objets de la catégorie des $R$-schémas en groupes $(L_{i'}')_{i'\in I'}\ra (L_i)_{i\in I}$.}
\end{enumerate}

\medskip

Il s'agit d'une extension de la théorie des modèles de Néron des groupes algébriques lisses sur $K$ à ceux dont la partie unipotente n'est pas totalement ployée. Ceci étant, on dispose d'un foncteur de prolongement des $K$-schémas en groupes de type fini lisses vers la ind-catégorie des $R$-schémas en groupes localement de type fini. Composant avec l'indisé du foncteur de Greenberg, on obtient finalement une projection vers la ind-catégorie des groupes proalgébriques sur $k$. On montre ensuite que si $f_K:L_{K,0}\ra L_{K,1}$ est une résolution de $G_K$ comme plus haut, alors le conoyau de la projection de $f_K$ se trouve dans l'image essentielle du plongement $\fQf\hra\cI\cP\hra\Ind(\cP)$ ; d'où une structure de groupe $\cQ$-admissible sur le groupe ordinaire $H^1(K,G_K)$. 

A ce stade, la structure indproalgébrique sur $H^1(K,G_K)$ dépend cependant \emph{a priori} du choix de la résolution du coefficient $G_K$. Contrairement au cas quasi-algébrique, ou plus généralement proalgébrique, il existe en effet des morphismes de groupes indproalgébriques induisant des isomorphismes sur les groupes ordinaires sous-jacents, mais qui ne sont pas des isomorphismes. Pour remédier à ce phénomène, on utilise l'équivalence entre groupes quasi-algébriques et \emph{schémas en groupes parfaits algébriques}. Celle-ci permet d'associer aux groupes quasi-algébriques des \emph{points génériques}, puis de définir un foncteur \emph{fibre générique sur la catégorie exacte $\cI\cP$}, consistant à évaluer un groupe indproalgébrique sur la catégorie de \emph{toutes} les extensions algébriquement closes de $k$, et non pas à lui associer seulement le groupe ordinaire de ses $k$-points. Sous une hypothèse de finitude convenable, on peut alors démontrer un résultat de pleine fidélité de ce foncteur fibre générique, en ayant recours à un procédé d'extension d'applications rationnelles, analogue pour les groupes indproalgébriques d'un lemme d'extension de Weil pour les applications rationnelles de $k$-groupes algébriques lisses qui sont génériquement des homomorphismes. Plus précisément, un morphisme entre les fibres génériques de deux groupes indproalgébriques induit par Galois équivariance des morphismes entre les groupes de points à valeurs dans les extensions \emph{parfaites} de $k$, en particulier sur les points génériques du ind-pro-schéma en groupes parfait algébrique associé à la source, auxquels les arguments de Weil s'appliquent. Maintenant, il est possible de choisir des résolutions $f_K:L_{K,0}\ra L_{K,1}$ de $G_K$, à composantes $\Gamma(K\widehat{\otimes}\kappa,\cdot)$-acycliques pour toute extension algébriquement close $\kappa/k$ ; la structure indproalgébrique sur $H^1(K,G_K)$ associée à l'une quelconque de ces résolutions a donc pour $\kappa$-points le groupe ordinaire $H^1(K\widehat{\otimes}\kappa,G_K)$. La géométrie des modèles de Néron généralisés montre en outre que ces structures satisfont à l'hypothèse du résultat de pleine fidélité précédent ; elles sont par conséquent isomorphes. On note 
$\cH^1(G_K)$ le groupe $\cQ$-admissible ainsi obtenu. On montre en outre que sa composante neutre est unipotente strictement connexe. 

Passons à la dualité. D'après ce qui précède, le groupe $\cQ$-admissible $\cH^1(G_K)$ peut être vu comme un \emph{espace de modules grossier sur $k$ des $G_K$-torseurs}. Grâce à cette interprétation, on peut dériver de la dualité de Cartier $G_K\times G_{K}^*\ra\bbG_{m,K}$ une flèche de dualité
$$
\xymatrix{
\cH^1(G_{K}^*)^0(\kappa)\ar[r] & (\cH^1(G_K)^0)^{\vee_S}(\kappa) 
}
$$
pour toute extension algébriquement close $\kappa/k$, de façon usuelle. En revanche, si $X$ est un $k$-schéma autre que le spectre d'un corps algébriquement clos extension de $k$, on ne dispose d'aucune interprétation géométrique du groupe $\cH^1(G_K)(X)$. Pour obtenir une flèche entre les espaces de modules $\cH^1(G_{K}^*)^0$ et $(\cH^1(G_K)^0)^{\vee_S}$ induisant les flèches précédentes sur les points géométriques, on fait donc à nouveau appel à la pleine fidélité du foncteur fibre générique sur les groupes indproalgébriques considérés. On obtient ainsi une flèche de dualité \emph{dans $\fQf$}
$$
\xymatrix{
\cH^1(G_{K}^*)^0\ar[r] &(\cH^1(G_K)^0)^{\vee_S}.
} 
$$

On montre ensuite que cette flèche est un isomorphisme sur les points géométriques. Quitte à appliquer le changement de base $k\ra\kappa$, il suffit de le montrer pour les $k$-points. Puis à l'aide d'un dévissage, on se ramène à $G_K=\alpha_p$, $\bbZ/p$ ou $\mu_p$. On transforme alors le problème de nature fppf en un problème de nature \emph{cristalline}. D'une part, via l'identification $\cH^1(G_{K}^*)(k)=\Ext_{\bbZ/p}^1(\bbZ/p,G_{K}^*)$, on peut << saturer >> un $G_{K}^*$-torseur en un $K$-schéma en groupes fini, extension de $\bbZ/p$ par $G_{K}^*$, et prendre son cristal de Dieudonné. D'autre part, via  l'identification $(\cH^1(G_K)^0)^{\vee_S}(k)=\Ext_{\cI\cP(1)}^1(\cH^1(G_K)^0,\bbZ/p)$ et une application différentielle logarithmique attachée au coefficient $G_K$, on peut expliciter les $k$-points de 
$(\cH^1(G_K)^0)^{\vee_S}$ en termes différentiels. On peut alors identifier la flèche de dualité sur les $k$-points à un isomorphisme classique de la théorie du corps de classes de $K$, composé d'un accouplement différentiel et d'une application résidu. Dans les cas $G_K=\alpha_p$ et $G_K=\mu_p$, les $G_{K}^*$-schémas en groupes finis qui apparaissent par saturation sont annulés par Verschiebung, leur cristal de Dieudonné s'identifie donc à la $p$-algèbre de Lie de leur dual de Cartier et l'interprétation du $H^1$ de $G_K$ en termes différentiels est alors due à Artin-Milne \cite{AM}. Dans le cas $G_K=\bbZ/p$, les saturés sont des Barsotti-Tate tronqués d'échelon 1, et le passage des données fppf aux données cristallines est alors assuré par le complexe syntomique à coefficients de Trihan-Vauclair \cite{TV}. La question d'étendre ces formes explicites de l'isomorphisme de dualité à tout $K$-schéma en groupes fini sera abordée dans un travail ultérieur \cite{PV}.

On obtient enfin que la flèche de dualité est un isomorphisme dans $\fQf$, à nouveau par pleine fidélité ou par un critère générique d'isomorphisme direct.  

\bigskip

\textbf{Remerciements.} J'exprime toute ma gratitude à Michel Raynaud pour m'avoir introduit au sujet, et pour les échanges stimulants qui ont suivi. Je remercie Luc Illusie pour avoir attiré mon attention sur le site syntomique. Je remercie Jilong Tong pour un groupe de travail sur l'article de Lucile Bégueri \cite{Beg} qui m'a été très utile. Je remercie enfin l'équipe d'algèbre de l'Université de Leuven pour m'avoir accueilli en post-doctorat, m'ayant ainsi permis de commencer ce travail.

\tableofcontents

\bigskip

\bigskip

\no\textbf{Conventions :} 
\begin{itemize}

\medskip

\item Les groupes considérés sont supposés \emph{commutatifs}, à de rares exceptions près mentionnées explicitement. 

\medskip 

\item Les systèmes inductifs (resp. projectifs) considérés sont supposés indexés par des \emph{petits ensembles préordonnés filtrants essentiellement dénombrables}.

\medskip

\item La cohomologie des schémas en groupes est calculée pour la topologie \emph{fppf}.
\end{itemize}

\chapter{Groupes indproalgébriques et proindalgébriques}

Soit $k$ un corps algébriquement clos de caractéristique $p>0$. 

\section{Groupes parfaits algébriques} \label{groupesparfaitsalgébriques}

\begin{Def}
Un \emph{groupe parfait algébrique (\emph{resp.} localement algébrique) sur $k$} est un $k$-schéma en groupes isomorphe au perfectisé d'un $k$-schéma en groupes de type fini (resp. localement de type fini). 
\end{Def}

\begin{Lem}\label{reppf}
Soient $\cG_1$ et $\cG_2$ des $k$-schémas en groupes de type fini. Soit $\varphi:\cG_1^{\pf}\ra\cG_2^{\pf}$ un morphisme de $k$-schémas en groupes. Alors il existe $n\in\bbN$ et un morphisme de $k$-schémas en groupes $f:\cG_1^{(p^{-n})}\ra\cG_2$ tel que $f^{\pf}=\varphi$.
\end{Lem}

\begin{proof}
Notons $\pi_2:\cG_{2}^{\pf}\ra\cG_2$ le morphisme canonique. Comme $\cG_2$ est de type fini sur $k$, il existe $n\in\bbN$, un recouvrement affine fini $\{U_i\}_{i\in I}$ de $\cG_1^{(p^{-n})}$ et des morphismes de $k$-schémas $f_i:U_i\ra\cG_2$ tels que $f_{i}^{\pf}=\pi_2\circ(\varphi|_{U_{i}^{\pf}}):U_{i}^{\pf}\ra\cG_2$. En particulier, pour tous $i,j\in I$, $f_i$ et $f_j$ coïncident  sur $U_i\cap U_j$ ensemblistement, puis quitte à augmenter $n$, schématiquement. Ils se recollent alors en un morphisme de $k$-schémas $f:\cG_1^{(p^{-n})}\ra\cG_2$ tel que $f^{\pf}=\varphi$. Notant $m$ les lois de groupes, le diagramme 
$$
\xymatrix{
\cG_1^{(p^{-n})}\times_k\cG_1^{(p^{-n})} \ar[r]^<<<<<<{m_{\cG_1^{(p^{-n})}}} \ar[d]_{f\times f} & \cG_1^{(p^{-n})} \ar[d]^f \\
\cG_2\times_k\cG_2 \ar[r]^{m_{\cG_2}} & \cG_2 
}
$$
commute alors ensemblistement, puis quitte à augmenter à nouveau $n$, schématiquement, et $f$ est alors un morphisme de $k$-schémas en groupes tel que $f^{\pf}=\varphi$.
\end{proof}

\begin{Def}
\begin{itemize}

\medskip

\item Une suite de $k$-schémas en groupes localement de type fini est \emph{exacte} si elle l'est dans la catégorie des faisceaux abéliens sur le gros site fppf de $k$.

\medskip

\item Une suite de $k$-schémas en groupes parfaits localement algébriques est \emph{exacte} si elle l'est dans la catégorie des faisceaux abéliens sur le site parfait de $k$ muni de la topologie étale.
\end{itemize}
\end{Def}

\begin{Lem} \emph{(cf. Berthelot \cite{Ber} 2.3)} \label{pfexact}
Soit 
$$
\xymatrix{
\cG' \ar[r] & \cG \ar[r] & \cG''
}
$$
une suite exacte de $k$-schémas en groupes localement de type fini. Alors la suite de $k$-schémas en groupes parfaits localement algébriques
$$
\xymatrix{
(\cG')^{\pf} \ar[r] & \cG^{\pf} \ar[r] & (\cG'')^{\pf}
}
$$
 est exacte.
\end{Lem}

\begin{proof}
Soit $\beta$ le morphisme canonique du gros topos fppf de $k$ dans le gros topos étale de $k$, $\alpha$ le morphisme canonique du gros topos étale de $k$ dans le gros topos parfait étale de $k$, et $\epsilon:=\alpha\circ\beta$. Le foncteur 
$\alpha_*$ est exact puisqu'un schéma étale sur un schéma parfait est parfait. Si $\cG$ est un $k$-schéma en groupes lisse, on a $R\beta_*\cG=\cG$ d'après \cite{}, donc $R\epsilon_*\cG=R\alpha_*\cG=\cG^{\pf}$. Si $\cG$ est un $k$-schéma en groupes radiciel, sa composante neutre possède une suite de composition dont les quotients successifs sont isomorphes à $\mu_p$ ou $\alpha_p$, donc $R\epsilon_*\cG=0$ d'après ce qui précède et les suites exactes
$$
\xymatrix{
0 \ar[r] & \mu_p \ar[r] & \bbG_m \ar[r]^F & \bbG_m \ar[r] & 0,
}
$$
$$
\xymatrix{
0 \ar[r] & \alpha_p \ar[r] & \bbG_a \ar[r]^F & \bbG_a \ar[r] & 0.
}
$$
Enfin si $\cG$ est un $k$-schéma en groupes localement de type fini quelconque, il est extension d'un $k$-schéma en groupes radiciel par le $k$-schéma en groupes lisse $\cG_{\red}$, d'où $R\epsilon_*\cG=\cG^{\pf}=(\cG_{\red})^{\pf}$ ; en particulier $R^1\epsilon_*\cG=0$.

Si l'on considère maintenant une suite exacte de $k$-schémas en groupes localement de type fini
$$
\xymatrix{
\cG' \ar[r]^f & \cG \ar[r]^g & \cG'',
}
$$
le faisceau fppf image de $f$ est représentable par le sous-$k$-schéma en groupes fermé noyau de $g$, et l'on a donc une suite exacte de $k$-schémas en groupes localement de type fini
$$
\xymatrix{
0 \ar[r] & \Ker f \ar[r] & \cG' \ar[r] & \Ima f \ar[r] & 0.
}
$$
Comme $R^1\epsilon_*(\Ker f)=0$, la suite
$$
\xymatrix{
0 \ar[r] & (\Ker f)^{\pf} \ar[r] & (\cG')^{\pf} \ar[r] & (\Ima f)^{\pf} \ar[r] & 0
}
$$
qui s'en déduit est donc exacte. En particulier, $f^{\pf}$ induit un épimorphisme de $(\cG')^{\pf}$ sur $(\Ima f)^{\pf}=(\Ker g)^{\pf}$, i.e. $\Ima f^{\pf}=(\Ker g)^{\pf}$. Comme par ailleurs $(\Ker g)^{\pf}=\Ker g^{\pf}$, la suite
$$
\xymatrix{
(\cG')^{\pf} \ar[r]^{f^{\pf}} & \cG^{\pf} \ar[r]^{g^{\pf}} & (\cG'')^{\pf}
}
$$
est donc exacte.
\end{proof}

\begin{Cor}\label{gppfab}
La catégorie des $k$-schémas en groupes parfaits algébriques est abélienne, et le plongement 
$$
\xymatrix{
\{\textrm{\emph{Groupes parfaits algébriques}}\} \ar@{^{(}->}[r] & \{\textrm{\emph{Faisceaux abéliens sur $(\parf/k)_{\Et}$}}\}
}
$$
est exact.
\end{Cor}

\begin{proof}
D'après \ref{reppf}, tout morphisme de $k$-schémas en groupes parfaits algébriques peut être réalisé comme le perfectisé d'un morphisme de $k$-schémas en groupes de type fini. Comme le foncteur de perfectisation de la catégorie abélienne des $k$-schémas en groupes de type fini dans celle des faisceaux abéliens sur $(\parf/k)_{\Et}$ est exact d'après \ref{pfexact}, et que son image essentielle est contenue dans la catégorie des groupes parfaits algébriques, on obtient le corollaire.
\end{proof}

\section{Groupes quasi-algébriques} \label{groupesquasialgébriques}

\begin{Def}
Soit $G$ un groupe ordinaire. 

\medskip

Une \emph{structure de groupe algébrique sur $G$} est un couple $(\cG,\iota)$ où $\cG$ est un $k$-schéma en groupes de type fini et $\iota$ un isomorphisme de groupes $\cG(k)\xrightarrow{\sim} G$. 

\medskip

Un \emph{morphisme de structures de groupe algébrique sur $G$} $(\cG_1,\iota_1)\ra (\cG_2,\iota_2)$ est un morphisme de $k$-schémas en groupes $f:\cG_1\ra \cG_2$ tel que le diagramme
$$
\xymatrix{
\cG_1(k) \ar[dr]_{\rotatebox{135}{$\sim$}}^{\iota_1} \ar[rr]^{f(k)}  && \cG_2(k) \ar[dl]^{\rotatebox{45}{$\sim$}}_{\iota_2} \\
& G & 
}
$$
commute.
\end{Def}

\begin{Prop} \emph{(cf. \cite{S2}  1.1 Proposition 2).} \label{condequivQ} Soient $(\cG_1,\iota_1)$ et $(\cG_2,\iota_2)$ deux structures de groupe algébrique sur un groupe ordinaire $G$. Les conditions suivantes sont équivalentes :

\begin{itemize}

\medskip

\item[(i)] Il existe une structure de groupe algébrique $(\cG_3,\iota_3)$ sur $G$ avec des morphismes 
$(\cG_3,\iota_3)\ra(\cG_1,\iota_1)$, $(\cG_3,\iota_3)\ra(\cG_2,\iota_2)$.

\medskip

\item[(ii)] Il existe une structure de groupe algébrique $(\cG_4,\iota_4)$ sur $G$ avec des morphismes 
$(\cG_1,\iota_1)\ra(\cG_4,\iota_4)$, $(\cG_2,\iota_2)\ra(\cG_4,\iota_4)$.

\medskip

\item[(iii)] Il existe $n\in\bbN$ avec un morphisme $(\cG_{1}^{(p^{-n})},\iota_1\circ F_{\cG_{1}^{(p^{-n})}/k}^n(k))\ra (\cG_2,\iota_2)$.

\medskip

\item[(iv)] Il existe $n\in\bbN$ avec un morphisme $(\cG_1,\iota_1)\ra (\cG_{2}^{(p^n)},\iota_2\circ (F_{\cG_2/k}^n(k))^{-1})$.

\medskip

\item[(v)] Il existe un isomorphisme de $k$-schémas en groupes $f:\cG_{1}^{\pf}\xrightarrow{\sim}\cG_{2}^{\pf}$ tel que le diagramme
$$
\xymatrix{
\cG_{1}^{\pf}(k)  \ar[d]^{\rotatebox{90}{$\sim$}} \ar[rr]^{f(k)}  && \cG_{2}^{\pf}(k) \ar[d]^{\rotatebox{90}{$\sim$}} \\
\cG_1(k) \ar[dr]_{\rotatebox{135}{$\sim$}}^{\iota_1}  && \cG_2(k) \ar[dl]^{\rotatebox{45}{$\sim$}}_{\iota_2} \\
& G & 
}
$$
commute.
\end{itemize}
\end{Prop}

\begin{proof}
$(iii)$ et $(iv)$ sont évidemment équivalents et impliquent $(i)$ et $(ii)$.

\medskip

$(i) \Rightarrow (v).$ Notant $f_{3,1}$ le morphisme $\cG_3\ra\cG_1$, on a une suite exacte 
$$
\xymatrix{
0\ar[r] & \Ker f_{3,1} \ar[r] & \cG_3 \ar[r]^{f_{3,1}} &\cG_1 \ar[r] & \Coker f_{3,1} \ar[r] & 0.
}
$$
Comme l'application $f_{3,1}(k)$ est bijective, $\Ker f_{3,1}$ et $\Coker f_{3,1}$ sont des $k$-schémas en groupes radiciels ; leurs prefectisés sont donc nuls. Il résulte alors de \ref{pfexact} que $f_{3,1}^{\pf}:\cG_{3}^{\pf}\ra\cG_{1}^{\pf}$ est un isomorphisme. De même $f_{3,2}^{\pf}:\cG_{3}^{\pf}\ra\cG_{2}^{\pf}$ est un isomorphisme, et $f:=f_{3,2}^{\pf}\circ(f_{3,1}^{\pf})^{-1}:\cG_{1}^{\pf}\ra \cG_{2}^{\pf}$ est alors un isomorphisme convenable.

\medskip

$(ii) \Rightarrow (v).$ Analogue à $(i) \Rightarrow (iv).$

\medskip

$(v) \Rightarrow (iii).$ Résulte de \ref{reppf}.
\end{proof}

\begin{Def}
Soit $G$ un groupe ordinaire. Deux structures de groupe algébrique sur $G$ sont \emph{équivalentes} si elles vérifient les conditions de \ref{condequivQ}.
\end{Def}

\begin{Def}
Un groupe \emph{quasi-algébrique} est un groupe ordinaire $G$ muni d'une classe d'équivalence $\cC_G$ de structures de groupe algébrique sur $G$.
\end{Def}

\begin{Pt}\label{foncQpf1}
Soit $G$ un groupe quasi-algébrique. D'après \ref{condequivQ} $(v)$, et puisque l'ensemble des $k$-points d'un schéma en groupes parfait algébrique est schématiquement dense dans celui-ci, le $k$-schéma en groupes $\cG^{\pf}$ ne dépend pas du choix de $(\cG^{\pf},\iota)\in\cC_G$ à isomorphisme canonique près ; on note $G^{\pf}$ le $k$-schéma en groupes parfait algébrique ainsi défini. 
\end{Pt}

\begin{Def}
Un \emph{morphisme de groupes quasi-algébriques} $(G_1,\cC_{G_1})\ra (G_2,\cC_{G_2})$ est un morphisme de groupes ordinaires $f:G_1\ra G_2$ tel qu'il existe $(\cG_1,\iota_1)\in\cC_{G_1}$, $(\cG_2,\iota_2)\in\cC_{G_2}$ et un morphisme de $k$-schémas en groupes $\varphi:\cG_1\ra\cG_2$ tels que le diagramme 
$$
\xymatrix{
\cG_{1}(k)  \ar[d]^{\rotatebox{90}{$\sim$}}_{\iota_1} \ar[r]^{\varphi(k)}  & \cG_{2}(k) \ar[d]^{\rotatebox{90}{$\sim$}}_{\iota_2} \\
G_1 \ar[r]^f & G_2 
}
$$
commute.
\end{Def}

\begin{Pt}\label{equivmorphQ}
De façon équivalente, d'après \ref{reppf}, un morphisme de groupes quasi-algébriques $(G_1,\cC_{G_1})\ra (G_2,\cC_{G_2})$ est un morphisme de groupes ordinaires $f:G_1\ra G_2$ tel qu'il existe un morphisme de $k$-schémas en groupes 
$\varphi:\cG_{1}^{\pf}\ra\cG_{2}^{\pf}$ tel que le diagramme 
$$
\xymatrix{
G_{1}^{\pf}(k)  \ar[d]^{\rotatebox{90}{$\sim$}} \ar[r]^{\varphi(k)}  & G_{2}^{\pf}(k) \ar[d]^{\rotatebox{90}{$\sim$}} \\
G_1 \ar[r]^f & G_2,
}
$$
où les flèches verticales sont les applications canoniques, commute. En particulier, le composé de deux morphismes de groupes quasi-algébriques est un morphisme de groupes quasi-algébriques (ce qui peut aussi se déduire de la définition et de \ref{condequivQ} $(iv)$).
\end{Pt}

\begin{Not}
On note $\cQ$ la catégorie des groupes quasi-algébriques. 
\end{Not}

\begin{Pt}\label{foncQpf2}
Comme l'ensemble des $k$-points d'un schéma en groupes parfait algébrique est schématiquement dense dans celui-ci, le morphisme $\varphi$ de \ref{equivmorphQ} est uniquement déterminé par $f$, et l'association $G\ra G^{\pf}$ est donc fonctorielle.
\end{Pt}

\begin{Prop}\label{Qgppf}
Le foncteur $(\cdot)^{\pf}$ \ref{foncQpf1} et \ref{foncQpf2} est une équivalence de catégories :
$$
\xymatrix{
(\cdot)^{\pf}:\cQ \ar[r]^<<<<{\sim} & \{\textrm{\emph{Groupes parfaits algébriques}}\}.
}
$$
\end{Prop}

\begin{proof}
Ce foncteur est évidemment fidèle, et il est plein d'après \ref{reppf}. Il est essentiellement surjectif par définition.
\end{proof}

\begin{Cor}
La catégorie $\cQ$ est abélienne. 
\end{Cor}

\begin{Pt}
Comme $k$ est algébriquement clos, le foncteur des sections globales sur la catégorie des faisceaux abéliens sur le site parfait étale de $k$ est exact. Réciproquement, une suite de schémas en groupes parfaits algébriques qui est exacte sur les $k$-points est exacte, à nouveau parce que l'ensemble des $k$-points d'un schéma en groupes parfait algébrique est schématiquement dense dans celui-ci. Le foncteur d'oubli de $\cQ$ dans la catégorie des groupes qui a un groupe quasi-algébrique associe le groupe ordinaire sous-jacent est donc exact, et pour qu'une suite de groupes quasi-algébriques soit exacte il suffit que la suite induite sur les groupes ordinaires sous-jacents le soit. En particulier, le foncteur d'oubli en question est conservatif.
\end{Pt}

Notons également la description suivante des suites exactes courtes de $\cQ$ :

\begin{Lem}\label{descrseQ}
Toute suite exacte courte de groupes quasi-algébriques est induite par une suite exacte courte de $k$-schémas en groupes lisses de type fini.
\end{Lem}

\begin{proof}
Soit 
$$
\xymatrix{
0 \ar[r] & G' \ar[r] & G \ar[r] & G'' \ar[r] & 0
}
$$
une suite exacte de $\cQ$. Par définition, $G'\ra G$ est induit par un morphisme de $k$-schémas en groupes de type fini $\iota:\cG'\ra\cG$. Quitte à passer aux réduits, on peut supposer $\cG'$ et $\cG$ lisses. Quitte à quotienter $\cG'$ par le noyau radiciel de $\iota$, on peut supposer que 
$\iota$ est une immersion fermée. Alors 
$$
\xymatrix{
0 \ar[r] & \cG' \ar[r]^{\iota} & \cG \ar[r] & \Coker\iota \ar[r] & 0
}
$$
est une suite exacte de $k$-schémas en groupes lisses de type fini induisant la suite exacte de $\cQ$ donnée.
\end{proof}

\section{Groupes proalgébriques} \label{groupesproalgébriques}

\begin{Not}\label{NotQ}
On note $\cQ$ la catégorie abélienne des groupes quasi-algébriques (sur $k$) : Serre \cite{S2} 1.2 Définition 2 et Proposition 5.
\end{Not}

Rappelons la définition d'un groupe proalgébrique : Serre \emph{loc. cit.} 2.1 Définition 1.

\begin{Def}\label{defP}
Un \emph{groupe proalgébrique} est un groupe ordinaire, muni d'une famille non vide $S$ de sous-groupes et, pour tout $H\in S$, d'une structure de groupe quasi-algébrique sur $G/H$, ces données vérifiant les axiomes suivants :
\begin{itemize}

\medskip

\item (P-1) $H,H'\in S\Rightarrow \exists H''\in S : H''\subset H\cap H'$.

\medskip

\item (P-2) Si $H\in S$, les sous-groupes $H'$ contenant $H$ qui appartiennent à $S$ sont les images réciproques des sous-groupes fermés de $G/H$.

\medskip

\item (P-3) Si $H,H'\in S$, et si $H\subset H'$, l'application canonique $G/H\ra G/H'$ est un morphisme de groupes quasi-algébriques.

\medskip

\item (P-4) L'application canonique $G\ra\lp_{H\in S} G/H$ est un isomorphisme de groupes.

\medskip

\end{itemize}

L'ensemble $S$ est l'\emph{ensemble complet de définition} du groupe proalgébrique $(G,S)$, et les sous-groupes $H\in S$ sont les \emph{sous-groupes de définition} de $(G,S)$.

\medskip

Un \emph{morphisme de groupes proalgébriques} $(G_1,S_1)\ra (G_2,S_2)$ est un morphisme de groupes $f:G_1\ra G_2$ vérifiant la condition suivante : pour tout $H_2\in S_2$, on a $f^{-1}(H_2)\in S_1$, et l'application $G_1/f^{-1}(H_2)\ra G_2/H_2$ induite par $f$ est un morphisme de groupes quasi-algébriques.
\end{Def}

\begin{Rem} \label{P2'}
Les axiomes (P-1) et (P-2) impliquent que $S$ est stable par intersection et par somme. Si $H,H'\in S$, il existe $H''\in S$ contenu dans $H\cap H'$, $H/H''$ et $H'/H''$ sont alors fermés dans $G/H''$ ; or l'intersection (resp. la somme) de deux sous-groupes fermés d'un groupe quasi-algébrique est un sous-groupe fermé, donc $(H/H'')\cap(H'/H'')=(H\cap H')/H''$ (resp. $(H/H'')+(H'/H'')=(H+H')/H'')$ est fermé dans $G/H''$, et donc $H\cap H'\in S$ (resp. $H+H'\in S$).
\end{Rem}

\begin{Rem} (cf. \emph{loc. cit.} 2.1 Remarque 2). 
Soit $G$ un groupe, muni d'une famille $S'$ de sous-groupes vérifiant (P-1), (P-3) et (P-4). Alors il existe une plus petite famille de sous-groupes de $G$ contenant $S'$ et vérifiant (P-1), (P-2), (P-3) et (P-4). On l'appellera l'ensemble complet de définition de $G$ \emph{engendré par $S'$}, et on le notera $<S'>$. On dira aussi (comme dans \emph{loc. cit.}) que $S'$ est un \emph{ensemble de définition} du groupe proalgébrique $(G,<S'>)$. Explicitement, $<S'>$ est l'ensemble des sous-groupes $H\subset G$ tels qu'il existe $H'\in S'$ vérifiant : 

\medskip

\begin{itemize}
\item $H'\subset H$ ;
\item l'inclusion $H/H'\subset G/H'$ est une immersion fermée quasi-algébrique.
\end{itemize}

\medskip

En particulier, $S'$ est coinitial dans $S$. Par ailleurs, $<S'>$ est l'unique ensemble complet de définition de $G$ contenant $S'$.
\end{Rem}

\begin{Not}\label{NotP}
On note $\cP$ la catégorie des groupes proalgébriques. On a un foncteur pleinement fidèle
\begin{eqnarray*}
\cQ & \ra & \cP \\
G & \mapsto & (G,<\{0\}>).
\end{eqnarray*}
\end{Not}

\begin{Th} \emph{(Serre \cite{S2} 2.6 Proposition 12).} \label{Pab}
Le foncteur de $\cP$ dans $\Pro(\cQ)$ qui a $(G,S)$ associe $(G/H)_{H\in S}$ est une équivalence de catégories :
$$
\xymatrix{
\cP \ar[r]^<<<<{\sim} & \Pro(\cQ).
}
$$
\end{Th}

\begin{Pt}\label{PidemQ}
Comme toute pro-catégorie d'une catégorie abélienne, la catégorie $\cP$ est abélienne. De plus, le foncteur d'oubli de $\cP$ dans la catégorie des groupes, qui à un groupe proalgébrique $(G,S)$ associe le groupe sous-jacent $G$, est exact. Par conséquent, ce foncteur est \emph{conservatif} : pour qu'un morphisme de groupes proalgébriques $(G_1,S_1)\ra (G_2,S_2)$ soit un isomorphisme, il faut et il suffit que le morphisme de groupes sous-jacent $G_1\ra G_2$ soit un isomorphisme. Il en résulte plus généralement que, pour qu'une suite de groupes proalgébriques
$$ 
\xymatrix{
(G_1,S_1)\ar[r] &  (G_2,S_2) \ar[r] & (G_3,S_3)
}
$$
soit exacte, il faut et il suffit que la suite des groupes sous-jacents
$$ 
\xymatrix{
G_1\ar[r] &  G_2 \ar[r] & G_3
}
$$
soit exacte.
\end{Pt}

Notons enfin la description suivante des suites exactes de $\cP$ (cf. \cite{S2} 2.4).

\begin{Lem}\label{descrseP}
Une suite de $\cP$
$$
\xymatrix{
0\ar[r] & (G',S') \ar[r] & (G,S) \ar[r]^f & (G'',S'') \ar[r] & 0
}
$$
est exacte si et seulement si la suite de $(\Groupes)$
$$
\xymatrix{
0\ar[r] & G' \ar[r] & G \ar[r]^f & G'' \ar[r] & 0
}
$$
est exacte, auquel cas
$$
S'=\{G'\cap H\ |\ H\in S\}\quad\textrm{et}\quad S''=\{f(H)\ |\ H\in S\},
$$
et pour tout $H\in S$, la suite de $\cQ$
$$
\xymatrix{
0\ar[r] & G'/G'\cap H \ar[r] & G/H \ar[r]^f & G''/f(H) \ar[r] & 0
}
$$
est exacte.
\end{Lem}

\section{Groupes parfaits proalgébriques affines}

\begin{Not}
On note $\cQ(\cA)$ la sous-catégorie pleine de $\cQ$ dont les objets sont les groupes quasi-algébriques affines (sur $k$).
\end{Not}

\begin{Pt}
La catégorie $\cQ(\cA)$ est abélienne et le plongement $\cQ(\cA)\ra \cQ$ est exact.
\end{Pt}

\begin{Def}
Un \emph{groupe proalgébrique affine} est un groupe proalgébrique $(G,S)$ tel $G/H\in\cQ(\cA)$ pour tout $H\in S$.
\end{Def}

\begin{Pt}
Un groupe proalgébrique est affine si et seulement s'il possède un ensemble de définition $S'$ tel que $G/H'\in\cQ(\cA)$ pour tout $H'\in S'$ : cela résulte de (P-3) puisque $\cQ(\cA)\subset\cQ$ est stable par quotient.
\end{Pt}

\begin{Not}
On note $\cP(\cA)$ la sous-catégorie pleine de $\cP$ dont les objets sont les groupes proalgébriques affines. 
\end{Not}

D'après \ref{Pab} et \ref{Qgppf}, on a : 

\begin{Th}\label{PAab}
\begin{itemize}
\medskip

\item Le foncteur de $\cP(\cA)$ dans $\Pro(\cQ(\cA))$ qui a $(G,S)$ associe $(G/H)_{H\in S}$ est une équivalence de catégories :
$$
\xymatrix{
\cP(\cA) \ar[r]^<<<<{\sim} & \Pro(\cQ(\cA)).
}
$$
En particulier, la catégorie $\cP(\cA)$ est abélienne et le plongement $\cP(\cA)\ra\cP$ est exact.

\medskip

\item Le foncteur $(\cdot)^{\pf}$ de $\cQ(\cA)$ dans la catégorie des groupes parfaits algébriques affines se prolonge en une équivalence de catégories :
$$
\xymatrix{
\Pro(\cQ(\cA)) \ar[r]^<<<<{\sim} & \Pro(\{\emph{\textrm{Groupes parfaits algébriques affines}}\}).
}
$$
\end{itemize}
\end{Th}

\begin{Def}
Un \emph{groupe parfait proalgébrique affine} est un $k$-schéma en groupes isomorphe à la limite d'un système projectif de $k$-schémas en groupes parfaits algébriques affines.
\end{Def}

\begin{Prop}
Le foncteur limite projective de la pro-catégorie des groupes parfaits algébriques affines dans la catégorie des groupes parfaits proalgébriques affines est une équivalence de catégories :
$$
\xymatrix{
\lp:\Pro(\{\textrm{\emph{Groupes parfaits algébriques affines}}\}) \ar[r]^<<<<<{\sim} & \{\textrm{\emph{Groupes parfaits proalgébriques affines}}\}
}
$$
En particulier, la catégorie des groupes parfaits proalgébriques affines est abélienne.
\end{Prop}

\begin{proof}
Il résulte des équivalences \ref{PAab} que le foncteur $\lp$ est fidèle, et on vérifie facilement qu'il est plein. Il est essentiellement surjectif par définition. 
\end{proof}

\begin{Not}\label{cPApf}
On note encore $(\cdot)^{\pf}$ l'équivalence composée
$$
\xymatrix{
\cP(\cA) \ar[r]^<<<<{\sim} \ar[drr]^{\sim}_{(\cdot)^{\pf}} & \cP(\cQ(\cA)) \ar[r]^<<<<{\sim} & \Pro(\{\textrm{Groupes parfaits algébriques affines}\}) \ar[d]^{\rotatebox{90}{$\sim$}} \\
& & \{\textrm{Groupes parfaits proalgébriques affines}\}.
}
$$
\end{Not}

\section{Groupes indalgébriques}

\begin{Def}\label{defI}
Un \emph{groupe indalgébrique} est un groupe ordinaire $G$ muni d'une famille non vide $S$ de sous-groupes et, pour tout $H\in S$, d'une structure de groupe quasi-algébrique sur $H$, ces données vérifiant les axiomes suivants :
\begin{itemize}

\medskip

\item (I-1) $H,H'\in S\Rightarrow \exists H''\in S : H+H'\subset H''$.

\medskip

\item (I-2) Si $H\in S$, les sous-groupes $H'$ contenus dans $H$ qui appartiennent à $S$ sont les sous-groupes fermés de $H$.

\medskip

\item (I-3) L'application canonique $\li_{H\in S} H\ra G$ est un isomorphisme de groupes.

\medskip

\end{itemize}

L'ensemble $S$ est l'\emph{ensemble complet de définition} du groupe proalgébrique $(G,S)$, et les sous-groupes $H\in S$ sont les \emph{sous-groupes de définition} de $(G,S)$.

\medskip

Un \emph{morphisme de groupes indalgébriques} $(G_1,S_1)\ra (G_2,S_2)$ est un morphisme de groupes $f:G_1\ra G_2$ vérifiant la condition suivante : pour tout $H_1\in S_1$, on a $f(H_1)\in S_2$, et l'application $H_1\ra f(H_1)$ induite par $f$ est un morphisme de groupes quasi-algébriques.
\end{Def}

\begin{Rem}\label{I2'}
L'axiome (I-2) implique l'axiome suivant :

\medskip

\begin{itemize}
\item \emph{(I-2)' Si $H,H'\in S$ avec $H$ contenu dans $H'$, alors l'inclusion $H\subset H'$ est une immersion fermée quasi-algébrique}.
\end{itemize}

\medskip

De plus, les axiomes (I-1) et (I-2) impliquent que $S$ est stable par somme et par intersection. Si $H,H'\in S$, il existe $H''\in S$ contenant $H+H'$, $H$ et $H'$ sont alors fermés dans $H''$, donc $H+H'$ (resp. $H\cap H'$) est fermé dans $H''$, et donc $H+H'\in S$ (resp. $H\cap H'\in S$).
\end{Rem}

\begin{Rem}\label{completengI}
Soit $G$ un groupe, muni d'une famille $S'$ de sous-groupes vérifiant (I-1), (I-2)' et (I-3). Alors il existe une plus petite famille de sous-groupes de $G$ contenant $S'$ et vérifiant (I-1), (I-2) et (I-3). On l'appellera l'ensemble complet de définition de $G$ \emph{engendré par $S'$}, et on le notera $<S'>$. On dira aussi que $S'$ est un \emph{ensemble de définition} du groupe indalgébrique $(G,<S'>)$. Explicitement, $<S'>$ est l'ensemble des sous-groupes $H\subset G$ tels qu'il existe $H'\in S'$ vérifiant : 

\medskip

\begin{itemize}
\item $H\subset H'$ ;
\item l'inclusion $H\subset H'$ est une immersion fermée quasi-algébrique.
\end{itemize}

\medskip

En particulier, $S'$ est cofinal dans $S$. En revanche, il peut y avoir deux ensembles complets de définition $S_1\subset S_2$ distincts sur un même groupe ordinaire $G$, donnant lieu à des groupes indalgébriques $(G,S_1)$ et $(G,S_2)$ non isomorphes. On peut par exemple munir le groupe $G=(k,+)$ des deux ensembles complets de définition : 
$$
S_1:=\{\bbG_{a,k}(k)\}\cup \{F\ |\ F\ \textrm{sous-groupe fini de}\ \bbG_{a,k}(k)\},
$$
$$
S_2:=\{F\ |\ F\ \textrm{sous-groupe fini de}\ \bbG_{a,k}(k)\}.
$$
Le groupe indalgébrique $(G,S_1)$ n'est autre que le groupe additif quasi-algébrique, alors que $(G,S_2)$ n'est pas isomorphe à un groupe quasi-algébrique, puisque le groupe ordinaire $G$ n'appartient pas à $S_2$.
\end{Rem}

\begin{Not}\label{NotI}
On note $\cI$ la catégorie des groupes indalgébriques. On a un foncteur pleinement fidèle 
\begin{eqnarray*}
\cQ & \ra & \cI \\
G & \mapsto & (G,< \{G\}>).
\end{eqnarray*}
\end{Not}

\begin{Prop}\label{Iexacte}
Le foncteur de $\cI$ dans $\Ind(\cQ)$ qui à $(G,S)$ associe $(H)_{H\in S}$ est pleinement fidèle : 
$$
\xymatrix{
\cI \ar@{^{(}->}[r] & \Ind(\cQ),
}
$$
et son image essentielle est la sous-catégorie épaisse de la catégorie abélienne $\Ind(\cQ)$, formée des ind-objets de $\cQ$ \emph{isomorphes à un ind-objet à transitions injectives}. En particulier, la catégorie $\cI$ est exacte. 
\end{Prop}

On a de plus le lemme suivant :

\begin{Lem}\label{condisoinjI}
Soit $((H_i,S_i))_{i\in I}\in\Ind(\cQ)$. Alors, pour que $((H_i,S_i))_{i\in I}$ soit isomorphe à un ind-objet à transitions injectives, il faut et il suffit qu'il vérifie la condition suivante :

\medskip

$(*)\quad$ Pour tout $i\in I$, il existe $j_i\geq i$ tel que $\Ker (H_{i}\ra\li_{i\in I} H_i) = \Ker (H_{i}\ra H_{j_i})$.
\end{Lem}

\begin{proof}
Supposons qu'il existe un isomorphisme $\iota:((H_i,S_i))_{i\in I}\xrightarrow{\sim}((H_{i'}',S_{i'}'))_{i'\in I'}$, avec 
$((H_{i'}',S_{i'}'))_{i'\in I'}$ à transitions injectives, et soit $i\in I$. Il existe alors $i'\in I'$, $j\geq i\in I$, et un diagramme 
$$
\xymatrix{
(H_i,S_i) \ar[r]^{\iota_{i,i'}} \ar[d]^{\tau_{i,j}} & (H_{i'}',S_{i'}') \ar[dl]^{(\iota^{-1})_{i',j}} \\
(H_j,S_j).
}
$$
Alors $(\iota^{-1})_{i',j}\circ \iota_{i,i'}$ et $\tau_{i,j}$ représentent chacun la composante d'indice $i$ de l'identité dans
$$
\Hom_{\Ind(\cQ)}(((H_i,S_i))_{i\in I},((H_i,S_i))_{i\in I})=\lp_i(\li_j\Hom_{\cQ}((H_i,S_i),(H_j,S_j))).
$$ 
Quitte à augmenter $j$, on peut donc supposer que le diagramme précédent est commutatif ; on pose alors $j_i:=j$. Comme le morphisme canonique $H_{i'}'\ra\li_{i'}H_{i'}'$ est injectif, on en déduit que  
$$
\Ker (H_i\ra \li_{i\in I} H_i)=\Ker (H_i\xrightarrow{\iota_{i,i'}} H_{i'}) \subset \Ker (H_i\xrightarrow{\tau_{i,j_i}} H_{j_i}).
$$ 

Réciproquement, supposons que $((H_i,S_i))_{i\in I}$ vérifie la condition $(*)$. Remarquons alors que pour tout $i\in I$ et $i'\geq j_i\in I$,
$$
\Ker((H_i,S_i)\ra (H_{i'},S_{i'}))=\Ker((H_i,S_i)\ra (H_{j_i},S_{j_i}))=:(K_i,S_{K_i}).
$$
En particulier, la transition $\tau_{i,i'}:(H_i,S_i)\ra (H_{i'},S_{i'})$ induit un morphisme injectif
$$
\xymatrix{
\overline{\tau_{i,i'}}:\overline{(H_i,S_i)}:=(H_i,S_i)/(K_i,S_{K_i}) \ar[r] & (H_{i'},S_{i'}).
}
$$
Composant avec la projection canonique $ (H_{i'},S_{i'})\ra (H_{i'},S_{i'})/(K_{i'},S_{K_{i'}})$, on obtient un morphisme
$$
\xymatrix{
\overline{\overline{\tau_{i,i'}}}:\overline{(H_i,S_i)}:=(H_i,S_i)/(K_i,S_{K_i}) \ar[r] & (H_{i'},S_{i'})/(K_{i'},S_{K_{i'}})=:\overline{(H_{i'},S_{i'})}.
}
$$
Celui-ci est encore injectif puisque $K_i=\Ker(H_i\ra\li_{i\in I} H_i)$. Ceci étant, soit $(\overline{I},\leq_{\overline{I}})$ le petit ensemble préordonné filtrant défini par $\overline{I}:=I$ et 
$$i\leq_{\overline{I}} i' \Leftrightarrow j_i\leq_{I} i'.$$
Alors le ind-objet $(\overline{(H_i,S_i)})_{i\in\overline{I}}$ dont les transitions sont données par les 
$\overline{\overline{\tau_{i,i'}}}$ est isomorphe à $((H_i,S_i))_{i\in I}$. Les projections canoniques
$$
\xymatrix{
(H_i,S_i) \ar[r] & \overline{(H_i,S_i)}, & i\in I,
}
$$
définissent en effet un morphisme $((H_i,S_i))_{i\in I}\ra (\overline{(H_i,S_i)})_{i\in \overline{I}}$, ayant pour inverse le morphisme défini par les
$$
\xymatrix{
\overline{(H_i,S_i)} \ar[r]^{\overline{\tau_{i,j_i}}} & (H_{j_i},S_{j_i}), & i\in \overline{I}.
}
$$
\end{proof}

\begin{proof}[Démonstration de la proposition \ref{Iexacte}]
Par définition des morphismes dans $\cI$, si $(G,S)$ et $(G',S')$ sont des objets de $\cI$, alors on a une application naturelle
$$
\xymatrix{
\Hom_{\cI}((G,S),(G',S')) \ar[r] & \lp_{H\in S}\li_{H'\in S'}\Hom_{\cQ}(H,H'),
}
$$
qui est évidemment bijective. Le foncteur $\cI\ra\Ind(\cQ)$ que l'on considère est donc bien pleinement fidèle. Les ind-objets de son image sont à transitions injectives par construction. Inversement, soit $((H_i,S_i))_{i\in I}\in\Ind(\cQ)$ à transitions injectives. Posant $G:=\li_{i\in I}H_i$, les $H_i$ s'identifient à des sous-groupes de $G$, et la famille $S:=(H_i)_{i\in I}$ vérifie (I-1), (I-2)' et (I-3). D'où un groupe indalgébrique $(G,<S>)$, dont l'image dans $\Ind(\cQ)$ est isomorphe à $((H_i,S_i))_{i\in I}$.

Il reste à vérifier que la sous-catégorie pleine de $\Ind(\cQ)$ formée des ind-objets isomorphes à un ind-objet à transitions injectives est épaisse. Soit 
\begin{equation}\label{seIepaisse}
\xymatrix{
0\ar[r] & H' \ar[r] & H \ar[r] & H'' \ar[r] & 0
}
\end{equation}
une suite exacte dans $\Ind(\cQ)$, avec $H'$ et $H''$ isomorphes à des ind-objets à transitions injectives. D'après \cite{KS} 8.6.6 (a), il existe un petit ensemble préordonné filtrant $I$, et un système de suites exactes courtes 
$$
\xymatrix{
0 \ar[r] & (H_{i}',S_{i}') \ar[r] & (H_i,S_i) \ar[r] & (H_{i}'',S_{i}'') \ar[r] & 0, & i\in I,
}
$$
tel que la suite 
$$
\xymatrix{
0 \ar[r] & ((H_{i}',S_{i}'))_{i\in I} \ar[r] & ((H_i,S_i))_{i\in I} \ar[r] & ((H_{i}'',S_{i}''))_{i\in I} \ar[r] & 0
}
$$
est isomorphe à (\ref{seIepaisse}). D'après \ref{condisoinjI}, les ind-objets $((H_{i}',S_{i}'))_{i\in I}$ et $((H_{i}'',S_{i}''))_{i\in I}$ vérifient la condition $(*)$, et il suffit de montrer que $((H_i,S_i))_{i\in I}$ la vérifie également. Soit $i\in I$. Choisissons $i_1\geq i$ tel que $\Ker (H_{i}''\ra\li_{i\in I} H_{i}'') = \Ker (H_{i}''\ra H_{i_1}'')$ puis $i_2\geq i_1$ tel que $\Ker (H_{i_1}'\ra\li_{i\in I} H_{i}') = \Ker (H_{i_1}'\ra H_{i_2}')$. Alors le diagramme commutatif exact suivant montre $\Ker (H_{i}\ra\li_{i\in I} H_{i}) = \Ker (H_{i}\ra H_{i_2})$ :
$$
\xymatrix{
0 \ar[r] & H_{i}' \ar[r] \ar[d] & H_i \ar[r] \ar[d] & H_{i}'' \ar[r] \ar[d] & 0 \\ 
0 \ar[r] & H_{i_1}' \ar[r] \ar[d] & H_{i_1} \ar[r] \ar[d] & H_{i_1}'' \ar[r] \ar[d] & 0 \\
0 \ar[r] & H_{i_2}' \ar[r] \ar[d] & H_{i_2} \ar[r] \ar[d] & H_{i_2}'' \ar[r] \ar[d] & 0 \\
0 \ar[r] & \li_{i\in I} H_{i}' \ar[r] & \li_{i\in I} H_i \ar[r] & \li_{i\in I} H_{i}'' \ar[r] & 0.
}
$$
\end{proof}

\begin{Pt}
Le foncteur $\cI\ra\Ind(\cQ)$ de \ref{Iexacte} n'est pas essentiellement surjectif. Par exemple, le ind-objet de $\cQ$ dont les transitions sont les surjections canoniques
$$
\xymatrix{
\bbG_{a,k}/F_i \ar[r] & \bbG_{a,k}/F_j,& F_i\subset F_j,
}
$$
où $F_i$ parcourt l'ensemble des sous-groupes finis de $k$, n'est pas dans l'image essentielle de $\cI$. Pour le vérifier, on peut utiliser le critère \ref{condisoinjI}. Plus généralement, la catégorie $\cI$ n'est pas abélienne, puisque comme on l'a vu en \ref{completengI}, le foncteur d'oubli de $\cI$ dans la catégorie des groupes, qui a un groupe indalgébrique $(G,S)$ associe le groupe sous-jacent $G$, n'est pas conservatif, bien qu'il soit exact (cf. plus généralement \ref{foncfibgenIexact}). On introduit donc un foncteur plus fin.
\end{Pt}

\begin{Def}\label{deffibgenI}
\begin{itemize}

\medskip

\item Soit $(\parf/k)$ la catégorie des $k$-schémas parfaits. Le \emph{foncteur de points} d'un groupe indalgébrique $(G,S)$ est le foncteur 
\begin{eqnarray*}
h_{(G,S)}:(\parf/k)^{\circ} & \ra & (\Groupes) \\
X  & \ra &  \li_{H\in S} H(X).
\end{eqnarray*}

\medskip

\item Soit $\cT$ la sous-catégorie pleine de $(\parf/k)$, dont les objets sont les morphismes $\Spec(\kappa)\ra\Spec(k)$ induits par les morphismes de corps algébriquement clos $k\ra\kappa$. La \emph{fibre générique} d'un groupe indalgébrique $(G,S)$ est la restriction du foncteur de points de $(G,S)$ à la catégorie $\cT$ :
\begin{eqnarray*}
(G,S)_{\eta}:\cT^{\circ} & \ra & (\Groupes) \\
\kappa  & \ra &  \li_{H\in S} H(\kappa).
\end{eqnarray*}

\medskip

\item Le \emph{foncteur fibre générique} de la catégorie des groupes indalgébriques dans la catégorie des préfaisceaux en groupes sur $\cT$ est le foncteur
\begin{eqnarray*}
(\cdot)_{\eta}:\cI & \ra & \widehat{\cT} \\
(G,S)  & \ra &  (G,S)_{\eta}.
\end{eqnarray*}
\end{itemize}
\end{Def}

\begin{Rem}\label{cbQ}
Soit $\kappa\in\cT^{\circ}$, et notons $\cQ_{\kappa}$ la catégorie des groupes quasi-algébriques sur $\kappa$. Le foncteur \emph{changement de base} de la catégorie des groupes parfaits algébriques sur $k$ dans la catégorie des groupes parfaits algébriques sur $\kappa$ définit un foncteur $\cQ\ra \cQ_{\kappa}$, que l'on appellera également foncteur \emph{changement de base}. Explicitement, si $H$ est un groupe quasi-algébrique sur $k$, induit par un groupe algébrique $H_{a}$ sur $k$, alors $H_{\kappa}$ est le groupe quasi-algébrique sur $\kappa$ induit par $(H_{a})_{\otimes_k}\kappa$. En particulier, l'identification $H_a(\kappa)=((H_{a})_{\otimes_k}\kappa)(\kappa)$ induit une identification $H(\kappa)=H_{\kappa}(\kappa)$. Le groupe des $\kappa$-points du groupe quasi-algébrique $H$ s'identifie donc au groupe ordinaire sous-jacent au groupe quasi-algébrique $H_{\kappa}$ sur $\kappa$. On notera par ailleurs que le foncteur changement de base 
$\cQ\ra\cQ_{\kappa}$ est exact, puisque le foncteur changement de base de la catégorie des groupes parfaits algébriques sur $k$ dans celle des groupes parfaits algébriques sur $\kappa$ est exact.
\end{Rem}

\begin{Lem}\label{foncfibgenIexact}
Le foncteur fibre générique $(\cdot)_{\eta}:\cI\ra\widehat{\cT}$  est exact.
\end{Lem}

\begin{proof} 
Le foncteur $\cQ\ni H\mapsto H(k)\in (\Groupes)$ coïncide avec le foncteur qui a un groupe quasi-algébrique associe le groupe ordinaire sous-jacent ; ce foncteur est donc exact. Plus généralement, il résulte de \ref{cbQ} que pour tout $\kappa\in\cT^{\circ}$, le foncteur $\cQ\ni H\mapsto H(\kappa)\in (\Groupes)$ est exact. Par conséquent, le foncteur 
\begin{eqnarray*}
(\cdot)_{\eta}:\cQ & \ra & \widehat{\cT} \\
H  & \ra &  H_{\eta}
\end{eqnarray*}
est exact. Il en résulte que le foncteur 
\begin{eqnarray*}
\Ind(\cQ) & \ra & \widehat{\cT} \\
(H_i)_{i\in I}  & \ra &  \li_i (H_i)_{\eta}
\end{eqnarray*}
est également exact, d'après \cite{KS} 8.6.6 (a), et le fait que les limites inductives filtrantes dans la catégorie des groupes sont exactes. D'où le lemme puisque le foncteur $\cI\ra\Ind(\cQ)$ de \ref{Iexacte} est exact par définition.  
\end{proof}

\begin{Prop}\label{foncfibgenIcons}
Le foncteur fibre générique $(\cdot)_{\eta}:\cI\ra \widehat{\cT}$ est conservatif : pour qu'un morphisme de groupes indalgébriques $(G_1,S_1)\ra (G_2,S_2)$ soit un isomorphisme, il faut et il suffit que le morphisme $(G_1,S_1)_{\eta}\ra (G_2,S_2)_{\eta}$ induit sur les fibres génériques soit un isomorphisme.
\end{Prop}

\begin{proof}
Soit $f:(G_1,S_1)\ra (G_2,S_2)$ tel que $f_{\eta}:(G_1,S_1)_{\eta}\ra (G_2,S_2)_{\eta}$ est un isomorphisme. Soit $H_1\in S_1$. Par hypothèse, $f(H_1)\in S_2$ et le morphisme $H_1\ra f(H_1)$ est quasi-algébrique. Comme de plus $f(k)$ est injectif, $H_1\ra f(H_1)$ est un isomorphisme quasi-algébrique. Nous allons montrer que 
$$
f(S_1):=\{ f(H_1)\ |\ H_1\in S_1 \}
$$
est égal à $S_2$. Les morphismes quasi-algébriques inverses des $H_1\xrightarrow{\sim} f(H_1)$ définiront donc un morphisme indalgébrique $(G_2,S_2)\ra (G_1,S_1)$, inverse de $f$.

On a $f(S_1)\subset S_2$ par définition, il s'agit donc de montrer $S_2\subset f(S_1)$. Il suffit pour cela de montrer que pour tout $H_2\in S_2$, il existe $H_1\in S_1$ tel que $H_2\subset f(H_1)$ : d'après (I-2) pour $S_2$, l'inclusion $H_2\subset f(H_1)$ est alors une immersion fermée quasi-algébrique, puis d'après (I-2) pour $S_1$ (et l'isomorphisme $H_1\xrightarrow{\sim} f(H_1)$), $H_2\in f(S_1)$. 

Soit $H_2\in S_2$. Soient $\eta_1,\ldots,\eta_n$ les points génériques de $H_{2}^{\pf}$. Choisissons un corps algébriquement clos $\kappa$ contenant les corps résiduels $k(\eta_1),\ldots,k(\eta_n)$, et notons 
$\overline{\eta_1},\ldots,\overline{\eta_n}$ les $\kappa$-points de $H_{2}^{\pf}$ induits par les $\eta_i$. Comme par hypothèse
l'application
$$
\xymatrix{
f(\kappa):\li_{H_1\in S_1} H_{1}(\kappa) \ar[r] & \li_{H_2\in S_2} H_{2}(\kappa)
}
$$
est bijective, il existe $H_1\in S_1$ tel que $\overline{\eta_1},\ldots,\overline{\eta_n}\in f(\kappa)(H_{1}(\kappa))$. Choisissons alors $H_{2}'\in S_2$ contenant $f(H_1)$ et $H_{2}$ ((I-1) pour $S_2$). Les inclusions $f(H_1)\subset H_{2}'$ et $H_{2}\subset H_{2}'$ étant des immersions fermées quasi-algébriques ((I-2) pour $S_2$), elles donnent lieu à des immersions fermées de schémas $(f(H_{1}))^{\pf}\ra (H_{2}')^{\pf}$ et $H_{2}^{\pf}\ra (H_{2}')^{\pf}$. De plus, comme $\cQ \ni H \mapsto H_{\eta} \in \widehat{\cT}$ est exact (\ref{foncfibgenIexact}), $(f(H_1))^{\pf}(\kappa)=(f(H_1))(\kappa)=f(\kappa)(H_{1}(\kappa))$. Par conséquent, $(f(H_{1}))^{\pf}$ majore l'adhérence schématique des $\overline{\eta_i}$ dans $(H_{2}')^{\pf}$, qui est égale à $H_{2}^{\pf}$. En prenant les $k$-points, on obtient $f(H_1)\supset H_2$.
\end{proof}

\begin{Rem}
Le critère \label{foncfibgenIcons} explique le choix de la notation $\cT$ : c'est une catégorie sur laquelle on peut \emph{tester} si un morphisme de groupes indalgébriques est un isomorphisme. 
\end{Rem}

\begin{Pt}
La catégorie exacte $\cI$ des groupes indalgébriques est donc munie de deux foncteurs exacts vers des catégories abéliennes : le plongement \ref{Iexacte}, et le foncteur fibre générique \ref{foncfibgenIexact} qui est conservatif,
$$
\xymatrix{
\cI \ar@{^{(}->}[r] \ar[d]_{(\cdot)_{\eta}} & \Ind(\cQ) \\
\widehat{\cT}.
}
$$
 \end{Pt}

\begin{Lem}\label{descrseI}
\begin{itemize}

\medskip

\item Une suite de $\cI$
$$
\xymatrix{
0\ar[r] & (G',S') \ar[r] & (G,S) \ar[r]^f & (G'',S'') \ar[r] & 0
}
$$
est exacte si et seulement si $G'\ra G$ est injectif,
$$
S'=\{\Ker(f|_{H}^{f(H)})\ |\ H\in S\}\quad\textrm{et}\quad S''=\{f(H)\ |\ H\in S\},
$$
auquel cas $G'=\Ker f$, $f:G\ra G''$ est surjectif, $\Ker(f|_{H}^{f(H)})=G'\cap H$ et $f(H)=H/G'\cap H$ pour tout $H\in S$.

\medskip

\item Un monomorphisme de $\cI$ 
$$
\xymatrix{
(G',S')\ar[r] & (G,S)
}
$$ 
est strict si et seulement si $G'\ra G$ est injectif et
$$
S'\supset\{G'\cap H\ |\ H\in S\},
$$
auquel cas les deux ensembles sont égaux.

\medskip

\item Un épimorphisme de $\cI$ 
$$
\xymatrix{
(G,S)\ar[r] & (G'',S'')
}
$$
est strict si et seulement si
$$
S''\subset\{f(H)\ |\ H\in S\},
$$
auquel cas les deux ensembles sont égaux et $G\ra G''$ est surjectif.
\end{itemize}
\end{Lem}

\begin{proof}
Soit 
$$
\xymatrix{
0\ar[r] & (G',S') \ar[r] & (G,S) \ar[r]^f & (G'',S'') \ar[r] & 0
}
$$
une suite exacte de $\cI$, c'est-à-dire telle que la suite
\begin{equation}\label{secIndQ}
\xymatrix{
0\ar[r] & (H')_{H'\in S'} \ar[r] & (H)_{H\in S} \ar[r]^f & (H'')_{H''\in S''} \ar[r] & 0
}
\end{equation}
de $\Ind(\cQ)$ est exacte. Nous allons utiliser la démonstration (et pas seulement l'énoncé) de \cite{KS} 8.6.6 (a). On pose  
$$
I:=\{ (H,H'')\in S\times S'' \ |\ f(H)\subset H''\}
$$
muni du préordre filtrant $(H_1,H_{1}'')\leq (H_2,H_{2}'')\Leftrightarrow H_1\subset H_2 \textrm{ et } H_{1}''\subset H_{2}''$. Les projections sur le premier et le second facteur définissent deux applications préordonnées cofinales $I\ra S$ et $I\ra S''$. Pour $i=(H,H'')\in I$, on pose $Q_i:=H$, $R_i:=H''$, $f_i:Q_i\ra R_{i}:=f|_{H}^{H''}$ et $Q_{i}':=\Ker(f|_{H}^{H''})$. Alors $(f_i)_{i\in I}:(Q_i)_{i\in I}\ra (R_i)_{i\in I}$ s'identifie à $f:(H)_{H\in S} \ra  (H'')_{H''\in S''}$. Il en résulte que $(Q_{i}')_{i\in I}$ est isomorphe  à $(H')_{H'\in S'}$. Maintenant, comme $(H'')_{H''\in S''}$ est à transitions injectives, $(Q_{i}')_{i\in I}$ est aussi isomorphe à $(\Ker f|_{H}^{f(H)})_{H\in S}$. Or on vérifie facilement que l'ensemble de sous-groupes de $G'$ donné par $\{\Ker f|_{H}^{f(H)}\ |\ H\in S \}$ vérifie (I-1), (I-2) et (I-3) (noter que $\Ker f=G'$ d'après \ref{foncfibgenIexact}), définissant donc un groupe indalgébrique $(G',\{\Ker f|_{H}^{f(H)}\}_{H\in S})$. Le morphisme $(G',S')\ra (G',\{\Ker f|_{H}^{f(H)}\}_{H\in S})$ induit par l'identité de $G'$ est donc une égalité par pleine fidélité de $\cI\ra\Ind(\cQ)$ \ref{Iexacte}. Ceci étant, le morphisme $(\Ker f|_{H}^{f(H)})_{H\in S}\ra (H)_{H\in S}$ s'identifie à $(H')_{H'\in S'}\ra (H)_{H\in S}$. Il en résulte que $(f(H))_{H\in S}$ s'identifie à $(H'')_{H''\in S''}$. Or on vérifie facilement que l'ensemble de sous-groupes de $G''$ donné par $\{ f(H)\ |\ H\in S \}$ vérifie (I-1), (I-2) et (I-3) (noter que $\Ima f=G''$ d'après \ref{foncfibgenIexact}), définissant donc un groupe indalgébrique $(G'',\{f(H)\}_{H\in S})$. Le morphisme $(G'',\{f(H)\}_{H\in S})\ra (G'',S'')$ induit par l'identité de $G''$ est donc une égalité par pleine fidélité de $\cI\ra\Ind(\cQ)$ \ref{Iexacte}. 

Réciproquement, une suite de $\cI$ 
$$
\xymatrix{
0\ar[r] & (G',S') \ar[r] & (G,S) \ar[r]^f & (G'',S'') \ar[r] & 0
}
$$
telle que $G'\ra G$ est injectif, $S'=\{\Ker(f|_{H}^{f(H)})\ |\ H\in S\}$ et $S''=\{f(H)\ |\ H\in S\}$ est exacte, puisque les suites de $\cQ$
$$
\xymatrix{
0\ar[r] & \Ker(f|_{H}^{f(H)}) \ar[r] & H \ar[r]^f & f(H) \ar[r] & 0, & H\in S,
}
$$
le sont.

Les autres assertions du lemme sont maintenant claires.
\end{proof}

\section{Groupes indproalgébriques}

\begin{Def}\label{defIP}
Un \emph{groupe indproalgébrique} est un groupe ordinaire $G$ muni d'une famille non vide $S$ de sous-groupes et, pour tout $H\in S$, d'une structure de groupe proalgébrique sur $H$, ces données vérifiant les axiomes suivants :
\begin{itemize}

\medskip

\item (IP-1) $H,H'\in S\Rightarrow \exists H''\in S : H+H'\subset H''$.

\medskip

\item (IP-2) Si $H\in S$, les sous-groupes $H'$ contenus dans $H$ qui appartiennent à $S$ sont les sous-groupes fermés de $H$.

\medskip

\item (IP-3) L'application canonique $\li_{H\in S} H\ra G$ est un isomorphisme de groupes.

\medskip

\end{itemize}

L'ensemble $S$ est l'\emph{ensemble complet de définition} du groupe indproalgébrique $(G,S)$, et les sous-groupes 
$H\in S$ sont les \emph{sous-groupes de définition} de $(G,S)$.

\medskip

Un \emph{morphisme de groupes indproalgébriques} $(G_1,S_1)\ra (G_2,S_2)$ est un morphisme de groupes $f:G_1\ra G_2$ vérifiant la condition suivante : pour tout $H_1\in S_1$, on a $f(H_1)\in S_2$, et l'application $H_1\ra f(H_1)$ induite par $f$ est un morphisme de groupes proalgébriques.
\end{Def}

\begin{Rem}\label{IP2'}
L'axiome (IP-2) implique l'axiome suivant :

\medskip

\begin{itemize}
\item \emph{(IP-2)' Si $H,H'\in S$ avec $H$ contenu dans $H'$, alors l'inclusion $H\subset H'$ est une immersion fermée proalgébrique}.
\end{itemize}

\medskip

De plus, les axiomes (IP-1) et (IP-2) impliquent que $S$ est stable par somme et par intersection. Raisonnant comme en \ref{I2'}, on obtient tout de suite la stabilité par intersection, et pour la stabilité par somme, il suffit de montrer que la somme de deux sous-groupes fermés d'un groupe proalgébrique est un sous-groupe fermé. Comme l'image d'un sous-groupe fermé par un morphisme de groupes proalgébriques est un sous-groupe fermé (\cite{S2} 2.4 Proposition 6), il suffit de montrer la loi de groupe $m$ d'un groupe proalgébrique $(G,S)$ est un morphisme de groupes proalgébriques $(G\times G,<S\times S>)\ra (G,S)$. Or ceci est bien clair : pour tout $H\in S$, $m^{-1}(H)/H\times H$ est le noyau du morphisme de groupes quasi-algébriques $\om:G/H\times G/H\ra G/H$ ; donc $m^{-1}(H)/H\times H$ est fermé dans $G\times G/H\times H$, $m^{-1}(H)\in\ <S\times S>$ et le morphisme $G\times G/m^{-1}(H)\ra G/H$ induit par $m$ est quasi-algébrique.
\end{Rem}

\begin{Rem}\label{completengIP}
Soit $G$ un groupe, muni d'une famille $S'$ de sous-groupes vérifiant (IP-1), (IP-2)' et (IP-3). Alors il existe une plus petite famille de sous-groupes de $G$ contenant $S'$ et vérifiant (IP-1), (IP-2) et (IP-3). On l'appellera l'ensemble complet de définition de $G$ \emph{engendré par $S'$}, et on le notera $<S'>$. On dira aussi que $S'$ est un \emph{ensemble de définition} du groupe indproalgébrique $(G,<S'>)$. Explicitement, $<S'>$ est l'ensemble des sous-groupes $H\subset G$ tels qu'il existe $H'\in S'$ vérifiant : 

\medskip

\begin{itemize}
\item $H\subset H'$ ;
\item l'inclusion $H\subset H'$ est une immersion fermée proalgébrique.
\end{itemize}

\medskip

En particulier, $S'$ est cofinal dans $S$. En revanche, comme en \ref{completengI}, il peut y avoir deux ensembles complets de définition $S_1\subset S_2$ distincts sur un même groupe ordinaire $G$, donnant lieu à des groupes indproalgébriques 
$(G,S_1)$ et $(G,S_2)$ non isomorphes. 
\end{Rem}

\begin{Not}\label{NotIP}
On note $\cI\cP$ la catégorie des groupes indproalgébriques. On a des foncteurs pleinement fidèles 
\begin{eqnarray*}
\cP&\ra&\cI\cP\\
(G,S) & \mapsto & (G,<\{(G,S)\}>) ;
\end{eqnarray*}
\begin{eqnarray*}
\cI&\ra&\cI\cP\\
(G,S) & \mapsto & (G,S),
\end{eqnarray*}
où les groupes quasi-algébriques $H\in S$ sont vus comme des groupes proalgébriques via \ref{NotP}.
\end{Not}

\begin{Prop}\label{IPexacte}
Le foncteur de $\cI\cP$ dans $\Ind(\cP)$ qui à $(G,S)$ associe $(H)_{H\in S}$ est pleinement fidèle : 
$$
\xymatrix{
\cI\cP \ar@{^{(}->}[r] & \Ind(\cP),
}
$$
et son image essentielle est la sous-catégorie épaisse de la catégorie abélienne $\Ind(\cP)$, formée des ind-objets de $\cP$ \emph{isomorphes à un ind-objet à transitions injectives}. En particulier, la catégorie $\cI\cP$ est exacte. 
\end{Prop}

On a de plus le lemme suivant :

\begin{Lem}\label{condisoinjIP}
Soit $((H_i,S_i))_{i\in I}\in\Ind(\cP)$. Alors, pour que $((H_i,S_i))_{i\in I}$ soit isomorphe à un ind-objet à transitions injectives, il faut et il suffit qu'il vérifie la condition suivante :

\medskip

$(*)\quad$ Pour tout $i\in I$, il existe $j_i\geq i$ tel que $\Ker (H_{i}\ra\li_{i\in I} H_i) = \Ker (H_{i}\ra H_{j_i})$.
\end{Lem}

\begin{proof}[Démonstrations du lemme \ref{condisoinjIP} et de la proposition \ref{IPexacte}]
Identiques à celles du lemme \ref{condisoinjI} et de la proposition \ref{Iexacte}, en y remplaçant la catégorie $\cQ$ par la catégorie $\cP$ et la catégorie $\cI$ par la catégorie $\cI\cP$.
\end{proof}

\begin{Def}\label{deffibgenIP}
On conserve les notations générales de la définition \ref{deffibgenI}.
\begin{itemize}

\medskip

\item Le \emph{foncteur de points} d'un groupe proalgébrique $(G,S)$ est le foncteur 
\begin{eqnarray*}
h_{(G,S)}:(\parf/k)^{\circ} & \ra & (\Groupes) \\
X  & \ra &  \lp_{H\in S} (G/H)(X).
\end{eqnarray*}

\medskip

\item Le \emph{foncteur de points} d'un groupe indproalgébrique $(G,S)$ est le foncteur 
\begin{eqnarray*}
h_{(G,S)}:(\parf/k)^{\circ} & \ra & (\Groupes) \\
X  & \ra &  \li_{H\in S} H(X).
\end{eqnarray*}

\medskip

\item La \emph{fibre générique} d'un groupe indproalgébrique $(G,S)$ est la restriction du foncteur de points de $(G,S)$ à la catégorie $\cT$ :
\begin{eqnarray*}
(G,S)_{\eta}:\cT^{\circ} & \ra & (\Groupes) \\
\kappa  & \ra &  \li_{H\in S} H(\kappa).
\end{eqnarray*}

\medskip

\item Le \emph{foncteur fibre générique} de la catégorie des groupes indproalgébriques dans la catégorie des préfaisceaux en groupes sur $\cT$ est le foncteur
\begin{eqnarray*}
(\cdot)_{\eta}:\cI\cP & \ra & \widehat{\cT} \\
(G,S)  & \ra &  (G,S)_{\eta}.
\end{eqnarray*}
\end{itemize}
\end{Def}

\begin{Rem}\label{cbP}
Soit $\kappa\in\cT^{\circ}$, et notons $\cP_{\kappa}$ la catégorie des groupes proalgébriques sur $\kappa$. Le foncteur changement de base $\cQ\ra \cQ_{\kappa}$ \ref{cbQ} induit un foncteur $\Pro(\cQ)\ra\Pro(\cQ_k)$, et donc, d'après \ref{Pab}, un foncteur $\cP\ra \cP_{\kappa}$, encore appelé \emph{changement de base}. Explicitement, si $(G,S)$ est un groupe proalgébrique sur $k$, alors $(G,S)_{\kappa}$ est le groupe proalgébrique sur $\kappa$ de groupe ordinaire sous-jacent $\lp_{H\in S} (G/H)(\kappa)$ et d'ensemble complet de définition engendré par
$$
\xymatrix{
S_{\kappa}:=\{ \Ker\big(\lp_{H\in S} (G/H)(\kappa)\ra (G/H)(\kappa)\big)\ |\ H\in S\},
}
$$
le quotient $\lp_{H\in S} (G/H)(\kappa)/\Ker\big(\lp_{H\in S} (G/H)(\kappa)\ra (G/H)(\kappa)\big)=(G/H)(\kappa)$ étant muni de la structure quasi-algébrique $(G/H)_{\kappa}$ ; on notera que l'application $\lp_{H\in S} (G/H)_{\kappa}(\kappa)\ra (G/H)_{\kappa}(\kappa)$ est bien surjective d'après \cite{S2} 2.3 Proposition 2 b). En particulier, le groupe des $\kappa$-points du groupe proalgébrique $(G,S)$ s'identifie au groupe ordinaire sous-jacent au groupe proalgébrique $(G,S)_{\kappa}$ sur $\kappa$. On notera par ailleurs que le foncteur changement de base $\cP\ra \cP_{\kappa}$ est exact puisque $\cQ\ra\cQ_{\kappa}$ l'est (cf. \cite{KS} 8.6.6 (a)).
\end{Rem}

\begin{Lem}\label{foncfibgenIPexact}
Le foncteur fibre générique $(\cdot)_{\eta}:\cI\cP\ra\widehat{\cT}$  est exact.
\end{Lem}

\begin{proof}
Par restriction de $(\cdot)_{\eta}$ le long du plongement canonique $\cP\ra\cI\cP$ \ref{NotIP}, on obtient un foncteur 
\begin{eqnarray*}
(\cdot)_{\eta}:\cP & \ra & \widehat{\cT} \\
H  & \ra &  H_{\eta}.
\end{eqnarray*}
Celui-ci est exact d'après \ref{PidemQ} et \ref{cbP}. Par conséquent, le foncteur 
\begin{eqnarray*}
\Ind(\cP) & \ra & \widehat{\cT} \\
(H_i)_{i\in I}  & \ra &  \li_i (H_i)_{\eta}
\end{eqnarray*}
est aussi exact, d'après \cite{KS} 8.6.6 (a), et le fait que les limites inductives filtrantes dans la catégorie des groupes sont exactes. D'où le lemme puisque le foncteur $\cI\cP\ra\Ind(\cP)$ de \ref{IPexacte} est exact par définition.  
\end{proof}

\begin{Pt}\label{cefibgenIPcons}
Le foncteur fibre générique $(\cdot)_{\eta}:\cI\cP\ra\widehat{\cT}$ n'est pas conservatif. Soit en effet $G_2$ le produit dans 
$\cP$ indexé par $\bbN$ de copies du groupe fini $\bbZ/p$, et soit $G_1$ le groupe $(\bbZ/p)^{\bbN}$ muni de la structure indproalgébrique $S_1$ constituée des sous-groupes finis de $(\bbZ/p)^{\bbN}$. Alors l'identité de 
$(\bbZ/p)^{\bbN}$ définit un morphisme de groupes indproalgébriques $(G_1,S_1)\ra (G_2,S_2)$ qui n'est pas un isomorphisme.
\end{Pt}

\begin{Prop}\label{critfoncfibgenIPcons}
Pour qu'un morphisme de groupes indproalgébriques $f:(G_1,S_1)\ra (G_2,S_2)$ soit un isomorphisme, il faut et il suffit que les deux conditions suivantes soient vérifiées :
\begin{enumerate}
\item Le morphisme $f_{\eta}:(G_1,S_1)_{\eta}\ra (G_2,S_2)_{\eta}$ induit sur les fibres génériques est un isomorphisme.
\item Il existe un ensemble de définition $S_{2}'$ de $(G,S_2)$ tel que pour tout $H_{2}'\in S_{2}'$, il existe $H_1\in S_1$ avec $f(H_1)\subset H_{2}'$ et $H_{2}'/f(H_1)$ quasi-algébrique.
\end{enumerate}
\end{Prop}

\begin{Rem}
Par définition d'un morphisme de groupes indproalgébriques, on a $f(H_1)\in S_2$, donc $f(H_1)\subset H_{2}'$ si et seulement si $f(H_1)$ est un sous-groupe fermé de $H_{2}'$ (d'après (IP-2)'). Dans ce cas, le groupe proalgébrique  $H_{2}'/f(H_1)$ est quasi-algébrique si et seulement si $f(H_1)$ est un sous-groupe de définition du groupe proalgébrique $H_{2}'$ (définition \ref{defP}) : \cite{S2} 2.4 Proposition 5.2.
\end{Rem}

\begin{proof}[Démonstration du lemme \ref{critfoncfibgenIPcons}]
Soit $f:(G_1,S_1)\ra (G_2,S_2)$ vérifiant les conditions \emph{1.} et \emph{2.} du lemme. Comme en \ref{foncfibgenIcons}, $$
f(S_1):=\{ f(H_1)\ |\ H_1\in S_1 \}
$$
est un ensemble de groupes proalgébriques contenu dans $S_2$, et il s'agit de montrer qu'il est égal à $S_2$. De plus, il suffit pour cela de montrer que pour tout $H_{2}'\in S_{2}'$, il existe $H_{1}'\in S_1$ tel que $H_{2}'\subset f(H_{1}')$ : d'après (IP-2) pour $S_2$, l'inclusion $H_{2}'\subset f(H_{1}')$ est alors une immersion fermée proalgébrique, puis d'après (IP-2) pour $S_1$ (et l'isomorphisme $H_{1}'\xrightarrow{\sim} f(H_{1}')$), $H_{2}'\in f(S_1)$ ; donc $S_2=\ <S_{2}'>\ \subset f(S_1)$, puisque $f(S_1)$ vérifie (IP-1), (IP-2) et (IP-3).

Soit $H_{2}'\in S_{2}'$. Soit $H_1\in S_1$ avec $f(H_1)\subset H_{2}'$ et $H_{2}'/f(H_1)$ quasi-algébrique. Soient 
$\eta_1,\ldots,\eta_n$ les points génériques de $(H_{2}'/f(H_1))^{\pf}$. Choisissons un corps algébriquement clos $\kappa$ contenant les corps résiduels $k(\eta_1),\ldots,k(\eta_n)$, et notons $\overline{\eta_1},\ldots,\overline{\eta_n}$ les $\kappa$-points de $(H_{2}'/f(H_1))^{\pf}$ induits par les $\eta_i$. Le foncteur $\cP \ni H \mapsto H_{\eta} \in \widehat{\cT}$ étant exact (\ref{foncfibgenIPexact}), l'application $H_{2}'(\kappa)\ra (H_{2}'/f(H_1))(\kappa)$ est surjective, et on peut donc choisir des relèvements $\widetilde{\overline{\eta_1}},\ldots,\widetilde{\overline{\eta_n}}$ des $\overline{\eta_i}$ dans $H_{2}'(\kappa)$. Comme par hypothèse l'application
$$
\xymatrix{
f(\kappa):\li_{H_1\in S_1} H_{1}(\kappa) \ar[r] & \li_{H_2\in S_2} H_{2}(\kappa)
}
$$
est bijective, il existe $H_{1}'\in S_1$ tel que $\widetilde{\overline{\eta_1}},\ldots,\widetilde{\overline{\eta_n}}\in f(\kappa)(H_{1}'(\kappa))$. D'après (IP-1) pour $S_1$, on peut en outre supposer que $H_{1}'$ contient $H_1$. Considérons alors le sous-groupe $I':=f(H_{1}')\cap H_{2}'$ de $G_2$, et choisissons $H_{2}''\in S_2$ contenant $f(H_{1}')$ et $H_{2}'$ ((IP-1) pour $S_2$). Alors $I'\in S_2$, et réalise le produit fibré dans $\cP$ de $f(H_{1}')$ et $H_{2}'$ au-dessus de $H_{2}''$ (\ref{IP2'}, stabilité de $S_2$ par intersection). Comme $\cP \ni H \mapsto H_{\eta} \in \widehat{\cT}$ est exact (\ref{foncfibgenIPexact}), on en déduit que $I'(\kappa)=f(\kappa)(H_{1}'(\kappa))\cap H_{2}'(\kappa)\subset H_{2}''(\kappa)$. Par conséquent, l'image du groupe proalgébrique $I'$ dans le groupe quasi-algébrique $H_{2}'/f(H_1)$ est un sous-groupe fermé $F$ tel que $\overline{\eta_1},\ldots,\overline{\eta_n}\in F(\kappa)$. Mais alors $F^{\pf}$ est fermé et schématiquement dense dans $(H_{2}'/f(H_1))^{\pf}$, donc lui est égal. En prenant les $k$-points, on obtient $F=H_{2}'/f(H_1)$, d'où $I'=H_{2}'$ puisque $I'$ contient $f(H_1)$. On a ainsi obtenu $f(H_{1}')\cap H_{2}'=H_{2}'$, i.e. $H_{2}'\subset f(H_{1}')$.
\end{proof}

\begin{Pt}
La catégorie exacte $\cI\cP$ des groupes indproalgébriques est donc munie de deux foncteurs exacts vers des catégories abéliennes : le plongement \ref{IPexacte}, et le foncteur fibre générique \ref{foncfibgenIPexact},
$$
\xymatrix{
\cI\cP \ar@{^{(}->}[r] \ar[d]_{(\cdot)_{\eta}} & \Ind(\cP) \\
\widehat{\cT},
}
$$
celui-ci donnant lieu au critère d'isomorphisme \ref{critfoncfibgenIPcons}.
\end{Pt}
 
\begin{Lem}\label{descrseIP}
\begin{itemize}

\medskip

\item Une suite de $\cI\cP$
$$
\xymatrix{
0\ar[r] & (G',S') \ar[r] & (G,S) \ar[r]^f & (G'',S'') \ar[r] & 0
}
$$
est exacte si et seulement si $G'\ra G$ est injectif,
$$
S'=\{\Ker(f|_{H}^{f(H)})\ |\ H\in S\}\quad\textrm{et}\quad S''=\{f(H)\ |\ H\in S\},
$$
auquel cas $G'=\Ker f$, $f:G\ra G''$ est surjectif, $\Ker(f|_{H}^{f(H)})=G'\cap H$ et $f(H)=H/G'\cap H$ pour tout $H\in S$.

\medskip

\item Un monomorphisme de $\cI\cP$ 
$$
\xymatrix{
(G',S')\ar[r] & (G,S)
}
$$ 
est strict si et seulement si $G'\ra G$ est injectif et
$$
S'\supset\{G'\cap H\ |\ H\in S\},
$$
auquel cas les deux ensembles sont égaux.

\medskip

\item Un épimorphisme de $\cI\cP$ 
$$
\xymatrix{
(G,S)\ar[r] & (G'',S'')
}
$$
est strict si et seulement si
$$
S''\subset\{f(H)\ |\ H\in S\},
$$
auquel cas les deux ensembles sont égaux et $G\ra G''$ est surjectif.
\end{itemize}
\end{Lem}

\begin{proof}
Identique à celle du lemme \ref{descrseIP}, en y remplaçant la catégorie $\cQ$ par la catégorie $\cP$.
\end{proof}

\section{Groupes proindalgébriques}

\begin{Def} \label{defPI}
Un \emph{groupe proindalgébrique} est un groupe ordinaire $G$ muni d'une famille non vide $S$ de sous-groupes et, pour tout $H\in S$, d'une structure de groupe indalgébrique sur $G/H$, ces données vérifiant les axiomes suivants :
\begin{itemize}

\medskip

\item (PI-1) $H,H'\in S\Rightarrow \exists H''\in S : H''\subset H\cap H'$.

\medskip

\item (PI-2) Si $H\in S$, les sous-groupes $H'$ contenant $H$ qui appartiennent à $S$ sont les images réciproques des sous-groupes indalgébriques stricts de $G/H$.

\medskip

\item (PI-3) Si $H,H'\in S$, et si $H\subset H'$, l'application canonique $G/H\ra G/H'$ est un épimorphisme strict de groupes indalgébriques.

\medskip

\item (PI-4) L'application canonique $G\ra\lp_{H\in S} G/H$ est un isomorphisme de groupes.

\medskip

\end{itemize}

L'ensemble $S$ est l'\emph{ensemble complet de définition} du groupe proindalgébrique $(G,S)$, et les sous-groupes $H\in S$ sont les \emph{sous-groupes de définition} de $(G,S)$.

\medskip

Un \emph{morphisme de groupes proindalgébriques} $(G_1,S_1)\ra (G_2,S_2)$ est un morphisme de groupes $f:G_1\ra G_2$ vérifiant la condition suivante : pour tout $H_2\in S_2$, on a $f^{-1}(H_2)\in S_1$, et l'application $G_1/f^{-1}(H_2)\ra G_2/H_2$ induite par $f$ est un morphisme de groupes indalgébriques. \footnote{On notera que l'on ne demande pas que ce morphisme de groupes indalgébriques soit un monomorphisme strict, alors que la condition analogue pour les morphismes dans $\cP$, $\cI$ et $\cI\cP$ était automatique, et qu'elle ne l'est pas ici (prendre par exemple pour $f$ un monomorphisme non strict de groupes indalgébriques comme le morphisme $(G_1,S_1)\ra (G_2,S_2)$ de \ref{completengI}.) Cela permet d'avoir un foncteur naturel de $\cI$ dans $\cP\cI$ (cf. \ref{NotPI}), ainsi qu'un foncteur pleinement fidèle de $\cP\cI$ dans 
$\Pro(\cI)$ (cf. \ref{PIexacte}).}
\end{Def}

\begin{Rem} \label{PI2'} 
Les axiomes (PI-1) et (PI-2) impliquent que $S$ est stable par intersection. Si $H,H'\in S$, il existe $H''\in S$ contenu dans $H\cap H'$, $H/H''$ et $H'/H''$ sont alors des sous-groupes indalgébriques stricts de $G/H''$ ; or l'intersection de deux sous-groupes stricts d'un groupe indalgébrique est un sous-groupe strict (cf. \ref{descrseI}), donc $(H/H'')\cap(H'/H'')=(H\cap H')/H''$ est strict dans $G/H''$, et donc $H\cap H'\in S$.

L'axiome (PI-2) implique que si $H,H'\in S$ sont les images réciproques de sous-groupes quasi-algébriques de $G/H''\in\cI$ pour un $H''\in S$ inclus dans $H\cap H'$, alors $H+H'\in S$ : $H+H'$ est alors l'image réciproque du sous-groupe quasi-algébrique (\emph{a fortiori} strict) $H/H''+H'/H''$ de $G/H''$. 
\end{Rem}

\begin{Rem}\label{completengPI} 
Soit $G$ un groupe, muni d'une famille $S'$ de sous-groupes vérifiant (PI-1), (PI-3) et (PI-4). Alors il existe une plus petite famille de sous-groupes de $G$ contenant $S'$ et vérifiant (PI-1), (PI-2), (PI-3) et (PI-4). On l'appellera l'ensemble complet de définition de $G$ \emph{engendré par $S'$}, et on le notera $<S'>$. On dira aussi que $S'$ est un \emph{ensemble de définition} du groupe proindalgébrique $(G,<S'>)$. Explicitement, $<S'>$ est l'ensemble des sous-groupes $H\subset G$ tels qu'il existe $H'\in S'$ vérifiant : 

\medskip

\begin{itemize}
\item $H'\subset H$ ;
\item l'inclusion $H/H'\subset G/H'$ est un monomorphisme strict de groupes indalgébriques.
\end{itemize}

\medskip

En particulier, $S'$ est coinitial dans $S$. En revanche, comme en \ref{completengI}, il peut y avoir deux ensembles complets de définition $S_1\subset S_2$ distincts sur un même groupe ordinaire $G$, donnant lieu à des groupes proindalgébriques 
$(G,S_1)$ et $(G,S_2)$ non isomorphes. 
\end{Rem}

\begin{Not}\label{NotPI}
On note $\cP\cI$ la catégorie des groupes proindalgébriques. On a des foncteurs pleinement fidèles 
\begin{eqnarray*}
\cI&\ra&\cP\cI\\
(G,S) & \mapsto & (G,<\{0\}>) ; 
\end{eqnarray*}
\begin{eqnarray*}
\cP&\ra&\cP\cI\\
(G,S) & \mapsto & (G,S),
\end{eqnarray*}
où les groupes quasi-algébriques $(G/H)$, $H\in S$, sont vus comme des groupes indalgébriques via \ref{NotI}.
\end{Not}

\begin{Prop}\label{PIexacte}
Le foncteur de $\cP\cI$ dans $\Pro(\cI)$ qui à $(G,S)$ associe $(G/H)_{H\in S}$ est pleinement fidèle : 
$$
\xymatrix{
\cP\cI \ar@{^{(}->}[r] & \Pro(\cI),
}
$$
et son image essentielle est la sous-catégorie épaisse de la catégorie exacte $\Pro(\cI)$, formée des pro-objets de $\cI$ \emph{isomorphes à un pro-objet dont les transitions sont des épimorphismes stricts}. En particulier, la catégorie $\cP\cI$ est exacte. 
\end{Prop}

On a de plus le lemme suivant :

\begin{Lem}\label{condisoepistrictPI}
Soit $((I_n,S_n))_{n\in N}\in\Pro(\cI)$. Alors, pour que $((I_n,S_n))_{n\in N}$ soit isomorphe à un pro-objet dont les transitions sont des épimorphismes stricts, il faut et il suffit qu'il vérifie la condition suivante :

\medskip

$(*)\quad$ Pour tout $n\in N$, il existe $m_n\geq n$ tel que 
$$m\geq m_n\Rightarrow\Ima ((I_{m},S_{m})\ra (I_n,S_n)) = \Ima ((I_{m_n},S_{m_n})\ra (I_n,S_n)).$$
\end{Lem}

\begin{proof}
Supposons qu'il existe un isomorphisme $\iota:((I_{n'}',S_{n'}'))_{n'\in N'}\xrightarrow{\sim}((I_n,S_n))_{n\in N}$, avec 
$((I_{n'}',S_{n'}'))_{n'\in N'}$ ayant pour transitions des épimorphismes stricts, et soit $n\in N$. Il existe alors $n'\in n'$, $m\geq n\in I$, et un diagramme 
$$
\xymatrix{
                        & (I_m,S_m) \ar[d]^{\tau_{m,n}} \ar[dl]_{(\iota^{-1})_{m,n'}}\\
(I_{n'}',S_{n'}') \ar[r]^{\iota_{n',n}} & (I_n,S_n).
}
$$
Alors $\iota_{n',n}\circ(\iota^{-1})_{m,n'}$ et $\tau_{m,n}$ représentent chacun la composante d'indice $n$ de l'identité dans
$$
\Hom_{\Pro(\cI)}(((I_n,S_n))_{n\in N},((I_n,S_n))_{n\in N})=\lp_n(\li_{m}\Hom_{\cI}((I_{m},S_{m}),(I_n,S_n))).
$$ 
Quitte à augmenter $m$, on peut donc supposer que le diagramme précédent est commutatif ; on pose alors $m_n:=m$. Soit $m\geq m_n$. Soit $H_n\in S_n$ de la forme $H_n=\tau_{m_n,n}(H_{m_n})$ pour un $H_{m_n}\in S_{m_n}$. Alors 
$(\iota^{-1})_{m_n,n'}(H_{m_n})\in S_{n'}'$. Raisonnant comme précédemment, il existe $m'\geq n'$ assez grand tel que l'on a un diagramme commutatif
$$
\xymatrix{
                                                                     && (I_m,S_m) \ar[d]^{\tau_{m,m_n}} \ar@/^3pc/[dd]^{\tau_{m,n}} \\
                                                                    && (I_{m_n},S_{m_n}) \ar[d]^{\tau_{m_n,n}} \ar[dl]_{(\iota^{-1})_{m_n,n'}} \\
(I_{m'}',S_{m'}') \ar[r]^{\tau_{m',n'}} \ar[uurr]^{\iota_{m',m}}  & (I_{n'}',S_{n'}') \ar[r]^{\iota_{n',n}} & (I_n,S_n).
}
$$
Comme la transition $\tau_{m',n'}:(I_{m'}',S_{m'}')\ra (I_{n'}',S_{n'}')$ est un épimorphisme indalgébrique strict par hypothèse, il existe $H_{m'}'\in S_{m'}'$ tel que $\tau_{m',n'}(H_{m'}')=(\iota^{-1})_{m_n,n'}(H_{m_n})$, d'après \ref{descrseI}. Posant $H_m:=\iota_{m',m}(H_{m'}')\in S_m$, la commutativité du diagramme entraîne $\tau_{m,n}(H_m)=H_n$. On a ainsi montré
$$
\{\tau_{m,n}(H_m)\ |\ H_m\in S_m\}\supset \{\tau_{m_n,n}(H_{m_n})\ |\ H_{m_n}\in S_{m_n}\}.
$$
L'inclusion en sens opposé étant évidente, ces deux ensembles sont donc égaux, ce qui signifie que les deux groupes indalgébriques $\Ima ((I_{m},S_{m})\ra (I_n,S_n))$ et $\Ima ((I_{m_n},S_{m_n})\ra (I_n,S_n))$ sont égaux.

Réciproquement, supposons que $((I_n,S_n))_{n\in N}$ vérifie la condition $(*)$. Pour tout $n\in N$ et $m\geq m_n\in N$, posons  
$$
(J_n,S_{J_n}):=\Ima ((I_{m_n},S_{m_n})\ra (I_n,S_n)) = \Ima ((I_{m},S_{m})\ra (I_n,S_n)).
$$
La transition $\tau_{m,n}:(I_m,S_m)\ra (I_n,S_n)$ induit un épimorphisme indalgébrique strict 
$$
\xymatrix{
\underline{\tau_{m,n}}:(I_m,S_m) \ar[r] & (J_n,S_{J_n})
}
$$
(cf. \ref{descrseI}). Composant avec le morphisme canonique $(J_{m},S_{J_m})\ra (I_{m},S_{m})$, on obtient un morphisme
$$
\xymatrix{
\underline{\underline{\tau_{m,n}}}:(J_m,S_{J_m}) \ar[r] & (J_n,S_{J_n}).
}
$$
Celui-ci est encore un épimorphisme strict d'après la condition pour $((I_n,S_n))_{n\in N}$ et \ref{descrseI}. puisque $K_i=\Ker(H_i\ra\li_{i\in I} H_i)$. Ceci étant, soit $(\underline{N},\leq_{\underline{N}})$ le petit ensemble préordonné filtrant défini par $\underline{N}:=N$ et 
$$n\leq_{\underline{N}} m \Leftrightarrow m_n\leq_{N} m.$$
Alors le pro-objet $(J_n,S_{J_n})_{n\in\underline{N}}$ dont les transitions sont données par les
$\underline{\underline{\tau_{m,n}}}$ est isomorphe à $((I_n,S_n))_{n\in N}$. Les morphismes canoniques
$$
\xymatrix{
(J_n,S_{J_n}) \ar[r] & (I_n,S_n), & n\in N,
}
$$
définissent en effet un morphisme $((J_n,S_{J_n}))_{n\in \underline{N}}\ra ((I_n,S_n))_{n\in N}$, ayant pour inverse le morphisme défini par les
$$
\xymatrix{
(I_{m_n},S_{m_n}) \ar[r]^{\underline{\tau_{m_n,n}}} & (J_n,S_{J_n}), & n\in \underline{N}.
}
$$
\end{proof}

\begin{proof}[Démonstration de la proposition \ref{PIexacte}]
Par définition des morphismes dans $\cP\cI$, si $(G,S)$ et $(G',S')$ sont des objets de $\cP\cI$, alors on a une application 
$$
\xymatrix{
\Hom_{\cP\cI}((G,S),(G',S')) \ar[r] & \lp_{H'\in S'}\li_{H\in S}\Hom_{\cI}(G/H,G'/H'),
}
$$
qui a $f:(G,S)\ra (G',S')$ associe $(G/f^{-1}(H')\ra G'/H')_{H'\in S'}$. Cette application est évidemment injective. De plus, en utilisant (PI-2) et (PI-3), on vérifie que pour tout $H\in S$, la projection canonique $\pi_H:G\ra G/H$ est un morphisme de groupes proindalgébriques, puis en utilisant (PI-3) et la description des épimorphismes stricts de $\cI$ \ref{descrseI}, que si 
$$(f_{H}:G/H\ra G/H')_{H'\in S'}\in \lp_{H'\in S'}\li_{H\in S}\Hom_{\cI}(G/H,G'/H'),$$
alors $\lp_{H'\in S'}(\pi_H\circ f_{H})$ est un morphisme proindalgébrique $(G,S)\ra (G',S')$, d'image $(f_{H})_{H'\in S'}$ puisque les $\pi_H$ sont surjectives. Le foncteur $\cP\cI\ra\Pro(\cI)$ que l'on considère est donc bien pleinement fidèle. Les pro-objets de son image ont pour transitions des épimorphismes stricts d'après (PI-3).

Inversement, soit $(I_n,S_n)_{n\in N}\in\Pro(\cI)$ ayant pour transitions des épimorphismes stricts (en particulier, surjectifs, cf. \ref{descrseI}). Posons $G:=\lp_{n\in N}$, et pour $n\in\bbN$, notons $\pi_n:G\ra I_n$ l'application canonique. Comme $N$ est essentiellement dénombrable par hypothèse, l'application $\pi_n$ est surjective (\cite{BouE} III \S 7 N${}^{\circ}4$ Proposition 5). En particulier, posant $H_n:=\Ker\pi_n\subset G$, $G/H_n$ s'identifie à $I_n$, et se trouve ainsi muni d'une structure de groupe indalgébrique. De plus, la famille $S:=(H_n)_{n\in\bbN}$ vérifie (PI-1), (PI-3) et (PI-4). D'où un groupe proindalgébrique $(G,<S>)$ dont l'image dans $\Pro(\cI)$ est isomorphe à $(I_n,S_n)_{n\in N}$.

D'après \ref{Iexacte} et \cite{KS} 8.6.12, la catégorie $\Pro(\cI)$ est équivalente à une sous-catégorie épaisse de la catégorie abélienne $\Pro(\Ind(\cQ))$, donc est exacte. Il reste à vérifier que la sous-catégorie pleine de $\Pro(\cI)$ formée des pro-objets isomorphes à un pro-objet dont les  transitions sont des épimorphismes stricts est épaisse. Soit 
\begin{equation}\label{seIepaissePI}
\xymatrix{
0\ar[r] & I' \ar[r] & I \ar[r] & I'' \ar[r] & 0
}
\end{equation}
une suite exacte dans $\Pro(\cI)$, avec $I'$ et $I''$ isomorphes à des pro-objets dont les transitions sont des épimorphismes stricts. Comme dans la démonstration de \ref{descrseI}, il existe un petit ensemble préordonné filtrant $N$, et un système de suites exactes courtes 
$$
\xymatrix{
0 \ar[r] & (I_{n}',S_{n}') \ar[r] & (I_n,S_n) \ar[r]^{f_n} & (I_{n}'',S_{n}'') \ar[r] & 0, & n\in N,
}
$$
tel que la suite 
$$
\xymatrix{
0 \ar[r] & ((I_{n}',S_{n}'))_{n\in N} \ar[r] & ((I_n,S_n))_{n\in N} \ar[r]^{(f_n)_{n\in N}} & ((I_{n}'',S_{n}''))_{n\in N} \ar[r] & 0
}
$$
soit isomorphe à (\ref{seIepaissePI})\footnote{Dans la démonstration de \ref{descrseI}, on reprend la démonstration de \cite{KS} 8.6.6 (a) montrant qu'une suite exacte courte de $\Ind(\cQ)$ est isomorphe à un système inductif de suites exactes courtes de $\cQ$. Il est clair que les arguments restent valides quand on remplace $\Ind$ par $\Pro$. En revanche, le fait de remplacer $\cQ$ par $\cI$ n'est \emph{a priori} pas immédiat, puisque $\cI$ n'est pas abélienne mais seulement exacte. Cependant, le noyau et l'image d'un morphisme dans $\cI$, vu dans la catégorie abélienne $\Ind(\cQ)$ via le plongement \ref{Iexacte}, sont évidemment dans $\cI$. On voit alors que les arguments restent bien valides quand on remplace $\cQ$ par $\cI$.}. D'après \ref{condisoepistrictPI}, les pro-objets $((I_{n}',S_{n}'))_{n\in N}$ et 
$((I_{n}'',S_{n}''))_{n\in N}$ vérifient la condition $(*)$, et il suffit de montrer que $((I_n,S_n))_{n\in N}$ la vérifie également. Soit $n\in N$. Choisissons $n_1\geq n$ tel que $\Ima ((I_{m}',S_{m}')\ra (I_{n}',S_{n}')) = \Ima ((I_{n_1}',S_{n_1}')\ra (I_{n}',S_{n}'))$ pour tout $m\geq n_1$, puis $n_2\geq n_1$ tel que $\Ima ((I_{m}'',S_{m}'')\ra (I_{n_1}'',S_{n_1}'')) = \Ima ((I_{n_2}'',S_{n_2}'')\ra (I_{n_1}'',S_{n_1}''))$ pour tout $m\geq n_2$. Alors le diagramme commutatif exact suivant montre $\Ima ((I_{m},S_{m})\ra (I_{n},S_{n})) = \Ima ((I_{n_2},S_{n_2})\ra (I_{n},S_{n}))$ pour tout $m\geq n_2$ :
$$
\xymatrix{
0 \ar[r] & (I_{m}',S_{m}') \ar[r] \ar[d]^{\tau_{m,n_2}'} & (I_{m},S_{m}) \ar[r]^{f_m} \ar[d]^{\tau_{m,n_2}} & (I_{m}'',S_{m}'') \ar[r] \ar[d]^{\tau_{m,n_2}''} & 0 \\ 
0 \ar[r] & (I_{n_2}',S_{n_2}') \ar[r] \ar[d]^{\tau_{n_2,n_1}'} & (I_{n_2},S_{n_2}) \ar[r]^{f_{n_2}} \ar[d]^{\tau_{n_2,n_1}} & (I_{n_2}'',S_{n_2}'') \ar[r] \ar[d]^{\tau_{n_2,n_1}''} & 0 \\
0 \ar[r] & (I_{n_1}',S_{n_1}') \ar[r] \ar[d]^{\tau_{n_1,n}'} & (I_{n_1},S_{n_1}) \ar[r]^{f_{n_1}} \ar[d]^{\tau_{n_1,n}} & (I_{n_1}'',S_{n_1}'') \ar[r] \ar[d]^{\tau_{n_1,n}''} & 0 \\
0 \ar[r] & (I_{n}',S_{n}') \ar[r]  & (I_{n},S_{n}) \ar[r]^{f_n} & (I_{n}'',S_{n}'') \ar[r] & 0. 
}
$$
En effet, soit $H_n\in S_n$ de la forme $\tau_{n_2,n}(H_{n_2})$ avec $H_{n_2}\in S_{n_2}$, et montrons que $H_n$ est l'image par $\tau_{m,n}$ d'un élément de $S_m$. Posons $H_{n_1}:=\tau_{n_2,n_1}(H_{n_2})$. Comme $f_{n_1}(H_{n_1})=\tau_{n_2,n_1}''(f_{n_2}(H_{n_2}))$, il existe par définition de $n_2$ un élément $H_{m}''\in S_{m}''$ tel que $f_{n_1}(H_{n_1})=\tau_{m,n_1}''(H_{m}'')$. Comme $f_m$ est un épimorphisme strict de $\cI$, il existe alors $H_m\in S_m$ tel $f_m(H_m)=H_{m}''$ (\ref{descrseI}). Par construction, $f_{n_1}(\tau_{m,n_1}(H_m))=f_{n_1}(H_{n_1})$, 
d'où (en rappelant que $I_{n_1}'=\Ker f_{n_1}$ par \ref{descrseI})
$$
H_{n_1}\subset \tau_{m,n_1}(H_m) + I_{n_1}'\cap(H_{n_1}+\tau_{m,n_1}(H_m)),
$$
puis donc
$$
H_n\subset \tau_{m,n}(H_m) + \tau_{n_1,n}'(I_{n_1}'\cap(H_{n_1}+\tau_{m,n_1}(H_m))).
$$
Maintenant, comme $S_{n_1}$ est stable par somme (\ref{I2'}) et $(I_{n_1}',S_{n_1}') \ra (I_{n_1},S_{n_1})$ est un monomorphisme strict de $\cI$, on a $I_{n_1}'\cap(H_{n_1}+\tau_{m,n_1}(H_m))\in S_{n_1}'$ (\ref{descrseI}). Il existe donc par définition de $n_1$ un élément $H_{m}'\in S_{m}'$ tel que $\tau_{n_1,n}'(I_{n_1}'\cap(H_{n_1}+\tau_{m,n_1}(H_m)))=\tau_{m,n}'(H_{m}')(=\tau_{m,n}(H_{m}')$ en identifiant $H_{m}'$ à son image dans $I_m$).  Ainsi donc
$$
H_n\subset \tau_{m,n}(H_m) + \tau_{m,n}(H_{m}') = \tau_{m,n}(H_m+H_{m}').
$$
Comme $S_m$ est stable par somme (et que l'on a identifié $H_{m}'$ à un élément de $S_m$), on a $\tau_{m,n}(H_m+H_{m}')\in S_n$, si bien que $H_n$ est un sous-groupe fermé de $\tau_{m,n}(H_m+H_{m}')$ d'après (I-2) pour $S_n$. Son image réciproque par le morphisme quasi-algébrique $H_m+H_{m}'\ra\tau_{m,n}(H_m+H_{m}')$ induit par 
$\tau_{m,n}$ est donc un sous-groupe fermé de $H_m+H_{m}'$, donc un élément de $S_m$ par (I-2) pour $S_m$, dont l'image par $\tau_{m,n}$ est $H_n$.
\end{proof}

\begin{Def}\label{deffibgenPI}
On conserve les notations de la définition \ref{deffibgenI}.
\begin{itemize}

\medskip

\item Le \emph{foncteur de points} d'un groupe proindalgébrique $(G,S)$ est le foncteur 
\begin{eqnarray*}
h_{(G,S)}:(\parf/k)^{\circ} & \ra & (\Groupes) \\
X  & \ra &  \lp_{H\in S} (G/H)(X).
\end{eqnarray*}

\medskip

\item La \emph{fibre générique} d'un groupe proindalgébrique $(G,S)$ est la restriction du foncteur de points de $(G,S)$ à la catégorie $\cT$ :
\begin{eqnarray*}
(G,S)_{\eta}:\cT^{\circ} & \ra & (\Groupes) \\
\kappa  & \ra &  \lp_{H\in S} (G/H)(\kappa).
\end{eqnarray*}

\medskip

\item Le \emph{foncteur fibre générique} de la catégorie des groupes proindalgébriques dans la catégorie des préfaisceaux en groupes sur $\cT$ est le foncteur
\begin{eqnarray*}
(\cdot)_{\eta}:\cP\cI & \ra & \widehat{\cT} \\
(G,S)  & \ra &  (G,S)_{\eta}.
\end{eqnarray*}
\end{itemize}
\end{Def} 

\section{Groupes admissibles} \label{groupesadmissibles}

\begin{Def}\label{defIPf}
La catégorie $\cI\cP^{\ad}$ est la sous-catégorie pleine de la catégorie $\cI\cP$, dont les objets sont les groupes indproalgébriques $(G,S)$ vérifiant la condition suivante :

\medskip 

$(\ad)$ Il existe un ensemble de définition $S'$ de $(G,S)$ tel que pour tout $H_{1}', H_{2}'\in S'$ avec $H_{1}'\subset H_{2}'$, le groupe proalgébrique $H_{2}'/H_{1}'$ est quasi-algébrique.
\end{Def}

\begin{Rem}
Rappelons que d'après \cite{S2} 2.4 Proposition 5.2, le groupe proalgébrique $H_{2}'/H_{1}'$ est quasi-algébrique si et seulement si $H_{1}'$ est un sous-groupe de définition du groupe proalgébrique $H_{2}'$.
\end{Rem}

\begin{Rem}\label{IetPdansIPf}
Les plongements $\cP\ra\cI\cP$ et $\cI\ra\cI\cP$ \ref{NotIP} se factorisent par $\cI\cP^{\ad}$.
\end{Rem}

\begin{Def}
Une \emph{origine} d'un groupe indproalgébrique $(G,S)$ est un élément $H_0\in S$ tel que pour tout $H\in S$ avec $H_0\subset H$, le groupe proalgébrique $H/H_0$ est quasi-algébrique.
\end{Def}

\begin{Lem}\label{defequivIPf}
Si $(G,S)\in\cI\cP^{\ad}$ et si $S'$ est un ensemble de définition de $(G,S)$ vérifiant la condition $(\ad)$, alors tout élément de $S'$ est une origine de $(G,S)$. L'ensemble des origines de $(G,S)$ est alors un ensemble de définition de $(G,S)$. 

Réciproquement, si $(G,S)\in\cI\cP$ possède une origine $H_0$, alors l'ensemble des $H\in S$ contenant $H_0$ est un ensemble de définition de $(G,S)$ vérifiant la condition $(\ad)$, et en particulier $(G,S)\in\cI\cP^{\ad}$.
\end{Lem}

\begin{proof}
Soient $(G,S)\in\cI\cP^{\ad}$, $S'$ vérifiant $(\ad)$ et $H_0\in S'$. Soit $H\in S$ contenant $H_0$. Comme $S'$ est cofinal dans $S$ 
(\ref{IP2'}), il existe $H'\in S'$ contenant $H$. Le quotient $H/H_0$ est alors un sous-groupe fermé du groupe quasi-algébrique $H'/H_0$, donc est quasi-algébrique. Le sous-groupe $H_0$ est donc une origine de $(G,S)$. L'ensemble des origines de $(G,S)$ contient ainsi l'ensemble de définition $S'$, c'est donc lui aussi un ensemble de définition de $(G,S)$.

Réciproquement, soit $(G,S)\in\cI\cP$ possédant une origine $H_0$. Alors l'ensemble $S'$ des $H\in S$ contenant $H_0$ est évidemment un ensemble de définition de $(G,S)$, et il vérifie la condition $(\ad)$ : si $H_{1}', H_{2}'\in S'$ avec $H_{1}'\subset H_{2}'$, alors $H_{2}'/H_{1}'$ est quotient du groupe quasi-algébrique $H_{2}'/H_0$, donc est quasi-algébrique.
\end{proof}

\begin{Rem}\label{propS0IPf}
L'ensemble de définition formé des origines de $(G,S)\in\cI\cP^{\ad}$ est stable par somme, intersection, passage à un majorant dans $S$, et passage à un minorant sous-groupe de définition de la structure proalgébrique. En effet, la stabilité par passage à un majorant dans $S$ résulte de \ref{defequivIPf}, et entraîne la stabilité par somme d'après \ref{IP2'}. Passons à la stabilité par passage à un minorant sous-groupe de définition de la structure proalgébrique. Soient $H_0$ une origine de $(G,S)$ et $H\in S_{H_0}$. Si $H'\in S$ contient $H$, alors en choisissant $H''\in S$ contenant $H_0$ et $H'$, on a $H''/H$ quasi-algébrique car $H_0/H$ et $H''/H_0$ le sont, et $H'/H$ est fermé dans $H''/H$, donc quasi-algébrique. Justifions enfin la stabilité par intersection. Soient $H_0$, $H_{0}'$ des origines de $(G,S)$. D'après ce que l'on vient de voir, il suffit de montrer que le sous-groupe fermé $H_0\cap H_{0}'$ de $H_0$ est un sous-groupe de définition de $H_0$, i.e. que $H_0/H_0\cap H_{0}'$ est quasi-algébrique. Or en choisissant $H\in S$ contenant $H_0$ et $H_{0}'$, on obtient un morphisme proalgébrique injectif $H_0/H_0\cap H_{0}'\ra H/H_{0}'$, et donc que $H_0/H_0\cap H_{0}'$ est quasi-algébrique puisque $H/H_{0}'$ l'est. 
\end{Rem}

\begin{Pt}
Le plongement canonique $\cI\cP^{\ad}\ra \cI\cP$ n'est pas une équivalence. Par exemple, si $(G,S)$ est un groupe proalgébrique non quasi-algébrique, alors la somme directe d'une infinité de copies de $(G,S)$ est un groupe indproalgébrique ne possédant pas d'origine.
\end{Pt}

\begin{Prop}\label{IPfexacte}
La catégorie $\cI\cP^{\ad}$ est épaisse dans $\cI\cP$ ; en particulier, elle est exacte.
\end{Prop}

\begin{proof}
Soit
$$
\xymatrix{
0\ar[r] & (G',S') \ar[r] & (G,S) \ar[r] & (G'',S'') \ar[r] & 0
}
$$
une suite exacte de $\cI\cP$, avec $(G',S'),(G'',S'')\in\cI\cP^{\ad}$. Choisissons une origine $H_{0}''\in S''$ de 
$(G'',S'')$. D'après \ref{descrseIP}, il existe une suite exacte de $\cP$
$$
\xymatrix{
0\ar[r] & H_{1}' \ar[r] & H_1 \ar[r] & H_{0}'' \ar[r] & 0
}
$$
avec $H_{1}'\in S'$ et $H_{1}\in S$. Soit $H_{0}'$ une origine de $(G',S')$ contenant $H_{1}'$. Identifions $H_{0}'$ à un élément de $S$ et choisissons $H_2\in S$ contenant $H_{0}'$ et $H_1$. Alors $H_2$ est une origine de $(G,S)$. En effet, $G'\cap H_2\in S'$ et contient $H_{0}'$, donc est une origine de $(G',S')$. Puis $H_2/G'\cap H_2\in S''$ et contient $H_{0}''$, 
donc est une origine de $(G'',S'')$. Donc si $H\in S$ contient $H_2$, de sorte que $G'\cap H\in S'$ et $H/G'\cap H\in S''$, le lemme du serpent appliqué au diagramme commutatif exact de $\cP$
$$
\xymatrix{
0\ar[r] & G'\cap H_{2} \ar[r] \ar[d] & H_2 \ar[r] \ar[d] & H_2/G'\cap H_{2}  \ar[r] \ar[d] & 0 \\
0\ar[r] & G'\cap H \ar[r] & H  \ar[r] & H/G'\cap H \ar[r] & 0 
}
$$
montre que $H/H_2$ est quasi-algébrique, puisque $\cQ$ est épaisse dans $\cP$ (\cite{S2} 2.4 Corollaire 1).
\end{proof}

\begin{Def}\label{defPIf}
La catégorie $\cP\cI^{\ad}$ est la sous-catégorie pleine de la catégorie $\cP\cI$, dont les objets sont les groupes proindalgébriques $(G,S)$ vérifiant la condition suivante :

\medskip 

$(\ad)$ Il existe un ensemble de définition $S'$ de $(G,S)$ tel que pour tout $H_{1}', H_{2}'\in S'$ avec $H_{1}'\subset H_{2}'$, le groupe indalgébrique $H_{2}'/H_{1}'$ est quasi-algébrique.
\end{Def}

\begin{Rem}
Rappelons que par (PI-3), le groupe quotient $H_{2}'/H_{1}'$ est muni d'une structure de groupe indalgébrique, noyau de l'épimorphisme strict de groupes indalgébriques $G/H_{1}'\ra G/H_{2}'$.
\end{Rem}

\begin{Rem}\label{IetPdansPIf}
Les plongements $\cI\ra\cP\cI$ et $\cP\ra\cP\cI$ \ref{NotPI} se factorisent par $\cP\cI^{\ad}$.
\end{Rem}

\begin{Def}
Une \emph{origine} d'un groupe proindalgébrique $(G,S)$ est un élément $H_0\in S$ tel que pour tout $H\in S$ avec $H\subset H_0$, le groupe indalgébrique $H_0/H$ est quasi-algébrique.
\end{Def}

\begin{Lem}\label{defequivPIf}
Si $(G,S)\in\cP\cI^{\ad}$ et si $S'$ est un ensemble de définition de $(G,S)$ vérifiant la condition $(\ad)$, alors tout élément de $S'$ est une origine de $(G,S)$. L'ensemble des origines de $(G,S)$ est alors un ensemble de définition de $(G,S)$. 

Réciproquement, si $(G,S)\in\cP\cI$ possède une origine $H_0$, alors l'ensemble des $H\in S$ inclus dans $H_0$ est un ensemble de définition de $(G,S)$ vérifiant la condition $(\ad)$, et en particulier $(G,S)\in\cP\cI^{\ad}$.
\end{Lem}

\begin{proof}
Soient $(G,S)\in\cP\cI^{\ad}$, $S'$ vérifiant $(\ad)$ et $H_0\in S'$. Soit $H\in S$ inclus dans $H_0$. Comme $S'$ est coinitial dans $S$ (\ref{completengPI}), il existe $H'\in S'$ inclus dans $H$. On a alors le diagramme commutatif exact de $\cI$ (cf. (PI-3))
$$
\xymatrix{
0\ar[r] & H_0/H' \ar[r] \ar[d] & G/H' \ar[r] \ar[d] & G/H_0 \ar[r] \ar@{=}[d] & 0 \\
0\ar[r] & H_0/H \ar[r] & G/H \ar[r] & G/H_0 \ar[r] & 0.
}
$$
En particulier, le noyau de $H_0/H'\ra H_0/H$ dans $\Ind(\cQ)$ est isomorphe à celui de l'épimorphisme strict de groupes indalgébriques $G/H'\ra G/H$ (cf. PI-3), donc est dans l'image essentielle de $\cI\ra\Ind(\cQ)$. Par conséquent $H_0/H'\ra H_0/H$ est un épimorphisme strict de $\cI$, et comme $H_0/H'$ est quasi-algébrique, $H_0/H$ l'est donc également. Le sous-groupe $H_0$ est donc une origine de 
$(G,S)$. L'ensemble des origines de $(G,S)$ contient ainsi l'ensemble de définition $S'$, c'est donc lui aussi un ensemble de définition de $(G,S)$.

Réciproquement, soit $(G,S)\in\cP\cI$ possédant une origine $H_0$. Alors l'ensemble $S'$ des $H\in S$ inclus dans $H_0$ est évidemment un ensemble de définition de $(G,S)$, et il vérifie la condition $(\ad)$ : si $H_{1}', H_{2}'\in S'$ avec $H_{1}'\subset H_{2}'$, alors on a le diagramme commutatif exact de $\cI$ (cf. (PI-3)) 
$$
\xymatrix{
0\ar[r] & H_0/H_{1}' \ar[r] \ar[d] & G/H_{1}' \ar[r] \ar[d] & G/H_0 \ar[r] \ar@{=}[d] & 0 \\
0\ar[r] & H_0/H_{2}' \ar[r] & G/H_{2}' \ar[r] & G/H_0 \ar[r] & 0.
}
$$
En particulier, le groupe indalgébrique $H_{2}'/H_{1}'$ est le premier terme d'une suite exacte courte de $\Ind(\cQ)$ dont les deux autres termes sont dans $\cI$
$$
\xymatrix{
0 \ar[r] & H_{2}'/H_{1}' \ar[r] & H_0/H_{1}' \ar[r] & H_0/H_{2}' \ar[r] & 0.
}
$$
D'où un monomorphisme strict de $H_{2}'/H_{1}'$ dans le groupe quasi-algébrique $H_0/H_{1}'$, et par conséquent $H_{2}'/H_{1}'$ est quasi-algébrique.
\end{proof}

\begin{Rem}\label{propS0PIf}
L'ensemble de définition formé des origines de $(G,S)\in\cP\cI^{\ad}$ est stable par intersection, somme, passage à un minorant dans $S$, et passage à un majorant image réciproque d'un sous-groupe de définition de la structure indalgébrique quotient. En effet, la stabilité par passage à un minorant dans $S$ résulte de \ref{defequivPIf}, et entraîne la stabilité par intersection d'après \ref{PI2'}. Passons à la stabilité par passage à un majorant image réciproque d'un sous-groupe de définition de la structure indalgébrique quotient. Soient $H_0$ une origine de $(G,S)$ et $H\in S$ contenant $H_0$ tel que $H/H_0\in S_{G/H_0}$. Si $H'\in S$ est inclus dans $H$, alors en choisissant $H''\in S$ inclus dans $H_0$ et $H'$, on a $H/H''$ quasi-algébrique car $H/H_0$ et $H_0/H''$ le sont, et $H/H'$ est quotient strict de $H/H''$, donc quasi-algébrique. Justifions enfin la stabilité par somme. Soient $H_0$, $H_{0}'$ des origines de $(G,S)$. D'après ce que l'on vient de voir, il suffit de montrer que $H_0+H_{0}'\in S$ (cf. \ref{PI2'}) est tel que le sous-groupe indalgébrique strict 
$(H_0+H_{0}')/H_0\subset G/H_0$ est de définition, i.e. est quasi-algébrique. Or on a un morphisme indalgébrique $H_{0}'/(H_0\cap H_{0}')\ra (H_0+H_{0}')/H_0$ bijectif avec $H_{0}'/(H_0\cap H_{0}')$ quasi-algébrique, et par conséquent $(H_0+H_{0}')/H_0$ est bien quasi-algébrique.
\end{Rem}

\begin{Pt}
Le plongement canonique $\cP\cI^{\ad}\ra \cP\cI$ n'est pas une équivalence. Par exemple, si $(G,S)$ est un groupe indalgébrique non quasi-algébrique, alors le produit d'une infinité de copies de $(G,S)$ est un groupe proindalgébrique ne possédant pas d'origine.
\end{Pt}

\begin{Pt} \label{foncIPPI}
Soit $(G,S)\in\cI\cP^{\ad}$. Si $H_0$ est un groupe proalgébrique origine de $(G,S)$, alors son ensemble complet de définition $S_{H_0}$ définit une structure $\cP\cI^{\ad}$ sur $G$, de la façon suivante. Soit $H\in S_{H_0}$. Alors $H$ est une origine de $(G,S)$ d'après \ref{propS0IPf}. De plus, l'ensemble de sous-groupes de $G/H$
$$
S_{G/H}:=\{H'/H\ |\ H\subset H'\in S\}
$$
vérifie (I-1), (I-2) et (I-3), munissant donc $G/H$ d'une structure de groupe indalgébrique. Ceci étant, l'ensemble
$$
\{H\subset G\ |\ H\in S_{H_0}\}
$$
vérifie (PI-1), (PI-3) et (PI-4), et par conséquent munit le groupe $G$ d'une structure de groupe proindalgébrique, notée 
$(G,<S_{H_0}>)$ (cf. \ref{completengPI}). Enfin, pour tout $H,H'\in S_{H_0}$ avec $H\subset H'$, on a $H'/H\in\cQ$, si bien que $(G,<S_{H_0}>)\in\cP\cI^{\ad}$. 

La structure proindalgébrique ainsi obtenue sur $G$ ne dépend pas du choix de l'origine $H_0$. Vérifions en effet que si $H_{0}'$ est une autre origine de $(G,S)$, alors $(G,<S_{H_{0}'}>)=(G,<S_{H_0}>)$. On peut supposer $H_0\subset H_{0}'$.
Alors $S_{H_0}\subset S_{H_{0}'}$, \emph{a fortiori} $<S_{H_0}>\ \subset\ <S_{H_{0}'}>$. Inversement, si $H'\in S_{H_{0}'}$, il existe $H\in S_{H_{0}'}$ inclus dans $H_0\cap H'$ ((P-1) pour $S_{H_{0}'}$). Alors $H\in S_{H_0}$ (cf. \ref{propS0IPf}), $H\subset H'\in S$ et l'inclusion $H'/H\subset G/H$ est un monomorphisme strict de $\cI$ par définition de la structure indalgébrique de $G/H$. Par conséquent $H'\in\ <S_{H_0}>$, et finalement $<S_{H_{0}'}>\ \subset\ <S_{H_0}>$, donc $<S_{H_{0}'}>\ =\ <S_{H_0}>$. Par ailleurs, pour $H\in S_{H_0}$, le groupe indalgébrique $(G/H,S_{G/H})$ ne dépend pas de $H_0$ par définition. Comme $S_{H_0}$ est coinitial dans $<S_{H_0}>$, il résulte de (PI-3) que pour tout $H\in\ <S_{H_0}>$, le groupe indalgébrique $(G/H,S_{G/H})$ ne dépend pas de $H_0$.

Soit maintenant $f:(G,S)\ra (G',S')\in\cI\cP^{\ad}$. Choisissons une origine $H_0$ de $(G,S)$, une origine $H_{0}'$ de $(G',S')$ contenant $f(H_0)$, et montrons que $f:(G,<S_{H_0}>)\ra (G',<S_{H_{0}'}>)$ est proindalgébrique. Si $H'\in S_{H_{0}'}$, il s'agit de montrer que $f^{-1}(H')\in\ <S_{H_0}>$ et que l'application $G/f^{-1}(H')\ra G'/H'$ induite par $f$ est indalgébrique. Comme $(f|_{H_0}^{H_{0}'})^{-1}(H')\in S_{H_0}$ est une origine de $(G,S)$ et $H'$ une origine de $(G',S')$ contenant 
$f((f|_{H_0}^{H_{0}'})^{-1}(H'))$ (\ref{propS0IPf}), on peut supposer $H'=H_{0}'$. On a alors $H_0\subset f^{-1}(H_{0}')$ par construction, et pour tout $H\in S$ contenant $H_0$,
$$f^{-1}(H_{0}')/H_0\cap H/H_0=(f^{-1}(H_{0}')\cap H)/H_0=(f|_{H}^{f(H)})^{-1}(H_{0}'\cap f(H))/H_0.$$
Comme $f|_{H}^{f(H)}:H\ra f(H)$ est proalgébrique par définition des morphismes dans $\cI\cP$, et que $H_{0}'\cap f(H)$ est fermé dans $f(H)$ (\ref{IP2'} et (IP-2) pour $S'$), on obtient que $(f|_{H}^{f(H)})^{-1}(H_{0}'\cap f(H))$ est fermé dans $H$, et donc appartient à $S$ ((IP-2) pour $S$) ; son quotient par l'origine $H_0$ est donc un groupe quasi-algébrique. Ainsi donc, $f^{-1}(H_{0}')/H_0$ se trouve muni d'une structure de groupe indalgébrique faisant de l'inclusion $f^{-1}(H_{0}')/H_0\subset G/H_0$ un monomorphisme strict de $\cI$ (cf. \ref{descrseI}) ; d'où $f^{-1}(H_{0}')\in\ <S_{H_0}>$. Ceci étant, pour voir que $G/f^{-1}(H_{0}')\ra G'/H_{0}'$ est indalgébrique, il suffit de voir que $G/H_0\ra G'/H_{0}'$ l'est, puisque $G/H_0\ra G/f^{-1}(H_{0}')$ est un épimorphisme strict d'après (PI-3). Soit $H\in S$ contenant $H_0$. Choisissons $H_{1}'\in S$ contenant $H_{0}'$ et $f(H)$. Alors la composée $H\ra f(H) \ra H_{1}'\ra H_{1}'/H_{0}'$ est proalgébrique, si bien que $H/H_0\ra H_{1}'/H_{0}'$ est quasi-algébrique. Donc $G/H_0\ra G'/H_{0}'$ est bien indalgébrique.

On a ainsi défini un foncteur 
$$
\xymatrix{
\can:\cI\cP^{\ad} \ar[r] & \cP\cI^{\ad}.
}
$$
\end{Pt}

\begin{Pt}\label{foncPIIP}
Soit $(G,S)\in\cP\cI^{\ad}$. Si $H_0$ est un sous-groupe origine de $(G,S)$, alors l'ensemble complet de définition $S_{G/H_0}$ du groupe indalgébrique $G/H_0$ définit une structure $\cI\cP^{\ad}$ sur $G$, de la façon suivante. Soit $H\in S$ contenant $H_0$ tel que $H/H_0\in S_{G/H_0}$. Alors $H$ est une origine de $(G,S)$ d'après \ref{propS0PIf}. De plus, l'ensemble de sous-groupes de $H$
$$
S_H:=\{ H'\in S  \ |\ H'\subset H \}
$$
vérifie (P-1), (P-2), (P-3) et (P-4), munissant donc $H$ d'une structure de groupe proalgébrique. Ceci étant, l'ensemble
$$
\tilde{S}_{G/H_0}:=\{H\subset G\ |\ H_0\subset H\in S \textrm{ et } H/H_0\in S_{G/H_0}\}
$$
vérifie (IP-1), (IP-2)' et (IP-3), et par conséquent munit le groupe $G$ d'une structure de groupe indproalgébrique, notée 
$(G,<\tilde{S}_{G/H_0}>)$ (cf. \ref{completengIP}). Enfin, pour tout $H,H'\in \tilde{S}_{G/H_0}$, si $H\subset H'$, alors $H'/H\in\cQ$, si bien que $(G,<\tilde{S}_{G/H_0}>)\in\cI\cP^{\ad}$.

La structure indproalgébrique ainsi obtenue sur $G$ ne dépend pas du choix de l'origine $H_0$. Vérifions en effet que si $H_{0}'$ est une autre origine de $(G,S)$, alors $(G,<\tilde{S}_{G/H_{0}'}>)=(G,<\tilde{S}_{G/H_0}>)$. On peut supposer $H_{0}'\subset H_0$. Alors $<\tilde{S}_{G/H_0}>\ \subset\ <\tilde{S}_{G/H_{0}'}>$. Inversement, si $H'\in \tilde{S}_{G/H_{0}'}$, il existe $H\in S$ contenant $H_0+H'$ tel que $H/H_{0}'\in S_{G/H_{0}'}$ ((I-1) pour $S_{G/H_{0}'}$). Alors $H_0\subset H$, $H/H_0\in S_{G/H_0}$ (cf. \ref{propS0PIf}), $S\ni H'\subset H$ et cette inclusion est une immersion fermée de $\cP$ par définition de la structure proalgébrique de $H$. Par conséquent $H'\in\ <\tilde{S}_{G/H_0}>$, et finalement $<\tilde{S}_{G/H_{0}'}>\ \subset\ <\tilde{S}_{G/H_0}>$, donc $<\tilde{S}_{G/H_{0}'}>\ =\ <\tilde{S}_{G/H_0}>$. Par ailleurs, pour tout $H\in \tilde{S}_{G/H_0}$, le groupe proalgébrique 
$(H,S_H)$ ne dépend pas de $H_0$ par définition. Comme $\tilde{S}_{G/H_0}$ est cofinal dans $<\tilde{S}_{G/H_0}>$, il résulte de (IP-2)' que pour tout $H\in\ <\tilde{S}_{G/H_0}>$, le groupe proalgébrique $(H,S_H)$ ne dépend pas de $H_0$.

Soit maintenant $f:(G,S)\ra (G',S')\in\cP\cI^{\ad}$. Choisissons une origine $H_{0}'$ de $(G',S')$, une origine $H_{0}$ de $(G,S)$ incluse dans $f^{-1}(H_{0}')$, et montrons que $f:(G,<\tilde{S}_{G/H_0}>)\ra (G',<\tilde{S}_{G'/H_{0}'}>)$ est indproalgébrique. Si $H\in \tilde{S}_{G/H_0}$, il s'agit de montrer que $f(H)\in\ <\tilde{S}_{G'/H_{0}'}>$ et que l'application $H\ra f(H)$ induite par $f$ est proalgébrique. Comme $f(H)+H_{0}'$ est une origine de $(G',S')$ et $H$ une origine de $(G,S)$ incluse dans $f^{-1}(f(H)+H_{0}')$ (\ref{propS0PIf}), on peut supposer $H=H_0$. On a alors $f(H_0)\subset H_{0}'$ par construction, et pour tout $H'\in S'$ inclus dans $H_{0}'$,
$$(f|_{H_0}^{H_{0}'})^{-1}(H')=f^{-1}(H')\cap H_0\in S$$
puisque $S$ est stable par intersection (\ref{PI2'}). De plus, $H_0/(f^{-1}(H')\cap H_0)\ra H_{0}'/H'$ est quasi-algébrique : le morphisme indalgébrique $G/(f^{-1}(H')\cap H_0)\ra G/f^{-1}(H')$ induit une bijection $H_0/(f^{-1}(H')\cap H_0)\xrightarrow{\sim} (H_0+f^{-1}(H'))/f^{-1}(H')$, munissant le groupe $(H_0+f^{-1}(H'))/f^{-1}(H')$ d'une structure de sous-groupe de définition de $G/f^{-1}(H')$, de sorte que $H_0/(f^{-1}(H')\cap H_0)\ra H_{0}'/H'$ s'identifie à la composée quasi-algébrique $(H_0+f^{-1}(H'))/f^{-1}(H')\ra (f(H_0)+H')/H' \ra H_{0}'/H'$. L'application $f:H_0\ra H_{0}'$ est donc proalgébrique, si bien que $f(H_0)$ est fermé dans $H_{0}'$, et par conséquent appartient à $<\tilde{S}_{G'/H_{0}'}>$, et le morphisme $H_0\ra f(H_0)$ qui s'en déduit est proalgébrique.

On a ainsi défini un foncteur 
$$
\xymatrix{
\can:\cP\cI^{\ad} \ar[r] & \cI\cP^{\ad}.
}
$$
\end{Pt}

\begin{Th}\label{canequiv}
Les foncteurs $\can:\cI\cP^{\ad}\ra\cP\cI^{\ad}$ \ref{foncIPPI} et $\can:\cP\cI^{\ad}\ra\cI\cP^{\ad}$ \ref{foncPIIP} sont des équivalences de catégories inverses l'une de l'autre. Ils induisent des diagrammes commutatifs 
$$
\xymatrix{
\cI\cP^{\ad}\ar[r]^{\can} & \cP\cI^{\ad}  & & \cI\cP^{\ad}\ar[r]^{\can} & \cP\cI^{\ad} \\
\cI \ar[u]^{\ref{IetPdansIPf}} \ar@{=}[r] & \cI \ar[u]_{\ref{IetPdansPIf}} & & \cP \ar[u]^{\ref{IetPdansIPf}} \ar@{=}[r] & \cP \ar[u]_{\ref{IetPdansPIf}}
}
$$
$$
\xymatrix{
\cP\cI^{\ad}\ar[r]^{\can} & \cI\cP^{\ad}  & & \cP\cI^{\ad}\ar[r]^{\can} & \cI\cP^{\ad} \\
\cI \ar[u]^{\ref{IetPdansPIf}} \ar@{=}[r] & \cI \ar[u]_{\ref{IetPdansIPf}} & & \cP \ar[u]^{\ref{IetPdansPIf}} \ar@{=}[r] & \cP. \ar[u]_{\ref{IetPdansIPf}}
}
$$
\end{Th}

\begin{proof}
Soit $(G,S)\in\cI\cP^{\ad}$. Choisissons une origine $H_0$ de $(G,S)$, de sorte que $\can(G,S)=(G,<S_{H_0}>)\in\cP\cI^{\ad}$. Alors $H_0$ est une origine de $\can(G,S)$ : si $H_{-1}\in\ <S_{H_0}>$ est inclus dans $H_0$, il existe $H_{-2}\in S_{H_0}$ inclus dans $H_{-1}$ tel que $H_{-1}/H_{-2}\subset G/H_{-2}$  est un monomorphisme strict de $\cI$, et $H_0/H_{-1}\in\cI$ est alors quotient strict de $H_0/H_{-2}\in\cQ$, donc est dans $\cQ$. Par conséquent, $\can\circ\can(G,S)=(G,<\tilde{S}_{G/H_0}>)\in\cI\cP^{\ad}$. Montrons que $(G,<\tilde{S}_{G/H_0}>)=(G,S)$.

Par définition,
$$
\tilde{S}_{G/H_0}=\{H\subset G\ |\ H_0\subset H\in\ <S_{H_0}> \textrm{ et } H/H_0\in S_{G/H_0}\},
$$
où $S_{G/H_0}$ est la structure de groupe indalgébrique sur $G/H_0$ associée à $(G,<S_{H_0}>)\in\cP\cI^{\ad}$,
$$
S_{G/H_0}=\{H/H_0\ |\ H_0\subset H\in S\}.
$$
En particulier, $H_0\in \tilde{S}_{G/H_0}\subset\ <\tilde{S}_{G/H_0}>$. Mais comme $<\tilde{S}_{G/H_0}>$ ne dépend pas du choix de l'origine $H_0$ de $(G,S)$, on en déduit que l'ensemble $S_0$ des origines de $(G,S)$ est contenu dans $<\tilde{S}_{G/H_0}>$, et donc que $S=\ <S_0>\ \subset\ <\tilde{S}_{G/H_0}>$. Inversement, on lit sur les définitions ci-dessus que $\tilde{S}_{G/H_0}\subset S$, d'ou $<\tilde{S}_{G/H_0}>\ \subset S$, et finalement $<\tilde{S}_{G/H_0}>\ =S$. Pour montrer que $(G,<\tilde{S}_{G/H_0}>)=(G,S)$, il reste à vérifier que pour tout $H\in S$, la structure proalgébrique sur $H$ vu comme élément de $<\tilde{S}_{G/H_0}>$ est la même que la structure proalgébrique sur $H$ vu dans $S$. Comme $S_0$ est cofinal dans $<\tilde{S}_{G/H_0}>\ =S$, il suffit d'après (IP-2) de le vérifier lorsque $H$ est une origine de $(G,S)$, et comme $(G,<\tilde{S}_{G/H_0}>)$ ne dépend pas du choix de $H_0$, il suffit de traiter le cas $H=H_0$. Par définition, la structure proalgébrique sur $H_0\in\ <\tilde{S}_{G/H_0}>$ est l'ensemble de sous-groupes
$$
\{H\in\ <S_{H_0}>\ |\ H\subset H_0\},
$$
le quotient $H_0/H$ pour $H\in\ <S_{H_0}>$ étant muni de la structure quasi-algébrique noyau de l'épimorphisme indalgébrique strict
$$
\xymatrix{
0\ar[r] & H_0/H \ar[r] & G/H \ar[r] & G/H_0 \ar[r] & 0.
}
$$
Soit $H\in\ <S_{H_0}>$ inclus dans $H_0$. Par définition de $<S_{H_0}>$, il existe $H_{-1}\in S_{H_0}$ tel que $H/H_{-1}$ est un sous-groupe indalgébrique strict de $G/H_{-1}$. Pour $H=H_0$, il s'agit exactement du groupe quasi-algébrique $H_0/H_{-1}$ quotient du groupe proalgébrique $H_0\in S$, par définition de $(G/H_{-1},S_{G/H_{-1}})\in\cI$ associé à $(G,<S_{H_0}>)\ \in\cP\cI^{\ad}$. Considérons alors le diagramme commutatif exact de $\cI$
$$
\xymatrix{
0\ar[r] & H/H_{-1} \ar[r] \ar[d] & G/H_{-1} \ar[r] \ar@{=}[d] & G/H \ar[r] \ar[d] & 0 \\
0\ar[r] & H_0/H_{-1} \ar[r] & G/H_{-1} \ar[r] & G/H_0 \ar[r] & 0. 
}
$$
La première colonne s'insère dans une suite exacte de $\Ind(\cQ)$
$$
\xymatrix{
0\ar[r] & H/H_{-1} \ar[r] &  H_0/H_{-1} \ar[r] & H_0/H \ar[r] & 0,
}
$$
dont les trois termes sont dans $\cI$, i.e. une suite exacte de $\cI$. En particulier $H/H_{-1} \ra  H_0/H_{-1}$ est un monomorphisme indalgébrique strict, si bien que $H/H_{-1}$ est quasi-algébrique puisque $H_0/H_{-1}$ l'est, puis que $H$ appartient à l'ensemble complet de définition $S_{H_0}$ du groupe proalgébrique $H_0\in S$. De plus, la structure quasi-algébrique de $H_0/H$ est quotient de celle de $H_0/H_{-1}$, donc de la structure proalgébrique de $H_0\in S$. On a donc
l'égalité d'ensembles 
$$
\{H\in\ <S_{H_0}>\ |\ H\subset H_0\}=S_{H_0},
$$
et pour tout $H$ dans cet ensemble, l'égalité des structures quasi-algébriques associées sur le quotient $H_0/H$. Ceci achève de montrer que $\can\circ\can=\Id:\cI\cP^{\ad}\ra\cI\cP^{\ad}$ sur les objets, et c'est évident sur les morphismes.

Soit maintenant $(G,S)\in\cP\cI^{\ad}$. Choisissons une origine $H_0$ de $(G,S)$, de sorte que $\can(G,S)=(G,<\tilde{S}_{G/H_0}>)\in\cI\cP^{\ad}$. Alors $H_0$ est une origine de $\can(G,S)$ : si $H_{1}\in\ <\tilde{S}_{G/H_0}>$ contient $H_0$, il existe $H_{2}\in \tilde{S}_{G/H_0}$ avec une inclusion fermée $H_{1}\subset H_2$, et $H_1/H_0\in\cP$ est alors fermé dans $H_2/H_0\in\cQ$, donc est dans $\cQ$. Par conséquent, $\can\circ\can(G,S)=(G,<S_{H_0}>)\in\cP\cI^{\ad}$. Montrons que 
$(G,<S_{H_0}>)=(G,S)$.

Par définition,
$$
S_{H_0}=\{H\in\ <\tilde{S}_{G/H_0}>\ |\ H\subset H_0\},
$$
où
$$
\tilde{S}_{G/H_0}=\{H\subset G\ |\ H_0\subset H\in S \textrm{ et } H/H_0\in S_{G/H_0}\}.
$$
En particulier, $H_0\in S_{H_0}\subset\ <S_{H_0}>$. Mais comme $<S_{H_0}>$ ne dépend pas du choix de l'origine $H_0$ de $(G,S)$, on en déduit que l'ensemble $S_0$ des origines de $(G,S)$ est contenu dans $<S_{H_0}>$, et donc que $S=\ <S_0>\ \subset\ <S_{H_0}>$. Inversement, on lit sur les définitions ci-dessus que $S_{H_0}\subset S$, d'ou $<S_{H_0}>\ \subset S$, et finalement $<S_{H_0}>\ =S$. Pour montrer que $(G,<S_{H_0}>)=(G,S)$, il reste à vérifier que pour tout $H\in S$, la structure indalgébrique sur $G/H$ pour $H$ vu comme élément de $<S_{H_0}>$ est la même que pour $H$ vu dans $S$. Comme $S_0$ est cofinal dans $<S_{H_0}>\ =S$, il suffit d'après (PI-2) de le vérifier lorsque $H$ est une origine de $(G,S)$, et comme $(G,<S_{H_0}>)$ ne dépend pas du choix de $H_0$, il suffit de traiter le cas $H=H_0$. Par définition, la structure indalgébrique sur $G/H_0$ pour $H_0\in\ <S_{H_0}>$ est l'ensemble de sous-groupes
$$
\{H/H_0\ |\ H_0\subset H\in\ <\tilde{S}_{G/H_0}>\},
$$
le quotient $H/H_0$ pour $H_0\subset H\in\ <\tilde{S}_{G/H_0}>$ étant muni de la structure quasi-algébrique quotient de l'inclusion fermée de groupes proalgébriques $H_0\subset H$, où 
$$
S_H=\{H'\in S\ |\ H'\subset H\}.
$$
Soit $H\in\ <\tilde{S}_{G/H_0}>$ contenant $H_0$. Par définition de $<\tilde{S}_{G/H_0}>$, il existe $H_1\in \tilde{S}_{G/H_0}$ tel que $H$ est un sous-groupe proalgébrique fermé de $H_1$. Pour $H=H_0$, le groupe proalgébrique quotient $H_1/H_0$ est exactement le groupe quasi-algébrique $H_1/H_0\in S_{G/H_0}$ pour $H_0\in S$, par définition de 
$(H_1,S_{H_1})\in\cP$ associé à $(G,<\tilde{S}_{G/H_0}>)\in\cI\cP^{\ad}$. Mais alors $H/H_0$ est fermé dans $H_1/H_0$, donc appartient à l'ensemble complet de définition $S_{G/H_0}$ du groupe indalgébrique $G/H_0$ pour $H_0\in S$. De plus, la structure quasi-algébrique quotient $H/H_0$ est induite par celle de $H_1/H_0$, donc par celle de $G/H_0\in\cI$, $H_0\in S$. On a donc l'égalité d'ensembles 
$$
\{H/H_0\ |\ H_0\subset H\in\ <\tilde{S}_{G/H_0}>\}=S_{G/H_0},
$$ 
et pour tout $H/H_0$ dans cet ensemble, l'égalité des structures quasi-algébriques associées. Ceci achève de montrer que $\can\circ\can=\Id:\cP\cI^{\ad}\ra\cP\cI^{\ad}$ sur les objets, et c'est évident sur les morphismes.

Il est bien clair que les diagrammes de l'énoncé du théorème sont commutatifs.
\end{proof}

\begin{Rem}\label{canS0}
Soit $(G,S)\in\cI\cP^{\ad}$ (resp. $\cP\cI^{\ad}$). On a en particulier démontré qu'un sous-groupe de $G$ est une origine de $(G,S)$ si et seulement si c'est une origine de $\can(G,S)$.
\end{Rem}

\begin{Prop}\label{canexacts}
Le foncteur $\can:\cI\cP^{\ad}\ra\cP\cI^{\ad}$ \ref{foncIPPI} est exact.
\end{Prop}

Commençons par le lemme suivant.

\begin{Lem}\label{origexact}
Soit 
$$
\xymatrix{
0\ar[r] & (G',S') \ar[r] & (G,S) \ar[r]^{f} & (G'',S'') \ar[r] & 0
}
$$
une suite exacte de $\cI\cP^{\ad}$. Si $H_0$ une origine de $(G,S)$, alors $G'\cap H_0$ est une origine de $(G',S')$ et $f(H_0)$ est une origine de $(G'',S'')$. Réciproquement, si $H_{0}'$ est une origine de $(G',S')$, alors il existe une origine $H_0$ de $(G,S)$ telle que $H_{0}'=G'\cap H_0$, et si $H_{0}''$ est une origine de $(G'',S'')$, alors il existe une origine $H_0$ de $(G,S)$ telle que $H_{0}''=f(H_0)$.
\end{Lem}

\begin{proof}
Soit $H_0$ une origine de $(G,S)$, et soit $H'\in S'$ contenant $G'\cap H_{0}$. Choisissons $H_1\in S$ contenant $H'$ (vu dans $S$) et $H_0$. Alors $G'\cap H_1$ contient $H'$ et on a un diagramme commutatif exact de $\cP$
$$
\xymatrix{
0\ar[r] & G'\cap H_{0} \ar[r] \ar[d] & H_0 \ar[r] \ar[d] & f(H_0) \ar[r] \ar[d] & 0 \\
0\ar[r] & G'\cap H_1 \ar[r] & H_1 \ar[r] & f(H_1) \ar[r] & 0.
}
$$
D'après le lemme du serpent dans $\cP$, $G'\cap H_1/G'\cap H_0$ est un sous-groupe fermé du groupe quasi-algébrique $H_1/H_0$, donc est quasi-algébrique. \emph{A fortiori}, le sous-groupe fermé $H'/G'\cap H_0\subset G'\cap H_1/G'\cap H_0$ est quasi-algébrique, et $G'\cap H_0$ est donc une origine de $(G',S')$. Soit maintenant $H''\in S''$ contenant $f(H_0)$. Choisissons $H\in S$ tel que $f(H)=H''$ puis $H_1$ contenant $H$ et $H_0$. Alors $f(H_1)$ contient $H''$ et on a un diagramme commutatif exact de $\cP$ comme ci-dessus, et le lemme du serpent dans $\cP$ montre que $f(H_1)/f(H_0)$ est quotient du groupe quasi-algébrique $H_1/H_0$ par un sous-groupe fermé, donc est quasi-algébrique. \emph{A fortiori}, le sous-groupe fermé $H''/f(H_0)\subset f(H_1)/f(H_0)$ est quasi-algébrique, et $f(H_0)$ est donc une origine de $(G'',S'')$.

Réciproquement, soit $H_{0}'$ est une origine de $(G',S')$. Alors $\overline{G'}:=G'/H_{0}'\in\cI$ est un sous-objet strict de $\overline{G}:=G/H_{0}'\in\cI\cP^{\ad}$. Choisissons une origine $H_0$ de $(G,S)$ contenant $H_{0}'$ et notons $\overline{H_{0}}$ son image dans $\overline{G}$. Alors $\overline{H_{0}}$ est un sous-groupe de définition de $\overline{G}$
et $\overline{H_{0}}':=\overline{G'}\cap\overline{H_{0}}$ est un sous-groupe de définition de $\overline{G'}\in\cI$. Par conséquent $\overline{H_{0}}'$ est un sous-groupe quasi-algébrique fermé du groupe proalgébrique $\overline{H_{0}}$. Il existe alors un sous-groupe de définition $\overline{H_{-1}}$ de $\overline{H_0}$ tel que $\overline{H_{-1}}\cap \overline{H_{0}}'=0$ : en effet, l'intersection des sous-groupes de définition de $\overline{H_0}$ est nulle et l'ensemble des sous-groupes fermés de $\overline{H_0}'$ vérifie la condition d'Artin. Ceci étant, soit $H_{-1}\subset G$ l'image réciproque de $\overline{H_{-1}}\subset\overline{G}$. Alors $H_{-1}$ est un sous-groupe de définition de $H_0$, donc une origine de $(G,S)$ d'après \ref{propS0IPf}, et 
$$G'\cap H_{-1}=G'\cap H_0 \cap H_{-1}=H_{0}'.$$
Soit enfin $H_{0}''$ une origine de $(G'',S'')$. Choisissons $\tilde{H}_{0}''\in S$ tel que $f(\tilde{H}_{0}'')=H_{0}''$, puis une origine $H_0$ de $(G,S)$ contenant $\tilde{H}_{0}''$. Alors $f(H_0)\in S''$ contient $H_{0}''$, et le quotient $f(H_0)/H_{0}''$ est donc quasi-algébrique. Le noyau du morphisme proalgébrique composé $H_0\ra f(H_0)\ra f(H_0)/ H_{0}''$ est donc un sous-groupe de définition $H_{-1}$ de $H_0$, contenant $\tilde{H}_{0}''$. C'est donc une origine de $(G,S)$ d'après \ref{propS0IPf}, telle que $f(H_{-1})=H_{0}''$.
\end{proof}

\begin{proof}[Démonstration de la proposition \ref{canexacts}]
Soit 
$$
\xymatrix{
0\ar[r] & (G',S') \ar[r] & (G,S) \ar[r]^{f} & (G'',S'') \ar[r] & 0
}
$$
une suite exacte de $\cI\cP^{\ad}$. Choisissons une origine $H_0$ de $(G,S)$. Alors $H_{0}':=G'\cap H_0$ et $H_{0}'':=f(H_0)$ sont des origines de $(G',S')$ et $(G'',S'')$ respectivement (\ref{origexact}), et on a une suite exacte de $\cP$
$$
\xymatrix{
0\ar[r] & H_{0}' \ar[r] & H_0 \ar[r] & H_{0}'' \ar[r] & 0.
}
$$
Si $H\in S$ contient $H_0$, alors en posant $H':=G'\cap H$ et $H'':=f(H)$, on obtient une suite exacte de $\cP$
$$
\xymatrix{
0\ar[r] & H' \ar[r] & H \ar[r] & H'' \ar[r] & 0
}
$$
contenant la précédente. D'où une suite exacte de $\cQ$
$$
\xymatrix{
0\ar[r] & H'/H_{0}' \ar[r] & H/H_0 \ar[r] & H''/H_{0}'' \ar[r] & 0.
}
$$
D'après (IP-1) pour S, l'ensemble $\{G'\cap H\ |\ H_0\subset H\in S \}$ est cofinal dans $\{H'\in S'\ |\ H_{0}'\subset H'\}$ et $\{f(H)\ |\ H_0\subset H\in S \}$ est cofinal dans $\{H''\in S''\ |\ H_{0}''\subset H''\}$. Faisant varier $H$, on obtient donc une suite exacte de $\cI$
$$
\xymatrix{
0\ar[r] & G'/H_{0}' \ar[r] & G/H_0 \ar[r] & G''/H_{0}'' \ar[r] & 0.
}
$$
Maintenant, d'après \ref{origexact}, l'ensemble $\{G'\cap H_0\ |\ H_0\textrm{ origine de } (G,S) \}$ est égal à l'ensemble des origines de $(G',S')$ et $\{f(H_0)\ |\ H_0\textrm{ origine de } (G,S) \}$ est égal à l'ensemble des origines de $(G'',S'')$. Or les origines d'un objet de $\cI\cP^{\ad}$ et de son image par $\can:\cI\cP^{\ad}\ra\cP\cI^{\ad}$ sont les mêmes (\ref{canS0}), et l'ensemble des origines d'un objet de $\cP\cI^{\ad}$ est coinitial dans l'ensemble complet de définition de celui-ci (\ref{defequivPIf}). Faisant varier $H_0$, on obtient donc que la suite de $\cP\cI^{\ad}$
$$
\xymatrix{
0\ar[r] & \can(G',S') \ar[r] & \can(G,S) \ar[r]^{f} & \can(G'',S'') \ar[r] & 0
}
$$
est exacte.
\end{proof}

\begin{Def}\label{deffQf}
Un \emph{groupe $\cQ$-admissible} est un groupe ordinaire $G$ muni d'une famille non vide $S$ de sous-groupes et, pour tous $H,H'\in S$ avec $H\subset H'$, d'une structure de groupe quasi-algébrique sur $H'/H$, ces données vérifiant les axiomes suivants :
\begin{itemize}

\medskip

\item ($\cQ$-ad-1) $H,H'\in S \Rightarrow \exists H^+\in S : H+H'\subset H^+$ et $\exists H_{-}\in S : H_{-}\subset H\cap H'$.

\medskip

\item ($\cQ$-ad-2) Si $H,H'\in S$ avec $H\subset H'$, alors les images réciproques des sous-groupes fermés de $H'/H$  appartiennent à $S$.

\medskip

\item ($\cQ$-ad-3) Si $H_{-},H,H^+\in S$ avec $H_{-}\subset H\subset H^+$, alors les applications canoniques $H/H_{-}\ra H^+/H_{-}$ et $H^{+}/H_{-}\ra H^{+}/H$ sont des morphismes de groupes quasi-algébriques.

\medskip

\item ($\cQ$-ad-4) Pour tout $H\in S$, les applications canoniques $\li_{H\subset H^+\in S}H^+/H\ra G/H$ et  $H\ra \lp_{S\ni H_{-}\subset H}H/H_{-}$ sont des isomorphismes de groupes.

\medskip

\end{itemize}

L'ensemble $S$ est l'\emph{ensemble complet des origines} du groupe $\cQ$-admissible $(G,S)$, et les sous-groupes 
$H\in S$ sont les \emph{origines} de $(G,S)$.

\medskip

Un \emph{morphisme de groupes $\cQ$-admissibles} $(G_1,S_1)\ra (G_2,S_2)$ est un morphisme de groupes $f:G_1\ra G_2$ vérifiant les conditions suivantes : 
\begin{itemize}

\medskip

\item (Hom-$\cQ$-ad-1) Pour tout $H_1\in S_1$, il existe $H_2\in S_2$ contenant $f(H_1)$ tel que pour tout $H_{2-}\in S_2$ inclus dans $H_2$, $f^{-1}(H_{2-})\cap H_1\in S_1$ et l'application $H_1/f^{-1}(H_{2-})\cap H_1\ra H_2/H_{2-}$ induite par $f$ est quasi-algébrique.

\medskip

\item (Hom-$\cQ$-ad-2) Pour tout $H_2\in S_2$, il existe $H_1\in S_1$ inclus dans $f^{-1}(H_2)$ tel que pour tout $H_{1+}\in S_1$ contenant $H_1$, $f(H_{1+})+H_2\in S_2$ et l'application $H_{1+}/H_1\ra f(H_{1+})+H_2/H_2$ induite par $f$ est quasi-algébrique.
\end{itemize}
\end{Def}

\begin{Rem}  
Soit $G$ un groupe, muni d'une famille $S'$ de sous-groupes vérifiant ($\cQ$-ad-1), ($\cQ$-ad-3) et ($\cQ$-ad-4). Alors il existe une plus petite famille de sous-groupes de $G$ contenant $S'$ et vérifiant ($\cQ$-ad-1), ($\cQ$-ad-2), ($\cQ$-ad-3) et ($\cQ$-ad-4). On l'appellera l'ensemble complet des origines de $G$ \emph{engendré par $S'$}, et on le notera $<S'>$. On dira aussi que $S'$ est un \emph{ensemble de définition} du groupe $\cQ$-admissible $(G,<S'>)$. Explicitement, $<S'>$ est l'ensemble des sous-groupes $H\subset G$ tels qu'il existe $H_{-}',H_{+}'\in S'$ vérifiant : 

\medskip

\begin{itemize}
\item $H_{-}'\subset H\subset H_{+}'$ ;
\item l'inclusion $H/H_{-}'\subset H_{+}'/H_{-}'$ est une immersion fermée quasi-algébrique.
\end{itemize}

\medskip

En particulier, $S'$ est coinitial et cofinal dans $S$. Par ailleurs, $<S'>$ est l'unique ensemble complet d'origines de $G$ contenant $S'$.
\end{Rem}

\begin{Not}
On note $\fQf$ la catégorie des groupes $\cQ$-admissibles.
\end{Not}

\begin{Pt}\label{fQfIPf}
Soit $(G,S)\in\fQf$. Si $H\in S$, l'ensemble de sous-groupes de $H$
$$
S_H:=\{H_-\in S \ |\ H_-\subset H\}
$$
vérifie (P-1), (P-2), (P-3) et (P-4), et munit donc $H$ d'une structure de groupe proalgébrique $(H,S_H)$. L'ensemble $S$ de sous-groupes de $G$ vérifie alors (IP-1), (IP-2)' et (IP-3), et pour tout $H,H'\in S$ avec $H\subset H'$, le groupe proalgébrique $H'/H$ est quasi-algébrique par hypothèse. D'où un objet $(G,<S>)\in\cI\cP^{\ad}$. De plus, si $f:(G_1,S_1)\ra (G_2,S_2)\in\fQf$, alors d'après (Hom-$\cQ$-ad-1), pour tout $H_1\in S_1$, il existe $H_2\in S_2$ contenant $f(H_1)$ tel que l'application $H_1\ra H_2$ induite par $f$ est proalgébrique : 
$$
H_1\ra\lp_{S_2\ni H_{2-}\subset H_2}H_1/f^{-1}(H_{2-})\cap H_1 \ra \lp_{S_2\ni H_{2-}\subset H_2}H_2/H_{2-}=H_2 ;
$$
d'où $f:(G_1,<S_1>)\ra (G_2,<S_2>)\in\cI\cP^{\ad}$, et finalement un foncteur 
$$
\xymatrix{
\indpro:\fQf\ar[r] & \cI\cP^{\ad}.
}
$$

Réciproquement, soit $(G,S)\in\cI\cP^{\ad}$. Alors pour toutes origines $H_0\subset H_{0}'\in S$, on a $H_{0}'/H_0\in\cQ$, et l'ensemble $S_0$ des origines de $(G,S)$ vérifie ($\cQ$-ad-1), ($\cQ$-ad-2), ($\cQ$-ad-3) et ($\cQ$-ad-4) (cf. \ref{propS0IPf}), donc définit un objet $\ad(G,S):=(G,S_0)$ de $\fQf$. De plus, si $f:(G_1,S_1)\ra (G_2,S_2)\in\cI\cP^{\ad}$, alors pour tout $H_{1,0}\in S_{1,0}$, il existe $H_{2,0}\in S_{2,0}$ contenant $f(H_{1,0})$, et alors pour tout $H_{2,0-}\in S_{2,0}$ inclus dans $H_{2,0}$, $f^{-1}(H_{2,0-})\cap H_{1,0}=(f|_{H_{1,0}}^{H_{2,0}})^{-1}(H_{2,0-})\in S_{H_{1,0}}\subset S_{1,0}$ (cf. \ref{propS0IPf}), et l'application $H_{1,0}/f^{-1}(H_{2,0-})\cap H_{1,0}\ra H_{2,0}/H_{2,0-}$ induite par $f$ est quasi-algébrique. 
\end{Pt}

\begin{Pt}\label{fQfPIf}
Soit $(G,S)\in\fQf$. Si $H\in S$, l'ensemble de sous-groupes de $G/H$
$$
S_{G/H}:=\{H_+/H \ |\ H\subset H^+\in S\}
$$
vérifie (I-1), (I-2) et (I-3), et munit donc $G/H$ d'une structure de groupe indalgébrique $(G/H,S_{G/H})$. L'ensemble $S$ de sous-groupes de $G$ vérifie alors (PI-1), (PI-3) et (PI-4), et pour $H,H'\in S$ avec $H\subset H'$, le groupe indalgébrique $H'/H$ est quasi-algébrique par hypothèse. D'où un objet $(G,<S>)\in\cP\cI^{\ad}$. De plus, si $f:(G_1,S_1)\ra (G_2,S_2)\in\fQf$, alors d'après (Hom-$\cQ$-ad-2), pour tout $H_2\in S_2$, il existe $H_1\in S_1$ inclus dans $f^{-1}(H_2)$ tel que l'application $G/H_1\ra G/H_2$ induite par $f$ est indalgébrique :
$$
G/H_1=\li_{H_1\subset H_{1+}\in S_1} H_{1+}/H_1 \ra\li_{H_1\subset H_{1+}\in S_1} f(H_{1+})+H_2/H_2\ra G/H_2 ;
$$
d'où $f:(G_1,<S_1>)\ra (G_2,<S_2>)\in\cP\cI^{\ad}$, et finalement un foncteur 
$$
\xymatrix{
\proind:\fQf\ar[r] & \cP\cI^{\ad}.
}
$$

Réciproquement, soit $(G,S)\in\cP\cI^{\ad}$. Alors pour toutes origines $H_0\subset H_{0}'\in S$, on a $H_{0}'/H_0\in\cQ$, et l'ensemble $S_0$ des origines de $(G,S)$ vérifie ($\cQ$-ad-1), ($\cQ$-ad-2), ($\cQ$-ad-3) et ($\cQ$-ad-4) (cf. \ref{propS0PIf}), donc définit un objet $\ad(G,S):=(G,S_0)$ de $\fQf$. De plus, si $f:(G_1,S_1)\ra (G_2,S_2)\in\cP\cI^{\ad}$, alors pour tout $H_{2,0}\in S_{2,0}$, il existe $H_{1,0}\in S_{1,0}$ inclus dans $f^{-1}(H_{2,0})$, et alors pour tout $H_{1,0+}\in S_{1,0}$ contenant $H_{1,0}$, $f(H_{1,0+})+H_{2,0}\in S_{2,0}$ car image réciproque de l'image de $H_{1,0+}/H_{1,0}$ par le morphisme indalgébrique $G_1/H_{1,0}\ra G_2/H_{2,0}$ induit par $f$ (cf. \ref{propS0PIf}), et l'application $H_{1,0+}/H_{1,0}\ra f(H_{1,0+})+H_{2,0}/H_{2,0}$ qui s'en déduit est quasi-algébrique.
\end{Pt}

\begin{Pt}\label{finfoncfQf}
Soit $(G,S)\in\cI\cP^{\ad}$. On a vu en \ref{fQfIPf} que $(G,S)$ définissait un objet $\ad(G,S)\in\fQf$. Par ailleurs, $\can(G,S)\in\cP\cI^{\ad}$, donc définit un objet $\ad(\can(G,S))\in\fQf$ par \ref{fQfPIf}. D'après \ref{canS0} et la définition de 
$\can(G,S)$, on a $\ad(G,S)=\ad(\can(G,S))$. Soit maintenant $f:(G_1,S_1)\ra (G_2,S_2)\in\cI\cP^{\ad}$. On a vu en \ref{fQfIPf} que $f:\ad(G_1,S_1)\ra \ad(G_2,S_2)$ vérifiait l'axiome (Hom-$\cQ$-ad-1). De plus $\can(f):\can(G_1,S_1)\ra \can(G_2,S_2)\in\cP\cI^{\ad}$ est tel $\can(f):\ad(\can(G_1,S_1))\ra \ad(\can(G_2,S_2))$ vérifie (Hom-$\cQ$-ad-2) d'après \ref{fQfPIf}. D'après ce qui précède, on en tire que $f:\ad(G_1,S_1)\ra \ad(G_2,S_2)\in\fQf$. D'où un foncteur 
$$
\xymatrix{
\ad:\cI\cP^{\ad}\ar[r] & \fQf.
}
$$ 
On obtient de même un foncteur 
$$
\xymatrix{
\ad:\cP\cI^{\ad}\ar[r] & \fQf.
}
$$
\end{Pt}

\begin{Rem}
On notera que si $f:(G_1,S_1)\ra (G_2,S_2)\in\cI\cP^{\ad}$ et $H_{1,0}$ est une origine de $(G_1,S_1)$, alors $f(H_{1,0})$ n'est pas une origine de $(G_2,S_2)$ en général. On peut par exemple prendre pour $f$ l'endomorphisme nul d'un objet de $\cI\cP^{\ad}$ qui n'est pas dans $\cI$. De même, si $f:(G_1,S_1)\ra (G_2,S_2)\in\cP\cI^{\ad}$ et $H_{2,0}$ est une origine de 
$(G_2,S_2)$, alors $f^{-1}(H_{2,0})$ n'est pas une origine de $(G_1,S_1)$ en général. On peut par exemple prendre pour $f$ l'endomorphisme nul d'un objet de $\cP\cI^{\ad}$ qui n'est pas dans $\cP$. 
\end{Rem}

\begin{Th}\label{fQfIPfPIf}
Les foncteurs $\indpro:\fQf\ra\cI\cP^{\ad}$ \ref{fQfIPf} et $\ad:\cI\cP^{\ad}\ra\fQf$ \ref{finfoncfQf} sont des équivalences de catégories inverses l'une de l'autre. Les foncteurs $\proind:\fQf\ra\cP\cI^{\ad}$ \ref{fQfPIf} et $\ad:\cP\cI^{\ad}\ra\fQf$ \ref{finfoncfQf} sont des équivalences de catégories inverses l'une de l'autre. Le diagramme
$$
\xymatrix{
& & \fQf\ar@<-1ex>[ddll]_{\indpro} \ar@<1ex>[ddrr]^{\proind} & & \\
\\
\cI\cP^{\ad} \ar@<1ex>[rrrr]^{\can} \ar@<-1ex>[uurr]_{\ad} & & & & \cP\cI^{\ad}  \ar@<1ex>[llll]^{\can} \ar@<1ex>[uull]^{\ad}
}
$$
est commutatif. Les axiomes (Hom-$\cQ$-add-1) et (Hom-$\cQ$-add-2) sont équivalents.
\end{Th}

\begin{proof}
Soit $(G,S)\in\fQf$. On a $\ad\circ\indpro(G,S)=(G,<S>_0)\in\fQf$. Montrons que $(G,<S>_0)=(G,S)$. Soit $H_0\in S$. Si $H\in\ <S>$ contient $H_0$, alors il existe $H_1\in S$ contenant $H$, et le groupe proalgébrique $H/H_0$ est alors fermé dans $H_1/H_0\in\cQ$, donc est quasi-algébrique. Par conséquent $S\subset\ <S>_0$. Inversement, soit $H\in\ <S>_0$. Il existe $H_1\in S$ contenant $H$, et le quotient $H_1/H$ est alors quasi-algébrique, de sorte que $H$ est un sous-groupe de définition du groupe proalgébrique $H_1$. Par définition de celui-ci, on a donc $H\in S$. Ainsi $<S>_0\ =\ S$. De plus, si $H_-,H\in S$ avec $H_-\subset H$, alors le groupe proalgébrique quotient $H/H_-$ est le groupe quasi-algébrique $H/H_-$ associé à $(G,S)$ par définition de la structure proalgébrique de $H$ associée à $\indpro(G,S)$. D'ou $(G,<S>_0)=(G,S)$, ce qui montre que $\ad\circ\indpro=\Id:\fQf\ra\fQf$ sur les objets, et c'est évident sur les morphismes.

Soit $(G,S)\in\cI\cP^{\ad}$. On a $\indpro\circ\ad(G,S)=(G,<S_0>)\in\fQf$. Montrons que $(G,<S_0>)=(G,S)$. On a $<S_0>=S$ d'après \ref{defequivIPf}. Soit alors $H\in S$ et montrons que les structures proalgébriques sur $H$ associées à $(G,<S_0>)$ et à $(G,S)$ sont identiques. Comme $S_0$ est cofinal dans $<S_0>\ =S$, il suffit d'après (IP-2)' de traiter le cas où $H\in S_0$. Mais alors 
$$
\{H_-\in S_0\ |\ H_-\subset H\}
$$
est égal à l'ensemble complet de définition de $H\in S$ (cf. \ref{propS0IPf}), et pour $H_-$ dans cet ensemble, le groupe quasi-algébrique $H/H_-$ est par définition de $\ad$ le quotient du groupe proalgébrique $H\in S$ par le sous-groupe fermé $H_0$. D'ou $(G,<S_0>)=(G,S)$, ce qui montre que $\indpro\circ\ad=\Id:\cI\cP^{\ad}\ra\cI\cP^{\ad}$ sur les objets, et c'est évident sur les morphismes.

Soit $(G,S)\in\fQf$. On a $\ad\circ\proind(G,S)=(G,<S>_0)\in\fQf$. Montrons que $(G,<S>_0)=(G,S)$. Soit $H_0\in S$. Si $H\in\ <S>$ est inclus dans $H_0$, alors il existe $H_{-1}\in S$ inclus dans $H$, et le groupe indalgébrique $H_0/H$ est alors quotient strict de $H_0/H_{-1}\in\cQ$, donc est quasi-algébrique. Par conséquent $S\subset\ <S>_0$. Inversement, soit $H\in\ <S>_0$. Il existe $H_{-1}\in S$ inclus dans $H$, et le quotient $H/H_{-1}$ est alors quasi-algébrique, donc est un sous-groupe de définition du groupe indalgébrique $G/H_{-1}$. Par définition de celui-ci, on a donc $H\in S$. Ainsi $<S>_0\ =\ S$. De plus, si $H,H_+\in S$ avec $H\subset H_+$, alors le sous-groupe indalgébrique strict $H_+/H\subset G/H$ est le groupe quasi-algébrique $H_+/H$ associé à $(G,S)$ par définition de la structure indalgébrique de $G/H$ associée à $\proind(G,S)$. D'ou $(G,<S>_0)=(G,S)$, ce qui montre que $\ad\circ\proind=\Id:\fQf\ra\fQf$ sur les objets, et c'est évident sur les morphismes.

Soit $(G,S)\in\cP\cI^{\ad}$. On a $\proind\circ\ad(G,S)=(G,<S_0>)\in\fQf$. Montrons que $(G,<S_0>)=(G,S)$. On a $<S_0>=S$ d'après \ref{defequivPIf}. Soit alors $H\in S$ et montrons que les structures indalgébriques sur $G/H$ associées à $(G,<S_0>)$ et à $(G,S)$ sont identiques. Comme $S_0$ est coinitial dans $<S_0>\ =S$, il suffit d'après (PI-3) de traiter le cas où $H\in S_0$. Mais alors 
$$
\{H_+/H \ |\ H\subset H_+\in S_0\}
$$
est égal à l'ensemble complet de définition de $G/H$ pour $H\in S$ (cf. \ref{propS0PIf}), et pour $H\subset H_+\in S_0$, le groupe quasi-algébrique $H_+/H$ est par définition de $\ad$ le sous-groupe indalgébrique strict $H_+/H\subset G/H$, $H\in S$. D'ou $(G,<S_0>)=(G,S)$, ce qui montre que $\proind\circ\ad=\Id:\cP\cI^{\ad}\ra\cP\cI^{\ad}$ sur les objets, et c'est évident sur les morphismes.

On a vu en \ref{finfoncfQf} que $\ad\circ\can=\ad:\cI\cP^{\ad}\ra\cQ$ sur les objets, et c'est évident sur les morphismes. Autrement dit, le triangle intérieur du diagramme de l'énoncé du théorème est commutatif. Il en résulte que tout le diagramme est commutatif.

Enfin, si $(G_1,S_1),(G_2,S_2)\in\fQf$ et $f:G_1\ra G_2$ est un morphisme de groupes, alors on a vu en \ref{fQfIPf} que $f$ était un morphisme $\indpro(G_1,S_1)\ra\indpro(G_2,S_2)\in\cI\cP^{\ad}$ dès que $f$ vérifiait (Hom-$\cQ$-ad-1). Puis on a vu en \ref{finfoncfQf} que $f$ vérifiait alors (Hom-$\cQ$-ad-2). De même (Hom-$\cQ$-ad-2) implique (Hom-$\cQ$-ad-1).
\end{proof}

\begin{Pt}\label{fQfexacte}
On munit $\fQf$ de la structure de catégorie exacte induite par l'équivalence 
$$
\xymatrix{
\indpro:\fQf\ar[r]^<<<<<{\sim} & \cI\cP^{\ad}.
}
$$
D'après \ref{canexacts}, le plongement $\fQf\xrightarrow{\proind}\cP\cI^{\ad}\ra\cP\cI$ est alors exact.
\end{Pt}

\begin{Pt}
Par définition, les foncteurs $\can$ \ref{foncIPPI} et \ref{foncPIIP} induisent l'identité au niveau des groupes ordinaires sous-jacents. Comme de plus ils commutent à toute extension de corps algébriquement close $\kappa/k$ (\ref{cbQ}), on obtient en particulier que pour tout $(G,S)\in\fQf$, les fibres génériques $\indpro(G,S)_{\eta}$ et $\proind(G,S)_{\eta}$ sont égales. 
\end{Pt}

\begin{Def} \label{deffibgenad}
Avec les notations de \ref{deffibgenIP} et \ref{deffibgenPI}, la \emph{fibre générique} d'un groupe admissible $(G,S)$ est le foncteur
$$
(G,S)_{\eta}:=\indpro(G,S)_{\eta}=\proind(G,S)_{\eta}:\cT^{\circ} \ra (\Groupes).
$$
\end{Def}

\begin{Pt}\label{resumeplong}
En résumé : 
$$
\xymatrix{
                              & & \cI\cP^{\ad}  \ar@{^{(}->}[r] & \cI\cP \ar@{^{(}->}[r] & \Ind(\cP) \ar@{=}[r] & \Ind(\Pro(\cQ)) \\
\cQ \ar@{^{(}->}[r] & \cI,\cP \ar@{^{(}->}[dr] \ar@{^{(}->}[ur]   & \fQf \ar@{<->}[d] \ar@{<->}[u] \\
                              & & \cP\cI^{\ad}  \ar@{^{(}->}[r] &\cP\cI \ar@{^{(}->}[r] &  \Pro(\cI) \ar@{^{(}->}[r] & \Pro(\Ind(\cQ)).
}
$$
\end{Pt}

\section{Variantes unipotentes}

\begin{Not}
On note $\cU$ la sous-catégorie pleine de $\cQ$ dont les objets sont les groupes quasi-algébriques unipotents (sur $k$) : \cite{S2} 8.1. C'est une catégorie abélienne, le plongement $\cU\ra\cQ$ est exact et fait de $\cU$ une sous-catégorie épaisse de $\cQ$.
\end{Not}

\begin{Pt}\label{resumeplongU}
Les considérations des sections \ref{groupesproalgébriques} à \ref{groupesadmissibles} peuvent être menées à partir de la catégorie $\cU$ à la place de la catégorie $\cQ$. Les objets obtenus sont alors qualifiés d'\emph{unipotents}. En particulier, on définit des catégories exactes et des plongements exacts donnant lieu à un diagramme commutatif 
$$
\xymatrix{
                              & & \cI\cP^{\ad}(\cU)  \ar@{^{(}->}[r] & \cI\cP(\cU) \ar@{^{(}->}[r] & \Ind(\cP(\cU)) \ar@{=}[r] & \Ind(\Pro(\cU)) \\
\cU \ar@{^{(}->}[r] & \cI(\cU),\cP(\cU) \ar@{^{(}->}[dr] \ar@{^{(}->}[ur]   & \fUf \ar@{<->}[d] \ar@{<->}[u] \\
                              & & \cP\cI^{\ad}(\cU)  \ar@{^{(}->}[r] &\cP\cI(\cU) \ar@{^{(}->}[r] &  \Pro(\cI(\cU)) \ar@{^{(}->}[r] & \Pro(\Ind(\cU)),
}
$$
où $\cP(\cU)\xrightarrow{\sim}\Pro(\cU)$ est abélienne. De plus, si $\cC$ est l'une des catégories quasi-algébriques, alors le fait qu'un objet soit dans $\cC(\cU)$ peut se tester sur un ensemble de définition.
\end{Pt}

\begin{Not}
Soit $n\in\bbN$. On note $\cU(n)$ la sous-catégorie pleine de $\cU$ dont les objets sont les groupes quasi-algébriques unipotents annulés par $p^n$. C'est une catégorie abélienne, le plongement $\cU(n)\ra\cU$ est exact, mais $\cU(n)$ n'est pas épaisse dans $\cU$.
\end{Not}

\begin{Pt}\label{resumeplongU(n)}
On a encore des catégories exactes et des plongements exacts donnant lieu à un diagramme commutatif 
$$
\xymatrix{
                              & & \cI\cP^{\ad}(\cU(n))  \ar@{^{(}->}[r] & \cI\cP(\cU(n)) \ar@{^{(}->}[r] & \Ind(\cP(\cU(n))) \ar@{=}[r] & \Ind(\Pro(\cU(n))) \\
\cU(n) \ar@{^{(}->}[r] & \cI(\cU(n)),\cP(\cU(n)) \ar@{^{(}->}[dr] \ar@{^{(}->}[ur]   & \fUnf \ar@{<->}[d] \ar@{<->}[u] \\
                              & & \cP\cI^{\ad}(\cU(n))  \ar@{^{(}->}[r] &\cP\cI(\cU(n)) \ar@{^{(}->}[r] &  \Pro(\cI(\cU(n))) \ar@{^{(}->}[r] & \Pro(\Ind(\cU(n))),
}
$$
où $\cP(\cU(n))\xrightarrow{\sim}\Pro(\cU(n))$ est abélienne. De plus, si $\cC$ est l'une des catégories quasi-algébriques, alors le fait qu'un objet soit dans $\cC(\cU(n))$ peut se tester sur un ensemble de définition.
\end{Pt}

\begin{Not}
On note $\cV$ la catégorie abélienne des $k$-espaces vectoriels de dimension finie. Le foncteur $\cV\ra\cU(1)$ qui a un $k$-espace vectoriel associe le groupe quasi-algébrique vectoriel de sections $V$ est fidèle et exact, mais pas plein.
\end{Not}

\begin{Pt}\label{resumeplongV}
On a encore des catégories exactes et des plongements exacts donnant lieu à un diagramme commutatif 
$$
\xymatrix{
                              & & \cI\cP^{\ad}(\cV)  \ar@{^{(}->}[r] & \cI\cP(\cV) \ar@{^{(}->}[r] & \Ind(\cP(\cV)) \ar@{=}[r] & \Ind(\Pro(\cV)) \\
\cV \ar@{^{(}->}[r] & \cI(\cV),\cP(\cV) \ar@{^{(}->}[dr] \ar@{^{(}->}[ur]   & \fVf \ar@{<->}[d] \ar@{<->}[u] \\
                              & & \cP\cI^{\ad}(\cV)  \ar@{^{(}->}[r] &\cP\cI(\cV) \ar@{^{(}->}[r] &  \Pro(\cI(\cV)) \ar@{^{(}->}[r] & \Pro(\Ind(\cV)),
}
$$
où $\cP(\cV)\xrightarrow{\sim}\Pro(\cV)$ est abélienne. Chacune de ces catégories se plonge dans la $\cU(1)$-catégorie correspondante. 
\end{Pt}

\chapter{Connexité dans les groupes indpro et proindalgébriques} 

Soit $k$ un corps algébriquement clos de caractéristique $p>0$.

\section{Groupes indpro et proindalgébriques connexes}

\begin{Not}
On note $\cQ^0$ la sous-catégorie pleine de $\cQ$ dont les objets sont les groupes quasi-algébriques connexes (sur $k$). 
\end{Not}

\begin{Pt}\label{PtQQ0}
La catégorie $\cQ^0$ est épaisse dans $\cQ$, en particulier elle est exacte ; les monomorphismes sont les injections quasi-algébriques et sont tous stricts, les épimorphismes sont les surjections quasi-algébriques et les épimorphismes stricts sont les épimorphismes dont le noyau est connexe.
\end{Pt}

\begin{Def} \emph{(Serre \cite{S2} 5.1).} \label{DefP0}
Un \emph{groupe proalgébrique connexe} est un groupe proalgébrique $(G,S)$ tel que $G/H\in \cQ^0$ pour tout $H\in S$.
\end{Def}

\begin{Pt} \label{PtDefP0}
Un groupe proalgébrique est connexe si et seulement s'il possède un ensemble de définition $S'$ tel que $G/H'\in \cQ^0$ pour tout $H'\in S'$ : cela résulte de (P-3), puisque $\cQ^0\subset\cQ$ est stable par quotient.
\end{Pt}

\begin{Not}
On note $\cP^0$ la sous-catégorie pleine de $\cP$ dont les objets sont les groupes proalgébriques connexes. 
\end{Not}

\begin{Prop}\label{P0exacte}
La catégorie $\cP^0$ est épaisse dans $\cP$ ; en particulier, elle est exacte. Les monomorphismes sont les injections proalgébriques et sont tous stricts. Les épimorphismes sont les surjections proalgébriques et les épimorphismes stricts sont les épimorphismes dont le noyau est connexe.
\end{Prop}

\begin{proof}
La catégorie $\cP^0$ est épaisse dans $\cP$ d'après \ref{descrseP} et le fait que $\cQ^0$ est épaisse dans $\cQ$. De plus, le quotient d'un groupe proalgébrique connexe par un sous-groupe fermé est connexe d'après \ref{descrseP} et le fait que 
$\cQ^0\subset\cQ$ est stable par quotient. 
\end{proof}

\begin{Def} \emph{(Serre \emph{loc. cit.}).} \label{G0P}
\begin{itemize}

\medskip

\item La \emph{composante neutre d'un groupe proalgébrique} $(G,S)$ est le sous-groupe fermé 
$$G^0:=\big(\lp_{H\in S}(G/H)^0, \{\Ker (\lp_{H\in S}(G/H)^0\ra(G/H)^0)\ |\ H\in S\}\big).$$ 

\medskip

\item Le \emph{groupe des composantes connexes d'un groupe proalgébrique} $(G,S)$ est le groupe quotient 
$$\pi_0(G):=G/G^0=\big(\lp_{H\in S}\pi_0(G/H), \{\Ker (\lp_{H\in S}\pi_0(G/H)\ra\pi_0(G/H))\ |\ H\in S\}\big).$$
\end{itemize}
\end{Def}

\begin{Pt}\label{PtG0P}
Soit $(G,S)\in\cP$. Alors $G^0$ est un groupe proalgébrique connexe, $\pi_0(G)$ est un groupe profini, et $(G,S)\in\cP^0\Leftrightarrow G^0=G \Leftrightarrow \pi_0(G)=0$.
\end{Pt}

\begin{Def} \label{DefI0}
Un \emph{groupe indalgébrique connexe} est un groupe indalgébrique $(G,S)$ vérifiant la condition suivante :

\medskip

$(c)$ Il existe un ensemble de définition $S'$ de $(G,S)$ tel que pour tout $H'\in S'$, le groupe quasi-algébrique $H'$ est connexe.
\end{Def}

\begin{Pt}\label{PtDefI0}
Un groupe indalgébrique est connexe si et seulement si l'ensemble de ses sous-groupes de définition connexes est un sous-ensemble de définition. 
\end{Pt}

\begin{Not}
On note $\cI^0$ la sous-catégorie pleine de $\cI$ dont les objets sont les groupes indalgébriques connexes.
\end{Not}

\begin{Prop} \label{I0exacte}
La catégorie $\cI^0$ est épaisse dans $\cI$ ; en particulier, elle est exacte.
\end{Prop}

\begin{proof}
Soit 
$$
\xymatrix{
0\ar[r] & (G',S') \ar[r] & (G,S) \ar[r]^f & (G'',S'') \ar[r] & 0
}
$$
une suite exacte de $\cI$, avec $(G',S'),(G'',S'')\in\cI^0$. D'après \ref{descrseI}, une telle suite est ind-représentée par les suites exactes de $\cQ$
$$
\xymatrix{
0\ar[r] & H' \ar[r] & H \ar[r] & H'' \ar[r] & 0, & H\in S,
}
$$
où $H':=G'\cap H$ et $H'':=f(H)$. Le foncteur $\pi_0$ sur $\cQ$ étant exact à droite, les suites exactes ci-dessus induisent un système inductif filtrant de suites exactes de groupes finis 
$$
\xymatrix{
\pi_0(H') \ar[r] & \pi_0(H) \ar[r] & \pi_0(H'') \ar[r] & 0, & H\in S,
}
$$
puis, en passant à la limite sur $H\in S$, une suite exacte de groupes 
$$
\xymatrix{
\li_{H'\in S'}\pi_0(H') \ar[r] & \li_{H\in S} \pi_0(H) \ar[r] & \li_{H''\in S''} \pi_0(H'') \ar[r] & 0.
}
$$
Par hypothèse, les deux termes extrémaux sont nuls, donc le terme central également. Comme les $\pi_0(H)$, $H\in S$, sont finis, on obtient $(G,S)\in \cI^0$. 
\end{proof}

\begin{Pt}
Soit $(G,S)\in\cI$. Pour tout $H\in S$, le groupe des composantes connexes $\pi_0(H)$ est fini. Par conséquent $(\pi_0(H))_{H\in S}\in\Ind(\cQ)$ vérifie la condition $(*)$ \ref{condisoinjI}, et définit donc un objet de $\cI$, quotient strict de $(G,S)$. La définition suivante a donc un sens :
\end{Pt}

\begin{Def} \label{G0I}
\begin{itemize}

\medskip

\item La \emph{composante neutre d'un groupe indalgébrique} $(G,S)$ est le sous-groupe strict 
$$G^0:=\big(\li_{H\in S}H^0, <\{H^0\ |\ H\in S\}>\big).$$ 

\medskip

\item Le \emph{groupe des composantes connexes d'un groupe indalgébrique} $(G,S)$ est le groupe quotient strict 
$$\pi_0(G):=G/G^0=\big(\li_{H\in S}\pi_0(H), \{\Ima (\pi_0(H)\ra \li_{H\in S}\pi_0(H))\ |\ H\in S\}\big).$$
\end{itemize}
\end{Def}

\begin{Pt} \label{PtG0I}
Soit $(G,S)\in\cI$. Alors $G^0$ est un groupe indalgébrique connexe, $\pi_0(G)$ est un groupe indfini, et $(G,S)\in\cI^0\Leftrightarrow G^0=G \Leftrightarrow \pi_0(G)=0$.
\end{Pt}

\begin{Def} \label{DefIP0}
Un \emph{groupe indproalgébrique connexe} est un groupe indproalgébrique $(G,S)$ vérifiant la condition suivante :

\medskip

$(c)$ Il existe un ensemble de définition $S'$ de $(G,S)$ tel que pour tout $H'\in S'$, le groupe proalgébrique $H'$ est connexe.
\end{Def}

\begin{Pt}\label{PtDefIP0}
Un groupe indproalgébrique est connexe si et seulement si l'ensemble de ses sous-groupes de définition connexes est un sous-ensemble de définition. 
\end{Pt}

\begin{Not}
On note $\cI(\cP^0)$ la sous-catégorie pleine de $\cI\cP$ dont les objets sont les groupes indproalgébriques connexes.
\end{Not}

\begin{Prop}\label{IP0exacte}
La catégorie $\cI(\cP^0)$ est épaisse dans $\cI\cP$ ; en particulier, elle est exacte.
\end{Prop}

\begin{proof}
Soit 
$$
\xymatrix{
0\ar[r] & (G',S') \ar[r] & (G,S) \ar[r]^f & (G'',S'') \ar[r] & 0
}
$$
une suite exacte de $\cI\cP$, avec $(G',S'),(G'',S'')\in\cI(\cP^0)$. D'après \ref{descrseIP}, une telle suite est ind-représentée par les suites exactes de $\cP$
$$
\xymatrix{
0\ar[r] & H' \ar[r] & H \ar[r] & H'' \ar[r] & 0, & H\in S,
}
$$
où $H':=G'\cap H$ et $H'':=f(H)$. Soit $H\in S$. Par hypothèse, il existe $H_{1}''\in S''$ connexe contenant $H'':=f(H)$. Choisissons $H_1\in S$ relevant $H_{1}''$. Quitte à remplacer $H_1$ par $H_1+H$ (\ref{IP2'}), on peut supposer que $H_1$ contient $H$. Posons alors $H_{1}':=G'\cap H_1$. Par hypothèse, il existe $H_{2}'\in S'$ connexe contenant $H_{1}'$. Alors $H_2:=H_{2}'+H_1$ appartient à $S$ (\ref{IP2'}), contient $H$, et est connexe d'après \ref{P0exacte} car extension de $H_{1}''$ par $H_{2}'$. L'ensemble des sous-groupes de définition connexes de $(G,S)$ est donc un sous-ensemble de définition, i.e. $(G,S)\in\cI(\cP^0)$.
\end{proof}

\begin{Pt}
Soit $(G,S)\in\cI\cP$. Pour tout $H\in S$, le groupe des composantes connexes $\pi_0(H)$ est profini (\ref{PtG0P}). En général, $(\pi_0(H))_{H\in S}\in\Ind(\cP)$ ne vérifie pas la condition $(*)$ \ref{condisoinjIP}, et ne définit donc pas un objet de $\cI\cP$. Cependant, si $(G,S)$ possède un ensemble de définition dont les éléments ont un $\pi_0$ fini, alors $(\pi_0(H))_{H\in S}\in\Ind(\cP)$ vérifie la condition $(*)$ \ref{condisoinjIP}.
\end{Pt}

\begin{Def}\label{DefIPpizf}
La catégorie $\cI(\cP^{\pizf})$ est la sous-catégorie pleine de $\cI\cP$ dont les objets sont les groupes indproalgébriques $(G,S)$ vérifiant la condition suivante : 

\medskip

$(\pizf)$ Il existe un ensemble de définition $S'$ de $(G,S)$ tel que pour tout $H'\in S'$, le groupe profini $\pi_0(H')$ est fini.
\end{Def}

\begin{Pt}\label{PtDefIPpizf}
Un groupe indproalgébrique est un objet de $\cI(\cP^{\pizf})$ si et seulement si l'ensemble de ses sous-groupes de définition dont le $\pi_0$ est fini est un sous-ensemble de définition. 
\end{Pt}

\begin{Def} \label{G0IPpizf}
\begin{itemize}

\medskip

\item La \emph{composante neutre de $(G,S)\in\cI(\cP^{\pizf})$} est le sous-groupe indproalgébrique strict 
$$G^0:=\big(\li_{H\in S}H^0, <\{H^0\ |\ H\in S\}>\big).$$ 

\medskip

\item Le \emph{groupe des composantes connexes de $(G,S)\in\cI(\cP^{\pizf})$} est le groupe indproalgébrique quotient strict 
$$\pi_0(G):=G/G^0=\big(\li_{H\in S}\pi_0(H), \{\Ima (\pi_0(H)\ra \li_{H\in S}\pi_0(H))\ |\ H\in S\}\big).$$
\end{itemize}
\end{Def}

\begin{Pt}
Soit $(G,S)\in\cI(\cP^{\pizf})$. Alors $G^0$ est un groupe indproalgébrique connexe, $\pi_0(G)$ est un groupe indfini, et $(G,S)\in\cI^0\Leftrightarrow G^0=G \Leftrightarrow \pi_0(G)=0$.
\end{Pt}

\begin{Prop}\label{pi0exactd}
Le foncteur $\pi_0:\cI(\cP^{\pizf})\ra\cI$ est exact à droite.
\end{Prop}

\begin{proof}
Soit 
$$
\xymatrix{
0\ar[r] & (G',S') \ar[r] & (G,S) \ar[r]^f & (G'',S'') \ar[r] & 0
}
$$
une suite exacte de $\cI(\cP^{\pizf})$. D'après \ref{descrseIP}, une telle suite est ind-représentée par les suites exactes de 
$\cP$
$$
\xymatrix{
0\ar[r] & H' \ar[r] & H \ar[r] & H'' \ar[r] & 0, & H\in S,
}
$$
où $H':=G'\cap H$ et $H'':=f(H)$. Le foncteur $\pi_0$ sur $\cP$ étant exact à droite d'après \cite{S2} 5.1 Proposition 2, les suites exactes ci-dessus induisent un système inductif filtrant de suites exactes de groupes profinis
$$
\xymatrix{
\pi_0(H') \ar[r] & \pi_0(H) \ar[r] & \pi_0(H'') \ar[r] & 0, & H\in S,
}
$$
définissant une suite exacte de $\Ind(\cP)$
$$
\xymatrix{
(\pi_0(H'))_{H'\in S'} \ar[r] & (\pi_0(H))_{H\in S} \ar[r] & (\pi_0(H''))_{H''\in S''} \ar[r] & 0.
}
$$
Mais par définition de $\pi_0:\cI(\cP^{\pizf})\ra\cI$, cette suite est isomorphe à l'image de la suite
$$
\xymatrix{
\pi_0(G',S') \ar[r] & \pi_0(G,S) \ar[r]^{\pi_0(f)} & \pi_0(G'',S'') \ar[r] & 0
}
$$
par le plongement $\cI\cP\ra\Ind(\cP)$ \ref{IPexacte}. La suite ci-dessus est donc exacte dans $\cI$.
\end{proof}

\begin{Def} \label{DefPI0}
Un \emph{groupe proindalgébrique connexe} est un groupe proindalgébrique $(G,S)$ tel que $G/H\in \cI^0$ pour tout $H\in S$.
\end{Def}

\begin{Pt} \label{PtDefPI0}
Pour qu'un groupe proindalgébrique soit connexe, il faut et il suffit qu'il possède un ensemble de définition $S'$ tel que $G/H'\in \cI^0$ pour tout $H'\in S'$ : cela résulte de (PI-3), puisque $\cI^0\subset\cI$ est stable par quotient strict d'après \ref{descrseI} et le fait que $\cQ^0\subset\cQ$ est stable par quotient.
\end{Pt}

\begin{Not}
On note $\cP(\cI^0)$ la sous-catégorie pleine de $\cP\cI$ dont les objets sont les groupes proindalgébriques connexes. 
\end{Not}

\begin{Pt}
Soit $(G,S)\in\cP\cI$. Pour tout $H,H'\in S$ avec $H\subset H'$, l'application canonique $G/H\ra G/H'$ est un épimorphisme strict de $\cI$ d'après (PI-3). Il résulte alors de \ref{descrseI} et du fait que le foncteur $\pi_0$ sur $\cQ$ transforme épimorphisme en épimorphisme, que l'application $\pi_0(G/H)\ra\pi_0(G/H')$ est un épimorphisme strict de $\cI$. D'après \ref{PIexacte}, $(\pi_0(G/H))_{H\in S}\in\Pro(\cI)$ définit donc un objet de $\cP\cI$. 
\end{Pt}

\begin{Def} \label{G0PI}
\begin{itemize}

\medskip

\item La \emph{composante neutre d'un groupe proindalgébrique} $(G,S)$ est le sous-groupe strict 
$$G^0:=\big(\lp_{H\in S}(G/H)^0, \{\Ker (\lp_{H\in S}(G/H)^0\ra(G/H)^0)\ |\ H\in S\}\big).$$ 

\medskip

\item Le \emph{groupe des composantes connexes d'un groupe proindalgébrique} $(G,S)$ est le groupe quotient strict 
$$\pi_0(G):=G/G^0=\big(\lp_{H\in S}\pi_0(G/H), \{\Ker (\lp_{H\in S}\pi_0(G/H)\ra\pi_0(G/H))\ |\ H\in S\}\big).$$
\end{itemize}
\end{Def}

\begin{Pt}
Soit $(G,S)\in\cP\cI$. Alors $G^0$ est un groupe proindalgébrique connexe, $\pi_0(G)$ est un groupe proindfini, et $(G,S)\in\cP(\cI^0)\Leftrightarrow G^0=G \Leftrightarrow \pi_0(G)=0$.
\end{Pt}

\section{Groupes parfaits algébriques connexes et parfaits proalgébriques affines connexes}

\begin{Pt}
Soit $(G,\cC_G)\in\cQ$. Si $(\cG,\iota)\in\cC_G$, alors la topologie induite sur $G$ par l'injection 
$$
\xymatrix{
G \ar[r]^<<<<<<{\iota^{-1}}_<<<<{\sim} & \cG(k) \ar@{^{(}->}[r] & \cG
}
$$
ne dépend pas du choix de $(\cG,\iota)\in\cC_G$, et coïncide avec celle induite par l'injection canonique
$$
\xymatrix{
G \ar[r]^<<<<<{\sim} & G^{\pf}(k) \ar@{^{(}->}[r] & G^{\pf}.
}
$$
On notera que le groupe $G$ muni de cette topologie n'est un groupe topologique que si $G$ est fini. L'application $G\ra G$, $x\mapsto -x$, et les translations de $G$ sont cependant des homéomorphismes.
\end{Pt}

\begin{Def}
\begin{itemize}

\medskip

\item La \emph{composante neutre d'un groupe quasi-algébrique $(G,\cC_G)$} est le sous-groupe quasi-algébrique ouvert et fermé ayant pour groupe ordinaire sous-jacent la composante connexe de l'élément neutre de $G$ et pour classe de structures de groupe algébrique 
$$
\{(\cG^0,\iota|_{\cG^0})\ |\ (\cG,\iota)\in\cC_G\}.
$$

\medskip

\item Le \emph{groupe des composantes connexes d'un groupe quasi-algébrique $(G,\cC_G)$} est le groupe quasi-algébrique quotient $\pi_0(G):=G/G_0$, qui est fini.
\end{itemize}
\end{Def}

\begin{Pt}
D'après \ref{pfexact}, la suite exacte de $\cQ$
$$
\xymatrix{
0 \ar[r] & G^0 \ar[r] & G \ar[r] & \pi_0(G) \ar[r] & 0
}
$$
induit une suite exacte de groupes parfaits algébriques 
$$
\xymatrix{
0 \ar[r] & (G^0)^{\pf} \ar[r] & G^{\pf} \ar[r] & \pi_0(G) \ar[r] & 0.
}
$$
La composante neutre du $k$-schéma en groupes $G^{\pf}$ (\cite{SGA 3nv} VI${}_A$ 2.6.5) est donc un sous-schéma en groupes parfait algébrique ouvert et fermé de $G^{\pf}$ canoniquement isomorphe à $(G^0)^{\pf}$, et le quotient 
$\pi_0(G^{\pf}):=G^{\pf}/(G^{\pf})^0$ est canoniquement isomorphe au $k$-schéma en groupes constant fini $\pi_0(G)$. Ainsi 
$(G^{\pf})^0=(G^0)^{\pf}$ et $\pi_0(G^{\pf})=\pi_0(G)$.
\end{Pt}

\begin{Pt}\label{cPA0}
Soit $(G,S)\in\cP(\cA)$. La suite exacte de $\cP(\cA)$ 
$$
\xymatrix{
0 \ar[r] & G^0 \ar[r] & G \ar[r] & \pi_0(G) \ar[r] & 0
}
$$
induit via \ref{cPApf} une suite exacte de groupes parfaits proalgébriques affines
$$
\xymatrix{
0 \ar[r] & (G^0)^{\pf} \ar[r] & G^{\pf} \ar[r] & \pi_0(G)^{\pf} \ar[r] & 0,
}
$$
limite du système projectif de suites exactes de groupes parfaits algébriques
$$
\xymatrix{
0 \ar[r] & ((G/H)^0)^{\pf} \ar[r] & (G/H)^{\pf} \ar[r] & \pi_0(G/H) \ar[r] & 0, & H\in S.
}
$$
En particulier, $(G^0)^{\pf}$ est canoniquement isomorphe au sous-schéma en groupes fermé noyau du morphisme de $k$-schémas en groupes $G^{\pf} \ra \pi_0(G)^{\pf}$, et celui-ci est fidèlement plat affine. L'immersion fermée  $(G^0)^{\pf} \ra G^{\pf}$ est plate, le $k$-schéma en groupes $(G^0)^{\pf}$ est intègre et l'espace topologique sous-jacent à $\pi_0(G)^{\pf}$ est totalement discontinu, de sorte que $(G^0)^{\pf}$ est une composante irréductible de $G^{\pf}$ ; c'est donc la \emph{composante neutre de $G^{\pf}$} (\cite{SGA 3nv} VI${}_A$ 2.6.5). Par conséquent, pour tout groupe parfait proalgébrique affine $G$, la composante neutre $G^0$ de $G$ est un groupe proalgébrique parfait affine, et le morphisme canonique $G^0\ra G$ est un monomorphisme dans la catégorie abélienne des groupes parfaits proalgébriques affines ; on appellera \emph{groupe des composantes connexes de $G$} et on notera $\pi_0(G)$ le quotient de $G$ par $G^0$ dans cette catégorie. C'est un $k$-schéma en groupes profini d'après ce qui précède, qui est constant si et seulement s'il est fini. Ainsi 
$(G^{\pf})^0=(G^0)^{\pf}$ et $\pi_0(G^{\pf})=\pi_0(G)^{\pf}$.
\end{Pt}

\begin{Pt}
Soit $(G,S)\in\cP(\cA)$. La topologie induite sur $G$ par l'injection canonique
$$
\xymatrix{
G \ar[r]^<<<<<{\sim} & G^{\pf}(k) \ar@{^{(}->}[r] & G^{\pf}
}
$$
coïncide avec la topologie limite projective des espaces topologiques $G/H$, $H\in S$. En particulier, la composante connexe de l'élément neutre de $G$ pour cette topologie coïncide avec le groupe ordinaire sous-jacent au groupe proalgébrique $G^0$ composante neutre de $G$. 
\end{Pt}

Pour les groupes indproalgébriques possédant un ensemble de définition dont les éléments sont proalgébriques affines à 
$\pi_0$ fini, le foncteur fibre générique est conservatif. Plus généralement, on a la variante suivante de \ref{critfoncfibgenIPcons} :

\begin{Prop}\label{foncfibgenIPcons00}
Soit $f:(G_1,S_1)\ra (G_2,S_2)$ un morphisme de groupes indproalgébriques. Supposons que $(G_2,S_2)$ possède un ensemble de définition dont les éléments sont proalgébriques affines à $\pi_0$ fini. Alors pour que $f$ soit un isomorphisme, il faut et il suffit que le morphisme $f_{\eta}:(G_1,S_1)_{\eta}\ra (G_2,S_2)_{\eta}$ induit sur les fibres génériques soit un isomorphisme.
\end{Prop}

\begin{proof}
Soit $f:(G_1,S_1)\ra (G_2,S_2)$ tel que $f_{\eta}:(G_1,S_1)_{\eta}\ra (G_2,S_2)_{\eta}$ soit un isomorphisme. Soit $H_1\in S_1$. Par hypothèse, $f(H_1)\in S_2$ et le morphisme $H_1\ra f(H_1)$ est proalgébrique. Comme de plus $f(k)$ est injectif, $H_1\ra f(H_1)$ est un isomorphisme proalgébrique. Nous allons montrer que 
$$
f(S_1):=\{ f(H_1)\ |\ H_1\in S_1 \}
$$
est égal à $S_2$. Les morphismes proalgébriques inverses des $H_1\xrightarrow{\sim} f(H_1)$ définiront donc un morphisme indproalgébrique $(G_2,S_2)\ra (G_1,S_1)$, inverse de $f$.

On a $f(S_1)\subset S_2$ par définition, il s'agit donc de montrer $S_2\subset f(S_1)$. Notant $S_{2}'$ un ensemble de définition de $(G_2,S_2)$ dont les éléments sont affines à $\pi_0$ fini, il suffit de montrer que pour tout $H_{2}'\in S_{2}'$, il existe $H_1\in S_1$ tel que $H_{2}'\subset f(H_1)$ : d'après (IP-2) pour $S_2$, l'inclusion $H_{2}'\subset f(H_1)$ est alors une immersion fermée proalgébrique, puis d'après (IP-2) pour $S_1$ (et l'isomorphisme $H_1\xrightarrow{\sim} f(H_1)$), $H_{2}'\in f(S_1)$, et finalement $S_2=\ <S_{2}'>\ \subset f(S_1)$. 

Soit $H_{2}'\in S_{2}'$. Par hypothèse et d'après \ref{cPA0}, le $k$-schéma en goupes parfait $(H_{2}')^{\pf}$ posséde un nombre fini de composantes irréductibles ; soient $\eta_1,\ldots,\eta_n$ les points génériques de 
$(H_{2}')^{\pf}$. Choisissons un corps algébriquement clos $\kappa$ contenant les corps résiduels $k(\eta_1),\ldots,k(\eta_n)$, et notons $\overline{\eta_1},\ldots,\overline{\eta_n}$ les $\kappa$-points de $(H_{2}')^{\pf}$ induits par les $\eta_i$. 
Comme par hypothèse l'application
$$
\xymatrix{
f(\kappa):\li_{H_1\in S_1} H_{1}(\kappa) \ar[r] & \li_{H_2\in S_2} H_{2}(\kappa)
}
$$
est bijective, il existe $H_1\in S_1$ tel que $\overline{\eta_1},\ldots,\overline{\eta_n}\in f(\kappa)(H_{1}(\kappa))$. Choisissons alors $H_{2}\in S_2$ contenant $f(H_1)$ et $H_{2}'$ ((IP-1) pour $S_2$). Les inclusions $f(H_1)\subset H_{2}$ et $H_{2}'\subset H_{2}$ étant des immersions fermées proalgébriques ((IP-2) pour $S_2$), elles donnent lieu à des immersions fermées de schémas $(f(H_{1}))^{\pf}\ra H_{2}^{\pf}$ et $(H_{2}')^{\pf}\ra H_{2}^{\pf}$. De plus, comme $\cP \ni H \mapsto H_{\eta} \in \widehat{\cT}$ est exact (\ref{foncfibgenIPexact}), $(f(H_1))^{\pf}(\kappa)=(f(H_1))(\kappa)=f(\kappa)(H_{1}(\kappa))$. Par conséquent, $(f(H_{1}))^{\pf}$ majore l'adhérence schématique des $\overline{\eta_i}$ dans $H_{2}^{\pf}$, qui est égale à $(H_{2}')^{\pf}$. En prenant les $k$-points, on obtient $f(H_1)\supset H_{2}'$.
\end{proof}

\section{Applications rationnelles de groupes indproalgébriques}

\begin{Prop}\label{ratIndPro0}
Soient $G=(G_i)_{i\in I}$ et $H=(H_j)_{j\in J}$ deux ind-objets de la pro-catégorie des $k$-schémas en groupes. On suppose :
\begin{enumerate}
\item Pour tout $i\in I$, $G_i$ est représenté par un système projectif dont les objets sont des $k$-schémas en groupes parfaits algébriques, et dont les transitions sont surjectives et induisent des bijections sur les groupes de composantes connexes.
\item Pour tout $j\in J$, $H_j$ est représenté par un système dont les objets sont des $k$-schémas en groupes parfaits localement algébriques.
\end{enumerate}
Soit $\cT'$ la sous-catégorie pleine de $(\parf/k)$ dont les objets sont les morphismes $\Spec(k')\ra\Spec(k)$ induits par les morphismes de corps parfaits $k\ra k'$. Soient $G_{\cT'}$ et $H_{\cT'}$ les restrictions des foncteurs de points de $G$ et $H$ à la catégorie $\cT'$. Alors l'application canonique
$$
\xymatrix{
\Hom_{\Ind(\Pro(\kschgr))}(G,H) \ar[r] & \Hom_{\widehat{\cT'}}(G_{\cT'},H_{\cT'})
}
$$
est bijective.
\end{Prop}

\begin{proof}
Commençons par montrer l'injectivité de la flèche considérée. Pour cela, on peut supposer que $I$ est un singleton, i.e. que $G$ est un pro-objet $(G_n)_{n\in N}$ de la catégorie des $k$-schémas en groupes. Posons par ailleurs
$$H=(H_j)_{j\in J}=\big((H_{m_{j}})_{m_{j}\in M_{j}}\big)_{j\in J}.$$
Par hypothèse, les transitions du pro-objet $(\pi_0(G_n))_{n\in N}$ de la catégorie des groupes sont bijectives. Les $G_n$ ont donc le même nombre $r$ de composantes irréductibles, et s'organisent en $r$ systèmes projectifs 
$(G_{n}^c)_{n\in N}$, $c=0,\ldots,r-1$, dont $(G_n)_{n\in N}$ est l'union disjointe. De plus, pour $c$ fixé et $n$ variable, les corps résiduels des points génériques $\eta_{G_{n}^c}$ sont parfaits, et définissent un ind-objet $(k(\eta_{G_{n}^c}))_{n\in N}$ de la catégorie des corps. Soit $k(\eta_{G^c})$ le corps parfait $\li_{n\in N} k(\eta_{G_{n}^c})$ et soit 
$p_{G_{n}^c}:\Spec(k(\eta_{G^c}))\ra \Spec(k(\eta_{G_{n}^c}))$ le morphisme induit par l'inclusion canonique $k(\eta_{G_{n}^c})\subset k(\eta_{G^c})$. Les morphismes composés
$$
\xymatrix{
p_{G_{n}^c}^*(\eta_{G_{n}^c}):\Spec(k(\eta_{G^c}))\ar[r]^>>>>>>{p_{G_{n}^c}} & \Spec(k(\eta_{G_{n}^c})) \ar[r]^>>>>>{\eta_{G_{n}^c}} & G_{n}^c,\quad n\in N,
}
$$
définissent un élément $\eta_{G^c}\in\lp_{n\in N} G_{n}^c(k(\eta_{G^c}))=G^c(k(\eta_{G^c}))$. Pour $n$ fixé, le morphisme 
$p_{G_{n}^c}^*(\eta_{G_{n}^c})$ est schématiquement dense dans $G_{n}^c$. L'injectivité annoncée en résulte. Soit en effet $f:G\ra H$ tel que $f_{\cT'}=0$. Par définition des morphismes de ind-objets, et puisque $I$ est un singleton, il existe $j\in J$ tel que $f$ se factorise en morphisme $f^j:G\ra H_j$. Quitte à augmenter $j\in J$, on peut supposer que $f^j(\eta_{G^c})=0\in H_{j}(k(\eta_{G^c}))$ pour tout $c=0,\ldots,r-1$. Choisissant alors une représentation du morphisme de pro-objets $f^j$ par une famille compatible de morphismes de $k$-schémas en groupes $f_{\varphi(m_j),m_j}^{j}:G_{\varphi(m_j)}\ra H_{m_j}$, où $\varphi$ est une application de $M_j$ dans $N$, on a pour tout $c=0,\ldots,r-1$
$$
0=f^{j}(\eta_{G^c})=\big(f_{\varphi(m_j),m_j}^j(p_{G_{\varphi(m_j)}^c}^*(\eta_{G_{\varphi(m_j)}^c}))\big)_{m_j\in M_j}\ \in H_j(k(\eta_{G^c}))=\lp_{m_j\in M_j} H_{m_j}(k(\eta_{G^c})).
$$
D'où $f_{\varphi(m_j),m_j}^j=0$ pour tout $m_j\in M_j$, i.e. $f^j=0$, puis $f=0$.

Passons à la surjectivité. D'après l'injectivité que l'on vient d'établir, on peut à nouveau supposer que $I$ est un singleton. Traitons alors dans un premier temps le cas où les $G_n$ sont connexes ; on note $\eta_{G_n}$ le point générique de $G_n$ et $\eta_{G}\in G(k(\eta_G))$ le point construit plus haut. Soit $f_{\cT'}\in\Hom_{\widehat{\cT'}}(G_{\cT'},H_{\cT'})$.  Choisissons $j\in J$ tel que $f_{\cT'}(\eta_{G})\in H_j(k(\eta_G))=\lp_{m_j\in M_j} H_{m_j}(k(\eta_G))$. On peut alors représenter $f_{\cT'}(\eta_G)$ par un système compatible de morphismes de $k$-schémas
$$
\xymatrix{
f_{\cT'}(\eta_G)_{m_j}:\Spec(k(\eta_G)) \ar[r] & H_{m_j},\quad m_j\in M_j.
}
$$
Fixons $m_j\in M_j$. Comme $H_{m_j}$ est le perfectisé d'un $k$-schéma localement de type fini sur $k$, et $k(\eta_G)$ limite inductive filtrante des corps parfaits $k(\eta_{G_n})$, il existe $\varphi(m_j)\in N$ tel que $f_{T'}(\eta_G)_{m_j}$ se factorise en 
$$
\xymatrix{
f_{\cT'}(\eta_{G})_{m_j}:\Spec(k(\eta_G)) \ar[r]^>>>>>{p_{G_{\varphi(m_j)}}} & \Spec(k(\eta_{G_{\varphi(m_j)}})) \ar[rr]^<<<<<<<<<<{f_{\cT'}(\eta_{G})_{\varphi(m_j),m_j}} & & H_{m_j},
}
$$
puis un ouvert non vide $U_{\varphi(m_j)}$ de $G_{\varphi(m_j)}$ et un morphisme de $k$-schémas $f_{U_{\varphi(m_j)}}:U_{\varphi(m_j)}\ra H_{m_j}$ prolongeant $f_{\cT'}(\eta_G)_{\varphi(m_j),m_j}$, i.e. $f_{\cT'}(\eta_G)_{\varphi(m_j),m_j}=f_{U_{\varphi(m_j)}}(\eta_{G_{\varphi(m_j)}})$; en particulier, si l'on pose $n:=\varphi(m_j)$, on a
$$f_{\cT'}(\eta_G)_{m_j}=f_{U_n}(p_{G_n}^*(\eta_{G_n}))\ \in H_{m_j}(k(\eta_{G})).$$

Ceci étant, pour $n\in N$, notons
$$
\xymatrix{
m_{G_n}: G_n\times_k G_n \ar[r] & G_n
}
$$
la loi de groupe de $G_n$ ; notation analogue pour les lois de groupe des $H_{m_j}$, $m_j\in M_j$. Montrons que, quitte à augmenter $j$, les diagrammes
\begin{equation}\label{mgUnipparf}
\xymatrix{
U_n\times_k U_n\cap m_{G_n}^{-1}(U_n) \ar[r]^<<<<{m_{G_n}} \ar[d]_{f_{U_n}\times_k f_{U_n}} & U_n \ar[d]^{f_{U_n}} \\
H_{m_j}\times_k H_{m_j} \ar[r]^<<<<<<<{m_{H_{m_j}}} & H_{m_j},&\quad m_j\in M_j, 
}
\end{equation}
où l'on a posé $n:=\varphi(m_j)$, sont commutatifs. Pour $n\in N$, le point générique $\eta_{G_n\times_kG_n}$ de $G_n\times_k G_n$
est schématiquement dense dans $G_n\times_k G_n$. Il suffit donc de montrer que, quitte à augmenter $j$, le diagramme d'indice $m_j$ commute après restriction à  $\eta_{G_n\times_kG_n}$, et ceci pour tout $m_j\in M_j$. On note $k(\eta_{G_n\times_k G_n})$ le corps résiduel de $\eta_{G_n\times_kG_n}$, qui est parfait, et $p_{1,n},p_{2,n}:G_n\times_k G_n\ra G_n$ les projections sur le premier et le second facteur.  Etant fidèlement plates, elles donnent lieu à des morphismes $p_{1\eta,n},p_{2\eta,n}$, rendant commutatif le diagramme
$$
\xymatrix{
\Spec(k(\eta_{G_n\times_k G_n})) \ar[dr]^{\eta_{G_n\times_k G_n}} \ar[rr]^{p_{2\eta,n}}  \ar[dd]_{p_{1\eta,n}} & & \Spec(k(\eta_{G_n})) \ar[dr]^{\eta_{G_n}} \\
 & G_n\times_k G_n \ar[rr]^{p_{2,n}} \ar[dd]_{p_{1,n}} & & G_n \\
\Spec(k(\eta_{G_n})) \ar[dr]^{\eta_{G_n}} & &  \\
& G_n, &
}
$$
i.e. $\eta_{G_n\times_k G_n}=(p_{1\eta,n}^*(\eta_{G_n}),p_{2\eta,n}^*(\eta_{G_n}))$. Il s'agit donc de montrer que, quitte à augmenter $j\in J$, on a pour tout $m_j\in M_j$ l'égalité suivante dans $H_{m_j}(k(\eta_{G_n\times_k G_n}))$, où $n:=\varphi(m_j)$ :
\begin{equation}\label{mqnipparf}
f_{U_n}(p_{1\eta,n}^*(\eta_{G_n})+p_{2\eta,n}^*(\eta_{G_n}))=f_{U_n}(p_{1\eta,n}^*(\eta_{G_n}))+f_{U_n}(p_{2\eta,n}^*(\eta_{G_n})).
\end{equation}
Notons $k(\eta_{G\times_k G})$ le corps parfait $\li_{n\in N} k(\eta_{G_n\times_k G_n})$ et  
$$
p_{G_n\times_kG_n}:\Spec(k(\eta_{G\times G}))\ra \Spec(k(\eta_{G_n\times G_n}))
$$ 
le morphisme induit par l'inclusion canonique $k(\eta_{G_n\times_k G_n})\subset k(\eta_{G\times_k G})$. Les morphismes composés
$$
\xymatrix{
p_{G_n\times_kG_n}^*(\eta_{G_n\times G_n}):\Spec(k(\eta_{G\times G}))\ar[rr]^<<<<<<<<<<<{p_{G_n\times_kG_n}} & & \Spec(k(\eta_{G_n\times G_n})) \ar[rr]^<<<<<<<<<<{\eta_{G_n\times_k G_n}} & & G_n\times_k G_n,
}
$$
définissent un élément $\eta_{G\times_k G}\in\lp_{n\in N} (G_n\times_k G_n)(k(\eta_{G\times G}))=(G\times_kG)(k(\eta_{G\times G}))$. Pour tout $n\in N$, on a alors un diagramme commutatif 
$$
\xymatrix{
\Spec(k(\eta_{G\times_k G})) \ar[d]_{p_{1\eta}:=\lp_n p_{1\eta,n}} \ar[rr]^{p_{G_n\times_kG_n}} & & \Spec(k(\eta_{G_n\times_kG_n})) \ar[d]^{p_{1\eta,n}} \ar[rr]^>>>>>>>>>>{\eta_{G_n\times_kG_n}} & & G_n\times_k G_n \ar[d]^{p_{1,n}} \\
\Spec(k(\eta_G)) \ar[rr]^{p_{G_n}} & & \Spec(k(\eta_{G_n})) \ar[rr]^{\eta_{G_n}} & & G_n,
}
$$
et de même avec $p_2$ à la place de $p_1$. Pour montrer l'égalité (\ref{mqnipparf}) dans $H_{m_j}(k(\eta_{G_n\times_k G_n}))$, il suffit de la montrer dans $H_{m_j}(k(\eta_{G\times_k G}))$, c'est-à-dire, avec les notations précédentes, il suffit de montrer que 
\begin{equation}\label{mqnbisipparf}
f_{U_n}(p_{1\eta}^*(p_{G_n}^*(\eta_{G_n}))+p_{2\eta}^*(p_{G_n}^*(\eta_{G_n})))=f_{U_n}(p_{1\eta}^*(p_{G_n}^*(\eta_{G_n})))+f_{U_n}(p_{2\eta}^*(p_{G_n}^*(\eta_{G_n}))).
\end{equation}
Comme $f_{\cT'}:G_{\cT'}\ra H_{\cT'}$ est un morphisme de préfaisceaux, on a
$$
f_{\cT'}(p_{1\eta}^*(\eta_G))  =  p_{1\eta}^*(f_{\cT'}(\eta_G)) = p_{1\eta}^*((f_{U_n}(p_{G_n}^*(\eta_{G_n})))_{m_j\in M_j}) =  (f_{U_n}(p_{1\eta}^*(p_{G_n}^*(\eta_{G_n}))))_{m_j\in M_j}
$$
dans $H(k(\eta_{G\times_k G}))$. Autrement dit, $f_{\cT'}(p_{1\eta}^*(\eta_G))\in H(k(\eta_{G\times_k G}))$ est représentable par l'élément $(f_{\cT'}(p_{1\eta}^*(\eta_G))_{m_j})_{m_j\in M_j}\in H_j(k(\eta_{G\times_k G}))$ défini par
\begin{equation}\label{mqdnipparf}
\forall m_j\in M_j,\ f_{\cT'}(p_{1\eta}^*(\eta_G))_{m_j}:=f_{U_n}(p_{1\eta}^*(p_{G_n}^*(\eta_{G_n})))\ \in H_{m_j}(k(\eta_{G\times_k G}))
\end{equation}
(où l'on pose toujours $n:=\varphi(m_j)$). De même avec $p_2$ à la place de $p_1$. Par ailleurs, $m_{G_n}:G_n\times_kG_n\ra G_n$ étant fidèlement plate, elle donne lieu à un morphisme $m_{G_{\eta,n}}$ rendant commutatif le diagramme
$$
\xymatrix{
\Spec(k(\eta_{G_n\times_k G_n}))  \ar[rr]^<<<<<<<<<<{\eta_{G_n\times_k G_n}}  \ar[d]_{m_{G_{\eta,n}}} & & G_n\times_k G_n \ar[d]^{m_{G_n}} \\
\Spec(k(\eta_{G_n}))  \ar[rr]^<<<<<<<<<<<<<{\eta_{G_n}} & & G_n, \\
}
$$
i.e. $(p_{1\eta,n}^*(\eta_{G_n})+p_{2\eta,n}^*(\eta_{G_n})=)m_{G_n}(\eta_{G_n\times_k G_n})=m_{G_{\eta,n}}^*(\eta_{G_n})$. Pour tout $n\in N$, on a alors un diagramme commutatif
$$
\xymatrix{
\Spec(k(\eta_{G\times_k G})) \ar[d]_{m_{G\eta}:=\lp_n m_{G_{\eta,n}}} \ar[rr]^{p_{G_n\times_kG_n}} & & \Spec(k(\eta_{G_n\times_kG_n})) \ar[d]^{m_{G_{\eta,n}}} \ar[rr]^>>>>>>>>>>{\eta_{G_n\times_kG_n}} & & G_n\times_k G_n \ar[d]^{m_{G_n}} \\
\Spec(k(\eta_G)) \ar[rr]^{p_{G_n}} & & \Spec(k(\eta_{G_n})) \ar[rr]^{\eta_{G_n}} & & G_n, 
}
$$
i.e. $p_{1\eta}^*(p_{G_n}^*(\eta_{G_n}))+p_{2\eta}^*(p_{G_n}^*(\eta_{G_n}))=m_{G\eta}^*(p_{G_n}^*(\eta_{G_n}))$ pour tout $n\in N$, ou encore $p_{1\eta}^*(\eta_G)+p_{2\eta}^*(\eta_G)=m_{G\eta}^*(\eta_G)$. Comme $f_{\cT'}:G_{\cT'}\ra H_{\cT'}$ est un morphisme de préfaisceaux, on a
\begin{eqnarray*}
f_{\cT'}(m_{G\eta}^*(\eta_G)) & = & m_{G\eta}^*(f_{\cT'}(\eta_G)) \\
                                              & = &  m_{G\eta}^*((f_{U_n}(p_{G_n}^*(\eta_{G_n})))_{m_j\in M_j})  \\
                                              & = & (f_{U_n}(m_{G\eta}^*(p_{G_n}^*(\eta_{G_n}))))_{m_j\in M_j}, 
\end{eqnarray*}
où $n:=\varphi(m_j)$. Autrement dit, $f_{\cT'}(p_{1\eta}^*(\eta_G)+p_{2\eta}^*(\eta_G))\in H(k(\eta_{G\times_k G}))$ est représentable
par l'élément $(f_{\cT'}(p_{1\eta}^*(\eta_G)+p_{2\eta}^*(\eta_G))_{m_j})_{m_j\in M_j} \in H_j(k(\eta_{G\times_k G}))$ défini par
\begin{equation}\label{mqgnipparf}
\forall m_j\in M_j,\ f_{\cT'}(p_{1\eta}^*(\eta_G)+p_{2\eta}^*(\eta_G))_{m_j} := f_{U_n}(p_{1\eta}^*(p_{G_n}^*(\eta_{G_n}))+p_{2\eta}^*(p_{G_n}^*(\eta_{G_n})))\ \in H_{m_j}(k(\eta_{G\times_k G})).
\end{equation}
Maintenant, 
$$
\xymatrix{
f_{\cT'}:G(k(\eta_{G\times_k G}))\ar[r] & H(k(\eta_{G\times_k G}))
}
$$ 
étant un morphisme de groupes, on a 
$$
f_{\cT'}(p_{1\eta}^*(\eta_G)+p_{2\eta}^*(\eta_G))=f_{\cT'}(p_{1\eta}^*(\eta_G))+f_{\cT'}(p_{2\eta}^*(\eta_G))\ \in H(k(\eta_{G\times_k G})).
$$
Il résulte donc des égalités (\ref{mqdnipparf}) et (\ref{mqgnipparf}) que, quitte à augmenter $j\in J$, 
$$
f_{U_n}(p_{1\eta}^*(p_{G_n}^*(\eta_{G_n}))+p_{2\eta}^*(p_{G_n}^*(\eta_{G_n})))=f_{U_n}(p_{1\eta}^*(p_{G_n}^*(\eta_{G_n})))+f_{U_n}(p_{2\eta}^*(p_{G_n}^*(\eta_{G_n})))
$$
dans $H_{m_j}(k(\eta_{G\times_k G}))$ pour tout $m_j\in M_j$. On a ainsi montré (\ref{mqnbisipparf}), et par conséquent, quitte à augmenter $j\in J$, les diagrammes (\ref{mgUnipparf}) sont bien commutatifs.

Pour tout $n\in N$, $k$-schéma $X$, et $x\in G_n(X)$, notons $\tau_x:(G_n)_X\ra (G_n)_X$ la translation par $x$. Pour $x\in U_n(k)$, avec $n=\varphi(m_j)$, on définit $f_{\tau_x(U_n)}:\tau_x(U_n)\ra H_{m_j}$ comme la composée
$$
\xymatrix{
\tau_x(U_n) \ar[r]^{\tau_{-x}} & U_n \ar[r]^{f_{U_n}} & H_{m_j} \ar[r]^{\tau_{f_{U_n}(x)}} & H_{m_j}.
}
$$
La commutativité des diagrammes (\ref{mgUnipparf}) assure que pour tout $m_j\in M_j$, $f_{U_n}$ et $f_{\tau_x(U_n)}$ coïncident sur $U_n\cap \tau_x(U_n)$. De plus, $G_n$ étant le perfectisé d'un $k$-schéma en groupes de type fini et connexe sur le corps algébriquement clos $k$, les $\tau_x(U_n)$ pour $x$ variant dans $U_n(k)$ recouvrent $G_n$. On obtient donc pour tout $m_j\in M_j$ un morphisme de $k$-schémas $f_n:G_n\ra H_{m_j}$ prolongeant $f_{U_n}$. Comme $U_n\times_k U_n\cap m_{G_n}^{-1}(U_n)$ est schématiquement dense dans $G_n\times_k G_n$, la commutativité de (\ref{mgUnipparf}) assure que $f_n$ est un morphisme de groupes. Lorsque $m_j\in M_j$ varie, les $f_n$ forment un système compatible :
si $m_{j}'\geq m_j\in M_j$, on choisit $n\in N$ majorant $\varphi(m_j)$ et $\varphi(m_{j}')$, on a alors le diagramme 
$$
\xymatrix{
                                    &                           & G_{\varphi(m_{j}')} \ar[r]^{f_{\varphi(m_{j}')}} & H_{m_{j}'} \ar[dd] \\
\Spec(k(\eta_G)) \ar[r]^<<<<<<{(\eta_G)_n} \ar[urr]^{(\eta_G)_{\varphi(m_{j}')}} \ar[drr]_{(\eta_G)_{\varphi(m_{j})}} & G_n \ar[ur] \ar[dr] &        & \\
                                    &                           & G_{\varphi(m_{j})} \ar[r]^{f_{\varphi(m_{j})}} & H_{m_{j}}, 
}
$$
qui est commutatif puisque $(f_{\varphi(m_j)}((\eta_G)_{\varphi(m_j)}))_{m_j\in M_j}=f_{\cT'}(\eta_G)\in \lp_{m_j\in M_j} H_{m_j}(k(\eta_G))$, si bien que le pentagone de droite est commutatif puisque $(\eta_G)_n$ est schématiquement dense dans $G_n$. On pose 
$$
\xymatrix{
f:=G\ar[rr]^<<<<<<<<<<{(f_{\varphi(m_j)})_{m_j\in M_j}} & & H_j \ar[r] & H.
}
$$

On a $f_{\cT'}(\eta_G)=f(\eta_G)\in H(k(\eta_G))$ par construction. Montrons que cette relation entraîne que le morphisme de préfaisceaux en groupes $G_{\cT'}\ra H_{\cT'}$ induit par $f$ est $f_{\cT'}$. Soit donc $k'/k$ une extension parfaite du corps $k$, et montrons que les morphismes $G(k')\rra H(k')$ induits par $f$ et $f_{\cT'}$ coïncident. 

On commence par se ramener au cas où $k'=k$. Soient $\cT''$ la sous-catégorie pleine de la catégorie des $k'$-schémas dont les objets sont les extensions de corps $k''/k'$, avec $k''$ parfait, et $\widehat{\cT''}$ la catégorie des préfaisceaux en groupes sur $\cT''$. Posons $G':=G\times_k k'$ et $H':=H\times_k k'$, et soient $G_{\cT''}'$ et $H_{\cT''}'$ les restrictions des foncteurs de points de $G'$ et $H'$ à la catégorie $\cT''$. On a alors des identifications canoniques $G_{\cT''}'=G_{\cT'}|_{\cT''}$ et $H_{\cT''}'=H_{\cT'}|_{\cT''}$. En particulier, $f_{\cT'}|_{\cT''}$ s'identifie à un morphisme $f_{\cT''}':G_{\cT''}'\ra H_{\cT''}'$. Par ailleurs, $f$ induit un morphisme $f'=f\times_k k':G'\ra H'$. De plus, on a encore la relation 
$f_{\cT''}'(\eta_{G'})=f'(\eta_{G'})\in H'(k(\eta_{G'}))$. En effet, pour $n\in N$, notons $\pr_{G_n}:G_{n}'\ra G_n$ la projection canonique. Etant fidèlement plate, elle induit un morphisme $\pr_{G_{\eta,n}}$ rendant commutatif le diagramme 
$$
\xymatrix{
\Spec(k(\eta_{G_{n}'})) \ar[r]^<<<<{\eta_{G_{n}'}} \ar[d]_{\pr_{G_{\eta,n}}} & G_{n}' \ar[d]^{\pr_{G_n}} \\
\Spec(k(\eta_{G_n})) \ar[r]^<<<<{\eta_{G_n}} & G_n,
}
$$
i.e. $\pr_{G_n}(\eta_{G_{n}'})=\pr_{G_{\eta,n}}^*(\eta_{G_n})$. On a alors pour tout $n\in N$ un diagramme commutatif
$$
\xymatrix{
\Spec(k(\eta_{G'})) \ar[d]_{\pr_{G\eta}:=\lp_n\pr_{G_{\eta,n}}} \ar[rr]^{p_{G_{n}'}} & & \Spec(k(\eta_{G_{n}'})) \ar[d]^{\pr_{G_{\eta,n}}} \ar[rr]^>>>>>>>>>>{\eta_{G_{n}'}} & & G_{n}' \ar[d]^{\pr_{G_n}} \\
\Spec(k(\eta_{G})) \ar[rr]^{p_{G_n}} & & \Spec(k(\eta_{G_n})) \ar[rr]^>>>>>>>>>>{\eta_{G_n}} & & G_n.
}
$$
Posons $\pr_G:=(\pr_{G_n})_{n\in N}:G'\ra G$. Alors $\eta_{G'}\in G'(k(\eta_{G'}))$ est identifié à $\pr_{G}(\eta_{G'})=\pr_{G_{\eta}}^*(\eta_G)\in G(k(\eta_{G'}))$ 
et $f_{\cT''}'(\eta_{G'})\in H'(k(\eta_{G'}))$ à $f_{\cT'}(\pr_{G_{\eta}}^*(\eta_{G}))\in H(k(\eta_{G'}))$. Comme $f_{\cT'}:G_{\cT'}\ra H_{\cT'}$ est un morphisme de préfaisceaux,
$$
f_{\cT'}(\pr_{G_{\eta}}^*(\eta_G))= \pr_{G_{\eta}}^*(f_{\cT'}(\eta_G)) = \pr_{G_{\eta}}^*(f(\eta_G)) = f(\pr_{G_{\eta}}^*(\eta_G)).
$$
Mais par commutativité de 
$$
\xymatrix{
\Spec(k(\eta_{G'})) \ar[d]_{\pr_{G_{\eta}}} \ar[rr]^{p_{G_{n}'}} & & \Spec(k(\eta_{G_{n}'})) \ar[d]^{\pr_{G_{\eta,n}}} \ar[rr]^>>>>>>>>>>{\eta_{G_{n}'}} & & G_{n}' \ar[d]^{\pr_{G_n}} \ar[r]^{f_{n}'} & H_{m_j}' \ar[d]^{\pr_{H_{m_j}}}\\
\Spec(k(\eta_G)) \ar[rr]^{p_{G_n}} & & \Spec(k(\eta_{G_n})) \ar[rr]^>>>>>>>>>>{\eta_{G_n}} & & G_n \ar[r]^{f_n} & H_{m_j}
}
$$
pour tout $m_j\in M_j$, avec $n:=\varphi(m_j)$, on a $f(\pr_{G_{\eta}}^*(\eta_G))=\pr_H(f'(\eta_{G'}))\in H(k(\eta_{G'}))$, identifié à $f'(\eta_{G'})\in H'(k(\eta_{G'})).$ On a donc bien $f_{\cT''}'(\eta_{G'})=f'(\eta_{G'})\in H'(k(\eta_{G'}))$, et on est ainsi ramené à vérifier que les morphismes $G(k')\rra H(k')$ induits par $f$ et $f_{\cT'}$ coïncident lorsque $k'=k$.

Soit $x=(x_n)_{n\in N}\in G(k)=\lp_{n\in N}G_n(k)$. Posons $\tau_x:=(\tau_{x_n})_{n\in N}:G\ra G$, et notons $x_{k(\eta_G)}$ l'image de $x$ dans $G(k(\eta_G))$. On a alors $\tau_x(\eta_G)=x_{k(\eta_G)}+\eta_G$, d'où
$$
f(x)=f_{\cT'}(x)\Longleftrightarrow f(x_{k(\eta_G)})=f_{\cT'}(x_{k(\eta_G)})\Longleftrightarrow f(\tau_x(\eta_G))=f_{\cT'}
(\tau_x(\eta_G))
$$
(la flèche de restriction $H(k)\ra H(k(\eta_G))$ est injective, et $f$ et $f_{\cT'}$ sont des morphismes de groupes qui coïncident en $\eta_G$). Ceci étant, $\tau_{x_n}$ induit un morphisme $\tau_{x_n\eta}$ rendant commutatif le diagramme
$$
\xymatrix{
\Spec(k(\eta_{G_n})) \ar[d]_{\tau_{x_n\eta}} \ar[r]^<<<<{\eta_{G_n}} & G_n \ar[d]^{\tau_{x_n}} \\
\Spec(k(\eta_{G_n})) \ar[r]^<<<<{\eta_{G_n}} & G_n,
}
$$
i.e. $\tau_{x_n}(\eta_{G_n})=\tau_{x_n\eta}^*(\eta_{G_n})$. Pour tout $n\in N$, on a alors un diagramme commutatif
$$
\xymatrix{
\Spec(k(\eta_{G})) \ar[d]_{\tau_{x\eta}:=\lp_n\tau_{x_n\eta}} \ar[rr]^{p_{G_{n}}} & & \Spec(k(\eta_{G_{n}})) \ar[d]^{\tau_{x_n\eta}} \ar[rr]^>>>>>>>>>>{\eta_{G_{n}}} & & G_{n} \ar[d]^{\tau_{x_n}} \\
\Spec(k(\eta_{G})) \ar[rr]^{p_{G_n}} & & \Spec(k(\eta_{G_n})) \ar[rr]^>>>>>>>>>>{\eta_{G_n}} & & G_n,
}
$$
i.e. $\tau_{x_n}(p_{G_n}^*(\eta_{G_n}))=\tau_{x\eta}^*(p_{G_n}^*(\eta_{G_n}))$ pour tout $n\in N$, ou encore $\tau_x(\eta_G)=\tau_{x\eta}^*(\eta_G)$. Comme $f_{\cT'}:G_{\cT'}\ra H_{\cT'}$ est un morphisme de préfaisceaux, 
$$
f_{\cT'}(\tau_{x\eta}^*(\eta_G))=\tau_{x\eta}^*(f_{\cT'}(\eta_G))=\tau_{x\eta}^*(f(\eta_G))=f(\tau_{x\eta}^*(\eta_G)),
$$
ce qui montre que $f(x)=f_{\cT'}(x)\in H(k)$. 

Pour achever la démonstration de la proposition, il reste à prouver l'assertion de surjectivité lorsque $I$ est un singleton, mais sans supposer les $G_n$ connexes. Soit dans ce cas $f_{\cT'}\in\Hom_{\widehat{\cT'}}(G_{\cT'},H_{\cT'})$. Soit $G_{\cT'}^0$ la restriction à la catégorie $\cT'$ du foncteur de points de la composante neutre $G^0:=(G_{n}^0)_{n\in N}$ de $G$. Il s'agit d'un sous-foncteur en groupes de $G_{\cT'}$ et, d'après ce qui précède, $f_{\cT'}^0:=f_{\cT'}|_{G_{\cT'}^0}$ provient d'un morphisme de ind-pro-$k$-schémas en groupes $f^0:G^0\ra H$, induit par un morphisme de pro-$k$-schémas en groupes $f^{0,j}:G^0\ra H_j$, que l'on peut représenter par un système compatible de morphismes de $k$-schémas en groupes $(f_{\varphi(m_j),m_j}^{0,j}:G_{\varphi(m_j)}^0\ra H_{m_j})_{m_j\in M_j}$.

Comme $k$ est algébriquement clos et $N$ essentiellement dénombrable, on peut choisir une section ensembliste $\{x_c\}_{c=0,\ldots,r-1}=\{(x_{c,n})_{n\in N}\}_{0=1,\ldots,r-1}$ de la flèche naturelle
$$
\xymatrix{
G(k)=\lp_{n\in N} G_n(k) \ar[r] &  \lp_{n\in N}\pi_0(G_n).
}
$$ 
Pour tous $c,c'\in\{0,\ldots,r-1\}$, on a une relation
$$
x_{c}+x_{c'}=x_{c+c'}+\delta_{c,c'},
$$
avec $\delta_{c,c'}\in G^0(k)$. Quitte à augmenter $j\in J$, on peut donc supposer que pour tout $c=0,\ldots,r-1$, $f_{\cT'}(x_c)$ est représentable par un élément de $H_j(k)$, et que l'on a dans $H_j(k)$ la relation 
$$
f_{\cT'}(x_c)+f_{\cT'}(x_{c'})=f_{\cT'}(x_{c+c'})+f^{0,j}(\delta_{c,c'}).
$$
On définit alors un morphisme de pro-$k$-schémas $f^{c,j}=(f_{\varphi(m_j),m_j}^{c,j})_{m_j\in M_j}$ par la composée
$$
\xymatrix{
G^c \ar[r]^{\tau_{-x_{c}}} & G^0 \ar[r]^{f^{0,j}} & H_j \ar[r]^{\tau_{f_{\cT'}(x_c)}} & H_j,
}
$$
puis on pose
$$
\xymatrix{
f^j=(f_{\varphi(m_j),m_j}^j)_{m_j\in M_j}:=(\coprod_{c=1}^r f_{\varphi(m_j),m_j}^{c,j})_{m_j\in M_j}:G\ar[r] & H_j.
}
$$
Montrons que les morphismes de $k$-schémas $f_{\varphi(m_j),m_j}^j:G_{\varphi(m_j)}\ra H_{m_j}$, $m_j\in M_j$, sont des morphismes de groupes. Comme pour tout $m_j\in M_j$, $H_{m_j}\times_k H_{m_j}$ est séparé et $G_{\varphi(m_j)}\times_k G_{\varphi(m_j)}$ parfait algébrique sur le corps algébriquement clos $k$, il suffit de montrer que les $f_{\varphi(m_j),m_j}^j(k):G_{\varphi(m_j)}(k)\ra H_{m_j}(k)$, $m_j\in M_j$, sont des morphismes de groupes. De plus, $N$ étant essentiellement dénombrable, les projections canoniques $G(k)\ra G_n(k)$, $n\in N$, sont surjectives. D'après les diagrammes commutatifs 
$$
\xymatrix{
G(k)\ar[rr]^{f^j} \ar[d] && H_j(k) \ar[d] \\
G_{\varphi(m_j)}(k) \ar[rr]^{f_{\varphi(m_j),m_j}^j(k)} && H_{m_j}(k),& m_j\in M_j,
}
$$
il suffit donc de montrer que $f^j:G(k)\ra H_j(k)$ est un morphisme de groupes. Soient $x\in G^c(k)$ et $x'\in G^{c'}(k)$. Ecrivons
$$
x=x_c+\delta,\quad x'=x_{c'}+\delta',
$$
avec $\delta,\delta'\in G^0(k)$. Alors, dans $H_j(k)$ :
\begin{eqnarray*}
f^j(x+x') & = & f^j(x_c+x_{c'}+\delta+\delta') \\
             & = & f^{c+c',j}(x_{c+c'}+\delta_{c,c'}+\delta+\delta')  \\
             & = & f_{\cT'}(x_{c+c'})+f^{0,j}(\delta_{c,c'}+\delta+\delta') \\
             & = & f_{\cT'}(x_c)+f_{\cT'}(x_{c'})-f^{0,j}(\delta_{c,c'})+f^{0,j}(\delta_{c,c'})+f^{0,j}(\delta)+f^{0,j}(\delta') \\
             & = & f^{c,j}(x_c+\delta)+f^{c',j}(x_{c'}+\delta') \\
             & = & f^j(x)+f^j(x').
\end{eqnarray*}
Ainsi donc, $f^j:G\ra H_j$ est un morphisme de pro-$k$-schémas en groupes. On définit alors un morphisme de ind-pro-$k$-schémas en groupes $f:G\ra H$ par la composée
$$
\xymatrix{
G \ar[r]^{f^j} & H_j \ar[r] & H.
}
$$
Montrons que si $k'/k$ est une extension parfaite du corps $k$, alors les morphismes $G(k')\rra H(k')$ induits par $f$ et $f_{\cT'}$ coïncident. Soit $x\in G^c(k')$, et écrivons $x=x_c+\delta$ avec $\delta\in G^0(k')$. Alors, dans $H(k')$ :
\begin{eqnarray*}
f(x) & = & f^j(x_c+\delta) \\
             & = & f_{\cT'}(x_c)+f^{0,j}(\delta)  \\
             & = & f_{\cT'}(x_c)+f_{\cT'}(\delta)  \\
             & = &  f_{\cT'}(x_c+\delta) \\
             & = & f_{\cT'}(x).
\end{eqnarray*}
\end{proof}

\begin{Cor}\label{ratIP0}
Soient $(G_1,S_1)\in\cI(\cP^{\pizf})$ et $(G_2,S_2)\in\cI\cP$ deux groupes indproalgébriques. Alors l'application de passage à la fibre générique (\ref{deffibgenIP}) 
$$
\xymatrix{
\Hom_{\cI\cP}((G_1,S_1),(G_2,S_2)) \ar[r] & \Hom_{\widehat{\cT}}((G_1,S_1)_{\eta},(G_2,S_2)_{\eta})
}
$$
est bijective.
\end{Cor}

\begin{proof}
D'après \ref{IPexacte} et \ref{Pab}, on a
$$
\xymatrix{
\Hom_{\cI\cP}((G_1,S_1),(G_2,S_2)) \ar[r]^<<<<<{\sim} & \Hom_{\Ind(\Pro(\cQ))}((G_1,S_1),(G_2,S_2)).
}
$$
D'après \ref{ratIndPro0}, on a ensuite 
$$
\xymatrix{
\Hom_{\Ind(\Pro(\cQ))}((G_1,S_1),(G_2,S_2)) \ar[r]^{\sim} & \Hom_{\widehat{\cT'}}((G_1,S_1)_{\cT'},(G_2,S_2)_{\cT'}).
}
$$
Par ailleurs, l'application de restriction à $\cT$
$$
\xymatrix{
\Hom_{\widehat{\cT'}}((G_1,S_1)_{\cT'},(G_2,S_2)_{\cT'}) \ar[r]^{\sim} & \Hom_{\widehat{\cT}}((G_1,S_1)_{\eta},(G_2,S_2)_{\eta})
}
$$
est bijective. Soit en effet $f_{\eta}\in \Hom_{\widehat{\cT}}((G_1,S_1)_{\eta},(G_2,S_2)_{\eta})$, et soit $k'/k$ une extension de corps parfaite. Choisissons une clôture algébrique $\kappa/k'$. Comme le morphisme de groupes 
$$
\xymatrix{
f_{\eta}(\kappa):(G_1,S_1)(\kappa)\ar[r] & (G_2,S_2)(\kappa)
}
$$ 
est naturel par rapport aux $k$-endomorphismes de $\kappa$ par hypothèse, il est en particulier $\Gal(\kappa/k')$-équivariant, induisant donc un morphisme de groupes 
$$
\xymatrix{
f_{\eta}(k'):=f_{\eta}(\kappa)^{\Gal(\kappa/k')}:(G_1,S_1)(\kappa)^{\Gal(\kappa/k')}\ar[r] & (G_2,S_2)(\kappa)^{\Gal(\kappa/k')}.
}
$$
Comme $(G,S)(\kappa)^{\Gal(\kappa/k')}=(G,S)(k')$ pour tout $(G,S)\in\cI\cP$, on obtient que $f_{\eta}(k')$ est un morphismes de groupes $(G_1,S_1)(k')\ra (G_2,S_2)(k')$. Faisant varier $k'/k$, on obtient une transformation naturelle $f:(G_1,S_1)_{\cT'}\ra(G_2,S_2)_{\cT'}$, dont la restriction à $\cT$ égale $f_{\eta}$. De plus $f\in\Hom_{\widehat{\cT'}}((G_1,S_1)_{\cT'},(G_2,S_2)_{\cT'})$ est uniquement déterminé par sa restriction $f_{\eta}$ puisque pour toute clôture algébrique $\kappa/k'$ comme ci-dessus, les flèches verticales du diagramme commutatif 
$$
\xymatrix{
(G_1,S_1)(\kappa)\ar[r]^{f_{\eta}(\kappa)} & (G_2,S_2)(\kappa) \\
(G_1,S_1)(k') \ar[u] \ar[r]^{f(k')} \ar[r]^{f(k')} & (G_2,S_2)(k') \ar[u]
}
$$
sont injectives. D'où finalement le résultat.
\end{proof}

\begin{Rem}\label{RemratIP0}
Au cours de la démonstration de \ref{ratIP0}, on a montré que pour tous groupes indproalgébriques $(G_1,S_1)$, 
$(G_2,S_2)$, l'application de restriction à $\cT$
$$
\xymatrix{
\Hom_{\widehat{\cT'}}((G_1,S_1)_{\cT'},(G_2,S_2)_{\cT'}) \ar[r]^{\sim} & \Hom_{\widehat{\cT}}((G_1,S_1)_{\eta},(G_2,S_2)_{\eta})
}
$$
est bijective.
\end{Rem}

\section{Groupes strictement connexes}

\begin{Pt}
Il existe des groupes proalgébriques connexes $(G,S)$ tels que pour tout $H,H'\in S$ avec $H\varsubsetneq H'\varsubsetneq G$, l'épimorphisme $G/H\ra G/H'$ de $\cQ^0$ n'est pas strict. Par exemple, on peut considérer une suite indexée par $\bbN$ d'isogénies d'Artin-Schreier
$$
\xymatrix{
\cdots \ar[r] & G_{n+1}:=\bbG_{a,k}^{\pf} \ar[r] & G_n:=\bbG_{a,k}^{\pf} \ar[r] & \cdots \ar[r] & G_0:=\bbG_{a,k}^{\pf}, & n\in\bbN.
}
$$
Posant $G:=\lp_{n\in\bbN}G_n(k)$ et $K_n:=\Ker(G\ra G_n(k))$, on a $G/K_n\xrightarrow{\sim} G_n(k)$ (\cite{S2} 2.3 Proposition 2), d'où un groupe proalgébrique $(G,S)$ avec $S:=\ <K_n,\ n\in\bbN>$. Si $H,H'\in S$ avec $H\varsubsetneq H'\varsubsetneq G$, choisissons $n\in\bbN$ tel que $K_n\subset H\cap H'$. Alors $H/K_n$ et $H'/K_n$ sont deux sous-groupes fermés de $G/K_n\simeq \bbG_{a,k}^{\pf}(k)$ distincts et distincts de $G/K_n$, de sorte que $(H'/K_n)/(H/K_n)$ est un groupe fini non nul. Le noyau $H'/H=(H'/K_n)/(H/K_n)$ de l'épimorphisme $G/H\ra G/H'$ n'est donc pas connexe.
\end{Pt}

\begin{Def}\label{DefP00}
Un \emph{groupe proalgébrique strictement connexe} est un groupe proalgébrique connexe $(G,S)$ vérifiant la condition suivante :

\medskip

$(sc)$ Il existe un ensemble de définition $S'$ de $(G,S)$ tel que pour tout $H_{1}',H_{2}'\in S'$ avec $H_{1}'\subset H_{2}'$, l'épimorphisme canonique $G/H_{1}'\ra G/H_{2}'$ de $\cQ^0$ est strict.
\end{Def}

\begin{Pt}\label{PtDefP00}
L'épimorphisme canonique $G/H_{1}'\ra G/H_{2}'$ de $\cQ^0$ est strict si et seulement si le groupe quasi-algébrique noyau $H_{2}'/H_{1}'$ est connexe.
\end{Pt}

\begin{Lem}\label{PS^0}
Si $(G,S)$ est un groupe proalgébrique strictement connexe et $S'$ un ensemble de définition de $(G,S)$ vérifiant la condition $(sc)$, alors tout élément de $S'$ est un groupe proalgébrique connexe. L'ensemble des sous-groupes de définition connexes de $(G,S)$ est alors un ensemble de définition de $(G,S)$.

Réciproquement, si $(G,S)$ est un groupe proalgébrique connexe possédant un ensemble de définition $S'$ dont les éléments sont connexes, alors $S'$ vérifie la condition $(sc)$, et en particulier $(G,S)$ est strictement connexe.
\end{Lem}

\begin{proof}
Soient $(G,S)$ un groupe proalgébrique strictement connexe, $S'$ vérifiant $(sc)$ et $H\in S'$. L'ensemble $S_{H}':=\{H'\in S'\ |\ H'\subset H\}$ est un sous-ensemble coinitial de l'ensemble complet de définition du groupe proalgébrique $H$ ; c'est donc un ensemble de définition de celui-ci. Par hypothèse (\ref{PtDefP00}), le groupe quasi-algébrique $H/H'$ est connexe pour tout $H'\in S_{H}'$. Il en résulte que $H$ est connexe (\ref{PtDefP0}). L'ensemble des sous-groupes de définition connexes de $(G,S)$ contient ainsi l'ensemble de définition $S'$, c'est donc lui aussi un ensemble de définition de $(G,S)$.

Réciproquement, soit $(G,S)$ un groupe proalgébrique connexe possédant un ensemble de définition $S'$ dont les éléments sont connexes. Alors pour tout $H',H\in S'$ avec $H'\subset H$, $H'$ est de définition dans $H$ et $H/H'$ est donc un groupe quasi-algébrique connexe.
\end{proof}

\begin{Not}
On note $\cP^{00}$ la sous-catégorie pleine $\cP^0$ dont les objets sont les groupes proalgébriques strictement connexes.
\end{Not}

\begin{Def} \label{DefIP00}
Un \emph{groupe indproalgébrique strictement connexe} est un groupe indproalgébrique $(G,S)$ vérifiant la condition suivante :

\medskip

$(sc)$ Il existe un ensemble de définition $S'$ de $(G,S)$ tel que pour tout $H'\in S'$, le groupe proalgébrique $H'$ est strictement connexe.
\end{Def}

\begin{Pt}\label{PtDefIP00}
Un groupe indproalégbrique est strictement connexe si et seulement si l'ensemble de ses sous-groupes de définition strictement connexes est un ensemble de définition.  
\end{Pt}

\begin{Not}
On note $\cI(\cP^{00})$ la sous-catégorie pleine de $\cI(\cP^0)$ dont les objets sont les groupes indproalgébriques strictement connexes.
\end{Not}

\begin{Def}\label{DefPI00}
Un \emph{groupe proindalgébrique strictement connexe} est un groupe proindalgébrique connexe $(G,S)$ vérifiant la condition suivante :

\medskip

$(sc)$ Il existe un ensemble de définition $S'$ de $(G,S)$ tel que pour tout $H_{1}',H_{2}'\in S'$ avec $H_{1}'\subset H_{2}'$, l'épimorphisme canonique $G/H_{1}'\ra G/H_{2}'$ de $\cI^0$ est strict.
\end{Def}

\begin{Pt}\label{PtDefPI00}
L'épimorphisme canonique $G/H_{1}'\ra G/H_{2}'$ de $\cI$ est strict par (PI-3) ; il est strict dans $\cI^0$ si et seulement si le groupe indalgébrique noyau $H_{2}'/H_{1}'$ est connexe.
\end{Pt}

\begin{Lem}\label{PIS^0}
Si $(G,S)$ est un groupe proindalgébrique strictement connexe et si $S'$ est un ensemble de définition de $(G,S)$ vérifiant la condition $(sc)$, alors tout élément de $S'$ est un groupe proindalgébrique connexe. L'ensemble des sous-groupes de définition connexes de $(G,S)$ est alors un ensemble de définition de $(G,S)$.

Réciproquement, si $(G,S)$ est un groupe proindalgébrique connexe possédant un ensemble de définition $S'$ dont les éléments sont connexes, alors $S'$ vérifie la condition $(sc)$, et en particulier $(G,S)$ est strictement connexe.
\end{Lem}

\begin{proof}
Soient $(G,S)$ un groupe proindalgébrique strictement connexe, $S'$ vérifiant $(sc)$ et $H\in S'$. L'ensemble $S_{H}':=\{H'\in S'\ |\ H'\subset H\}$ est un sous-ensemble coinitial de l'ensemble complet de définition du groupe proindalgébrique $H$ ; c'est donc un ensemble de définition de celui-ci. Par hypothèse (\ref{PtDefPI00}), le groupe indalgébrique $H/H'$ est connexe pour tout $H'\in S_{H}'$. Il en résulte que $H$ est connexe (\ref{PtDefPI0}). L'ensemble des sous-groupes de définition connexes de $(G,S)$ contient ainsi l'ensemble de définition $S'$, c'est donc lui aussi un ensemble de définition de $(G,S)$.

Réciproquement, soit $(G,S)$ un groupe proindalgébrique connexe possédant un ensemble de définition $S'$ dont les éléments sont connexes. Alors pour tout $H',H\in S'$ avec $H'\subset H$, $H'$ est de définition dans $H$ et $H/H'$ est donc un groupe indalgébrique connexe.
\end{proof}

\begin{Not}
On note $\cP^{00}(\cI)$ la sous-catégorie pleine de $\cP^0(\cI)$ dont les objets sont les groupes proindalgébriques strictement connexes.
\end{Not}

\begin{Notss}
On note $\cI(\cP^{00})^{\ad}$ (resp. $\cP(\cI^{00})^{\ad})$ la sous-catégorie pleine de $\cI\cP^{\ad}$ (resp. $\cP\cI^{\ad}$) dont les objets sont ceux de $\cI\cP^{\ad}$ (resp. $\cP\cI^{\ad}$) qui sont strictement connexes.  
\end{Notss}

\begin{Pt} \label{foncIPPI00}
Soit $(G,S)\in \cI(\cP^{00})^{\ad}$. Comme $S^{00}$ est cofinal dans $S$ et l'ensemble des origines de $(G,S)$ stable par majorant, l'ensemble $S_{0}^{00}$ des origines strictement connexes est un ensemble de définition de $(G,S)$. Choisissons $H_0\in S_{0}^{00}$. Si $H\in S_{H_0}$, alors l'ensemble
$$
\{H'/H\ |\ H\subset H'\in S^{00}\}
$$
est de définition dans $(G/H,S_{G/H})\in\cI$ et ses éléments sont connexes, donc $(G/H,S_{G/H})\in\cI^0$. De plus, l'ensemble  des sous-groupes de définition connexes de $H_0$ est de définition dans $H_0$, donc aussi dans l'ensemble complet de définition $<S_{H_0}>$ de $\can(G,S)$. Le groupe proindalgébrique $\can(G,S)$ est donc strictement connexe. D'où un foncteur 
$$
\xymatrix{
\can:\cI(\cP^{00})^{\ad}\ar[r] & \cP^{00}(\cI)^{\ad}.
}
$$
\end{Pt}

\begin{Pt} \label{foncPIIP00}
Soit $(G,S)\in \cP^{00}(\cI)^{\ad}$. Comme $S^0$ est coinitial dans $S$ et l'ensemble des origines de $(G,S)$ stable par minorant, l'ensemble $S_{0}^0$ des origines connexes est un ensemble de définition de $(G,S)$. Choisissons $H_{0}\in S_{0}^0$. Si $H\in\tilde{S}_{G/H_0}$, alors l'ensemble
$$
\{H'\in S^0\ |\ H'\subset H\}
$$
est de définition dans $(H,S_H)$ et ses éléments sont connexes, donc $(H,S_H)\in\cP^{00}$ s'il est connexe. De plus, l'ensemble des éléments connexes de $\tilde{S}_{G/H_0}$ est cofinal dans $\tilde{S}_{G/H_0}$ : l'ensemble des sous-groupes de définition connexes de $(G/H_0,S_{G/H_0})\in\cI^0$ est de définition dans $(G/H_0,S_{G/H_0})$, et l'image réciproque dans $G$ d'un tel sous-groupe est connexe puisque $H_0$ est connexe (\ref{P0exacte}). Donc l'ensemble des éléments strictement connexes de $\tilde{S}_{G/H_0}$ est de définition dans l'ensemble complet de définition $<\tilde{S}_{G/H_0}>$ de $\can(G,S)$. Le groupe indproalgébrique $\can(G,S)$ est donc strictement connexe. D'où un foncteur 
$$
\xymatrix{
\can:\cP^{00}(\cI)^{\ad}\ar[r] & \cI(\cP^{00})^{\ad}.
}
$$
\end{Pt}

Ceci étant, le théorème \ref{canequiv} admet le corollaire suivant :

\begin{Cor}\label{canequiv00}
Les foncteurs $\can:\cI(\cP^{00})^{\ad}\ra\cP^{00}(\cI)^{\ad}$ \ref{foncIPPI00} et $\can:\cP^{00}(\cI)^{\ad}\ra\cI(\cP^{00})^{\ad}$ \ref{foncPIIP00} sont des équivalences de catégories inverses l'une de l'autre. Ils induisent des diagrammes commutatifs 
$$
\xymatrix{
\cI(\cP^{00})^{\ad}\ar[r]^{\can} & \cP^{00}(\cI)^{\ad} & & \cI(\cP^{00})^{\ad}\ar[r]^{\can} & \cP^{00}(\cI)^{\ad} \\
\cI^0 \ar[u]^{\ref{IetPdansIPf}} \ar@{=}[r] & \cI^0 \ar[u]_{\ref{IetPdansPIf}} & & \cP^{00} \ar[u]^{\ref{IetPdansIPf}} \ar@{=}[r] & \cP^{00} \ar[u]_{\ref{IetPdansPIf}}
}
$$
$$
\xymatrix{
\cP^{00}(\cI)^{\ad}\ar[r]^{\can} & \cI(\cP^{00})^{\ad} & & \cP^{00}(\cI)^{\ad}\ar[r]^{\can} & \cI(\cP^{00})^{\ad} \\
\cI^0 \ar[u]^{\ref{IetPdansPIf}} \ar@{=}[r] & \cI^0 \ar[u]_{\ref{IetPdansIPf}} & & \cP^{00} \ar[u]^{\ref{IetPdansPIf}} \ar@{=}[r] & \cP^{00}. \ar[u]_{\ref{IetPdansIPf}}
}
$$
\end{Cor}

\begin{Def}
Un groupe \emph{$\cQ$-admissible strictement connexe} est un groupe $\cQ$-admissible $(G,S)$ vérifiant la condition suivante :

\medskip

$(sc)$  Il existe un ensemble de définition $S'$ de $(G,S)$ tel que pour tout $H_{1}',H_{2}'\in S'$ avec $H_{1}'\subset H_{2}'$, le groupe quasi-algébrique $H_{2}'/H_{1}'$ est connexe.
\end{Def}

\begin{Not}
On note $\fQfzz$ la catégorie des groupes $\cQ$-admissibles strictement connexes.
\end{Not}

\begin{Pt} \label{fQfIPf00}
Soit $(G,S)\in\fQfzz$. Choisissons $S'$ vérifiant la condition $(sc)$. Pour tout $H'\in S'$, l'ensemble de sous-groupes de $H'$
$$
S_{H'}':=\{H_{-}'\in S' \ |\ H_{-}'\subset H'\}
$$
est de définition dans le groupe proalgébrique $(H',S_{H'})$. Celui-ci est donc connexe puisque $H'/H_{-}'$ est connexe pour tout $H_{-}'\in S_{H'}'$. En particulier, pour chaque $H_{-}'\in S_{H'}'$, le groupe proalgébrique $(H_{-}',S_{H_{-}'})$ est connexe. Par conséquent $(H',S_{H'})$ est strictement connexe. Comme $S'$ est cofinal dans $S$, on en déduit que $\indpro(G,S)=(G,<S>)\in\cI\cP^{\ad}$ est strictement connexe. D'où un foncteur 
$$
\xymatrix{
\indpro:\fQfzz\ar[r] & \cI(\cP^{00})^{\ad}.
}
$$
\end{Pt}

\begin{Pt} \label{fQfPIf00}
Soit $(G,S)\in\fQfzz$. Choisissons $S'$ vérifiant la condition $(sc)$. Pour tout $H'\in S'$, l'ensemble de sous-groupes de $G/H'$
$$
S_{G/H'}':=\{H_{+}'/H'\in S' \ |\ H'\subset H_{+}'\in S'\}
$$
est de définition dans le groupe indalgébrique $(G/H',S_{G/H'})$. Celui-ci est donc connexe puisque $H_{+}'/H'$ est connexe pour tout $H_{+}'\in S'$ contenant $H'$. De plus, chaque groupe proalgébrique $H'\in S'$ est connexe (cf. \ref{fQfIPf00}).  Comme $S'$ est coinitial dans $S$, on en déduit que $\proind(G,S)=(G,<S>)\in\cP\cI^{\ad}$ est strictement connexe. D'où un foncteur 
$$
\xymatrix{
\proind:\fQfzz\ar[r] & \cP^{00}(\cI)^{\ad}.
}
$$
\end{Pt}

\begin{Pt} \label{finfoncfQf00}
Soit $(G,S)\in\cI(\cP^{00})^{\ad}$. Alors pour toutes origines strictement connexes $H_0\subset H_{0}'\in S$, on a $H_{0}'/H_0\in\cQ^0$. L'ensemble $S_{0}^{00}$ des origines strictement connexes de $(G,S)$ est cofinal dans l'ensemble $S_0$ des origines de $(G,S)$. L'ensemble $S_{0}^0$ des origines connexes de $\can(G,S)\in\cP^{00}(\cI)^{\ad}$ est coinitial dans l'ensemble des origines de $\can(G,S)$ qui est $S_0$ (\ref{canS0}), et les éléments de $S_{0}^0$ sont strictement connexes (\ref{PIS^0}). Par conséquent $S_{0}^{00}=S_{0}^0$ est cofinal et coinitial dans $S_0$, et $\ad(G,S):=(G,S_0)=(G,<S_{0}^{00}>)\in\fQfzz$. D'où un foncteur  
$$
\xymatrix{
\ad:\cI(\cP^{00})^{\ad}\ar[r] & \fQfzz.
}
$$
On obtient de même un foncteur  
$$
\xymatrix{
\ad:\cP(\cI^{00})^{\ad}\ar[r] & \fQfzz.
}
$$
\end{Pt}

Ceci étant, le théorème \ref{fQfIPfPIf} admet le corollaire suivant :

\begin{Cor} \label{fQfIPfPIf00}
Les foncteurs $\indpro:\fQfzz\ra\cI(\cP^{00})^{\ad}$  \ref{fQfIPf00} et $\ad:\cI(\cP^{00})^{\ad}\ra\fQfzz$ \ref{finfoncfQf00} sont des équivalences de catégories inverses l'une de l'autre. Les foncteurs $\proind:\fQfzz\ra\cP^{00}(\cI)^{\ad}$ \ref{fQfPIf00} et 
$\ad:\cP^{00}(\cI)^{\ad}\ra\fQfzz$ \ref{finfoncfQf00} sont des équivalences de catégories inverses l'une de l'autre. Le diagramme
$$
\xymatrix{
& & \fQfzz\ar@<-1ex>[ddll]_{\indpro} \ar@<1ex>[ddrr]^{\proind} & & \\
\\
\cI(\cP^{00})^{\ad} \ar@<1ex>[rrrr]^{\can} \ar@<-1ex>[uurr]_{\ad} & & & & \cP(\cI^{00})^{\ad}  \ar@<1ex>[llll]^{\can} \ar@<1ex>[uull]^{\ad}
}
$$
est commutatif. 
\end{Cor}

\section{Dualité de Serre pour les groupes unipotents}

\begin{Th} \emph{(Bégueri \cite{Beg} Proposition 1.2.1).} \label{dualSU0}
 Soit $G\in\cU^0$. Le faisceau 
$$\uExt_{(\parf/k)_{\et}}^1(G^{\pf},\bbQ_p/\bbZ_p)$$ 
sur le site des $k$-schémas parfaits muni de la topologie étale est représentable par le perfectisé d'un $k$-schéma en groupes unipotent connexe ; l'objet $G^{\vee_S}\in\cU^0$ correspondant est le \emph{dual de Serre} de $G$. Le foncteur $(\cdot)^{\vee_S}:(\cU^0)^{\circ}\ra\cU^0$ ainsi défini est exact et involutif.
\end{Th}

\begin{Pt} \label{S}
Avec Bégueri \cite{Beg} 1.3, on note encore $(\cdot)^{\vee_S}$ le foncteur composé
$$
\xymatrix{
\cU^{\circ}\ar[r]^{(\cdot)^0} & (\cU^0)^{\circ} \ar[r]^{(\cdot)^{\vee_S}} \ar[r] & \cU. 
}
$$ 
Il est exact à droite, et la sous-catégorie $(\cU^0)^{\circ}\subset (\cU)^{\circ}$  est $(\cdot)^{\vee_S}$-projective (au sens de \cite{KS} 13.3.4) ; il admet donc un dérivé gauche $L(\cdot)^{\vee_S}:D^+(\cU)^{\circ}\ra D^-(\cU)$. Si $G\in\cU$, on a par construction :
$$
\left\{ \begin{array}{ll}
L^i(G)^{\vee_S}=0 & \textrm{si $i\neq -1,0$}, \\
L^0(G)^{\vee_S}=G^{\vee_S}=(G^0)^{\vee_S}, &  \\
L^{-1}(G)^{\vee_S}=\li_{n\in\bbN}\Hom_{\cU}(G,\bbZ/p^n)=\li_{n\in\bbN}\Hom_{\cU}(\pi_0(G),\bbZ/p^n)=:\pi_0(G)^{\vee_P}. & 
\end{array} \right.
$$
\end{Pt}

\begin{Pt} \label{IS}
D'après \cite{KS} 15.3.7, la sous-catégorie $\Ind((\cU^0)^{\circ})=\Pro(\cU^0)^{\circ}\subset \Ind(\cU^{\circ})=\Pro(\cU)^{\circ}$ est $\Ind(\cdot)^{\vee_S}$-projective, où
$$
\Ind(\cdot)^{\vee_S}:\Pro(\cU)^{\circ}=\Ind(\cU^{\circ})\lra\Ind(\cU).
$$
Ce foncteur admet donc un dérivé gauche
$$
\xymatrix{
L\Ind(\cdot)^{\vee_S}:D^+(\Pro(\cU))^{\circ}\ar[r] & D^-(\Ind(\cU)),
}
$$
rendant commutatif le diagramme
$$
\xymatrix{
D^+(\cU)^{\circ} \ar[rr]^{L(\cdot)^{\vee_S}} \ar[d] && D^-(\cU) \ar[d] \\
D^+(\Pro(\cU))^{\circ}\ar[rr]^<<<<<<<<{L\Ind(\cdot)^{\vee_S}} && D^-(\Ind(\cU)),
}
$$
et pour tout $G\in\Pro(\cU)$, on a :
$$
\left\{ \begin{array}{ll}
L^i\Ind(G)^{\vee_S}=0 & \textrm{si $i\neq -1,0$}, \\
L^0\Ind(G)^{\vee_S}=\Ind(G)^{\vee_S}=\Ind(G^0)^{\vee_S}, &  \\
L^{-1}\Ind(G)^{\vee_S}=\Ind\pi_0(G)^{\vee_P}. & 
\end{array} \right.
$$
\end{Pt}

\begin{Pt}\label{PS}
La sous-catégorie $\Pro((\cU^0)^{\circ})=\Ind(\cU^0)^{\circ}\subset \Pro(\cU^{\circ})=\Ind(\cU)^{\circ}$ est $\Pro(\cdot)^{\vee_S}$-projective, où
$$
\Pro(\cdot)^{\vee_S}:\Ind(\cU)^{\circ}=\Pro(\cU^{\circ})\lra\Pro(\cU).
$$
Ce foncteur admet donc un dérivé gauche
$$
\xymatrix{
L\Pro(\cdot)^{\vee_S}:D^+(\Ind(\cU))^{\circ} \ar[r] & D^-(\Pro(\cU)),
}
$$
rendant commutatif le diagramme
$$
\xymatrix{
D^+(\cU)^{\circ} \ar[rr]^{L(\cdot)^{\vee_S}} \ar[d] && D^-(\cU) \ar[d] \\
D^+(\Ind(\cU))^{\circ}\ar[rr]^<<<<<<<<{L\Pro(\cdot)^{\vee_S}} && D^-(\Pro(\cU)),
}
$$
et pour tout $G\in\Ind(\cU)$, on a :
$$
\left\{ \begin{array}{ll}
L^i\Pro(G)^{\vee_S}=0 & \textrm{si $i\neq -1,0$}, \\
L^0\Pro(G)^{\vee_S}=\Pro(G)^{\vee_S}=\Pro(G^0)^{\vee_S}, &  \\
L^{-1}\Pro(G)^{\vee_S}=\Pro\pi_0(G)^{\vee_P}. & 
\end{array} \right.
$$
\end{Pt}

\begin{Pt}\label{IPS}
La sous-catégorie $\Ind\Pro((\cU^0)^{\circ})=\Pro\Ind(\cU^0)^{\circ}\subset \Ind\Pro(\cU^{\circ})=\Pro\Ind(\cU)^{\circ}$ est $\Ind\Pro(\cdot)^{\vee_S}$-projective, où
$$
\Ind\Pro(\cdot)^{\vee_S}:\Pro\Ind(\cU)^{\circ}=\Ind\Pro(\cU^{\circ})\lra\Ind\Pro(\cU).
$$
Ce foncteur admet donc un dérivé gauche
$$
\xymatrix{
L\Ind\Pro(\cdot)^{\vee_S}:D^+(\Pro\Ind(\cU))^{\circ}\ar[r] & D^-(\Ind\Pro(\cU)),
}
$$
rendant commutatif le diagramme
$$
\xymatrix{
D^+(\cU)^{\circ} \ar[rr]^{L(\cdot)^{\vee_S}} \ar[d] && D^-(\cU) \ar[d] \\
D^+(\Ind(\cU))^{\circ}\ar[rr]^{L\Pro(\cdot)^{\vee_S}} \ar[d] && D^-(\Pro(\cU)) \ar[d] \\
D^+(\Pro\Ind(\cU))^{\circ}\ar[rr]^<<<<<<<<<{L\Ind\Pro(\cdot)^{\vee_S}} && D^-(\Ind\Pro(\cU)),
}
$$
et pour tout $G\in\Pro\Ind(\cU)$, on a :
$$
\left\{ \begin{array}{ll}
L^i\Ind\Pro(G)^{\vee_S}=0 & \textrm{si $i\neq -1,0$}, \\
L^0\Ind\Pro(G)^{\vee_S}=\Ind\Pro(G)^{\vee_S}=\Ind\Pro(G^0)^{\vee_S}, &  \\
L^{-1}\Ind\Pro(G)^{\vee_S}=\Ind\Pro\pi_0(G)^{\vee_P}. & 
\end{array} \right.
$$
\end{Pt}

\begin{Pt} \label{PIS}
La sous-catégorie $\Pro\Ind((\cU^0)^{\circ})=\Ind\Pro(\cU^0)^{\circ}\subset \Pro\Ind(\cU^{\circ})=\Ind\Pro(\cU)^{\circ}$ est $\Pro\Ind(\cdot)^{\vee_S}$-projective, où
$$
\Pro\Ind(\cdot)^{\vee_S}:\Ind\Pro(\cU)^{\circ}=\Pro\Ind(\cU^{\circ})\lra\Pro\Ind(\cU).
$$
Ce foncteur admet donc un dérivé gauche
$$
\xymatrix{
L\Pro\Ind(\cdot)^{\vee_S}:D^+(\Ind\Pro(\cU))^{\circ}\ar[r] & D^-(\Pro\Ind(\cU)),
}
$$
rendant commutatif le diagramme
$$
\xymatrix{
D^+(\cU)^{\circ} \ar[rr]^{L(\cdot)^{\vee_S}} \ar[d] && D^-(\cU) \ar[d] \\
D^+(\Pro(\cU))^{\circ} \ar[rr]^{L\Ind(\cdot)^{\vee_S}} \ar[d] && D^-(\Ind(\cU)) \ar[d] \\
D^+(\Ind\Pro(\cU))^{\circ}\ar[rr]^<<<<<<<<<{L\Pro\Ind(\cdot)^{\vee_S}} && D^-(\Pro\Ind(\cU)),
}
$$
et pour tout $G\in\Ind\Pro(\cU)$, on a :
$$
\left\{ \begin{array}{ll}
L^i\Pro\Ind(G)^{\vee_S}=0 & \textrm{si $i\neq -1,0$}, \\
L^0\Pro\Ind(G)^{\vee_S}=\Pro\Ind(G)^{\vee_S}=\Pro\Ind(G^0)^{\vee_S}, &  \\
L^{-1}\Pro\Ind(G)^{\vee_S}=\Pro\Ind\pi_0(G)^{\vee_P}. & 
\end{array} \right.
$$
\end{Pt}

\section{Dualité de Serre pour les groupes unipotents strictement connexes}

\begin{Pt} \label{foncP00UI0U} 
Soit $(G,S)\in\cP^{00}(\cU)$. Posons $S^0:=\{H\in S\ |\ H\textrm{ connexe}\}$ (cf. \ref{PS^0}). Pour tout $H,H'\in S^0$ avec $H\subset H'$, l'épimorphisme strict $G/H\ra G/H'$ de $\cU^0$ induit un monomorphisme $(G/H')^{\vee_S}\ra (G/H)^{\vee_S}$ de $\cU^0$ (\ref{dualSU0}). Posons $G^{\vee_S}:=\li_{H\in S^0}(G/H)^{\vee_S}$. Pour tout $H\in S^0$, l'application canonique $(G/H)^{\vee_S}\ra G^{\vee_S}$ est injective, identifiant le groupe sous-jacent au groupe quasi-algébrique unipotent connexe $(G/H)^{\vee_S}$ à un sous-groupe de $G^{\vee_S}$. L'ensemble de sous-groupes de $G^{\vee_S}$
$$
(S^0)^{\vee_S}:=\{ (G/H)^{\vee_S}\subset G^{\vee_S}\ |\ H\in S^0\}
$$
vérifie alors (I-1), (I-2)' et (I-3), munissant donc $G^{\vee_S}$ d'une structure de groupe indalgébrique 
$(G^{\vee_S},S^{\vee_S})$, avec $S^{\vee_S}:=\ <(S^0)^{\vee_S}>$. On a $(G^{\vee_S},S^{\vee_S})\in\cI^{0}(\cU)$ par construction. De plus si $f:(G_1,S_1)\ra (G_2,S_2)\in \cP^{00}(\cU)$, alors pour tout $H_2\in S_{2}^0$, il existe $H_1\in S_{1}^0$ inclus dans $f^{-1}(H_2)$ tel que l'application $G_1/H_1\ra G_2/H_2$ induite par $f$ soit quasi-algébrique, et définisse donc un morphisme quasi-algébrique $(G_2/H_2)^{\vee_S}\ra (G_1/H_1)^{\vee_S}$ ; d'où un morphisme $f^{\vee_S}:(G_{2}^{\vee_S},S_{2}^{\vee_S})\ra (G_{1}^{\vee_S},S_{1}^{\vee_S})$, et finalement un foncteur 
$$
\xymatrix{
(\cdot)^{\vee_S}:\cP^{00}(\cU)^{\circ}\ar[r] & \cI^{0}(\cU).
}
$$
\end{Pt}

\begin{Pt} \label{foncI0UP00U} 
Soit $(G,S)\in\cI^0(\cU)$. Posons $S^0:=\{H\in S\ |\ H\textrm{ connexe}\}$ (cf. \ref{PtDefI0}). Pour tout $H,H'\in S^0$ avec $H\subset H'$, le monomorphisme $H\subset H'$ de $\cU^0$ est strict, et induit donc un épimorphisme strict $(H')^{\vee_S}\ra H^{\vee_S}$ de $\cU^0$ (\ref{dualSU0}). Posons $G^{\vee_S}:=\lp_{H\in S^0}H^{\vee_S}$. Pour tout $H\in S^0$, l'application canonique $G^{\vee_S}\ra H^{\vee_S}$ est alors surjective (\cite{S2} 2.3 Proposition 2) ; si l'on note $K_{H^{\vee_S}}$ son noyau, le quotient $G^{\vee_S}/K_{H^{\vee_S}}$ se trouve donc muni d'une structure de groupe quasi-algébrique unipotent connexe, canoniquement isomorphe à $H^{\vee_S}$. L'ensemble de sous-groupes de $G^{\vee_S}$
$$
(S^0)^{\vee_S}:=\{ K_{H^{\vee_S}}\subset G^{\vee_S}\ |\ H\in S^0\}
$$
vérifie alors (P-1), (P-3) et (P-4), munissant donc $G^{\vee_S}$ d'une structure de groupe proalgébrique 
$(G^{\vee_S},S^{\vee_S})$, avec $S^{\vee_S}:=\ <(S^0)^{\vee_S}>$. On a $(G^{\vee_S},S^{\vee_S})\in\cP^{00}(\cU)$ par construction. De plus si $f:(G_1,S_1)\ra (G_2,S_2)\in \cI^0(\cU)$, alors pour tout $H_1\in S_{1}^0$, il existe $H_2\in S_{2}^0$ contenant $f(H_1)$ tel que l'application $H_1\ra H_2$ induite par $f$ soit quasi-algébrique (on peut prendre $H_2=f(H_1)$),
et définisse donc un morphisme quasi-algébrique $H_{2}^{\vee_S}\ra H_{1}^{\vee_S}$ ; d'où un morphisme $f^{\vee_S}:(G_{2}^{\vee_S},S_{2}^{\vee_S})\ra (G_{1}^{\vee_S},S_{1}^{\vee_S})$, et finalement un foncteur 
$$
\xymatrix{
(\cdot)^{\vee_S}:\cI^{0}(\cU)^{\circ}\ar[r] & \cP^{00}(\cU).
}
$$
\end{Pt}

\begin{Th}\label{dualSP00UI0U}
Les foncteurs $(\cdot)^{\vee_S}:\cP^{00}(\cU)^{\circ}:\ra\cI^0(\cU)$ \ref{foncP00UI0U} et $(\cdot)^{\vee_S}:\cI^0(\cU)^{\circ}\ra\cP^{00}(\cU)$ \ref{foncI0UP00U} sont des équivalences de catégories quasi-inverses l'une de l'autre. De plus, ils sont exacts.
\end{Th}

\begin{proof}
Soit $(G,S)\in\cP^{00}(\cU)$. Avec les notations de \ref{foncP00UI0U}, on a 
$$(G,S)^{\vee_S}=(\li_{H\in S^0}(G/H)^{\vee_S},\ <(S^0)^{\vee_S}>).$$
Montrons que l'ensemble des sous-groupes de définition connexes de $(G,S)^{\vee_S}$ est réduit à l'ensemble $\{(G/H)^{\vee_S}\ |\ H\in S^0\}$. Soit $C$ un tel sous-groupe de $(G,S)^{\vee_S}$. Il existe $H\in S^0$ tel que $C\subset (G/H)^{\vee_S}$. L'inclusion $C\subset (G/H)^{\vee_S}$ est alors un monomorphisme de $\cU^0$, induisant un épimorphisme strict $G/H\ra C^{\vee_S}$ (\ref{dualSU0}). Son noyau est de la forme $H^0/H\in\cU^0$ pour un $H^0\in S$. Pour tout $H'\in S^0$ inclus dans $H$, le groupe quasi-algébrique $H^0/H'$ est alors le noyau de l'épimorphisme strict 
$G/H'\ra C^{\vee_S}$ de $\cU^0$ dual de $C\subset (G/H')^{\vee_S}$, et est donc connexe.
Il en résulte que $H^0$ est connexe, i.e. $H^0\in S^0$, et $C=(G/H^0)^{\vee_S}$ par construction. On a donc montré $<(S^0)^{\vee_S}>^0=\{(G/H)^{\vee_S}\ |\ H\in S^0\}$. Ceci étant, il est clair que $(G,S)^{\vee_S\vee_S}=(G,S)$, et que $(\cdot)^{\vee_S}\circ(\cdot)^{\vee_S}=\Id:\cP^{00}(\cU)\ra\cP^{00}(\cU)$ sur les morphismes aussi.

Soit $(G,S)\in\cI^0(\cU)$. Avec les notations de \ref{foncI0UP00U}, on a 
$$(G,S)^{\vee_S}=(\lp_{H\in S^0}H^{\vee_S},\ <(S^0)^{\vee_S}>).$$
Montrons que l'ensemble des sous-groupes de définition connexes de $(G,S)^{\vee_S}$ est réduit à l'ensemble $\{K_{H^{\vee_S}}\ |\ H\in S^0\}$. Soit $C$ un tel sous-groupe de $(G,S)^{\vee_S}$. Il existe $H\in S^0$ tel que $K_{H^{\vee_S}}\subset C$ et $C/K_{H^{\vee_S}}$ fermé dans $G/K_{H^{\vee_S}}=H^{\vee_S}$. D'où un épimorphisme strict $H^{\vee_S}\ra G/C$ de $\cU^0$, induisant un monomorphisme $(G/C)^{\vee_S}\ra H$ (\ref{dualSU0}). Par conséquent le groupe connexe $H^0:=(G/C)^{\vee_S}$ est de définition dans $(G,S)$, i.e. $H^0\in S^0$, et $K_{(H^0)^{\vee_S}}=C$ par construction. On a donc montré $<(S^0)^{\vee_S}>^0=\{K_{H^{\vee_S}}\ |\ H\in S^0\}$. Ceci étant, il est clair que $(G,S)^{\vee_S\vee_S}=(G,S)$, et que $(\cdot)^{\vee_S}\circ(\cdot)^{\vee_S}=\Id:\cI^0(\cU)\ra\cI^0(\cU)$ sur les morphismes aussi.

Les foncteurs $(\cdot)^{\vee_S}:\cP^{00}(\cU)^{\circ}:\ra\cI^0(\cU)$ et $(\cdot)^{\vee_S}:\cI^0(\cU)^{\circ}\ra\cP^{00}(\cU)$ sont exacts d'après \ref{IS} et \ref{PS}.
\end{proof}

\begin{Rem}\label{RemdualS2} 
Soit $(G,S)\in\cP^{00}(\cU)$. Si dans la construction \ref{foncP00UI0U} de $(G,S)^{\vee_S}\in\cI^0(\cU)$, on remplace l'ensemble $S^0$ des sous-groupes de définition connexes de $(G,S)$ par un ensemble de définition $S'$ dont les éléments sont connexes, alors le groupe indalgébrique obtenu 
$$(\li_{H\in S'}(G/H)^{\vee_S},\ <(S')^{\vee_S}>)$$
est identique à $(G,S)^{\vee_S}$, puisque 
$$
\{ (G/H)^{\vee_S}\subset \li_{H\in S^0}(G/H)^{\vee_S}\ |\ H\in S'\}
$$
est cofinal dans
$$
(S^0)^{\vee_S}:=\{ (G/H)^{\vee_S}\subset \li_{H\in S^0}(G/H)^{\vee_S}\ |\ H\in S^0\}.
$$
De même, si $(G,S)\in\cI^0(\cU)$, alors dans la construction \ref{foncI0UP00U} de $(G,S)^{\vee_S}\in\cP^{00}(\cU)$, on peut remplacer l'ensemble $S^0$ des sous-groupes de définition connexes de $(G,S)$ par un ensemble de définition dont les éléments sont connexes.
\end{Rem}

\begin{Notss}
Pour tout $n\in\bbN$, on note $\cU(n)$ la catégorie abélienne des groupes quasi-algébriques unipotents annulés par $p^n$. On note $\cP^{00}(\cU(n))$ (resp. $\cI^0(\cU(n))$) la sous-catégorie pleine de $\cP^{00}(\cU)$ (resp. 
$\cI^0(\cU)$) dont les objets sont ceux de $\cP^{00}(\cU)$ (resp. $\cI^0(\cU)$) tels que $(G/H)\in\cU(n)$ (resp. $H\in\cU(n)$) pour tout $H\in S$.
\end{Notss}

\begin{Pt}\label{ZZpn}
Si l'on note $\widetilde{(\parf/k)}_{\et}(n)$ la catégorie des faisceaux abéliens sur ${(\parf/k)}_{\et}$ annulés par $p^n$, alors d'après \cite{Ber} 2.5 (i), on a pour tout $G\in\cU(n)$,
$$
\xymatrix{
G^{\vee_S}:=\uExt_{(\parf/k)_{\et}}^1(G^{\pf},\bbQ_p/\bbZ_p)=\uExt_{\widetilde{(\parf/k)}_{\et}(n)}^1(G^{\pf},\bbZ/p^n).
}
$$
\end{Pt}

\begin{Lem} \label{dualSexplicit2}
Soient $n\in\bbN$ et $\kappa/k$ une extension de corps algébriquement close de $k$. Alors pour tout $(G,S)\in \cP^{00}(\cU(n))$ (resp. $\cI^0(\cU(n))$), on a  
$$(G,S)^{\vee_S}(\kappa)=\Ext_{\cP(\cU_{\kappa}(n))}^1((G,S)_{\kappa},\bbZ/p^n) \quad \big(\textrm{resp. } (G,S)^{\vee_S}(\kappa)=\Ext_{\cI(\cU_{\kappa}(n))}^1((G,S)_{\kappa},\bbZ/p^n)\big).$$ 
\end{Lem}

\begin{proof}
Pour tout $(G,S)\in\cP^{00}(\cU(n))$, on a
$$
(G,S)^{\vee_S}(\kappa)=\li_{H\in S^0}((G/H)^{\vee_S})(\kappa)=\li_{H\in S^0}\Ext_{\cU_{\kappa}(n)}^1((G/H)_{\kappa},\bbZ/p^n)
$$
d'après \ref{dualSU0} et \ref{ZZpn}, et
$$
\li_{H\in S^0}\Ext_{\cU_{\kappa}(n)}^1((G/H)_{\kappa},\bbZ/p^n)=\Ext_{\cP(\cU_{\kappa}(n))}^1((G,S)_{\kappa},\bbZ/p^n)
$$
d'après \cite{S2} 3.4 Proposition 7 (ou encore, d'après \ref{descrseP}). Pour tout $(G,S)\in\cI^0(\cU(n))$, on a
$$
(G,S)^{\vee_S}(\kappa)=\lp_{H\in S^0}(H^{\vee_S})(\kappa)=\lp_{H\in S^0}\Ext_{\cU_{\kappa}(n)}^1(H_{\kappa},\bbZ/p^n)
$$
d'après \ref{dualSU0} et \ref{ZZpn}, 
$$
\lp_{H\in S^0}\Ext_{\cU_{\kappa}(n)}^1(H_{\kappa},\bbZ/p^n)=\Ext_{\Ind(\cU_{\kappa}(n))}^1((H_{\kappa})_{H\in S^0},\bbZ/p^n)
=\Ext_{\cI(\cU_{\kappa}(n))}^1((G,S)_{\kappa},\bbZ/p^n)
$$
d'après \ref{Iexacte} et \ref{descrseI}.
\end{proof}

\begin{Pt} \label{foncP00IUIP00U} 
Soit $(G,S)\in\cP^{00}(\cI)(\cU)$. Posons $S^0:=\{H\in S\ |\ H\textrm{ connexe}\}$ (cf. \ref{PIS^0}). Pour tout $H,H'\in S^0$ avec $H\subset H'$, l'épimorphisme strict $G/H\ra G/H'$ de $\cI^0(\cU)$ induit un monomorphisme $(G/H')^{\vee_S}\ra (G/H)^{\vee_S}$ de $\cP^{00}(\cU)$ (\ref{dualSP00UI0U}). Posons $G^{\vee_S}:=\li_{H\in S^0}(G/H)^{\vee_S}$. Pour tout $H\in S^0$, l'application canonique $(G/H)^{\vee_S}\ra G^{\vee_S}$ est injective, identifiant le groupe sous-jacent au groupe proalgébrique unipotent strictement connexe $(G/H)^{\vee_S}$ à un sous-groupe de $G^{\vee_S}$. L'ensemble de sous-groupes de $G^{\vee_S}$
$$
(S^0)^{\vee_S}:=\{ (G/H)^{\vee_S}\subset G^{\vee_S}\ |\ H\in S^0\}
$$
vérifie alors (IP-1), (IP-2)' et (IP-3), munissant donc $G^{\vee_S}$ d'une structure de groupe indproalgébrique 
$(G^{\vee_S},S^{\vee_S})$, avec $S^{\vee_S}:=\ <(S^0)^{\vee_S}>$. On a $(G^{\vee_S},S^{\vee_S})\in\cI(\cP^{00})(\cU)$ par construction. De plus si $f:(G_1,S_1)\ra (G_2,S_2)\in \cP^{00}(\cI)(\cU)$, alors pour tout $H_2\in S_{2}^0$, il existe $H_1\in S_{1}^0$ inclus dans $f^{-1}(H_2)$ tel que l'application $G_1/H_1\ra G_2/H_2$ induite par $f$ soit dans $\cI^0(\cU)$, et définisse donc un morphisme $(G_2/H_2)^{\vee_S}\ra (G_1/H_1)^{\vee_S}\in\cP^{00}(\cU)$ ; d'où un morphisme $f^{\vee_S}:(G_{2}^{\vee_S},S_{2}^{\vee_S})\ra (G_{1}^{\vee_S},S_{1}^{\vee_S})$, et finalement un foncteur 
$$
\xymatrix{
(\cdot)^{\vee_S}:\cP^{00}(\cI)(\cU)^{\circ}\ar[r] & \cI(\cP^{00})(\cU).
}
$$
\end{Pt}

\begin{Pt} \label{foncIP00UP00IU} 
Soit $(G,S)\in\cI(\cP^{00})(\cU)$. Posons $S^{00}:=\{H\in S\ |\ H\textrm{ strictement connexe}\}$ (cf. \ref{PtDefIP00}). Pour tout $H,H'\in S^{00}$ avec $H\subset H'$, le monomorphisme $H\subset H'$ de $\cP^{00}(\cU)$ est strict, et induit donc un épimorphisme strict $(H')^{\vee_S}\ra H^{\vee_S}$ de $\cI^0(\cU)$ (\ref{dualSP00UI0U}). Posons $G^{\vee_S}:=\lp_{H\in S^{00}}H^{\vee_S}$. Pour tout $H\in S^{00}$, l'application canonique $G^{\vee_S}\ra H^{\vee_S}$ est alors surjective (\cite{BouE} III \S 7 N${}^{\circ}4$ Proposition 5) ; si l'on note $K_{H^{\vee_S}}$ son noyau, le quotient $G^{\vee_S}/K_{H^{\vee_S}}$ se trouve donc muni d'une structure de groupe indalgébrique unipotent connexe, canoniquement isomorphe à $H^{\vee_S}$. L'ensemble de sous-groupes de $G^{\vee_S}$
$$
(S^{00})^{\vee_S}:=\{ K_{H^{\vee_S}}\subset G^{\vee_S}\ |\ H\in S^{00}\}
$$
vérifie alors (PI-1), (PI-3) et (PI-4), munissant donc $G^{\vee_S}$ d'une structure de groupe proindalgébrique 
$(G^{\vee_S},S^{\vee_S})$, avec $S^{\vee_S}:=\ <(S^{00})^{\vee_S}>$. On a $(G^{\vee_S},S^{\vee_S})\in\cP^{00}(\cI)(\cU)$ par construction. De plus si $f:(G_1,S_1)\ra (G_2,S_2)\in \cI(\cP^{00})(\cU)$, alors pour tout $H_1\in S_{1}^{00}$, il existe $H_2\in S_{2}^{00}$ contenant $f(H_1)$ tel que l'application $H_1\ra H_2$ induite par $f$ soit dans $\cP^{00}(\cU)$ (on peut prendre $H_2=f(H_1)$), et définisse donc un morphisme $H_{2}^{\vee_S}\ra H_{1}^{\vee_S}\in\cI^0(\cU)$ ; d'où un morphisme $f^{\vee_S}:(G_{2}^{\vee_S},S_{2}^{\vee_S})\ra (G_{1}^{\vee_S},S_{1}^{\vee_S})$, et finalement un foncteur 
$$
\xymatrix{
(\cdot)^{\vee_S}:\cI(\cP^{00})(\cU)^{\circ}\ar[r] & \cP^{00}(\cI)(\cU).
}
$$
\end{Pt}

\begin{Th}\label{dualSP00IUIP00U}
Les foncteurs $(\cdot)^{\vee_S}:\cP^{00}(\cI)(\cU)^{\circ}\ra\cI(\cP^{00})(\cU)$ \ref{foncP00IUIP00U} et $(\cdot)^{\vee_S}:\cI(\cP^{00})(\cU)^{\circ}\ra\cP^{00}(\cI)(\cU)$ \ref{foncIP00UP00IU} sont des équivalences de catégories quasi-inverses l'une de l'autre. De plus, ils sont exacts.
\end{Th}

\begin{proof}
Soit $(G,S)\in\cP^{00}(\cI)(\cU)$. Avec les notations de \ref{foncP00IUIP00U}, on a 
$$(G,S)^{\vee_S}=(\li_{H\in S^0}(G/H)^{\vee_S},\ <(S^0)^{\vee_S}>).$$
Montrons que l'ensemble des sous-groupes de définition strictement connexes de $(G,S)^{\vee_S}$ est réduit à l'ensemble 
$\{(G/H)^{\vee_S}\ |\ H\in S^0\}$. Soit $C$ un tel sous-groupe de $(G,S)^{\vee_S}$. Il existe $H\in S^0$ tel que $C\subset (G/H)^{\vee_S}$. L'inclusion $C\subset (G/H)^{\vee_S}$ est alors un monomorphisme de $\cP^{00}(\cU)$, induisant un épimorphisme strict $G/H\ra C^{\vee_S}$ (\ref{dualSP00UI0U}). Son noyau est de la forme $H^0/H\in\cI^0(\cU)$ pour un $H^0\in S$. Pour tout $H'\in S^0$ inclus dans $H$, le groupe indalgébrique $H^0/H'$ est alors le noyau de l'épimorphisme strict 
$G/H'\ra C^{\vee_S}$ de $\cI^0(\cU)$ dual de $C\subset (G/H')^{\vee_S}$, et est donc connexe. Il en résulte que $H^0$ est connexe, i.e. $H^0\in S^0$, et $C=(G/H^0)^{\vee_S}$ par construction. On a donc montré $<(S^0)^{\vee_S}>^{00}=\{(G/H)^{\vee_S}\ |\ H\in S^0\}$. Ceci étant, il est clair que $(G,S)^{\vee_S\vee_S}=(G,S)$, et que $(\cdot)^{\vee_S}\circ(\cdot)^{\vee_S}=\Id:\cP^{00}(\cI)(\cU)\ra\cP^{00}(\cI)(\cU)$ sur les morphismes aussi.

Soit $(G,S)\in\cI(\cP^{00})(\cU)$. Avec les notations de \ref{foncIP00UP00IU}, on a 
$$(G,S)^{\vee_S}=(\lp_{H\in S^{00}}H^{\vee_S},\ <(S^{00})^{\vee_S}>).$$
Montrons que l'ensemble des sous-groupes de définition connexes de $(G,S)^{\vee_S}$ est réduit à l'ensemble $\{K_{H^{\vee_S}}\ |\ H\in S^{00}\}$. Soit $C$ un tel sous-groupe de $(G,S)^{\vee_S}$. Il existe $H\in S^{00}$ tel que $K_{H^{\vee_S}}\subset C$ et $C/K_{H^{\vee_S}}$ indalgébrique strict dans $G/K_{H^{\vee_S}}=H^{\vee_S}$. D'où un épimorphisme strict $H^{\vee_S}\ra G/C$ de $\cI^0(\cU)$, induisant un monomorphisme $(G/C)^{\vee_S}\ra H$ (\ref{dualSP00UI0U}). Par conséquent le groupe strictement connexe $H^{00}:=(G/C)^{\vee_S}$ est de définition dans $(G,S)$, i.e. $H^{00}\in S^{00}$, et $K_{(H^{00})^{\vee_S}}=C$ par construction. On a donc montré $<(S^{00})^{\vee_S}>^0=\{K_{H^{\vee_S}}\ |\ H\in S^{00}\}$. Ceci étant, il est clair que $(G,S)^{\vee_S\vee_S}=(G,S)$, et que $(\cdot)^{\vee_S}\circ(\cdot)^{\vee_S}=\Id:\cI(\cP^{00})(\cU)\ra\cI(\cP^{00})(\cU)$ sur les morphismes aussi.
 
Les foncteurs $(\cdot)^{\vee_S}:\cP^{00}(\cI)(\cU)^{\circ}\ra\cI(\cP^{00})(\cU)$ et $(\cdot)^{\vee_S}:\cI(\cP^{00})(\cU)^{\circ}\ra\cP^{00}(\cI)(\cU)$ sont exacts d'après \ref{IPS} et \ref{PIS}.
\end{proof}

\begin{Rem}\label{RemdualS3} 
Dans la construction \ref{foncP00IUIP00U} (resp. \ref{foncIP00UP00IU}), on peut remplacer l'ensemble $S^0$ (resp. $S^{00}$) des sous-groupes de définition connexes (resp. strictement connexes) de $(G,S)$ par un ensemble de définition dont les éléments sont connexes (resp. strictement connexes), cf. \ref{RemdualS2}.
\end{Rem}

\begin{Notss}
On note $\cP^{00}(\cI)(\cU(n))$ (resp. $\cI(\cP^{00})(\cU(n))$) la sous-catégorie pleine de $\cP^{00}(\cI)$ (resp. 
$\cI(\cP^{00})$) dont les objets sont ceux de $\cP^{00}(\cI)$ (resp. $\cI(\cP^{00})$) tels que $(G/H)\in\cI(\cU(n))$ (resp. $H\in\cP(\cU(n))$) pour tout $H\in S$.
\end{Notss}

\begin{Lem} \label{dualSexplicit3}
Soient $n\in\bbN$ et $\kappa/k$ une extension de corps algébriquement close de $k$. Alors pour tout $(G,S)\in\cI(\cP^{00})(\cU(n))$, on a  $$(G,S)^{\vee_S}(\kappa)=\Ext_{\cI\cP(\cU_{\kappa}(n))}^1((G,S)_{\kappa},\bbZ/p^n).$$ 
\end{Lem}

\begin{proof}
Pour tout $(G,S)\in\cI(\cP^{00})(\cU(n))$, on a
$$
(G,S)^{\vee_S}(\kappa)=\lp_{H\in S^{00}}(H^{\vee_S})(\kappa)=\lp_{H\in S^{00}}\Ext_{\cP(\cU_{\kappa}(n))}^1(H_{\kappa},\bbZ/p^n)
$$
d'après \ref{dualSexplicit2}, et
\begin{eqnarray*}
\lp_{H\in S^{00}}\Ext_{\cP(\cU_{\kappa}(n))}^1(H_{\kappa},\bbZ/p^n)&=&\Ext_{\Ind(\cP(\cU_{\kappa}(n)))}^1((H_{\kappa})_{H\in S^{00}},\bbZ/p^n)\\
&=&\Ext_{\cI\cP(\cU_{\kappa}(n))}^1((G,S)_{\kappa},\bbZ/p^n)
\end{eqnarray*}
d'après \ref{IPexacte} et \ref{descrseIP}.
\end{proof}

\begin{Pt} \label{foncP00IUIP00Uad} 
Soit $(G,S)\in\cP^{00}(\cI)(\cU)^{\ad}$. Comme $S^0$ est coinitial dans $S$ et l'ensemble des origines de $(G,S)$ stable par minorant, il existe une origine $H_0$ de $(G,S)$ dans $S^0$. Soit $H\in S^0$ tel que $(G/H_0)^{\vee_S}\subset (G/H)^{\vee_S}$. Alors $H\subset H_0$, d'où $H_0/H\in\cU^0$, puis $(H_0/H)^{\vee_S}\in\cU^0$. Comme $(H_0/H)^{\vee_S}$ est isomorphe à $(G/H)^{\vee_S}/(G/H_{0})^{\vee_S}$, on obtient en particulier que celui-ci est quasi-algébrique. L'ensemble $\{(G/H)^{\vee_S}\supset (G/H_0)^{\vee_S}\ |\ H\in S^0\}$ étant cofinal dans $S^{\vee_S}$, on en déduit que $(G/H_0)^{\vee_S}$ est une origine de $(G,S)^{\vee_S}$, et donc que $(G,S)^{\vee_S}\in\cI(\cP^{00})(\cU)^{\ad}.$ D'où un foncteur 
$$
\xymatrix{
(\cdot)^{\vee_S}:(\cP^{00}(\cI)(\cU)^{\ad})^{\circ}\ar[r] & \cI(\cP^{00})(\cU)^{\ad}.
}
$$
\end{Pt}

\begin{Pt} \label{foncIP00UP00IUad} 
Soit $(G,S)\in\cI(\cP^{00})(\cU)^{\ad}$. Comme $S^{00}$ est cofinal dans $S$ et l'ensemble des origines de $(G,S)$ stable par majorant, il existe une origine $H_0$ de $(G,S)$ dans $S^{00}$. Soit $H\in S^{00}$ tel que $K_{H^{\vee_S}}\subset K_{(H_0)^{\vee_S}}$. Alors $H_0\subset H$, d'où $H/H_0\in\cU^0$, puis $(H/H_0)^{\vee_S}\in\cU^0$. Comme $(H/H_0)^{\vee_S}$ est isomorphe à $K_{(H_0)^{\vee_S}}/K_{H^{\vee_S}}$, on obtient en particulier que celui-ci est quasi-algébrique. L'ensemble $\{K_{H^{\vee_S}}\subset K_{(H_0)^{\vee_S}}\ |\ H\in S^{00}\}$ étant coinitial dans $S^{\vee_S}$, on en déduit que $K_{(H_0)^{\vee_S}}$ est une origine de $(G,S)^{\vee_S}$, et donc que $(G,S)^{\vee_S}\in\cP^{00}(\cI)(\cU)^{\ad}.$ D'où un foncteur 
$$
\xymatrix{
(\cdot)^{\vee_S}:(\cI(\cP^{00})(\cU)^{\ad})^{\circ}\ra\cP^{00}(\cI)(\cU)^{\ad}.
}
$$
\end{Pt}

Ceci étant, le théorème \ref{dualSP00IUIP00U} admet le corollaire suivant :

\begin{Cor}\label{dualSP00IUIP00Uad}
Les foncteurs $(\cdot)^{\vee_S}:(\cP^{00}(\cI)(\cU)^{\ad})^{\circ}\ra\cI(\cP^{00})(\cU)^{\ad}$ \ref{foncP00IUIP00Uad} et 
$(\cdot)^{\vee_S}:(\cI(\cP^{00})(\cU)^{\ad})^{\circ}\ra\cP^{00}(\cI)(\cU)^{\ad}$ \ref{foncIP00UP00IUad} sont des équivalences de catégories exactes quasi-inverses l'une de l'autre. 
\end{Cor}

\begin{Pt}\label{egdualSerre}
Soit $(G,S)\in\fUfzz$. Posons 
$$
S^0:=\{H\in S\ |\ \forall  H'\in S, H'\subset H\Rightarrow H/H'\in\cU^0\}.
$$
Par construction, le groupe sous-jacent à $(\indpro(G,S))^{\vee_S}\in\cP^{00}(\cI)(\cU)^{\ad}$ (\ref{fQfIPfPIf00}) est   
$$
\lp_{H^+\in S^0}H_{+}^{\vee_S},
$$
avec pour ensemble de définition $\{K_{H_{+}^{\vee_S}}|H_+\in S^0\}$, et  
$$
(\lp_{H^+\in S^0}H_{+}^{\vee_S})/K_{H_{+}^{\vee_S}}=H_{+}^{\vee_S}=(\lp_{S^0\ni H_-\subset H_+}H_+/H_-)^{\vee_S}=\li_{S^0\ni H_-\subset H_+}(H_+/H_-)^{\vee_S}\in\cI^0(\cU).
$$
Pour tout $H_-,H_+\in S^0$ avec $H_-\subset H_+$, le diagramme commutatif exact de $\cP\cI$
$$
\xymatrix{
0\ar[r] & K_{H_{+}^{\vee_S}} \ar[r] \ar[d] & \lp_{H^+\in S^0}H_{+}^{\vee_S} \ar[r] \ar@{=}[d] & H_{+}^{\vee_S} \ar[r] \ar[d] & 0 \\
0\ar[r] & K_{H_{-}^{\vee_S}} \ar[r] & \lp_{H^+\in S^0}H_{+}^{\vee_S} \ar[r] & H_{-}^{\vee_S} \ar[r] & 0 
}
$$
montre que $K_{H_{-}^{\vee_S}}/K_{H_{+}^{\vee_S}} = \Ker (H_{+}^{\vee_S}\ra H_{-}^{\vee_S})=(H_+/H_-)^{\vee_S}$. On a donc 
$$
K_{H_{-}^{\vee_S}}=\lp_{S^0\ni H_+\supset H_-}(H_+/H_-)^{\vee_S}\in\cP^{00}(\cU),
$$
et $\{K_{H_{-}^{\vee_S}}|H_-\in S^0\}$ est un ensemble de définition de $\can((\indpro(G,S))^{\vee_S})\in\cI(\cP^{00})(\cU)^{\ad}$, dont le groupe sous-jacent est
$$
\li_{H_-\in S^0}K_{H_{-}^{\vee_S}}.
$$
Mais par construction, le groupe sous-jacent à $(\proind(G,S))^{\vee_S}\in\cI(\cP^{00})(\cU)^{\ad}$ est
$$
\li_{H_-\in S^0}(G/H_-)^{\vee_S},
$$
avec pour ensemble de définition $\{(G/H_-)^{\vee_S}\ |\ H_-\in S^0\}$, et 
$$
(G/H_-)^{\vee_S}=(\li_{S^0\ni H_+\supset H_-}H_+/H_-)^{\vee_S}=\lp_{S^0\ni H_+\supset H_-}(H_+/H_-)^{\vee_S}\in\cP^{00}(\cU).
$$
On a donc obtenu $(\can\circ(\cdot)^{\vee_S}\circ\indpro)(G,S)=((\cdot)^{\vee_S}\circ\proind)(G,S)\in\cI(\cP^{00})(\cU)^{\ad}$.
On en tire (\ref{fQfIPfPIf}) :
$$
\ad\circ(\cdot)^{\vee_S}\circ\indpro=\ad\circ(\cdot)^{\vee_S}\circ\proind:(\fUfzz)^{\circ}\ra\fUfzz.
$$
On notera encore $(\cdot)^{\vee_S}:(\fUfzz)^{\circ}\ra\fUfzz$ ce foncteur.
\end{Pt}

\begin{Def} \label{dualSfUfzz}
Le foncteur $(\cdot)^{\vee_S}:(\fUfzz)^{\circ}\ra\fUfzz$ \ref{egdualSerre} est le foncteur de \emph{dualité de Serre sur $\fUfzz$}.
\end{Def}

D'après \ref{dualSP00IUIP00Uad} et \ref{fQfIPfPIf00}, on a :

\begin{Cor} \label{dualSfUf00}
Le foncteur de dualité de Serre $(\cdot)^{\vee_S}:(\fUfzz)^{\circ}\ra\fUfzz$ est exact et involutif.
\end{Cor}

\begin{Lem} \label{dualSexplicit4}
Soient $n\in\bbN$ et $\kappa/k$ une extension de corps algébriquement close de $k$. Alors pour tout $(G,S)\in \fUnf^{00}$, on a  
$$(G,S)^{\vee_S}(\kappa)=\Ext_{\fUnf_{\kappa}}^1((G,S)_{\kappa},\bbZ/p^n).$$ 
\end{Lem}

\begin{proof} Pour tout $(G,S)\in \fUnf^{00}$, on a
\begin{eqnarray*}
(G,S)^{\vee_S}(\kappa) & = & \ad((\indpro(G,S))^{\vee_S})(\kappa) \\
                                      & = & (\indpro(G,S))^{\vee_S}(\kappa) \\
                                      & = & \Ext_{\cI\cP(\cU_{\kappa}(n))}^1(\indpro(G,S)_{\kappa},\bbZ/p^n) \\
                                      & = & \Ext_{\fUnf_{\kappa}}^1((G,S)_{\kappa},\bbZ/p^n),
\end{eqnarray*}
les deux premières égalités résultant des définitions, la troisième de \ref{dualSexplicit3} et la quatrième de la définition de la structure exacte de $\fUnf_{\kappa}$ (cf. \ref{fQfexacte}).
\end{proof}

\begin{Pt}\label{dualLV}
Le foncteur de \emph{dualité linéaire}
$$
\xymatrix{
\Hom_{\cV}(\cdot,k):\cV^{\circ} \ar[r] & \cV
}
$$
est exact et involutif. On le note $(\cdot)^{\vee_L}$. Pour tout $V\in\cV$, l'extension d'Artin-Schreier parfaite
$$
\xymatrix{
0 \ar[r] & \bbZ/p \ar[r] & \bbG_{a,k}^{\pf} \ar[r]^{1-F} & \bbG_{a,k}^{\pf} \ar[r] & 0 & =:\epsilon
}
$$
induit un morphisme de faisceaux abéliens sur $(\parf/k)_{\et}$
\begin{eqnarray*}
\uHom_{\klin}(V^{\pf},\bbG_{a,k}^{\pf}) & \lra & \uExt_{\bbZ/p}^1(V^{\pf},\bbZ/p) \\
f & \lmapsto & f^*\epsilon,
\end{eqnarray*}
qui est un isomorphisme d'après \cite{Beg} 1.1.1 b). Autrement dit, le diagramme
$$ 
\xymatrix{
\cV^{\circ}\ar[r] \ar[d]_{(\cdot)^{\vee_L}} & (\cU(1)^0)^{\circ} \ar[d]^{(\cdot)^{\vee_S}} \\
\cV \ar[r] & \cU(1)^0 
}
$$
est commutatif.
\end{Pt}

\begin{Pt} \label{foncPVIV} 
Soit $(V,S)\in\cP(\cV)$. Pour tout $W,W'\in S$ avec $W\subset W'$, l'épimorphisme $V/W\ra V/W'$ de $\cV$ induit un monomorphisme $(V/W')^{\vee_L}\ra (V/W)^{\vee_L}$. Posons $V^{\vee_L}:=\li_{W\in S}(V/W)^{\vee_L}$. Pour tout $W\in S$, l'application canonique $(V/W)^{\vee_L}\ra V^{\vee_L}$ est injective, identifiant $(V/W)^{\vee_L}$ à un sous-$k$-vectoriel de $V^{\vee_L}$. L'ensemble de sous-$k$-vectoriels de $V^{\vee_L}$
$$
S^{\vee_L}:=\{ (V/W)^{\vee_L}\subset V^{\vee_L}\ |\ W\in S\}
$$
vérifie alors (I$\cV$-1), (I$\cV$-2) et (I$\cV$-3), définissant donc un objet $(V^{\vee_L},S^{\vee_L})\in\cI(\cV)$. De plus si $f:(V_1,S_1)\ra (V_2,S_2)\in \cP(\cV)$, alors pour tout $W_2\in S_{2}$, $f^{-1}(W_2)\in S_1$ et l'application $k$-linéaire $V_1/f^{-1}(W_2)\ra V_2/W_2$ induite par $f$ définit une application $k$-linéaire $(V_2/W_2)^{\vee_L}\ra (V_1/f^{-1}(W_2))^{\vee_L}$ ; d'où un morphisme $f^{\vee_L}:(V_{2}^{\vee_L},S_{2}^{\vee_L})\ra (V_{1}^{\vee_L},S_{1}^{\vee_L})$, et finalement un foncteur 
$$
\xymatrix{
(\cdot)^{\vee_L}:\cP(\cV)^{\circ}\ar[r] & \cI(\cV).
}
$$
\end{Pt}

\begin{Pt} \label{foncIVPV} 
Soit $(V,S)\in\cI(\cV)$. Pour tout $W,W'\in S$ avec $W\subset W'$, le monomorphisme $W\ra W'$ de $\cV$ induit un épimorphisme $(W')^{\vee_L}\ra W^{\vee_L}$. Posons $V^{\vee_L}:=\lp_{W\in S}W^{\vee_L}$. Pour tout $W\in S$, l'application canonique $V^{\vee_L}\ra W^{\vee_L}$ est alors surjective (\cite{BouE} III \S 7 N${}^{\circ}4$ Proposition 5) ;
si l'on note $K_{W^{\vee_L}}$ son noyau, le quotient $V^{\vee_L}/K_{W^{\vee_L}}$ est donc de dimension finie, canoniquement isomorphe à $W^{\vee_L}$. L'ensemble de sous-$k$-vectoriels de $V^{\vee_L}$
$$
S^{\vee_L}:=\{ K_{W^{\vee_L}}\subset V^{\vee_L}\ |\ W\in S\}
$$
vérifie alors (P$\cV$-1), (P$\cV$-2), (P$\cV$-3) et (P$\cV$-4), définissant donc un objet $(V^{\vee_L},S^{\vee_L})\in\cP(\cV)$. De plus si $f:(V_1,S_1)\ra (V_2,S_2)\in \cI(\cV)$, alors pour tout $W_1\in S_1$, $f(W_1)\in S_2$ et l'application $k$-linéaire $W_1\ra f(W_1)$ induite par $f$ définit une application $k$-linéaire $f(W_1)^{\vee_L}\ra W_{1}^{\vee_L}$ ; d'où un morphisme $f^{\vee_L}:(V_{2}^{\vee_L},S_{2}^{\vee_L})\ra (V_{1}^{\vee_L},S_{1}^{\vee_L})$, et finalement un foncteur 
$$
\xymatrix{
(\cdot)^{\vee_L}:\cI(\cV)^{\circ}\ar[r] & \cP(\cV).
}
$$
\end{Pt}

\begin{Pt}\label{dualLPVIV}
Les foncteurs $\cP(\cV)^{\circ}\ra\cI(\cV)$ \ref{foncPVIV} et $\cI(\cV)^{\circ}\ra\cP(\cV)$ \ref{foncIVPV} sont des équivalences de catégories quasi-inverses l'une de l'autre. De plus, ils sont exacts. 
\end{Pt}

\begin{Pt}\label{dualLexplicit2}
Pour tout $(V,S)\in\cP(\cV)$ (resp. $\cI(\cV)$), on a
$$(V,S)^{\vee_L}=\Hom_{\cP(\cV)}((V,S),k) \quad \big(\textrm{resp. } (V,S)^{\vee_L}=\Hom_{\cI(\cV)}((V,S),k)\big).$$
\end{Pt}

\begin{Pt} \label{foncPIVIPV} 
Soit $(V,S)\in\cP\cI(\cV)$. Pour tout $W,W'\in S$ avec $W\subset W'$, l'épimorphisme strict $V/W\ra V/W'$ de $\cI(\cV)$ induit un monomorphisme $(V/W')^{\vee_L}\ra (V/W)^{\vee_L}$ de $\cP(\cV)$ (\ref{dualLPVIV}). Posons $V^{\vee_L}:=\li_{W\in S}(V/W)^{\vee_L}$. Pour tout $W\in S$, l'application canonique $(V/W)^{\vee_L}\ra V^{\vee_L}$ est injective, identifiant $(V/W)^{\vee_L}$ à un sous-$k$-vectoriel de $V^{\vee_L}$. L'ensemble de sous-$k$-vectoriels de $V^{\vee_L}$
$$
S^{\vee_L}:=\{ (V/W)^{\vee_L}\subset V^{\vee_L}\ |\ W\in S\}
$$
vérifie alors (IP$\cV$-1), (IP$\cV$-2) et (IP$\cV$-3), définissant donc un objet $(V^{\vee_L},S^{\vee_L})\in\cI\cP(\cV)$. De plus si $f:(V_1,S_1)\ra (V_2,S_2)\in \cP\cI(\cV)$, alors pour tout $W_2\in S_2$, $f^{-1}(W_2)\in S_1$ et l'application $V_1/f^{-1}(W_2)\ra V_2/W_2\in\cI(\cV)$ induite par $f$ définit un morphisme $(V_2/W_2)^{\vee_L}\ra (V_1/f^{-1}(W_2))^{\vee_L}\in\cP(\cV)$ ; d'où un morphisme $f^{\vee_L}:(V_{2}^{\vee_L},S_{2}^{\vee_L})\ra (V_{1}^{\vee_L},S_{1}^{\vee_L})$, et finalement un foncteur 
$$
\xymatrix{
(\cdot)^{\vee_L}:\cP\cI(\cV)^{\circ}\ar[r] & \cI\cP(\cV).
}
$$
\end{Pt}

\begin{Pt} \label{foncIPVPIV} 
Soit $(V,S)\in\cI\cP(\cV)$. Pour tout $W,W'\in S$ avec $W\subset W'$, le monomorphisme $W\ra W'$ de $\cP(\cV)$ est strict 
et induit donc un épimorphisme strict $(W')^{\vee_L}\ra W^{\vee_L}$ de $\cI(\cV)$ (\ref{dualLPVIV}). Posons $V^{\vee_L}:=\lp_{W\in S}W^{\vee_L}$. Pour tout $W\in S$, l'application canonique $V^{\vee_L}\ra W^{\vee_L}$ est alors surjective  (\cite{BouE} III \S 7 N${}^{\circ}4$ Proposition 5) ; si l'on note $K_{W^{\vee_L}}$ son noyau, le quotient $V^{\vee_L}/K_{W^{\vee_L}}$ se trouve donc muni d'une structure d'objet de $\cI(\cV)$, canoniquement isomorphe à $W^{\vee_L}$. L'ensemble de sous-$k$-vectoriels de $V^{\vee_L}$
$$
S^{\vee_L}:=\{ K_{W^{\vee_L}}\subset V^{\vee_L}\ |\ W\in S\}
$$
vérifie alors (PI$\cV$-1), (PI$\cV$-2), (PI$\cV$-3) et (PI$\cV$-4), définissant donc un objet  $(V^{\vee_L},S^{\vee_L})\in\cP\cI(\cV)$. De plus si $f:(V_1,S_1)\ra (V_2,S_2)\in \cI\cP(\cV)$, alors pour tout $W_1\in S_1$, $f(W_1)\in S_2$ et l'application $W_1\ra f(W_1)\in\cP(\cV)$ induite par $f$ définit un morphisme $f(W_1)^{\vee_L}\ra W_{1}^{\vee_L}\in\cI(\cV)$ ; d'où un morphisme $f^{\vee_L}:(V_{2}^{\vee_L},S_{2}^{\vee_L})\ra (V_{1}^{\vee_L},S_{1}^{\vee_L})$, et finalement un foncteur 
$$
\xymatrix{
(\cdot)^{\vee_L}:\cI\cP(\cV)^{\circ}\ar[r] & \cP\cI(\cV).
}
$$
\end{Pt}

\begin{Pt}\label{dualLPIVIPV}
Les foncteurs $\cP\cI(\cV)^{\circ}\ra\cI\cP(\cV)$ \ref{foncPIVIPV} et $\cI\cP(\cV)^{\circ}\ra\cP\cI(\cV)$ \ref{foncIPVPIV} sont des équivalences de catégories quasi-inverses l'une de l'autre. De plus, ils sont exacts. 
\end{Pt}

\begin{Pt}\label{dualLexplicit3}
Pour tout $(V,S)\in\cP\cI(\cV)$ (resp. $\cI\cP(\cV)$), on a
$$(V,S)^{\vee_L}=\Hom_{\cP\cI(\cV)}((V,S),k) \quad \big(\textrm{resp. } (V,S)^{\vee_L}=\Hom_{\cI\cP(\cV)}((V,S),k)\big).$$
\end{Pt}

\begin{Pt} \label{foncPIVIPVad} 
Soit $(V,S)\in\cP\cI(\cV)^{\ad}$. Soit $W_0$ une origine de $(V,S)$. Soit $W\in S$ tel que $(V/W_0)^{\vee_L}\subset (V/W)^{\vee_L}$. Alors $W\subset W_0$, d'où $W_0/W\in\cV$, puis $(W_0/W)^{\vee_L}\in\cV$. Comme $(W_0/W)^{\vee_L}$ est isomorphe à $(V/W)^{\vee_L}/(V/W_0)^{\vee_L}$, on obtient en particulier que celui-ci est de dimension finie. On en déduit que $(V/W_0)^{\vee_L}$ est une origine de $(V,S)^{\vee_L}$, et donc que $(V,S)^{\vee_L}\in\cI\cP(\cV)^{\ad}.$ D'où un foncteur 
$$
\xymatrix{
(\cdot)^{\vee_L}:(\cP\cI(\cV)^{\ad})^{\circ}\ar[r] & \cI\cP(\cV)^{\ad}.
}
$$
\end{Pt}

\begin{Pt} \label{foncIPVPIVad} 
Soit $(V,S)\in\cI\cP(\cV)^{\ad}$. Soit $W_0$ une origine de $(V,S)$. Soit $W\in S$ tel que $K_{W^{\vee_L}}\subset K_{(W_0)^{\vee_L}}$. Alors $W_0\subset W$, d'où $W/W_0\in\cV$, puis $(W/W_0)^{\vee_L}\in\cV$. Comme $(W/W_0)^{\vee_L}$ est isomorphe à $K_{(W_0)^{\vee_L}}/K_{W^{\vee_L}}$, on obtient en particulier que celui-ci est de dimension finie. On en déduit que $K_{(W_0)^{\vee_L}}$ est une origine de $(V,S)^{\vee_L}$, et donc que $(V,S)^{\vee_L}\in\cP\cI(\cV)^{\ad}.$ D'où un foncteur 
$$
\xymatrix{
(\cdot)^{\vee_L}:(\cI\cP(\cV)^{\ad})^{\circ}\ra\cP\cI(\cV)^{\ad}.
}
$$
\end{Pt}

\begin{Pt}\label{dualLPIVIPVad}
Les foncteurs $(\cP\cI(\cV)^{\ad})^{\circ}\ra\cI\cP(\cV)^{\ad}$ \ref{foncPIVIPVad} et $(\cI\cP(\cV)^{\ad})^{\circ}\ra\cP\cI(\cV)^{\ad}$ \ref{foncIPVPIVad} sont des équivalences de catégories quasi-inverses l'une de l'autre. De plus, ils sont exacts. 
\end{Pt}

\begin{Pt}\label{egdualLineaire}
Soit $(V,S)\in\fVf$. Par construction, le $k$-vectoriel sous-jacent à $(\indpro(V,S))^{\vee_L}\in\cP\cI(\cV)^{\ad}$ est   
$$
\lp_{W^+\in S}W_{+}^{\vee_L},
$$
avec pour ensemble de définition $\{K_{W_{+}^{\vee_L}}|W_+\in S\}$, et  
$$
(\lp_{W^+\in S}W_{+}^{\vee_L})/K_{W_{+}^{\vee_L}}=W_{+}^{\vee_L}=(\lp_{S\ni W_-\subset W_+}W_+/W_-)^{\vee_L}=\li_{S\ni W_-\subset W_+}(W_+/W_-)^{\vee_L}\in\cI(\cV).
$$
Pour tout $W_-,W_+\in S$ avec $W_-\subset W_+$, le diagramme commutatif exact de $\cP\cI(\cV)$
$$
\xymatrix{
0\ar[r] & K_{W_{+}^{\vee_L}} \ar[r] \ar[d] & \lp_{W^+\in S}W_{+}^{\vee_L} \ar[r] \ar@{=}[d] & W_{+}^{\vee_L} \ar[r] \ar[d] & 0 \\
0\ar[r] & K_{W_{-}^{\vee_L}} \ar[r] & \lp_{W^+\in S}W_{+}^{\vee_L} \ar[r] & W_{-}^{\vee_L} \ar[r] & 0 
}
$$
montre que $K_{W_{-}^{\vee_L}}/K_{W_{+}^{\vee_L}} = \Ker (W_{+}^{\vee_L}\ra W_{-}^{\vee_L})=(W_+/W_-)^{\vee_L}$. On a donc 
$$
K_{W_{-}^{\vee_L}}=\lp_{S\ni W_+\supset W_-}(W_+/W_-)^{\vee_L}\in\cP(\cV),
$$
et $\{K_{W_{-}^{\vee_L}}|W_-\in S\}$ est un ensemble de définition de $\can((\indpro(V,S))^{\vee_L})\in\cI\cP(\cV)^{\ad}$, dont le $k$-vectoriel sous-jacent est
$$
\li_{W_-\in S}K_{W_{-}^{\vee_L}}.
$$
Mais par construction, le $k$-vectoriel sous-jacent à $(\proind(V,S))^{\vee_L}\in\cI\cP(\cV)^{\ad}$ est
$$
\li_{W_-\in S}(V/W_-)^{\vee_L},
$$
avec pour ensemble de définition $\{(V/W_-)^{\vee_L}\ |\ W_-\in S\}$, et 
$$
(V/W_-)^{\vee_L}=(\li_{S\ni W_+\supset W_-}W_+/W_-)^{\vee_L}=\lp_{S\ni W_+\supset W_-}(W_+/W_-)^{\vee_L}\in\cP(\cV).
$$
On a donc obtenu $(\can\circ(\cdot)^{\vee_L}\circ\indpro)(V,S)=((\cdot)^{\vee_L}\circ\proind)(V,S)\in\cI\cP(\cV)^{\ad}$.
On en tire :
$$
\ad\circ(\cdot)^{\vee_L}\circ\indpro=\ad\circ(\cdot)^{\vee_L}\circ\proind:(\fVf)^{\circ}\ra\fVf.
$$
On notera encore $(\cdot)^{\vee_L}:(\fVf)^{\circ}\ra\fVf$ ce foncteur.
\end{Pt}

\begin{Def}
Le foncteur $(\cdot)^{\vee_L}:(\fVf)^{\circ}\ra\fVf$ \ref{egdualLineaire} est le foncteur de \emph{dualité linéaire} sur $\fVf$.
\end{Def}

On a alors :

\begin{Pt} \label{dualSfVf}
Le foncteur de dualité linéaire $(\cdot)^{\vee_L}:(\fVf)^{\circ}\ra\fVf$ est exact et involutif.
\end{Pt}

\begin{Lem}\label{dualLdualSad}
Les diagrammes
$$ 
\xymatrix{
\cP(\cV)^{\circ} \ar[r] \ar[d]_{(\cdot)^{\vee_L}} & \cP^{00}(\cU(1))^{\circ} \ar[d]^{(\cdot)^{\vee_S}} & & \cP\cI(\cV)^{\circ} \ar[r] \ar[d]_{(\cdot)^{\vee_L}} & \cP^{00}(\cI)(\cU(1))^{\circ} \ar[d]^{(\cdot)^{\vee_S}} \\
\cI(\cV)\ar[r] & \cI^0(\cU(1)) & & \cI\cP(\cV)\ar[r] & \cI(\cP^{00})(\cU(1))
}
$$
$$ 
\xymatrix{
 (\cP\cI(\cV)^{\ad})^{\circ} \ar[r] \ar[d]_{(\cdot)^{\vee_L}} & (\cP^{00}(\cI)(\cU(1))^{\ad})^{\circ} \ar[d]^{(\cdot)^{\vee_S}} & & (\fVf)^{\circ} \ar[r] \ar[d]_{(\cdot)^{\vee_L}} & (\fUuf^{00})^{\circ} \ar[d]^{(\cdot)^{\vee_S}} \\
\cI\cP(\cV)^{\ad} \ar[r] & \cI(\cP^{00})(\cU(1))^{\ad} & & \fVf \ar[r] & \fUuf^{00} 
}
$$
sont commutatifs.
\end{Lem}

\begin{proof}
Résulte de \ref{dualLV},  \ref{RemdualS2} et \ref{RemdualS3}.
\end{proof}

\chapter{Rappels sur les schémas en groupes de Greenberg  en égale caractéristique $p$}

Soit $R$ un anneau de valuation discrète complet de caractéristique $p>0$, de corps des fractions $K$ et de corps résiduel $k$. On note $\fm$ l'idéal maximal de $R$ et on pose $R_n:=R/\fm^n$ pour tout entier $n\geq 1$.

\section{Définition}

\begin{Pt}
Soit $X$ un $R$-schéma. Pour tout entier $n\geq 1$, on note 
$$\Gr_n(X):=\prod_{R_n/k}\ X\otimes_R R_n$$
la restriction de Weil de $X\otimes_R R_n$ le long du morphisme canonique $k\ra R_n$. Comme $k\ra R_n$ est fini, plat et radiciel, $\Gr_n(X)$ est un $k$-schéma (cf. \cite{BLR} 7.6/4 et sa preuve), le \emph{schéma de Greenberg de cran $n$ de $X$}. Pour $n=1$, il s'agit de la fibre spéciale $X_k$ du $R$-schéma $X$. Si $X/R$ est affine (resp. localement de type fini, resp. de type fini, resp. lisse), alors pour tout $n\geq 1$, $\Gr_n(X)$ aussi (cf. \cite{BLR} 7.6/4 et sa preuve, et 7.6/5).
\end{Pt}

\begin{Pt}\label{green}
Pour tout $n\geq 1$, la formation de $\Gr_{n}(X)$ est fonctorielle en le $k$-schéma $X$, et commute aux produits ; en particulier, $\Gr_n(\cdot)$ transforme $R$-schémas en groupes en $k$-schémas en groupes. De plus, pour tous $n'\geq n$, le morphisme canonique $\Gr_{n'}(X)\ra \Gr_n(X)$ est naturel en $X$. L'association 
$G\lmapsto (\Gr_n(G))_{n\in\bbN}$ définit donc un foncteur
$$
\xymatrix{
\Gr:(\textrm{$R$-schémas en groupes}) \ar[r] &\Pro(\textrm{$k$-schémas en groupes}),
}
$$
le \emph{foncteur de Greenberg}.
\end{Pt}

\begin{Not}\label{Kk'}
Si $k'/k$ une extension de corps, on pose 
$$
R_{k'}:=R\widehat{\otimes}_k k'=\lp_{n\in\bbN}(R/\fm^n\otimes_kk')\quad\textrm{et}\quad K_{k'}=\Frac(R_{k'}).
$$
Alors $R_{k'}$ est un anneau de valuation discrète complet, de corps des fractions $K_{k'}$ et de corps résiduel $k'$. Son idéal maximal est $\fm R_{k'}$.
\end{Not}

\begin{Pt}\label{greencb}
Pour toute extension de corps $k'/k$, pour tout $R$-schéma $X$, et pour tout entier $n\geq 1$, le morphisme canonique de $k'$-schémas 
$$
\xymatrix{
\Gr_n(X_{R_{k'}}) \ar[r] & \Gr_n(X)_{k'}
}
$$
est un isomorphisme. En particulier, le foncteur de Greenberg $\Gr$ commute à toute extension de corps $k'/k$.
\end{Pt}

\begin{Lem}\label{troncation}
Soit $G$ un $R$-schéma en groupes. Alors pour tout $m\geq n\geq 1$, on a une suite exacte 
$$
\xymatrix{
0\ar[r] & \Gr_n{V(\omega_G)} \ar[r] & \Gr_{n+m}(G)\ar[r] & \Gr_{m}(G),
}
$$
où l'on a noté $V(\omega_G)$ le $R$-groupe vectoriel défini par le $R$-module conormal $\omega_G$ de la section nulle de $G$.
\end{Lem} 

\begin{proof}  
Soit $\Lambda$ une $k$-algèbre. La réduction modulo $\fm^{m}$
$$
\xymatrix{
R_{n+m}\ar[r] & R_{m}
}
$$
induit une flèche
$$
\xymatrix{
r_{n+m,m}:G\otimes_R R_{n+m}(\Lambda\otimes_k R_{n+m})\ar[r] & G\otimes_R R_{n+m}(\Lambda\otimes_k R_{m}).
}
$$
Notons 
$$I:=\fm^{m}/\fm^{n+m}$$
le noyau de $R_{n+m}\ra R_{m}.$ Comme $m\geq n$, l'idéal $I$ de $R_{n+m}$ est de carré nul. Posant $\Lambda_n:=\Lambda\otimes_k R_{n}$, on a donc \begin{eqnarray*}
\Ker(r_{n+m,m})&=&\Hom_{\Lambda_{n+m}}\big((\omega_{G\otimes_R R_{n+m}})\otimes_{R_{n+m}}\Lambda_{n+m}\ ,\ I\Lambda_{n+m}\big)\\
&=&\Hom_{\Lambda_{n+m}}\big((\omega_{G}\otimes_R \Lambda_{n+m}\ ,\ I\Lambda_{n+m}\big)\\
&=&\Hom_{R}\big(\omega_G\ ,\ (I\Lambda_{n+m})_R\big),
\end{eqnarray*}
où le $R$-module $(I\Lambda_{n+m})_R$ est le groupe $I\Lambda_{n+m}$ muni de l'action de $R$ déduite du morphisme composé $R\ra R_{n+m}\ra \Lambda_{n+m}$. Maintenant, 
$$I\Lambda_{n+m}=\fm^m\Lambda_{n+m}$$
est un $\Lambda_{n+m}$-module isomorphe à $\Lambda_{n}$, celui-ci étant muni de la structure de $\Lambda_{n+m}$-module induite par $\Lambda_{n+m}\ra \Lambda_{n}$. Notant $(\Lambda_n)_R$ le $R$-module qui s'en déduit via
$$R\ra\Lambda_{n+m}\ra\Lambda_{n}=R\ra R_{n}\ra\Lambda_n,$$
on trouve donc
\begin{eqnarray*}
\Ker(r_{n+m,m})&=&\Hom_{R}\big(\omega_G\ ,\ (\Lambda_n)_R\big)\\
&=&\Gr_n V(\omega_G)(\Lambda).
\end{eqnarray*}
\end{proof}

\section{Quelques propriétés dans le cas lisse}

\begin{Lem}\label{troncationlisse}
Soit $L$ un $R$-schéma en groupes lisse. 
\begin{enumerate}

\medskip

\item Pour tout $m\geq n\geq 1$, la suite de $k$-schéma en groupes
$$
\xymatrix{
0\ar[r] & \Gr_n V(\omega_L) \ar[r] & \Gr_{n+m}(L)\ar[r] & \Gr_{m}(L) \ar[r] & 0
}
$$
est exacte sur le gros site Zariski des $k$-schémas $(\Sch/k)_{\ZAR}$. En particulier, on a la suite exacte
$$
\xymatrix{
0\ar[r] & V(\omega_{L_k}) \ar[r] & \Gr_{n+1}(L)\ar[r] & \Gr_n(L) \ar[r] & 0.
}
$$

\medskip

\item Pour tout $n\geq 1$, on a une suite de $k$-schémas en groupes  
$$
\xymatrix{
0 \ar[r] & U_n \ar[r] & \Gr_{n}(L) \ar[r] & L_k \ar[r] & 0
}
$$
exacte sur $(\Sch/k)_{\ZAR}$, où $U_n$ est un $k$-schéma en groupes lisse connexe unipotent de dimension $(n-1)\dim(L/R)$. En particulier, le $k$-schéma en groupes $\Gr_n(L)$ est lisse de dimension $n\dim(L/R)$.

\medskip

\item Pour tout $n\geq 1$, l'inclusion $L^0\ra L$ induit un isomorphisme 
$$
\xymatrix{
\Gr_n(L^0) \ar[r]^{\sim} & \Gr_n(L)^0,
}
$$
et la transition $\Gr_{n+1}(L)\ra\Gr_n(L)$ induit une suite exacte 
$$
\xymatrix{
0\ar[r] & V(\omega_{L_k}) \ar[r] & \Gr_{n+1}(L)^0 \ar[r] & \Gr_n(L)^0 \ar[r] & 0
}
$$
et un isomorphisme
$$
\xymatrix{
\pi_0(\Gr_{n+1}(L)) \ar[r]^{\sim} & \pi_0(\Gr_n(L)). 
}
$$
\end{enumerate}
\end{Lem}

\begin{proof}  
\emph{1.} Avec les notations de la démonstration de \ref{troncation}, $L\otimes_R\Lambda_{n+m}$ est alors lisse sur $\Lambda_{n+m}$, et comme l'idéal $I\Lambda_{n+m}$ est de carré nul, la flèche $r_{n+m,m}$ est surjective dans ce cas. On obtient donc une suite exacte
$$
\xymatrix{
0\ar[r] & \Gr_nV(\omega_L)\ar[r] & \Gr_{n+m}(L)\ar[r] & \Gr_{m}(L)\ar[r] & 0
}
$$
sur $(\Sch/k)_{\ZAR}$, en particulier une suite exacte
$$
\xymatrix{
0\ar[r] & V(\omega_{L_k})\ar[r] & \Gr_{n+1}(L)\ar[r] & \Gr_n(L)\ar[r] & 0
}
$$
puisque $\Gr_1V(\omega_L):=V(\omega_L)_k=V(\omega_{L_k})$. 

\medskip

\emph{2.} Le $k$-schéma en groupes $V(\omega_{L_k})$ est isomorphe à $\bbG_{a,k}^d$, où $d:=\dim L_k$. Pour $n\in\bbN^*$ fixé, on a donc une suite exacte
$$
\xymatrix{
0 \ar[r] & U_n \ar[r] & \Gr_{n}(L) \ar[r] & L_k \ar[r] & 0,
}
$$
où $U_n$ possède une suite de composition de longueur $n-1$ dont les quotients successifs sont isomorphes à $\bbG_{a,k}^d$. 

\medskip

\emph{3.} Comme $V(\omega_{L_k})$ est égal à $V(\omega_{L_{k}^0})$ et est connexe, on a un diagramme commutatif exact
$$
\xymatrix{
0\ar[r] & V(\omega_{L_k})\ar[r] \ar@{=}[d] & \Gr_{n+1}(L^0)\ar[r] \ar[d] & \Gr_n(L^0)\ar[r] \ar[d] & 0 \\
0\ar[r] & V(\omega_{L_k})\ar[r]  & \Gr_{n+1}(L)\ar[r] & \Gr_n(L)\ar[r] & 0 \\
0\ar[r] & V(\omega_{L_k})\ar[r] \ar@{=}[u] & \Gr_{n+1}(L)^0\ar[r] \ar[u] & \Gr_n(L)^0\ar[r] \ar[u] & 0. 
}
$$
On en déduit par récurrence sur $n\geq 1$ que $\Gr_n(L^0)\xrightarrow{\sim}\Gr_n(L)^0$. La connexité de $V(\omega_{L_k})$ entraîne également l'isomorphisme $\pi_0(\Gr_{n+1}(L)) \xrightarrow{\sim} \pi_0(\Gr_n(L))$ à partir de la suite exacte du milieu.
\end{proof}

\begin{Lem}\label{Grsurjlisse}
Soit $f:L\ra L'$ un morphisme surjectif lisse de $R$-schéma en groupes lisses. Alors pour tout $n\geq 1$, le morphisme 
$\Gr_n(f):\Gr_n(L)\ra\Gr_n(L')$ est surjectif lisse.
\end{Lem}

\begin{proof}
Pour tout $n\geq 1$, on a d'après \ref{troncationlisse} le diagramme commutatif exact
$$
\xymatrix{
0\ar[r] & V(\omega_{L_k}) \ar[r] \ar[d]^{V(f_{k}^*)} & \Gr_{n+1}(L) \ar[r] \ar[d]^{\Gr_{n+1}(f)} & \Gr_n(L) \ar[r] \ar[d]^{\Gr_n(f)} & 0 \\
0\ar[r] & V(\omega_{L_{k}'}) \ar[r] & \Gr_{n+1}(L') \ar[r] & \Gr_n(L') \ar[r] & 0.
}
$$
Comme $f$ est surjectif, $\Gr_1(f)=f_k$ est surjectif. Comme $f$, et donc $f_k$, est lisse, $\Lie(f_k):\Lie(L_k)\ra\Lie(L_{k'})$ est surjectif, donc
$V(f_{k}^*)$ est surjectif. D'où la surjectivité de $\Gr_n(f)$ par récurrence sur $n$. Notant $N$ le noyau de $f$, on obtient donc pour tout $n\geq 1$ la suite exacte
$$
\xymatrix{
0\ar[r] & \Gr_n(N) \ar[r] & \Gr_n(L) \ar[r]^{\Gr_n(f)} & \Gr_n(L') \ar[r] & 0.
}
$$
Comme $\Gr_n(N)$ est lisse sur $k$ d'après \ref{troncationlisse}, l'épimorphisme $\Gr_n(f)$ est lisse.
\end{proof}

\begin{Prop} \label{imrelisse}
Soit $\alpha:L\ra L'$ un morphisme de $R$-schémas en groupes lisses séparés. Supposons que le morphisme $\alpha_K:L_K\ra L_{K}'$ induit sur les fibres génériques soit lisse, et que pour tout $n\geq 1$, le conoyau du morphisme $\Gr_n(\alpha):\Gr_n(L)\ra \Gr_n(L')$ soit représentable. Alors pour tout $n$ assez grand, le morphisme canonique 
$$
\xymatrix{
\Coker \Gr_{n+1}(\alpha)\ar[r] & \Coker \Gr_n(\alpha)
}
$$
est un isomorphisme.
\end{Prop}

Pour démontrer cette proposition, nous allons utiliser le résultat topologique suivant, cf. \cite{MB} 2.2.1 (i) :

\begin{Lem}\label{lemimrelisse}
Soit $\alpha_K:X_K\ra Y_K$ un morphisme lisse de $K$-schémas localement de type fini. Alors $\alpha_K(K):X_K(K)\ra Y_K(K)$ est ouvert pour la topologie forte sur $X_K(K)$ et $Y_K(K)$.
\end{Lem}

\begin{proof}
Le lemme est vrai si $f_K$ est le morphisme structural d'un espace affine standard au-dessus de $Y_K$, ou si $f_K$ est étale  d'après le théorème d'inversion local. Le cas général en résulte.
\end{proof}

\begin{proof}[Démonstration de la proposition \ref{imrelisse}.]
Comme $\alpha_K:L_{K}\ra L_{K}'$ est lisse par hypothèse, le morphisme $\alpha_K(K):L_{K}(K)\ra L_{K}'(K)$ est ouvert d'après \ref{lemimrelisse}. En particulier, l'image du morphisme $\alpha(R):L(R)\ra L'(R)$ est un voisinage ouvert de $0$ dans $L'(R)$. Il existe donc un entier $n_0\geq 1$ tel que $\Ima \alpha(R) \supset \cV_{n_0}(L'):=\Ker(L'(R)\ra L'(R/\fm^{n_0}))$. Soit $n\geq n_0$ et considérons le diagramme commutatif exact
$$
\xymatrix{
\Gr_n(L')\ar[r]^<<<<<{=:q_n} \ar[d]^{\tau_{n,n_0}} & \Coker \Gr_n(\alpha) \ar[r] \ar[d]^{\overline{\tau_{n,n_0}}}  & 0 \\
\Gr_{n_0}(L') \ar[r]^<<<<<{=:q_{n_0}} & \Coker \Gr_{n_0}(\alpha) \ar[r] & 0.
}
$$
D'après \ref{troncationlisse} \emph{1.}, $\tau_{n,n_0}$ est un épimorphisme. Nous allons montrer que la projection canonique $q_n$ se factorise par $\tau_{n,n_0}$. Comme le noyau de $\tau_{n,n_0}$ est lisse d'après \emph{loc. cit.} et $k$ est algébriquement clos, il suffit de montrer que 
$$
\xymatrix{
q_n(k):L'(R/\fm^n) \ar[r] & L'(R/\fm^n)/\Ima \alpha(R/\fm^n) 
}
$$
s'annule sur le noyau de $\tau_{n,n_0}(k):L'(R/\fm^n)\ra L'(R/\fm^{n_0})$. Comme $L'(R)\ra L'(R/\fm^n)$ est surjectif (\ref{Hensel}), ce noyau coïncide avec l'image du morphisme composé canonique
$$
\xymatrix{
\cV_{n_0}(L') \ar[r] & L'(R) \ar[r] & L'(R/\fm^n). 
}
$$
Par définition de $n_0$, le diagramme commutatif
$$
\xymatrix{
& & \cV_{n_0}(L') \ar[d] \\
L(R) \ar[rr]^{\alpha(R)} \ar[d] & &  L'(R) \ar[d] \\
L(R/\fm^n) \ar[rr]^{\alpha(R/\fm^n)} & & L'(R/\fm^n)
}
$$
montre que $\Ima \alpha(R/\fm^n) \supset \Ima (\cV_{n_0}(L')\ra L'(R/\fm^n))$, i.e. que $q_n(k)$ s'annule bien sur $\Ima (\cV_{n_0}(L')\ra L'(R/\fm^n))$.

On dispose ainsi d'un système projectif de suites exactes
\begin{equation} \label{systprojmn}
\xymatrix{
0\ar[r] & \Ker(\overline{q_m}) \ar[r] \ar[d] &\Gr_{n_0}(L') \ar[r]^<<<<<{\overline{q_m}} \ar@{=}[d] &  \Coker \Gr_m(\alpha) \ar[r] \ar[d] & 0 \\
0\ar[r] & \Ker(\overline{q_n}) \ar[r] &\Gr_{n_0}(L') \ar[r]^<<<<<{\overline{q_n}} &  \Coker \Gr_n(\alpha) \ar[r] & 0, & 
m\geq n\geq n_0.
}
\end{equation}
Comme la composante neutre de $\Gr_{n_0}(L')$ est un schéma noethérien, la suite 
$$
\big(\Ker(\overline{q_n})\cap\Gr_{n_0}(L')^0\big)_{n\geq n_0}
$$ 
est stationnaire. Par ailleurs, pour tout $n\geq 1$, le diagramme commutatif exact (\ref{troncationlisse} \emph{1.})
$$
\xymatrix{
0\ar[r] & V(\omega_{L_k}) \ar[r] \ar[d]^{V(\alpha_{k}^*)} & \Gr_{n+1}(L) \ar[r] \ar[d]^{\Gr_{n+1}(\alpha)} & \Gr_n(L) \ar[r] \ar[d]^{\Gr_n(\alpha)} & 0 \\
0\ar[r] & V(\omega_{L_{k}'}) \ar[r] & \Gr_{n+1}(L') \ar[r] & \Gr_n(L') \ar[r] & 0 
}
$$
montre, avec le lemme du serpent, que $\Coker(\Gr_{n+1}(\alpha))\ra\Coker(\Gr_n(\alpha))$ est un épimorphisme à noyau connexe, et donc que l'application $\pi_0(\Coker(\Gr_{n+1}(\alpha)))\ra\pi_0(\Coker(\Gr_n(\alpha)))$ qui s'en déduit est bijective. En appliquant $\pi_0$ au diagramme (\ref{systprojmn}), on obtient donc un système projectif de suites exactes
$$
\xymatrix{
0\ar[r] & \pi_0(\Ker(\overline{q_m})\cap\Gr_{n_0}(L')^0) \ar[r] \ar[d] & \pi_0(\Ker(\overline{q_m})) \ar[r] \ar[d] &\pi_0(\Gr_{n_0}(L'))  \ar@{=}[d]  \\
0\ar[r] & \pi_0(\Ker(\overline{q_n})\cap\Gr_{n_0}(L')^0) \ar[r] & \pi_0(\Ker(\overline{q_n})) \ar[r] & \pi_0(\Gr_{n_0}(L'))  \\
& \pi_0(\Gr_{n_0}(L')) \ar[r]^<<<<<<<<{\pi_0(\overline{q_m})} \ar@{=}[d] &  \pi_0(\Coker \Gr_m(\alpha)) \ar[r] \ar[d]^{\rotatebox{90}{$\sim$}}  & 0 \\
& \pi_0(\Gr_{n_0}(L')) \ar[r]^<<<<<<<<{\pi_0(\overline{q_n})} &  \pi_0(\Coker\Gr_n(\alpha)) \ar[r] & 0,
}
$$
$m\geq n\geq n_0$, dont la quatrième flèche verticale est bijective et la première bijective pour tout $n$ assez grand ; la seconde flèche verticale est donc bijective pour tout $n$ assez grand. On en déduit finalement que la suite $(\Ker(\overline{q_n}))_{n\geq n_0}$ est stationnaire, et donc que $\Coker(\Gr_{n+1}(\alpha))\ra\Coker(\Gr_n(\alpha))$ est un isomorphisme pour tout $n$ assez grand.
\end{proof}

\section{Le lemme de Hensel}

\begin{Pt}
Soit $G$ un $R$-schéma en groupes fidèlement plat localement de type fini. Comme $\Omega_{G/R}^1$ est invariant par les translations de $G$, ses idéaux de Fitting aussi ; ceux-ci proviennent donc de $R$.
\end{Pt}

\begin{Def}\label{Defd}
Soit $G$ un $R$-schéma en groupes fidèlement plat localement de type fini, de dimension relative $r$. La \emph{différente absolue} de $G$ sur $R$ est l'idéal de $R$ engendrant le $r$-ième idéal de Fitting de $\Omega_{G/R}^1$. 
\end{Def}

\begin{Pt}\label{calculd}
Si de plus $G$ est fini, la différente absolue de $G$ coïncide avec celle de \cite{R} Appendice, Définition 8. La différente absolue de $G$ est non nulle si et seulement si la fibre générique de $G$ est lisse, auquel cas sa valuation est égale au défaut de lissité de $G$ le long de chacune de ses sections (\cite{BLR} 3.3/2). En particulier, si $G$ est lisse sur $R$, sa différente absolue est l'idéal unité de $R$.
\end{Pt}

\begin{Lem}\label{Hensel}
Soit $G$ un $R$-schéma en groupes fidèlement plat localement de type fini, à fibre générique lisse. Soit $d\in\bbN$ la valuation de la différente absolue de $G$ sur $R$. Alors pour tout $n\geq d+1$, on a
$$
\Ima\big(G(R)\lra G(R/\fm^n)\big)=\Ima\big(G(R/\fm^{n+d})\lra G(R/\fm^n)\big).
$$
En particulier, si $G$ est lisse sur $R$, alors pour tout $n\geq 1$, l'application de réduction
$$
\xymatrix{
G(R)\ar[r] & G(R/\fm^n)
}
$$ 
est surjective.
\end{Lem}

\begin{proof}
Il s'agit d'un cas particulier du lemme de Hensel pour les intersections complètes locales \cite{CLNS} 4.2.1. En effet, notant $D$ le sous-schéma fermé de $G$ défini par l'idéal $\fm^d\cO_G$, l'ensemble $D(R/\fm^{d+1})$ est vide. Il résulte donc de \emph{loc. cit.} que pour tout $n\geq d+1$, l'image dans $G(R/\fm^n)$ d'un élément de $G(R/\fm^{n+d})$ peut être relevée dans $G(R)$. Lorsque $G$ est lisse sur $R$, le lemme résulte aussi directement de la définition de la lissité formelle et du fait que $R$ est complet.
\end{proof}

\begin{Rem}\label{cteGreen} 
En particulier, pour tout $n\geq 1$, on a
$$
\Ima\big(G(R)\lra G(R/\fm^n)\big)=\Ima\big(G(R/\fm^{n+2d})\lra G(R/\fm^n)\big).
$$
La constante de Greenberg de $G$ (cf. \cite{CLNS} 4.2.8) est donc égale à $1$.
\end{Rem}

\chapter{La cohomologie des $R$-schémas en groupes finis plats}

Soient $R$ un anneau de valuation discrète complet de caractéristique $p>0$ et $k$ son corps résiduel, que l'on suppose algébriquement clos. Pour tout $n\geq 1$, on pose $R_n:=R/\fm^n$.

\section{Nullité du $H^i$ pour $i>1$}

\begin{Lem}\label{dimcohlR} \label{dimcohlRn}
Soit $L$ un $R$-schéma en groupes lisse. Alors 
\begin{enumerate}
\item $\forall n\geq 1,\quad\forall i>0,\quad H^i(R_n,L)=0 ;$
\item $\forall i>0,\quad H^i(R,L)=0.$
\end{enumerate}
\end{Lem}

\begin{proof}
La cohomologie d'un schéma en groupes lisse sur une base $S$ calculée sur le gros site fppf de $S$ coïncide avec la cohomologie de l'image directe de $L$ sur le petit site étale de $S$, et les $R_n$, $\geq 1$, ainsi que $R$, sont strictement hensélien par hypothèse.
\end{proof}

\begin{Cor} \label{dimcohfR} \label{dimcohfRn}
Soit $G$ un $R$-schéma en groupes fini plat. Alors 
\begin{enumerate}
\item $\forall n\geq 1,\quad\forall i>1,\quad H^i(R_n,G)=0 ;$
\item $\forall i>1,\quad H^i(R,G)=0.$
\end{enumerate}
\end{Cor}

\begin{proof}
Un $R$-schéma en groupes fini plat $G$ possède des résolutions lisses de longueur $1$, par exemple sa résolution lisse standard (\cite{Beg} 2.2.1), et les composantes d'une telle résolution sont $\Gamma(R_n,\cdot)$ et $\Gamma(R,\cdot)$-acycliques d'après \ref{dimcohlR}. 
\end{proof}

\section{Structure quasi-algébrique sur le $H^0$}

\begin{Def}\label{DefcH0R}
Soit $G$ un $R$-schéma en groupes fini plat. La \emph{structure de groupe quasi-algébrique sur $H^0(R,G)$} est l'unique groupe quasi-algébrique fini sur $k$ dont le groupe ordinaire sous-jacent est $H^0(R,G)$. 
\end{Def}

\begin{Not}\label{NotcH0R}
On note $\cH^0(G)$ la structure de groupe quasi-algébrique sur $H^0(R,G)$.
\end{Not}

\begin{Pt} \label{cH0foncR}
On a ainsi un \emph{foncteur}
$$
\xymatrix{
\cH^0:\textrm{($R$-schémas en groupes finis plats)} \ar[r] & \cQ.
}
$$
\end{Pt}

\begin{Lem}\label{schfiniplatpt}
Soient $X$ un $R$-schéma séparé à fibre générique finie, et $k'/k$ une extension de corps. Alors l'application canonique
$$
\xymatrix{
X(R)\ar[r] & X_{R_{k'}}(R_{k'})
}
$$
(notation \ref{Kk'}) est bijective.
\end{Lem}

\begin{proof}
D'après \ref{schfinipt} l'application canonique $X_K(K)\ra X_{K_{k'}}(K_{k'})$ est bijective. Comme $X$ est séparé sur 
$R$, on en déduit que l'application considérée ici est injective. Comme de plus $R_{k'}\cap K=R$, on en déduit aussi qu'elle est surjective.
\end{proof}

\begin{Cor}\label{cH0cbR}
Si $k'/k$ est une extension de corps parfaite de $k$, alors on a un isomorphisme canonique
$$
\xymatrix{
\cH^0(G)_{k'}\ar[r] & \cH^0(G_{R_{k'}}).
}
$$
\end{Cor}

\section{Structures quasi-algébrique sur le $H^1$ mod $\fm^n$ et proalgébrique sur le $H^1$}

\begin{Pt}
On a défini en \ref{green} le foncteur de Greenberg de cran $n\geq 1$
$$
\xymatrix{
\textrm{($R$-schémas en groupes localement de type fini)} \ar[d]^{\Gr_n} \\
\textrm{($k$-schémas en groupes localement de type fini)},
}
$$
et on dispose par ailleurs du foncteur de perfectisation
$$
\xymatrix{
\textrm{($k$-schémas en groupes localement de type fini)} \ar[d]^{(\cdot)^{\pf}} \\
\textrm{($k$-schémas en groupes parfaits localement algébriques)}.
}
$$
En sous-entendant les plongements de catégories évidents, on peut donc considérer le foncteur composé
$$
\xymatrix{
(\textrm{$R$-schémas en groupes lisses}) \ar[d]^{(\cdot)^{\pf}\circ\Gr_n} \\
\textrm{($k$-schémas en groupes parfaits localement algébriques)},
}
$$
que l'on notera $\pi_{n*}$.
\end{Pt} 

\begin{Pt}\label{pi*kk'Rn}
D'après \ref{greencb}, le foncteur $\pi_{n*}$ commute à toute extension de corps parfaite $k'/k$.
\end{Pt}

\begin{Lem}\label{proppi*Rn}
Soit $L$ un $R$-schéma en groupes lisse. On a $(\pi_{n*}L)(k)=L(R_n)$. De plus, 
$$
\big(\Gr_n(L)^{\pf}(k)\big)_{n\geq 1}=\big(\Gr_n(L)(k)\big)_{n\geq 1}\in\Pro(\textrm{\emph{Groupes}})
$$ 
est à transitions surjectives, et pour tous $m\geq n\geq 1$, le noyau de la transition 
$$\Gr_m(L)^{\pf}\ra \Gr_n(L)^{\pf}$$
est un $k$-schéma en groupes parfait connexe unipotent. En particulier, pour tous $m\geq n\geq 1$, les applications canoniques 
$$\pi_0(\Gr_m(L))\ra \pi_0(\Gr_n(L)) \ra \pi_0(L_k)$$
sont bijectives.
\end{Lem}

\begin{proof}
Résulte de \ref{troncationlisse} et \ref{pfexact}.
\end{proof}

\begin{Prop}\label{cH1adRn}
Soit $G$ un $R$-schéma en groupes fini plat. Alors pour tout $n\geq 1$, il existe un unique groupe quasi-algébrique 
$\cH^1(G_{R_n})$ de fibre générique (\ref{deffibgenI})
\begin{eqnarray*}
\cT^{\circ} & \ra & (\textrm{\emph{Groupes}}) \\
\kappa & \ra & H^1((R_{\kappa})_n,G) \ ;
\end{eqnarray*}
il est connexe unipotent. Le morphisme $\cH^1(G_{R_{n+1}})_{\eta}\ra \cH^1(G_{R_n})_{\eta}$ provient d'un unique morphisme
$$
\xymatrix{
\cH^1(G_{R_{n+1}})\ar[r] & \cH^1(G_{R_n}),
}
$$
qui est un épimorphisme dont le noyau est un groupe quasi-algébrique vectoriel.

Si $f:L_0\ra L_1$ est un morphisme fidèlement plat de $R$-schémas en groupes lisses de noyau $G$, alors l'image de 
$\cH^1(G_{R_n})$ dans 
$$\textrm{\emph{\big(Faisceaux abéliens sur $(\parf/k)_{\et}$\big)}}$$
réalise le conoyau de
$$
\pi_{n*}f:\pi_{n*}L_0\ra \pi_{*n}L_1.
$$
\end{Prop}

\begin{proof}
L'unicité résulte de \ref{ratIP0}. Pour l'existence, choisissons un morphisme fidèlement plat de $R$-schémas en groupes lisses $f:L_0\ra L_1$ de noyau $G$ (par exemple, la résolution lisse standard). Pour tout $n\geq 1$, $\Gr_n(f):\Gr_n(L_0)\ra\Gr_n(L_1)$ est affine puisque $f$ l'est, de sorte que $\Coker(\Gr_n(f))$ est représentable par un $k$-schéma en groupes localement de type fini. Pour toute extension de corps algébriquement close $\kappa/k$, $\Coker(\Gr_n(f))(\kappa)$ est le conoyau de $f((R_{\kappa})_n):L_0((R_{\kappa})_n)\ra L_1((R_{\kappa})_n)$, qui s'identifie à $H^1((R_{\kappa})_n,G)$ d'après \ref{dimcohlRn}. Pour tout $n\geq 1$, on a d'après \ref{troncationlisse} le diagramme commutatif 
$$
\xymatrix{
0\ar[r] & V(\omega_{L_{0,k}}) \ar[r] \ar[d] & \Gr_{n+1}(L_0)\ar[r] \ar[d]^{\Gr_{n+1}(f)} & \Gr_n(L_0) \ar[r] \ar[d]^{\Gr_n(f)} & 0 \\ 
0\ar[r] & V(\omega_{L_{1,k}}) \ar[r] & \Gr_{n+1}(L_1) \ar[r] & \Gr_n(L_1) \ar[r] & 0,
}
$$
qui montre via le lemme du serpent que $\Coker(\Gr_{n+1}(f))$ est extension de $\Coker(\Gr_n(f))$ par un $k$-groupe vectoriel. Comme $f_k:L_{0,k}\ra L_{1,k}'$ est un épimorphisme, i.e. $\Coker(\Gr_1(f))=0$, on en déduit que pour tout $n\geq 1$, $\Coker(\Gr_n(f))$ est (lisse) connexe unipotent. D'où le lemme d'après \ref{pfexact} et \ref{Qgppf}.
\end{proof}

\begin{Def}\label{DefcH1Rn}
Soit $G$ un $R$-schéma en groupes fini plat. La \emph{structure de groupe quasi-algébrique sur $H^1(R_n,G)$} est le groupe quasi-agébrique sur $k$ défini par \ref{cH1adRn}.
\end{Def}

\begin{Not}\label{NotcH1Rn}
On note $\cH^1(G_{R_n})$ la structure de groupe quasi-algébrique sur $H^1(R_n,G)$. Le groupe ordinaire sous-jacent à 
$\cH^1(G_{R_n})$ est $H^1(R_n,G)$, et plus généralement $\cH^1(G_{R_n})(\kappa)=H^1((R_{\kappa})_n,G)$ pour toute extension de corps algébriquement close $\kappa/k$.
\end{Not}

\begin{Pt}\label{cH1cbRn}
D'après \ref{cH1adRn} et \ref{RemratIP0} (ou \ref{pi*kk'Rn}), pour toute extension de corps parfaite $k'/k$, on a un isomorphisme canonique
$$
\xymatrix{
\cH^1(G_{R_n})_{k'}\ar[r]^{\sim} & \cH^1(G_{(R_{k'})_n}).
}
$$
\end{Pt}

\begin{Pt}\label{cH1foncQRn}
D'après \ref{cH1adRn} et \ref{ratIP0}, la formation de $\cH^1(G_{R_n})$ est fonctorielle en $G$. On a ainsi défini un \emph{foncteur}
$$
\xymatrix{
\cH^1((\cdot)_{R_n}):\textrm{($R$-schémas en groupes finis plats)} \ar[r] & \cU^0,
}
$$
induisant sur les points à valeurs dans une extension algébriquement close $\kappa/k$ arbitraire le foncteur
$$
\xymatrix{
H^1((R_{\kappa})_n,\cdot):\textrm{($R$-schémas en groupes finis plats)} \ar[r] & \textrm{(Groupes)}.
}
$$
\end{Pt}

\begin{Pt}\label{preDefcH1R}
Soit $G$ un $R$-schéma en groupes fini plat. D'après \ref{cH1adRn}, les transitions du pro-objet 
$$
\big(\cH^1(G_{R_n})\big)_{n\geq 1}\in\Pro(\cU^0)
$$
sont surjectives à noyau connexe. Passant à la limite projective dans $\cP$, on obtient donc un groupe 
proalgébrique unipotent strictement connexe
$$
\lp_{n\geq 1}\cH^1(G_{R_n}) \in \cP^{00}(\cU).
$$ 
Il a pour fibre générique le foncteur  
\begin{eqnarray*}
\cT^{\circ} & \ra & (\textrm{Groupes}) \\
\kappa & \ra & H^1(R_{\kappa},G) \ ;
\end{eqnarray*}
en effet si $f:L_0\ra L_1$ est une résolution lisse de $G$ de longueur 1, alors les transitions du pro-groupe $(L_0((R_{\kappa})_n))_{n\geq 1}$ sont surjectives, de sorte que la suite
$$
\xymatrix{
\lp_{n\geq 1}\Gr_n(L_0)(\kappa) \ar[rr]^{\Gr(f)(\kappa)} && \lp_{n\geq 1}\Gr_n(L_1)(\kappa) \ar[r] & \lp_{n\geq 1}\Coker(\Gr_n(f))(\kappa) \ar[r] & 0
}
$$
est exacte, identifiant donc le groupe $\lp_{n\geq 1}\Coker(\Gr_n(f))(\kappa)$ au conoyau de $f(R_{\kappa}):L_0(R_{\kappa})\ra L_1(R_{\kappa})$, qui est canoniquement isomorphe à $H^1(R_{\kappa},G)$ d'après \ref{dimcohlR}.
\end{Pt}

\begin{Def}\label{DefcH1R}
Soit $G$ un $R$-schéma en groupes fini plat. La \emph{structure de groupe proalgébrique sur $H^1(R,G)$} est le groupe proalgébrique sur $k$ $\lp_{n\geq 1}\cH^1(G_{R_n})$.
\end{Def}

\begin{Not}\label{NotcH1R}
On note $\cH^1(G)$ la structure de groupe proalgébrique sur $H^1(R,G)$. Le groupe ordinaire sous-jacent à $\cH^1(G)$ est $H^1(R,G)$, et plus généralement $\cH^1(G)(\kappa)=H^1(R_{\kappa},G)$ pour toute extension de corps algébriquement close $\kappa/k$.
\end{Not}

\begin{Pt}\label{cH1cbR}
D'après \ref{cH1cbRn}, si $k'/k$ est une extension de corps parfaite, alors on a un isomorphisme canonique
$$
\xymatrix{
\cH^1(G_R)_{k'}\ar[r]^{\sim} & \cH^1(G_{R_{k'}}).
}
$$
\end{Pt}

\begin{Pt}\label{cH1foncQR}
D'après \ref{cH1foncQRn}, les $\cH^1(G)$ pour $G$ variable définissent un \emph{foncteur}
$$
\xymatrix{
\cH^1:\textrm{($R$-schémas en groupes finis plats)} \ar[r] & \cP^{00}(\cU),
}
$$
induisant sur les points à valeurs dans une extension algébriquement close $\kappa/k$ arbitraire le foncteur
$$
\xymatrix{
H^1(R_{\kappa},\cdot):\textrm{($R$-schémas en groupes finis plats)} \ar[r] & \textrm{(Groupes)}.
}
$$
\end{Pt}

\begin{Pt}\label{cH1foncQRbis}
Le foncteur $\cH^1$ \ref{cH1foncQR} peut se calculer de manière analogue au foncteur $\cH^1((\cdot)_{R_n})$ \ref{cH1foncQRn}. On a en effet défini en \ref{green} le foncteur de Greenberg
$$
\xymatrix{
\textrm{($R$-schémas en groupes localement de type fini)} \ar[d]^{\Gr} \\
\Pro(\textrm{$k$-schémas en groupes localement de type fini)},
}
$$
et le foncteur de perfectisation se prolonge de manière unique en un foncteur
$$
\xymatrix{
\Pro\textrm{($k$-schémas en groupes de localement type fini)}) \ar[d]^{\Pro((\cdot)^{\pf})} \\
\Pro\textrm{($k$-schémas en groupes parfaits localement algébriques)}.
}
$$
En sous-entendant les plongements de catégories évidents, on peut donc considérer le foncteur composé
$$
\xymatrix{
(\textrm{$R$-schémas en groupes lisses}) \ar[d]^{\Pro((\cdot)^{\pf})\circ\Gr} \\
\Pro\textrm{($k$-schémas en groupes parfaits localement algébriques)}),
}
$$
que l'on notera $\pi_*$. Alors si $G$ est un $R$-schéma en groupes fini plat, et $L_0\ra L_1$ une résolution lisse de $G$ de longueur 1, le conoyau de
$$
\pi_*f:\pi_*L_0\ra \pi_*L_1
$$
dans $\Pro\textrm{(Faisceaux abéliens sur $(\parf/k)_{\et}$)}$ est représenté par l'image de $\cH^1(G)$.
\end{Pt}

\begin{Exs}\label{ExsR}

\medskip

\emph{1. $G=\mu_{p,R}$.} 
On peut choisir pour résolution de $\mu_{p,R}$ la suite exacte de Kummer 
$$
\xymatrix{
0\ar[r] & \mu_{p,R} \ar[r] & \bbG_{m,R} \ar[r]^p & \bbG_{m,R} \ar[r] & 0.
}
$$
Pour tout $n\geq 1$, le diagramme commutatif exact
$$
\xymatrix{
0 \ar[r] & V(\omega_{\bbG_{m,k}}) \ar[r] \ar[d]^0 & \Gr_{n+1}(\bbG_{m,R}) \ar[r]^{\tau_{n+1,n}} \ar[d]^p & \Gr_n(\bbG_{m,R}) \ar[r] \ar[d]^p & 0\\
0 \ar[r] & V(\omega_{\bbG_{m,k}}) \ar[r] & \Gr_{n+1}(\bbG_{m,R}) \ar[r]^{\tau_{n+1,n}} & \Gr_n(\bbG_{m,R}) \ar[r] & 0
}
$$
induit un diagramme commutatif exact
$$
\xymatrix{
           &    & 0                                                                     \ar[d] & 0                                        \ar[d] & \\
           &   & p\Gr_{n+1}(\bbG_{m,R}) \ar[r]^{\tau_{n+1,n}} \ar[d] & p\Gr_n(\bbG_{m,R}) \ar[r] \ar[d] & 0\\
0 \ar[r] & V(\omega_{\bbG_{m,k}}) \ar[r] \ar[d] & \Gr_{n+1}(\bbG_{m,R}) \ar[r]^{\tau_{n+1,n}} \ar[d] & \Gr_n(\bbG_{m,R}) \ar[r] \ar[d] & 0\\
0 \ar[r] & \Ker(\overline{\tau_{n+1,n}}) \ar[r] \ar[d] & \Gr_{n+1}(\bbG_{m,R})/p \ar[r]^{\overline{\tau_{n+1,n}}} \ar[d] & \Gr_n(\bbG_{m,R})/p \ar[r] \ar[d] & 0 \\
           & 0                                                              & 0                                                                                                & 0 &
}
$$
Si $p$ ne divise pas $n$, le morphisme de réduction
$$
\tau_{n+1,n}(k):p\Gr_{n+1}(\bbG_{m,R})(k)=\big((R/\fm^{n+1})^*\big)^p \ra p\Gr_n(\bbG_{m,R})(k)=\big((R/\fm^n)^*\big)^p 
$$
est injectif, et $\cH^1(\mu_{p,R_{n+1}})$ est donc extension de $\cH^1(\mu_{p,R_n})$ par $\bbG_{a,k}^{\pf}$. Si $p$ divise $n$, $p\Gr_{n+1}(\bbG_{m,R})$ contient l'image de $V(\omega_{\bbG_{m,k}})$, et le morphisme canonique $\cH^1(\mu_{p,R_{n+1}})\ra\cH^1(\mu_{p,R_n})$ est alors un isomorphisme. Par conséquent, 
$\cH^1(\mu_{p,R})$ est un $k$-groupe provectoriel de dimension infinie.

\medskip

\emph{2. $G=\alpha_{p,R}$.} 
On peut choisir pour résolution de $\alpha_{p,R}$ le Frobenius relatif de 
$\bbG_{a,R}$ 
$$
\xymatrix{
0\ar[r] & \alpha_{p,R} \ar[r] & \bbG_{a,R} \ar[r]^{F} & \bbG_{a,R} \ar[r] & 0.
}
$$
Pour tout $n\geq 1$, le diagramme commutatif exact
$$
\xymatrix{
0 \ar[r] & V(\omega_{\bbG_{a,k}}) \ar[r] \ar[d]^0 & \Gr_{n+1}(\bbG_{a,R}) \ar[r] \ar[d]^{\Gr_{n+1}(F)} & \Gr_n(\bbG_{a,R}) \ar[r] \ar[d]^{\Gr_n(F)} & 0 \\
0 \ar[r] & V(\omega_{\bbG_{a,k}}) \ar[r] & \Gr_{n+1}(\bbG_{a,R}) \ar[r] & \Gr_n(\bbG_{a,R}) \ar[r] & 0
}
$$
s'identifie au diagramme
$$
\xymatrix{
0 \ar[r] & \bbG_{a,k} \ar[r] \ar[d]^0& \bbG_{a,k}^{n+1} \ar[r] \ar[d]^{F_{n+1}} & \bbG_{a,k}^n \ar[r] \ar[d]^{F_n} & 0 \\
0 \ar[r] & \bbG_{a,k} \ar[r] & \bbG_{a,k}^{n+1} \ar[r] & \bbG_{a,k}^n \ar[r] & 0
}
$$
où 
$$
F_n=F\times \underbrace{0\times\ldots \times 0}_{p-1 \textrm{ fois}}\times F\times \underbrace{0\times\ldots \times 0}_{p-1 \textrm{ fois}}\times F\ldots F\times\underbrace{0\times\ldots\times0}_{n-p\lfloor\frac{n-1}{p}\rfloor -1 \textrm{ fois}},
$$
le Frobenius relatif $F$ de $\bbG_{a,k}$ apparaissant $\lfloor (n-1)/p \rfloor +1$ fois. Pour tout $n\geq 1$, $\cH^1(\alpha_{p,R_n})$ est donc un $k$-groupe quasi-algébrique vectoriel de dimension 
$n-\lfloor (n-1)/p \rfloor -1$, et le morphisme canonique $\cH^1(\alpha_{p,R_{n+1}})\ra \cH^1(\alpha_{p,R_n})$ est un épimorphisme de noyau $\bbG_{a,k}^{\pf}$ si $p$ ne divise pas $n$, un isomorphisme si $p$ divise $n$. Par conséquent 
$\cH^1(\alpha_{p,R})$ est un $k$-groupe provectoriel de dimension infinie.
\end{Exs}

Lorsque que le $R$-schéma en groupes fini plat $G$ est à fibre générique étale, on a comme en inégales caractéristiques \cite{Beg} 4.2.2 :

\begin{Prop} \label{fibgenetQ}
Soit $G$ un $R$-schéma en groupes fini plat, à fibre générique étale. Soit $d\in\bbN$ la valuation de la différente absolue de $G$ (\ref{Defd}). Alors la projection canonique
$$
\xymatrix{
\cH^1(G)\ar[r] & \cH^1(G_{R_{2d+1}})
}
$$
est un isomorphisme. La dimension de $\cH^1(G_{R_{2d+1}})$ est $d$.
\end{Prop}

\begin{proof}
Soit $f:L\ra L'$ une résolution lisse de $G$ de longueur 1. D'après \ref{troncation} et \ref{troncationlisse}, on a pour tout 
$n\geq d$ un diagramme commutatif exact 
$$
\xymatrix{
               & 0 \ar[d]                                        & 0 \ar[d]                                &   0 \ar[d]                              & & \\
0 \ar[r]    & \Gr_d(V(\omega_G)) \ar[d]\ar[r]         & \Gr_{n+d}(G) \ar[d]\ar[r] &   \Gr_n(G) \ar[d]                   & & \\
0 \ar[r]   & \Gr_d(V(\omega_L)) \ar[r] \ar[d] & \Gr_{n+d}(L) \ar[r] \ar[d]^{\Gr_{n+d}(f)} & \Gr_n(L) \ar[r] \ar[d]^{\Gr_n(f)} & 0 & \\
0 \ar[r]   & \Gr_d(V(\omega_{L'}))  \ar[r]       & \Gr_{n+d}(L')\ar[r] \ar[d]     & \Gr_n(L')  \ar[r] \ar[d]                 &  0 & \\
     &                                             &  \Coker(\Gr_{n+d}(f))  \ar[d] \ar[r]              & \Coker(\Gr_n(f))  \ar[d] \ar[r] & 0 & \\
     &                                             &  0                                     & 0.                        & &                                
}
$$
D'après \ref{Hensel}, pour tout $n\geq d+1$,
$$
\Ima\big(\Gr_{n+d}(G)(k)\lra \Gr_n(G)(k)\big)=\Ima\big(\Gr(G)(k)\lra \Gr_n(G)(k)\big), 
$$
qui est un groupe fini puisque $\Gr(G)(k)=G(R)$ est fini. Par conséquent,
$$
\dim \Gr_d(V(\omega_G)) = \dim \Gr_{n+d}(G) = \dim \Coker(\Gr_{n+d}(f)), 
$$
la seconde égalité résultant du calcul de dimension \ref{troncationlisse} \emph{2}. En particulier, pour tout $n\geq 2d+1$, on a
$$
\dim \Coker(\Gr_{n+1}(f))=\dim \Coker(\Gr_n(f)).
$$
Comme le noyau de l'épimorphisme $\Coker(\Gr_{n+1}(f))\ra \Coker(\Gr_n(f))$ est lisse et connexe (cf. \ref{cH1adRn}), cet épimorphisme est donc un isomorphisme pour tout $n\geq 2d+1$. 

Montrons enfin que $\Gr_d(V(\omega_G))$ est de dimension $d$. Le groupe vectoriel $\Gr_d(V(\omega_G))$ a pour dimension la longueur du $R$-module $\Gr_d(V(\omega_G))(k)=\Hom_R(\omega_G,R/t^d)$. Or d'après \ref{calculd}, on a 
$$\lo_R \Hom_R(\omega_G,R/t^d)=\lo_R \omega_G = d.$$
\end{proof}

\chapter{Modèles de Néron généralisés}

Soient $K$ un corps discrètement valué, $R$ son anneau d'entiers et $k$ son corps résiduel.

\section{Modèle de Néron généralisé d'un groupe algébrique lisse sur $K$}

\subsection{Définition}

\begin{Def*}\label{Defmodel}
Soit $G_K$ un $K$-schéma en groupes. 
\begin{enumerate}
\item Un \emph{modèle} de $G_K$ sur $R$ est un couple $(G,\iota)$ formé d'un $R$-schéma en groupes $G$ et d'un isomorphisme de $K$-schémas en groupes
$$
\xymatrix{
\iota: G\otimes_R K \ar[r]^<<<<{\sim} & G_K.
}
$$
\item Un \emph{morphisme de modèles} $(G,\iota)\ra(G',\iota)$ est un morphisme de $R$-schémas en groupes $f:G\ra G'$ tel que le diagramme
$$
\xymatrix{
G\otimes_R K \ar[dr]_{\rotatebox{135}{$\sim$}}^{\iota} \ar[rr]^{f_K}  && G'\otimes_R K \ar[dl]^{\rotatebox{45}{$\sim$}}_{\iota'} \\
& G_K & 
}
$$
commute.
\end{enumerate}
\end{Def*}

\begin{Rem*}\label{axiome3mn}
Si $(G,\iota)$ et $(G',\iota')$ sont des modèles de $G_K$, avec $G$ plat sur $R$ et $G'$ séparé sur $R$, alors il existe au plus un morphisme $(G,\iota)\ra(G',\iota')$.
\end{Rem*}

\begin{Pt*}
On sous-entendra souvent l'isomorphisme $\iota$ du couple $(G,\iota)$ définissant un modèle de $G_K$. 
\end{Pt*}

\begin{Def*} \label{Defmn}
Soit $L_K$ un $K$-schéma en groupes lisse. Un \emph{modèle de Néron} de $L_K$ sur $R$ est un ind-objet $L=(L_i,\iota_i)_{i\in I}$ de la catégorie des modèles de $L_K$, avec $L_i$ lisse et séparé sur $R$ pour tout $i\in I$, vérifiant la propriété suivante : 
 
\medskip
 
Pour tout $R$-schéma lisse de type fini $Y$, l'application induite par passage à la fibre générique
$$
\xymatrix{
\li_{i\in I}\Hom_{\SchR}(Y,L_i)\ar[r] &\Hom_{\SchK}(Y_K,L_K)
}
$$
est bijective. 
\end{Def*}

\begin{Rem*}
Le ind-objet $L\otimes_R K:=(L_i\otimes_R K)_{i\in I}$ de la catégorie des $K$-schémas en groupes déduit de $(L_i,\iota_i)_{i\in I}$ par oubli des isomorphismes $\iota_i$ et passage à la fibre générique, est isomorphe au ind-objet constant $L_K$ :
$$
\xymatrix{
(\iota_i)_{i\in I}:L\otimes_R K \ar[r]^<<<<<{\sim} & L_K.
}
$$
\end{Rem*}

\begin{Pt*}
Lorsque l'ensemble $I$ est un singleton, l'application
$$
\xymatrix{
\Hom_{\SchR}(Y,L)\ar[r] &\Hom_{\SchK}(Y_K,L_K)
}
$$
est encore bijective pour tout $R$-schéma lisse $Y$ qui n'est pas de type fini. La définition \ref{Defmn} généralise donc la définition des modèles de Néron ltf de \cite{BLR} 10.1/1 dans le cas des $K$-schémas en groupes (commutatifs).
\end{Pt*}

\subsection{Unicité}

\begin{Pt*} \label{unitf}
Soit $L_K$ un $K$-schéma en groupes lisse de type fini. Alors $L_K$ possède \emph{au plus un} modèle de Néron $L=(L_i)_{i\in I}$ \emph{avec $L_i$ de type fini sur $R$ pour tout $i\in I$}, à isomorphisme unique près de ind-objets de la catégorie des modèles de $L_K$. Cela résulte directement de la propriété universelle d'un tel modèle de Néron.
Cette propriété d'unicité peut être généralisée, lorsque $R$ est strictement hensélien.
\end{Pt*}

\begin{Not*} \label{Notpi0}
Soit $X$ un $k$-schéma en groupes localement de type fini. On note $\pi_0(X)$ le $k$-schéma en groupes étale de ses composantes connexes, identifié à un groupe ordinaire lorsque $k$ est séparablement clos.
\end{Not*}

\begin{Prop*}\label{conduni}
Supposons $R$ strictement hensélien. Soient $L_K$ un $K$-schéma en groupes lisse et $L=(L_i)_{i\in I}$ un modèle de Néron de $L_K$. Soit $Y$ un $R$-schéma en groupes lisse et séparé, à fibre générique de type fini, et tel que le groupe ordinaire 
$\pi_0(Y_k)$ soit de type fini. Alors pour tout morphisme de $K$-schéma en groupes  
$$
\xymatrix{
f_K:Y_K \ar[r] & L_K,
}
$$
il existe $i\in I$ et un morphisme de $R$-schémas en groupes
$$
\xymatrix{
f:Y \ar[r] & L_i
}
$$
prolongeant $f_K$.
\end{Prop*}

\begin{proof}
Par hypothèse, il existe un entier $r\in\bbN$ et un groupe fini $\pi_0(Y_k)_{\tors}$ tel que 
$$\pi_0(Y_k)\simeq \bbZ^r\times \pi_0(Y_k)_{\tors}.$$
Notons $Y_{\tors}$ le sous-$R$-schéma en groupes ouvert de $Y$ ayant pour fibre générique $Y_K$ et pour fibre spéciale la réunion des composantes connexes de $Y_k$ appartenant à $\pi_0(Y_k)_{\tors}$. Notons 
$Y_{\tf}$ le sous-schéma ouvert de $Y$, réunion de $Y_{\tors}$ et des composantes connexes $e_1,\ldots,e_r$ de $Y_k$ formant la base canonique $e_1,\ldots,e_r$ de $\bbZ^r\subset \pi_0(Y_k)$. Comme le schéma en groupes $Y$ est lisse sur $R$, à fibre générique de type fini, le schéma $Y_{\tf}$ est de type fini sur $R$ (cf. \cite{SGA 3}${}_I$ VI${}_B$ 3.10 et 3.6).
Par définition du modèle de Néron $L$, il existe donc $i\in I$ et un morphisme de $R$-schémas
$$
\xymatrix{
f:Y_{\tf} \ar[r] & L_i
}
$$
prolongeant le morphisme de $K$-schémas en groupes $f_K:Y_K\ra L_K$. En particulier, la restriction de $f$ à $Y_{\tors}$ est un morphisme de $R$-schémas en groupes.

Dans cette situation, $f$ s'étend naturellement à $Y$ tout entier. Pour le voir, choisissons des sections $s_1,\ldots,s_r$ de $Y$ se spécialisant respectivement dans  $e_1,\ldots,e_r$. Soit $c\in\pi_0(Y_k)$. Notons $Y_c$ le sous-schéma ouvert de $Y$ réunion de $Y_K$ et de la composante $c$ de $Y_k$, et fixons une section $s_c$ de $Y_c$. La translation par 
$-s_c$ définit un isomorphisme
$$
\xymatrix{
-s_c:Y_c \ar[r]^{\sim} & Y_0,
}
$$
où $Y_0$ est la réunion de $Y_K$ et de la composante neutre de $Y_k$. De plus, si la projection de $c$ sur $\bbZ^r\subset\pi_0(Y_k)$ s'écrit
$$n_1e_1+\ldots+n_re_r,$$
alors 
$$s_c=n_1s_1+\ldots+n_rs_r+s_{c,\tors}$$
avec $s_{c,\tors}\in Y_{\tors}(R)$, et on peut considérer la translation dans $L_i$, notée $\tau_{f(s_c)}$, définie par la section 
$$n_1f(s_1)+\ldots+n_rf(s_r)+f(s_{c,\tors}).$$ 
On définit alors la composée
$$
\xymatrix{
f_{s_c}:Y_c \ar[r]^<<<<{-s_c} & Y_0 \ar[r]^{f|_{Y_0}} & L_i \ar[r]^{\tau_{f(s_c)}} & L_i.
}
$$
Comme la restriction de $f$ à la fibre générique $Y_K$ est un morphisme de groupes, les morphismes $f_{s_c}$ et $f$ coïncident sur $Y_K$. En particulier $f_{s_c}$ ne dépend pas du choix de la section $s_c$ de $Y_c$, et on peut poser $f_{Y_c}:=f_{s_c}$. Si de plus $Y_c$ est contenu dans $Y_{\tf}$, alors $f_{Y_c}$ coïncide avec la restriction de $f$ à $Y_c$. Comme les $Y_c$, $c\in\pi_0(Y_k)$, forment un recouvrement ouvert de $Y$, on peut donc bien étendre $f$ à $Y$ au moyen des $f_{Y_c}$.
\end{proof}

\begin{Cor*}\label{unicitemn}
Supposons $R$ strictement hensélien. Soient $L_K$ un $K$-schéma en groupes lisse de type fini et $M_K$ un $K$-schéma en groupes lisse. Soit $L=(L_i,\iota_i)_{i\in I}$ un ind-objet de la catégorie des modèles de $L_K$, tel que pour tout $i$, $L_i$ soit lisse et séparé sur $R$ et $\pi_0(L_{i,k})$ de type fini. Soit $M=(M_j,\iota_j)_{j\in J}$ un modèle de Néron de $M_K$. Alors l'application induite par passage à la fibre générique
$$
\xymatrix{
\Hom_{\Ind(\SchgrR)}((L_i)_{i\in I},(M_j)_{i\in J})\ar[r] & \Hom_{\SchgrK}(L_K,M_K)
}
$$
est bijective. En particulier, un $K$-schéma en groupes lisse de type fini $L_K$ possède au plus un modèle de Néron $L=(L_i,\iota_i)_{i\in I}$ tel que $\pi_0(L_{i,k})$ soit de type fini pour tout $i\in I$, à isomorphisme unique près de ind-objets de la catégorie des modèles de $L_K$. De plus, le ind-objet de la catégorie des $R$-schémas en groupes sous-jacent à $L$ dépend fonctoriellement de $L_K$.
\end{Cor*}

\begin{proof}
Résulte de \ref{conduni}, et du fait qu'un morphisme d'un $R$-schéma plat dans un $R$-schéma séparé est uniquement déterminé par sa fibre générique.
\end{proof}

\begin{Rem*}
Dans la proposition \ref{conduni}, le fait que $f_K$ soit un morphisme de $K$-schémas \emph{en groupes} est indispensable : voir \ref{nonprol} pour un contre-exemple si l'on retire cette hypothèse. 
\end{Rem*}

\begin{Pt*} \label{HN}
Soit $L_K$ un $K$-schéma en groupes lisse de type fini, possédant un modèle de Néron ltf $(L,\iota)$ au sens de \cite{BLR} 10.1/1. Notons $k^s$ une clôture séparable de $k$. D'après \cite{HN} 3.5, le groupe des composantes connexes de $L\otimes_R k^s$ est de type fini. D'après \ref{unicitemn}, lorsque $R$ est strictement hensélien, tout modèle de Néron $(L_i,\iota_i)_{i\in I}$ de $L_K$ au sens de \ref{Defmn} tel que $\pi_0(L_{i,k})$ est de type fini pour tout $i\in I$ est donc canoniquement isomorphe à $(L,\iota)$.
\end{Pt*}

\subsection{Compatibilité aux extensions formellement lisses} 

\begin{Pt*}
Rappelons qu'un morphisme local d'anneaux de valuation discrète $R\ra R'$ est \emph{essentiellement lisse} si $R'$ est $R$-isomorphe au localisé d'un $R$-schéma lisse. Un tel morphisme est en particulier formellement lisse pour les topologies $\fm$-préadique sur $R$ et $\fm'$-préadique sur $R'$, $\fm$ et $\fm'$ désignant les idéaux maximaux de $R$ et $R'$ (\cite{EGA IV}${}_4$ 17.5.3). Dans la suite, on emploiera toujours la terminologie << formellement lisse >> relativement à ces topologies, sans qu'on le mentionne explicitement. Elle est équivalente à la terminologie << avoir indice de ramification 1 >> de \cite{BLR} 3.6/1 (\cite{EGA IV}${}_1$ 0.19.7.1, 0.19.6.5 et 0.19.6.1).
\end{Pt*}

\begin{Lem*}\label{Ncaractesslisse}
Soient $L_K$ un $K$-schéma en groupes lisse et $L=(L_i)_{i\in I}$ un ind-objet de la catégorie des modèles de $L_K$, avec $L_i$ lisse et séparé sur $R$ pour tout $i\in I$. Alors $L$ est un modèle de Néron de $L_K$ si et seulement si pour tout morphisme local essentiellement lisse d'anneaux de valuation discrète $R\rightarrow R'$, l'application induite par passage à la fibre générique
$$
\xymatrix{
\li_{i\in I}\Hom_{\SchR}(\Spec(R'),L_i)\ar[r] &\Hom_{\SchK}(\Spec(K'),L_K)
}
$$
est surjective. 
\end{Lem*}

\begin{proof}
La condition est évidemment nécessaire, et il s'agit de voir qu'elle est suffisante. Soit donc $Y$ un $R$-schéma lisse de type fini, et considérons l'application de passage à la fibre générique de la définition \ref{Defmn}. Cette application est injective puisqu'un morphisme d'un $R$-schéma plat dans un $R$-schéma séparé est uniquement déterminé par sa fibre générique. 
Montrons que sous la condition du lemme, elle est également surjective. Soit $Y_K\ra L_K$ un $K$-morphisme. Soit $R'$ le localisé de $Y$ en un point maximal de sa fibre spéciale. Par hypothèse, le morphisme composé $\Spec(K')\ra Y_K\ra L_K$ se prolonge en un $R$-morphisme $\Spec(R')\ra L_i$ pour un certain $i\in I$. Comme la fibre spéciale de $Y/R$ ne compte qu'un nombre fini de points maximaux, on peut supposer que $i$ ne dépend pas du point maximal considéré. Le morphisme $Y_K\ra L_K$ se prolonge alors en une application $R$-rationnelle de $Y$ dans $L_i$. Mais comme $L_i$ est un $R$-schéma en groupes séparé, cette application est nécessairement un morphisme d'après le théorème d'extension de Weil \cite{BLR} 4.4/1, et le $R$-morphisme $Y\ra L_i$ ainsi obtenu prolonge $Y_K\ra L_K$. 
\end{proof}

\begin{Lem*}\label{Npropflisse}
Soient $L_K$ un $K$-schéma en groupes lisse et $L=(L_i)_{i\in I}$ un modèle de Néron de $L_K$. Soit $R\ra R'$ un morphisme local formellement lisse d'anneaux de valuation discrète. Alors l'application induite par passage à la fibre générique
$$
\xymatrix{
\li_{i\in I}\Hom_{\SchR}(\Spec(R'),L_i)\ar[r] &\Hom_{\SchK}(\Spec(K'),L_K)
}
$$
est bijective.
\end{Lem*}

\begin{proof}
L'application considérée est injective puisque les $L_i$ sont séparés sur $R$. Réciproquement, soit $a_{K'}\in L_K(K')$. Choisissons un ouvert affine $U_K$ de $L_K$ contenant $a_{K'}$. En chassant les dénominateurs, on peut trouver un $R$-schéma affine de type fini $U$ de fibre générique isomorphe à $U_K$ tel que $a_{K'}$ se prolonge en un élément $a\in U(R')$. D'après l'hypothèse sur le morphisme $R\ra R'$, et puisque $U_K$ est lisse sur $K$, il existe alors un $R$-morphisme 
$U'\ra U$ composé d'un nombre fini de dilatations centrées dans les fibres fermées, tel que $a$ se relève en un $R'$-point  
$\widetilde{a}$ de l'ouvert de lissité $\widetilde{U}$ du $R$-schéma $U'$ : \cite{BLR} 3.6/4. Comme $\widetilde{U}$ est de type fini sur $R$ par construction, et $L$ est un modèle de Néron de $L_K$, l'immersion ouverte $\widetilde{U}_K\ra L_K$ se prolonge en un morphisme de $R$-schémas $\widetilde{U} \ra L_i$ pour un certain $i\in I$. L'image de $\widetilde{a}$ par ce morphisme est alors un $R'$-point de $L_i$ prolongeant $a_{K'}$.
\end{proof}

\begin{Cor*}\label{Neroncb}
Soient $L_K$ un $K$-schéma en groupes lisse et $L=(L_i)_{i\in I}$ un ind-objet de la catégorie des modèles de $L_K$, avec $L_i$ lisse et séparé sur $R$ pour tout $i\in I$. Alors pour tout morphisme local formellement lisse d'anneaux de valuation discrète $R\rightarrow R'$, $L\otimes_R R':=(L_i\otimes_R R')_{i\in I}$ est un modèle de Néron de $L_K\otimes_K K'$.
\end{Cor*}

\section{Existence}

\begin{Pt} \label{defGa}
On définit un ind-objet $\cG_a=(\cG_{a,i},\iota_i)_{i\in\bbN}$ de la catégorie des modèles du groupe additif $\bbG_{a,K}$ de la façon suivante. On choisit une uniformisante $t$ de $R$. Pour tout $i\in \bbN$, on pose
$$
\xymatrix{
\cG_{a,i}:=\bbG_{a,R} & \textrm{et} & \iota_i:\cG_{a,i}\otimes_R K=\bbG_{a,K}\ar[r]^<<<<<{t^{-i}} & \bbG_{a,K}.
}
$$
Pour tous $i\leq j$, on prend pour transition $\tau_{i,j}:\cG_{a,i}\ra\cG_{a,j}$ l'homothétie de rapport $t^{j-i}$
$$
\xymatrix{
\bbG_{a,R} \ar[r]^{t^{j-i}} & \bbG_{a,R}.
}
$$
Alors $\tau_{i,j}$ est un morphisme de modèles $(\cG_{a,i},\iota_i)\ra(\cG_{a,j},\iota_j)$.
\end{Pt}

\begin{Lem} \label{Ga}
$\cG_a=(\cG_{a,i})_{i\in\bbN}$ est un modèle de Néron de $\bbG_{a,K}$.
\end{Lem}

\begin{proof}
L'application 
$$
\xymatrix{
\li_{i\in \bbN}\Hom_{\SchR}(\Spec(R),\cG_{a,i})\ar[r] &\Hom_{\SchK}(\Spec(K),\bbG_{a,K})
}
$$
induite par passage à la fibre générique est bijective par construction. Comme la formation de $\cG_a$ commute aux extensions d'anneaux de valuation discrète qui envoient uniformisante sur uniformisante, il résulte de \ref{Ncaractesslisse} que $\cG_a$ est un modèle de Néron de $\bbG_{a,K}$.
\end{proof}

\begin{Rem}\label{nonprol}
Soit $Y$ le $R$-schéma lisse sous-jacent au modèle de Néron du groupe multiplicatif $\bbG_{m,K}$. Alors pour tout $i\in\bbN$, la composée de l'immersion ouverte canonique $Y_K\ra \bbG_{a,K}$ et de l'homothétie $\iota_i:\bbG_{a,K}\ra\bbG_{a,K}$ ne se prolonge pas en un morphisme de $R$-schémas $Y\ra\cG_{a,i}$ : si un tel morphisme existait, l'ensemble $\bbG_{a,K}(K)\setminus\{0\}=K^{*}$ se prolongerait en un ensemble de sections de $\cG_{a,i}=\bbG_{a,R}$, ce qui n'est pas le cas.
\end{Rem}

\begin{Lem} \label{unipd}
Soit $U_K$ un $K$-schéma en groupes  lisse unipotent connexe déployé. Alors $U_K$ possède un modèle de Néron 
$(U_i)_{i\in \bbN}$, tel que pour tout $i\in\bbN$, le $R$-schéma en groupes $U_i$ admette une suite de composition dont les quotients successifs sont isomorphes à $\bbG_{a,R}$, et pour tous $i<j\in\bbN$, le morphisme de transition $U_i\ra U_j$ soit nul sur les fibres spéciales.
\end{Lem}

\begin{proof}
Raisonnant par récurrence sur la longueur d'une suite de composition de $U_K$, on peut supposer que $U_K$ est extension de $\bbG_{a,K}$ par un groupe $V_K$ possédant un modèle de Néron $V=(V_{\alpha})_{\alpha\in\bbN}$, tel que pour tout $\alpha\in\bbN$, le $R$-schéma en groupes $V_{\alpha}$ admette une suite de composition dont les quotients successifs sont isomorphes à $\bbG_{a,R}$, et pour tous $\alpha<\beta\in\bbN$, le morphisme de transition $V_{\alpha}\ra V_{\beta}$ soit nul sur les fibres spéciales : 
$$
\xymatrix{
0\ar[r] & V_K \ar[r] & U_K \ar[r] & \bbG_{a,K} \ar[r] & 0.
}
$$
Comme $H^1(\bbG_{a,K},V_K)=0$, cette extension définit un $2$-cocycle symétrique
$$
\xymatrix{
c:\bbG_{a,K}\times_K \bbG_{a,K} \ar[r] & V_K.
}
$$ 
Par définition de $V$, il existe $\alpha_0\in\bbN$ tel que $c$ se prolonge en un $R$-morphisme
$$
\xymatrix{
\cG_{a,0}\times_R \cG_{a,0} \ar[r] & V_{\alpha_0}.
}
$$
De même, il existe $\alpha_1\in\bbN$, que l'on peut choisir strictement supérieur à $\alpha_0$, tel que 
\begin{displaymath}
\xymatrix{
\bbG_{a,K}\times_K \bbG_{a,K} \ar[r]_{\sim}^{t^{-1}\times t^{-1}} & \bbG_{a,K}\times_K \bbG_{a,K} \ar[r]^>>>>>{c} & V_K
}
\end{displaymath}
se prolonge en un $R$-morphisme
$$
\xymatrix{
\cG_{a,1}\times_R \cG_{a,1} \ar[r] & V_{\alpha_1}.
}
$$
D'où un cube commutatif
\begin{displaymath}
\xymatrix{
\bbG_{a,K}\times_K\bbG_{a,K} \ar[rr]^<<<<<<<<<<<<<<<<<<<{c} \ar@{^{(}->}[dr] \ar[dd]^>>>>>>>{t\times t}_>>>>>>>{\rotatebox{90}{$\sim$}} &                 & V_K\simeq V_{\alpha_0}\otimes_R K \ar[dd]_>>>>>>{\rotatebox{90}{$\sim$}} \ar@{^{(}->}[drr]&  \\
                         & \cG_{a,0}\times_R\cG_{a,0} \ar[rrr] \ar[dd]^<<<<<<<{t\times t} &                &  &V_{\alpha_0} \ar[dd]^<<<<<<<{\tau_{\alpha_0,\alpha_1}}\\
\bbG_{a,K}\times_K\bbG_{a,K} \ar[rr] \ar@{^{(}->}[dr]       &                 & V_K\simeq V_{\alpha_1}\otimes_R K \ar@{^{(}->}[drr]& \\
                         & \cG_{a,1}\times_R\cG_{a,1} \ar[rrr]           &                 &  & V_{\alpha_1}
}
\end{displaymath}
et par conséquent un morphisme d'extensions sur $R$
\begin{displaymath}
\xymatrix{
0\ar[r] & V_{\alpha_0} \ar[r] \ar[d]^{\tau_{\alpha_0,\alpha_1}} & U_0 \ar[r] \ar[d]^{(t,\tau_{\alpha_0,\alpha_1})} & \cG_{a,0} \ar[r] \ar[d]^t & 0\\
0\ar[r] & V_{\alpha_1} \ar[r] & U_1 \ar[r] & \cG_{a,1} \ar[r] & 0,
}
\end{displaymath}
nul sur les fibres spéciales. Plus généralement, on construit un système inductif d'extensions indexé par $i\in\bbN$
\begin{displaymath}
\xymatrix{
0\ar[r] & V_{\alpha_i} \ar[r] \ar[d]^{\tau_{\alpha_i,\alpha_{i+1}}} & U_i \ar[r] \ar[d]^{(t,\tau_{\alpha_i,\alpha_{i+1}})} & \cG_{a,i} \ar[r] \ar[d]^t & 0\\
0\ar[r] & V_{\alpha_{i+1}} \ar[r] & U_{i+1} \ar[r] & \cG_{a,i+1} \ar[r] & 0,
}
\end{displaymath}
avec $\alpha_{i+1}>\alpha_i$, nul sur les fibres spéciales. Montrons que $U:=(U_i)_{i\in\bbN}$ est un modèle de Néron de $U_K$, et pour cela, utilisons le critère \ref{Ncaractesslisse}. Comme la formation de $U$ commute aux extensions essentiellement lisses d'anneaux de valuation discrète, il suffit de voir que l'application $\li_{i\in\bbN}U_i(R)\ra U_K(K)$ est surjective. Soit $u_K\in U_K(K)=U_0(K)$. Son image dans $\cG_{a,0}(K)$ se prolonge en une section $a_i$ de $\cG_{a,i}$ pour un certain $i\in\bbN$.  Comme $H^1(R,V_{\alpha_i})=0$ puisque $V_{\alpha_i}$ se dévisse en $\bbG_{a,R}$ par hypothèse, il existe $u_i\in U_i(R)$ relevant $a_i$. La différence $\delta_{i,K}:=u_K-u_{i,K}$ appartient au sous-groupe $V_{\alpha_i}(K)$, et par conséquent se prolonge en un $\delta_j\in V_{\alpha_j}(R)$ pour un certain $j\geq i$. Notant $u_j$ l'image de $u_i$ dans $U_j(R)$, la section $\delta_j+u_j\in U_j(R)$ prolonge $u_K$. Ainsi $U$ est bien un modèle de Néron de $U_K$. Par construction, chaque $U_i$ admet une suite de composition dont les quotients successifs sont isomorphes à $\bbG_{a,R}$, et les transitions non triviales sont nulles sur les fibres spéciales.
\end{proof}

\begin{Rem} \label{structuremnunipd}
Par construction, le ind-$R$-schéma en groupes $(U_i)_{i\in\bbN}$ admet une suite de composition dont les quotients successifs sont isomorphes au ind-$R$-schéma en groupes $\cG_a$ (ceci dans la ind-catégorie abélienne des faisceaux sur le gros site Zariski de $\Spec(R)$).
\end{Rem}

\begin{Ex} 
Soit $\bbW_{2,K}$ le schéma en groupes sur $K$ des vecteurs de Witt de longueur $2$. Le $2$-cocycle symétrique de Lazard
\begin{eqnarray*}
L:\bbG_{a,K}\times_K\bbG_{a,K} & \lra & \bbG_{a,K} \\
(x,y) & \lmapsto & \frac{1}{p}(x^p+y^p-(x+y)^p)
\end{eqnarray*}
définit la suite exacte
\begin{displaymath}
\xymatrix{
0\ar[r] & \bbG_{a,K} \ar[r]^{V} & \bbW_{2,K} \ar[r]^{R} & \bbG_{a,K} \ar[r] & 0 
}
\end{displaymath}
où $V(x)=(0,x)$ et $R(x_0,x_1)=x_0$ (\cite{S1} Chap. VII § 2.7). Le cocycle $L$ est un polynôme homogène de degré $p$ en les variables $x$ et $y$, à coefficients dans $\bbZ$. Le cube
\begin{displaymath}
\xymatrix{
\bbG_{a,K}\times_K\bbG_{a,K} \ar[rr]^<<<<<<<<<<<<<<<<<<<<<<<{L} \ar@{^{(}->}[dr] \ar[dd]^>>>>>>>{t\times t}_>>>>>>>{\rotatebox{90}{$\sim$}} &                 & \bbG_{a,K} \ar[dd]^>>>>>>{t^p}_>>>>>>{\rotatebox{90}{$\sim$}} \ar@{^{(}->}[drr]&  \\
                         & \cG_{a,0}\times_R\cG_{a,0} \ar[rrr] \ar[dd]^<<<<<<<{t} &                &  &\cG_{a,0} \ar[dd]^<<<<<<<{t^p}\\
\bbG_{a,K}\times_K\bbG_{a,K} \ar[rr]^<<<<<<<<<<<<<<<<<<<<<<<{L} \ar@{^{(}->}[dr]       &                 & \bbG_{a,K} \ar@{^{(}->}[drr]& \\
                         & \cG_{a,1}\times_R\cG_{a,1} \ar[rrr]           &                 &  & \cG_{a,p}
}
\end{displaymath}
est donc commutatif, et correspond à un morphisme d'extensions sur $R$
\begin{displaymath}
\xymatrix{
0\ar[r] & \cG_{a,0} \ar[r] \ar[d]^{t^p} & \cW_{2,0} \ar[r] \ar[d]^{(t,t^p)} & \cG_{a,0} \ar[r] \ar[d]^t & 0\\
0\ar[r] & \cG_{a,p} \ar[r] & \cW_{2,1} \ar[r] & \cG_{a,1} \ar[r] & 0.
}
\end{displaymath}
Plus généralement, on construit un système inductif d'extensions indexé par $i\in\bbN$
\begin{displaymath}
\xymatrix{
0\ar[r] & \cG_{a,ip} \ar[r] \ar[d]^{t^p} & \cW_{2,i} \ar[r] \ar[d]^{(t,t^p)} & \cG_{a,i} \ar[r] \ar[d]^t & 0\\
0\ar[r] & \cG_{a,(i+1)p} \ar[r] & \cW_{2,i+1} \ar[r] & \cG_{a,i+1} \ar[r] & 0.
}
\end{displaymath}
On obtient ainsi une extension (au sens de la remarque de \ref{structuremnunipd})
\begin{displaymath}
\xymatrix{
0\ar[r] & \cG_{a} \ar[r]  & \cW_{2} \ar[r] & \cG_{a} \ar[r] & 0
}
\end{displaymath}
où $\cW_2:=(\cW_{2,i})_{i\in\bbN}$ est un modèle de Néron de $\bbW_{2,K}$.
\end{Ex}

\begin{Ex}\label{defWn}
Dans le cas des $K$-schémas en groupes de vecteurs de Witt tronqués, on peut cependant donner une construction directe 
du modèle de Néron fourni par la construction \ref{unipd}. Soit $n$ un entier $\geq 1$ et soit $\bbW_{n,K}$ le schéma en groupes sur $K$ des vecteurs de Witt de longueur $n$. Notant toujours $t$ une uniformisante de $R$, on pose, pour tout $i\in \bbN$, 
$$
\xymatrix{
\cW_{n,i}:=\bbW_{n,R} & \textrm{et} & \iota_i:\cW_{n,i}\otimes_R K=\bbW_{n,K}\ar[r]^<<<<<{[t^{-i}]} & \bbW_{n,K},
}
$$
où $[\cdot]:K\ra\bbW_{n,K}(K)$ désigne le représentant de Teichmüller. Pour tous $i\leq j$, on prend pour transition $\tau_{i,j}:\cW_{n,i}\ra\cW_{n,j}$ l'homothétie de rapport $[t^{j-i}]\in\bbW_{n,R}(R)$
$$
\xymatrix{
\bbW_{n,R} \ar[r]^{[t^{j-i}]} & \bbW_{n,R}.
}
$$
Alors $\tau_{i,j}$ est un morphisme de modèles $(\cW_{n,i},\iota_i)\ra(\cW_{n,j},\iota_j)$. Le ind-objet $\cW_n:=(\cW_{n,i},\iota_i)_{i\in\bbN}$ de la catégorie des modèles de $\bbW_{n,K}$ est alors celui fourni par \ref{unipd}. On peut retrouver directement que c'est un modèle de Néron de $\bbW_{n,K}$ grâce au critère \ref{Ncaractesslisse} et à la formule
$$
[t](a_0,\ldots,a_{n-1})=(ta_0,\ldots,t^{p^{n-1}}a_{n-1}) \quad \textrm{pour tout } (a_0,\ldots,a_{n-1})\in\bbW_{n,K}(K).
$$
\end{Ex}

\begin{Th}\label{casgenmn}
Supposons $R$ complet et $k$ séparablement clos. Alors tout $K$-schéma en groupes lisse de type fini possède un modèle de Néron $(L_i)_{i\in\bbN}$, avec $\pi_0(L_{i,k})$ de type fini pour tout $i\in\bbN$.
\end{Th}

\begin{Lem}\label{Tits}
Soient $K$ un corps et $L_K$ un $K$-schéma en groupes  lisse de type fini. Alors il existe une extension
$$
\xymatrix{
0\ar[r] & T_K\times_K U_K \ar[r] & L_K \ar[r] & B_K \ar[r] & 0 
} 
$$
où $T_K$ est un tore déployé, $U_K$ est un groupe unipotent connexe lisse déployé, et $B_K$ ne possède pas de sous-groupe isomorphe à $\bbG_{m,K}$ ou $\bbG_{a,K}$.
\end{Lem}

\begin{proof}
Commençons par le cas où le $K$-schéma en groupes lisse $L_K$ est affine et connexe. Alors $L_K$ est extension d'un groupe unipotent lisse connexe $L_{K,u}$ par un tore $L_{K,t}$ (\cite{SGA 3} XVII 7.2). D'après Tits, le groupe $L_{K,u}$ est lui-même extension d'un groupe unipotent totalement ployé $L_{K,u,p}$ par un groupe unipotent déployé $L_{K,u,d}$ (cf. \cite{O} paragraphe V.5 ou \cite{CGP} Appendice B pour une démonstration). L'image réciproque de l'extension
$$ 
\xymatrix{
0\ar[r] & L_{K,t} \ar[r] & L_K \ar[r] & L_{K,u} \ar[r] & 0
}
$$
par l'inclusion $L_{K,u,d}\ra L_{K,u}$ est scindée (\cite{SGA 3} XVII 6.1.1 A) ii)). On note encore $L_{K,u,d}$ l'unique relèvement de $L_{K,u,d}$ dans $L_K$. Par ailleurs, le groupe $L_{K,t}$ possède un plus grand sous-tore déployé $L_{K,t,d}$, de sorte que le quotient $L_{K,t}/L_{K,t,d}$ est un tore anisotrope $L_{K,t,a}$. On a alors le diagramme commutatif exact suivant :
\begin{displaymath}
\xymatrix{
            & 0             \ar[d]           & 0   \ar[d]                                                         & 0  \ar[d] \\               
0 \ar[r] & L_{K,t,d} \ar[d] \ar[r] &  L_{K,t,d}\times_K L_{K,u,d} \ar[d] \ar[r]        & L_{K,u,d}  \ar[d] \ar[r] & 0 \\
0 \ar[r] & L_{K,t}    \ar[d] \ar[r] & L_K \ar[d] \ar[r]                                               & L_{K,u} \ar[r] \ar[d]      &  0 \\
0 \ar[r] & L_{K,t,a} \ar[d] \ar[r] & L_K/(L_{K,t,d}\times_K L_{K,u,d}) \ar[r] \ar[d]  & L_{K,u,p} \ar[r] \ar[d] & 0 \\
            & 0                               & 0                                                                & 0.
}
\end{displaymath}
D'où le lemme dans le cas lisse affine connexe. Dans le cas général, $L_K$ contient un plus grand sous-groupe lisse affine connexe $L_{K}'$ ; le quotient  $L_{K}'':=L_K/L_{K}'$ ne contient alors pas de sous-groupe lisse affine connexe non nul.
Le lemme appliqué à $L_{K}'$ fournit des groupes $T_K$, $U_K$ et $B_K$, puis le diagramme commutatif exact :
\begin{displaymath}
\xymatrix{
            & 0             \ar[d]           & 0   \ar[d]                                                         & \\               
            & T_K\times_K U_K \ar[d] \ar@{=}[r] &  T_K \times_K U_K \ar[d]         \\
0 \ar[r] & L_{K}'    \ar[d] \ar[r] & L_K \ar[d] \ar[r]                                               & L_{K}'' \ar[r] \ar@{=}[d]      &  0 \\
0 \ar[r] & B_{K} \ar[d] \ar[r]       & L_K/(T_K\times_K U_K) \ar[r] \ar[d]               & L_{K}'' \ar[r] & 0 \\
            & 0                               & 0.                                                               
}
\end{displaymath}
D'où lemme.
\end{proof}

\begin{Rem}
Si $K$ est de caractéristique zéro, le $K$-schéma en groupes $L_{K}''$ qui apparaît dans la démonstration précédente est une variété abélienne, d'après le théorème de Chevalley. Mais si $K$ est de caractéristique $>0$, donc ici imparfait, $L_{K}''$ n'est pas une variété abélienne en général ; lorsque $L_K$ est connexe, $L_{K}''$ est une variété \emph{pseudo-abélienne} dans la terminologie de Totaro \cite{To}.
\end{Rem}

\begin{Nots}
Soient $L$ et $M$ des $R$-schémas en groupes lisses, séparés, fidèlement plats et de type fini. On note $\Ext_K(L_K,M_K)_{\rat}$ (resp. $\Ext_R(L,M)_{\Rrat}$) le groupe des classes d'isomorphisme d'extensions de $L_K$ par $M_K$ qui possèdent une section $K$-rationnelle (resp. de $L$ par $M$ qui possèdent une section $R$-rationnelle, c'est-à-dire une pseudo-section relativement à $R$ dans la terminologie de \cite{EGA IV} 20.5.). On note $H_{\rat}^2(L_K,M_K)_{\sym}$ (resp. $H_{\Rrat}^2(L,M)_{\sym}$) le groupe des $2$-cocycles symétriques $K$-rationnels  $L_K\times_K L_K\dashrightarrow M_K$ (resp. $R$-rationnels $L\times_R L\dashrightarrow M$), modulo les cobords.
\end{Nots}

\begin{Lem}\label{rat}
On a un diagramme commutatif 
\begin{displaymath}
\xymatrix{
\Ext_R(L,M)_{\Rrat} \ar[r]^>>>>>>>{\sim} \ar@{^{(}->}[d] & H_{\Rrat}^2(L,M)_{\sym} \ar@{^{(}->}[d] \\
\Ext_K(L_K,M_K)_{\rat} \ar[r]^>>>>>>{\sim} & H_{\rat}^2(L_K,M_K)_{\sym}
}
\end{displaymath}
où les flèches horizontales sont bijectives et les flèches verticales sont injectives.
\end{Lem}

\begin{proof}
La démonstration est la même que celle de Serre \cite{S1} Chap. VII §1.4 Prop. 4 b), en remplaçant la référence au théorème de Weil d'extension des lois de groupes $K$-birationnelles par la référence \cite{BLR} 5.1/5.
\end{proof}

\begin{proof}[Démonstration du théorème \ref{casgenmn}]
Le lemme \ref{Tits} fournit une extension
$$
\xymatrix{
0\ar[r] & T_K\times_K U_K \ar[r] & L_K \ar[r] & B_K \ar[r] & 0. 
} 
$$
Considérons le diagramme commutatif
\begin{displaymath}
\xymatrix{
            & 0        \ar[d]           & 0   \ar[d] \\               
           & T_K \ar[d] \ar@{=}[r] & T_K \ar[d]   \\
0 \ar[r] & T_K\times_K U_K \ar[d] \ar[r] & L_K \ar[d] \ar[r] & B_K \ar[r] \ar@{=}[d]&  0 \\
0 \ar[r] & U_K \ar[d] \ar[r] & H_K:=L_K/T_K \ar[d] \ar[r] & B_K  \ar[r] & 0 \\
            & 0                      & 0.                  
}
\end{displaymath}
Posant $\overline{L_K}:=L_K/U_K$, la deuxième suite verticale se réalise comme image réciproque de $\overline{L_K}$ :
\begin{displaymath}
\xymatrix{
              0   \ar[d] & 0 \ar[d]\\               
             T_K \ar[d] \ar@{=}[r]& T_K \ar[d] \\
  L_K \ar[d] \ar[r] & \overline{L_K} \ar[d]  \\
  H_K \ar[d] \ar[r] & B_K  \ar[d] \\
             0            & 0.                  
}
\end{displaymath}
Les hypothèses sur $R$ et $k$ assurent d'une part que $B_K$ possède un modèle de Néron $B$, qui est un $R$-schéma de type fini (\cite{BLR} 10.2/1), et d'autre part que $\overline{L_K}$ possède un modèle de Néron $\overline{L}$, qui est un $R$-schéma localement de type fini (\cite{BLR} 10.2/2). Montrons maintenant que $H_K$ possède un modèle de Néron $H$ au sens de \ref{Defmn}. Pour cela, considérons l'extension
\begin{equation}\label{extBK}
\xymatrix{
0 \ar[r] & U_K \ar[r] & H_K \ar[r] & B_K \ar[r] & 0.
}
\end{equation}
Comme $U_K$ admet une suite de composition dont les quotients successifs sont isomorphes à $\bbG_{a,K}$,  l'extension (\ref{extBK}) possède une section $K$-rationnelle
(\cite{S1} Chap. VI §1.6 Prop. 6 et Lemme 2). Elle peut donc être définie à partir d'un $2$-cocycle symétrique $K$-rationnel 
$$
\xymatrix{
c_K:B_K \times_K B_K \ar@{-->}[r] & U_K.
}
$$
D'après \ref{unipd}, le $K$-groupe $U_K$ possède un modèle de Néron $U=(U_i)_{i\in\bbN}$ où chaque $U_i$ admet une suite de composition dont les quotients successifs sont isomorphes à $\bbG_{a,R}$. Comme $B$ est de type fini sur $R$, la propriété de prolongement pour $U$ fournit un $i_0\in \bbN$ et un $2$-cocycle symétrique $R$-rationnel 
$$
\xymatrix{
c_{i_0}:B\times_R B \ar@{-->}[r] & U_{i_0} 
}
$$
de fibre générique $c_K$. Quitte à remplacer $U=(U_i)_{i\in\bbN}$ par $(U_i)_{i\geq i_0}$, on peut supposer que l'application $c:B\times_R B \dashrightarrow U$ définie par $c_{i_0}$ se factorise par $U_i$ pour tout $i\in\bbN$. D'après \ref{rat}, les cocycles $c_i$, $i\in\bbN$, définissent un système inductif d'extensions
\begin{displaymath}
\xymatrix{
           &  0            \ar[r]       & U_i   \ar[d] \ar[r] & H_i \ar[d] \ar[r] & B \ar@{=}[d] \ar[r] & 0 \\               
            &  0            \ar[r]        & U_j          \ar[r]  & H_j \ar[r]           & B  \ar[r] & 0  & \quad    (j\geq i)    
}
\end{displaymath}
qui prolongent (\ref{extBK}) au-dessus de $R$. Alors $H:=(H_i)_{i\in\bbN}$ est un modèle de Néron de $H_K$. D'après \ref{Ncaractesslisse}, il suffit en effet pour cela que pour toute extension essentiellement lisse d'anneaux de valuation discrète $R'/R$, l'application $\li_{i\in\bbN}H_i(R')\ra H_K(K')$ soit surjective. Or c'est bien le cas puisque $B$ et $U$ sont des modèles de Néron et $H^1(R',(U_0)_{R'})=0$.

Pour terminer la démonstration, on remarque que les morphismes $H_i\ra B$ sont lisses et séparés, et que le $K$-morphisme $\overline{L_K}\ra B_K$ se prolonge en un $R$-morphisme $\overline{L}\ra B$. Alors, en utilisant à nouveau \ref{Ncaractesslisse}, on vérifie que $L:=H\times_B \overline{L}=(H_i\times_B \overline{L})_{i\in\bbN}$ est un modèle de Néron de $L_K$. Le groupe des composantes connexes de la fibre spéciale de $H_i\times_B \overline{L}$ est isomorphe à celui de $\overline{L}$ d'après la suite exacte
$$
\xymatrix{
0\ar[r] & U_i \ar[r] & H_i\times_B \overline{L} \ar[r] & \overline{L} \ar[r] & 0,
}
$$
et par conséquent est de type fini d'après \cite{HN} 3.5. 
\end{proof}

\begin{Pt}\label{structuremn}
Soit $L_K$ un $K$-schéma en groupes lisse de type fini. La structure du (ind-$R$-schéma en groupes sous-jacent au) modèle de Néron $L$ de $L_K$ fourni par la construction \ref{casgenmn} est la suivante. Notant $U_K$ le plus grand sous-groupe unipotent connexe lisse déployé de $L_K$, et posant $\overline{L_K}:=L_K/U_K$, $L$ est extension du modèle de Néron usuel localement de type fini $\overline{L}$ de $\overline{L_K}$, par le modèle de Néron $U$ \ref{unipd} de $U_K$ :
\begin{equation}\label{extGbar}
\xymatrix{
0\ar[r] & U \ar[r] & L \ar[r] & \overline{L} \ar[r] & 0.
}
\end{equation}
Le modèle de Néron $U$ admet une suite de composition dont les quotients successifs sont isomorphes au modèle de Néron $\cG_a$ de $\bbG_{a,K}$ défini en \ref{defGa} (\ref{structuremnunipd}). D'autre part, soit $T_K$ le plus grand sous-tore déployé de $L_K$. Posons $H_K:=L_K/T_K$ et $B_K:=L_K/(T_K\times_K U_K)$. Le modèle de Néron usuel $B$ de $B_K$ est un $R$-schéma en groupes de type fini, et le modèle de Néron $H$ de $H_K$ est un ind-$R$-schéma en groupes de type fini, extension de $B$ par $U$ : 
\begin{equation}\label{extB}
\xymatrix{
0\ar[r] & U \ar[r] & H \ar[r] & B \ar[r] & 0.
}
\end{equation}
L'extension (\ref{extGbar}) est l'image réciproque de l'extension (\ref{extB}) par la flèche canonique $\overline{L}\ra B$ :
\begin{equation}\label{LoLHB}
\xymatrix{
           &  0            \ar[r]       & U   \ar@{=}[d] \ar[r] & L \ar[d] \ar[r] & \overline{L} \ar[d] \ar[r] & 0 \\               
            &  0            \ar[r]        & U          \ar[r]  & H \ar[r]           & B  \ar[r] & 0.      
}
\end{equation}
\end{Pt}

\begin{Pt}\label{structuremncompl}
Si $L_K$ est connexe, on peut compléter cette description de la structure du modèle de Néron $L$ de \ref{casgenmn}. Notons 
$T$ le modèle de Néron usuel localement de type fini du $K$-tore déployé $T_K$. D'après le lemme \ref{lemTotaro} qui suit, et la démonstration de \cite{BLR} 10.1/7, la suite exacte 
$$
\xymatrix{
0\ar[r] & T_K \ar[r] & \overline{L_K} \ar[r] & B_K \ar[r] & 0
}
$$
se prolonge en une suite exacte
$$
\xymatrix{
0\ar[r] & T \ar[r] & \overline{L} \ar[r] & B \ar[r] & 0.
}
$$
En prenant son image réciproque par le morphisme canonique $H\ra B$, on obtient un diagramme commutatif exact complétant (\ref{LoLHB})
\begin{equation}\label{LoLHBT}
\xymatrix{
            &                               &                                & 0  \ar[d]                & 0 \ar[d] & \\
           &                                &                               & T \ar[d] \ar@{=}[r] & T \ar[d] & \\ 
           &  0            \ar[r]       & U   \ar@{=}[d] \ar[r] & L \ar[d] \ar[r] & \overline{L} \ar[d] \ar[r] & 0 \\               
            &  0            \ar[r]        & U          \ar[r]  & H \ar[r] \ar[d]          & B  \ar[r] \ar[d] & 0 \\
            &                               &                       &  0                    & 0. &      
}
\end{equation}
\end{Pt}

\begin{Lem}\label{lemTotaro}
Soient $K$ un corps et $H_K$ un $K$-schéma en groupes lisse connexe. Supposons que $H_K$ ne possède pas de sous-groupe isomorphe à $\bbG_{m,K}$. Alors $\Hom_{\SchgrK}(H_K,\bbG_{m,K})=0$.
\end{Lem}

\begin{proof}
Soit $h_K\in\Hom_{\SchgrK}(H_K,\bbG_{m,K})$. Notons $H_{K}'$ le plus grand sous-groupe lisse affine connexe de $H_K$. 
Le tore maximal $H_{K,m}'$ de $H_{K}'$ étant anisotrope par hypothèse, la restriction de $h_K$ à $H_{K,m}'$ est nulle. Le quotient  $H_{K}'/H_{K,m}'$ étant unipotent, on en déduit que la restriction de $h_K$ à $H_{K}'$ est nulle, et passe donc au quotient en un morphisme 
$\overline{h_K}:H_{K}'':=H_{K}/H_{K}'\ra\bbG_{m,K}$. Par construction, le $K$-schéma en groupes lisse connexe $H_{K}''$ est une variété pseudo-abélienne au sens de \cite{To}, donc est extension d'un groupe lisse unipotent connexe par une variété abélienne d'après \emph{loc. cit.} 2.1. Il en résulte que $\overline{h_K}$ est nul.
\end{proof}

\begin{Pt}\label{prefoncNer}
D'après \ref{unicitemn}, le modèle de Néron de $L_K$ fourni par la construction \ref{casgenmn} est l'unique modèle de Néron de $L_K$ composé de $R$-schémas en groupes dont le groupe ordinaire des composantes connexes de la fibre spéciale est de type fini, et le ind-$R$-schéma en groupes sous-jacent à ce modèle de Néron dépend fonctoriellement de $L_K$. 
\end{Pt}

\begin{Def} \label{foncNer}
Supposons $R$ complet et $k$ séparablement clos. Le \emph{foncteur de Néron} est le foncteur 
$$
\xymatrix{
(\textrm{\emph{$K$-schémas en groupes lisses de type fini}}) \ar[d]^{\N} \\
\Ind(\textrm{\emph{$R$-schémas en groupes lisses séparés}})
}
$$
défini par \ref{prefoncNer}.
\end{Def}

\section{Composante neutre et groupe des composantes connexes}

\begin{Pt} \label{precnpi0}
D'après \cite{SGA 3} $\VIB$ 3.10, on dispose d'un foncteur \emph{composante neutre} 
$$
\xymatrix{
(\cdot)^{0}:(\textrm{$R$-schémas en groupes lisses}) \ar[r] & (\textrm{$R$-schémas en groupes lisses de type fini}).
}
$$
Par ailleurs, en composant le foncteur \emph{fibre spéciale} et le foncteur $\pi_0$ \ref{Notpi0} \emph{schéma en groupes des composantes connexes}, on obtient un foncteur \emph{schéma en groupes des composantes connexes de la fibre spéciale}
$$
\xymatrix{
\pi_0((\cdot)_k):(\textrm{$R$-schémas en groupes localement de type fini}) \ar[r] & (\textrm{$k$-schémas en groupes étales})
}
$$
Ces foncteurs s'étendent aux ind-catégories (\cite{SGA 4} I 8.6) :
\end{Pt}

\begin{Def} \label{cnpi0}
\begin{itemize}

\medskip

\item Le foncteur \emph{composante neutre} 
$$
\xymatrix{
\Ind(\textrm{\emph{$R$-schémas en groupes lisses}}) \ar[d]^{(\cdot)^{0}} \\
 \Ind(\textrm{\emph{$R$-schémas en groupes lisses de type fini}})
}
$$
est l'indisé du foncteur composante neutre \ref{precnpi0}.

\medskip

\item Le foncteur \emph{schéma en groupes des composantes connexes de la fibre spéciale}
$$
\xymatrix{
\Ind(\textrm{\emph{$R$-schémas en groupes localement de type fini}}) \ar[d]^{\pi_0((\cdot)_k)} \\
\Ind(\textrm{\emph{$k$-schémas en groupes étales}})
}
$$
est l'indisé du foncteur schéma en groupes des composantes connexes de la fibre spéciale \ref{precnpi0}.
\end{itemize}
\end{Def}

\begin{Prop}\label{gctf}
Supposons $R$ complet et $k$ séparablement clos. Soient $L_K$ un $K$-schéma en groupes lisse de type fini, et $L=(L_i)_{i\in I}$ un modèle de Néron de $L_K$ avec $\pi_0(L_{i,k})$ de type fini pour tout $i\in I$. Alors le ind-objet de la catégorie des groupes $\pi_0(L_k)$ est isomorphe à un groupe de type fini.
\end{Prop}

\begin{proof}
D'après \ref{prefoncNer}, on peut supposer $L=\N(L_K)$. Alors, en reprenant les notations de \ref{structuremn}, le système inductif des $L_i$ s'insère dans un système inductif d'extensions
$$
\xymatrix{
0  \ar[r]   & U_i   \ar[r] & L_i \ar[r] & \overline{L} \ar[r] & 0,\quad i\in \bbN,                 
}
$$
où les $U_i$ sont à fibres connexes, et $\pi_0(\overline{L}_k)$ est un groupe ordinaire de type fini. Le ind-groupe 
$\pi_0(L_k)=(\pi_0(L_{i,k}))_{i\in\bbN}$ est alors isomorphe à $\pi_0(\overline{L}_k)$.
\end{proof}

\section{Modèles de Néron et épimorphismes}

\subsection{Préliminaires sur les tores maximaux et maximaux déployés}

\begin{Lem*}\label{preDeftoremax}
Soient $K$ un corps et $L_K$ un $K$-schéma en groupes lisse connexe. Soient : 
\begin{itemize}

\medskip

\item $L_{K}'$ le plus grand sous-groupe lisse affine connexe de $L_K$, et $L_{K,m}'$ le tore maximal de $L_{K}'$ ;

\medskip

\item $C_K$ un sous-groupe affine connexe de $L_K$ tel que le quotient $L_K/C_K$ soit une variété abélienne (cf. \cite{BLR} 9.2/1), et $C_{K,m}$ le plus grand sous-groupe de type multiplicatif de $C_K$.

\medskip

\end{itemize}

Alors $L_{K,m}'=(C_{K,m})_{\red}^0$.
\end{Lem*}

\begin{proof}
Comme $(C_{K,m})_{\red}^0$ est un tore, il est contenu dans $L_{K}'$ puis dans $L_{K,m}'$. Réciproquement, l'image du tore $L_{K,m}'$ dans la variété abélienne $L_K/C_K$ est nulle, donc $L_{K,m}'$ est contenu dans $C_K$, puis dans 
$(C_{K,m})_{\red}^0$.
\end{proof}

\begin{Def*}\label{Deftoremax}
Soient $K$ un corps et $L_K$ un $K$-schéma en groupes lisse connexe. Le \emph{tore maximal} de $L_K$ est le tore 
$$
L_{K,m}:=L_{K,m}'=(C_{K,m})_{\red}^0
$$
de \ref{preDeftoremax}.
\end{Def*}

\begin{Lem*} \label{surjtoresmax}
Soient $K$ un corps et $f_K:L_{K,0}\ra L_{K,1}$ un morphisme surjectif de $K$-schémas en groupes lisses connexes. Alors le morphisme induit par $f_K$ sur les tores maximaux est surjectif.
\end{Lem*}

\begin{proof}
Soit $C_{K,0}$ un sous-groupe affine connexe de $L_{K,0}$ tel que $L_{K,0}/C_{K,0}$ soit une variété abélienne. Posant $C_{K,1}:=f_K(C_{K,0})$, on obtient un morphisme fidèlement plat $C_{K,0}\ra C_{K,1}$ et un diagramme commutatif exact
$$
\xymatrix{
0 \ar[r] & C_{K,0} \ar[r] \ar[d] & L_{K,0} \ar[r] \ar[d]^{f_K} & L_{K,0}/C_{K,0} \ar[r] \ar[d] & 0 \\
0 \ar[r] & C_{K,1} \ar[r] & L_{K,1} \ar[r] & L_{K,1}/C_{K,1} \ar[r] & 0,
}
$$
qui montre que $L_{K,1}/C_{K,1}$ est une variété abélienne. Pour $i=0,1$, le quotient de $C_{K,i}$ par son plus grand sous-groupe de type multiplicatif $C_{K,i,m}$ est unipotent. Le lemme du serpent appliqué au diagramme commutatif exact
$$
\xymatrix{
0 \ar[r] & C_{K,0,m} \ar[r] \ar[d] & C_{K,0} \ar[r] \ar[d] & C_{K,0}/C_{K,0,m} \ar[r] \ar[d] & 0 \\
0 \ar[r] & C_{K,1,m} \ar[r] & C_{K,1} \ar[r] & C_{K,1}/C_{K,1,m} \ar[r] & 0,
}
$$
montre donc que le morphisme $C_{K,0,m}\ra C_{K,1,m}$ induit par $f_K$ est fidèlement plat. Le morphisme $(C_{K,0,m})_{\red}^0\ra (C_{K,1,m})_{\red}^0$ qui s'en déduit est donc surjectif.
\end{proof}

\begin{Lem*} \label{surjtoresdepmax}
Supposons $R$ complet et $k$ algébriquement clos. Soit $f_K:L_{K,0}\ra L_{K,1}$ un morphisme surjectif de $K$-schémas en groupes lisses connexes. Alors le morphisme induit par $f_K$ sur les plus grands sous-tores déployés est surjectif.
\end{Lem*}

\begin{proof}
D'après \ref{surjtoresmax}, le morphisme $L_{K,0,m}\ra L_{K,1,m}$ induit par $f_K$ sur les tores maximaux est surjectif. Notant $N_K$ son noyau, on a donc un diagramme  commutatif exact
$$
\xymatrix{
0 \ar[r] & (N_K)_{\red}^0 \ar[r] \ar[d] & L_{K,0,m} \ar[r] \ar@{=}[d] & \widetilde{L_{K,1,m}} \ar[r] \ar[d] & 0 \\
0 \ar[r] & N_K \ar[r] & L_{K,0,m} \ar[r] & L_{K,1,m} \ar[r] & 0. 
}
$$
Le morphisme $\widetilde{L_{K,1,m}}\ra L_{K,1,m}$ ainsi défini est une isogénie de $K$-tores ; le morphisme $\widetilde{T_{K,1}}\ra T_{K,1}$ qui s'en déduit sur les plus grands sous-tores déployés est donc également une isogénie. Par ailleurs, la première ligne du diagramme est une suite exacte de $K$-tores. D'après les hypothèses sur $K$, ceux-ci sont cohomologiquement triviaux : \cite{Beg} 3.2.1 (la preuve de \emph{loc. cit.} est valable pour $K$ de valuation discrète complet à corps résiduel algébriquement clos, en inégales comme en égales caractéristiques). La suite obtenue sur les plus grands sous-tores déployés est donc exacte, en particulier $T_{K,0}\ra \widetilde{T_{K,1}}$ est surjectif. Par composition, on obtient que $T_{K,0}\ra T_{K,1}$ est surjectif.
\end{proof}

\begin{Rem*}\label{passurjunip}
Le morphisme induit par $f_K$ sur les plus grands sous-groupes lisses unipotents connexes déployés n'est lui pas surjectif en général. Par exemple, supposons $R$ de caractéristique $p>0$ et considérons le sous-groupe $L_{K,0}$ de $\bbG_{a,K}^2$ d'équation
$$
x+x^p+ty^p=0
$$
où $t$ est une uniformisante de $R$. La projection sur le premier (resp. second) facteur de $\bbG_{a,K}^2$ induit une isogénie $f_K:L_{K,0}\ra \bbG_{a,K}=:L_{K,1}$. Mais comme $L_{K,0}$ est une forme non triviale de $\bbG_{a,K}$, le morphisme induit sur les plus grands sous-groupes lisses unipotents connexes déployés est $0\ra\bbG_{a,K}$. 
\end{Rem*}

\subsection{Effet d'un épimorphisme sur les modèles de Néron}

\begin{Lem*} \label{CokerGrreptoresdep}
Soit $K$ un corps de valuation discrète. Soit $f_K:T_{K,0}\ra T_{K,1}$ un morphisme surjectif de $K$-tores déployés. Soit $f:T_0\ra T_1$ le prolongement de $f_K$ aux modèles de Néron. Alors pour tout $n\geq 1$, le conoyau de $\Gr_n(f):\Gr_n(T_0)\ra\Gr_n(T_1)$ dans la catégorie des faisceaux fppf sur le corps résiduel $k$ est représentable. Le $k$-schéma en groupes $\Coker(f_k)$ est constant fini, et pour tout $n\geq 1$, $\Coker(\Gr_{n+1}(f))$ est extension de $\Coker(\Gr_n(f))$ par un $k$-groupe vectoriel.

Si de plus $f_K$ est fini, alors $f_k$ est fini, et pour tout $n\geq 1$, $\Gr_n(f):\Gr_n(T_0)\ra\Gr_n(T_1)$ est affine.
\end{Lem*}

\begin{proof}
Notant $N_K$ le noyau de $f_K$, on a un diagramme commutatif exact
$$
\xymatrix{
0 \ar[r] & (N_K)_{\red}^0 \ar[r] \ar[d] & T_{K,0} \ar[r]^{\widetilde{f_K}} \ar@{=}[d] & \widetilde{T_{K,1}} \ar[r] \ar[d]^{i_K} & 0 \\
0 \ar[r] & N_K \ar[r] & T_{K,0} \ar[r]^{f_K} & T_{K,1} \ar[r] & 0. 
}
$$
Comme $T_{K,0}$ est un $K$-tore déployé, les $K$-schémas en groupes $(N_K)_{\red}^0$ et $\widetilde{T_{K,1}}$ sont des $K$-tores déployés, et la première ligne est donc exacte scindée. Le prolongement $\widetilde{f}:T_0\ra \widetilde{T_1}$ de $\widetilde{f_K}$ aux modèles de Néron est donc scindé, puis, pour tout $n\geq 1$, $\Gr_n(\widetilde{f}):\Gr_n(T_0)\ra\Gr_n(\widetilde{T_1})$ est scindé, en particulier est un épimorphisme. Le conoyau de $\Gr_n(f):\Gr_n(T_0)\ra \Gr_n(T_1)$ coïncide donc avec celui de $\Gr_n(i):\Gr_n(\widetilde{T_1})\ra \Gr_n(T_1)$, où $i$ est le prolongement de $i_K$ aux modèles de Néron. 

Pour $n=1$, considérons le diagramme commutatif exact
$$
\xymatrix{
0\ar[r] & (\widetilde{T_1})_{k}^0 \ar[r] \ar[d]^{i_{k}^0} & (\widetilde{T_1})_k \ar[r] \ar[d]^{i_k} & \pi_0((\widetilde{T_1})_k) \ar[r] \ar[d]^{\pi_0(i_k)} & 0 \\
0\ar[r] & (T_1)_{k}^0 \ar[r] & (T_1)_k \ar[r]  & \pi_0((T_1)_k) \ar[r]  & 0,
}
$$
où, par construction du modèle de Néron d'un tore déployé, $(\widetilde{T_1})_{k}^0$ et $ (T_1)_{k}^0$ sont des $k$-tores (déployés) de dimension $d=\dim \widetilde{T_{K,1}}=\dim T_{K,1}$, et $\pi_0((\widetilde{T_1})_k)$ et $\pi_0((T_1)_k)$ sont des $\bbZ$-modules libres de rang $d$. Comme $i_K$ est une isogénie de $K$-tores par construction, $i_{k}^0$ est une isogénie de $k$-tores et $\pi_0(i_k)$ est injectif à conoyau fini. En particulier $i_k$ est quasi-compact, puis fini, et son conoyau est isomorphe à un $k$-schéma en groupes constant fini. 

Pour tout $n\geq 1$, on a d'après \ref{troncationlisse} un diagramme commutatif exact
$$
\xymatrix{
0\ar[r] & V(\omega_{(\widetilde{T_1})_{k}}) \ar[r] \ar[d]^{V(i_{k}^*)} & \Gr_{n+1}(\widetilde{T_1}) \ar[r] \ar[d]^{\Gr_{n+1}(i)} & \Gr_n(\widetilde{T_1}) \ar[r] \ar[d]^{\Gr_{n}(i)} & 0 \\
0\ar[r] & V(\omega_{(T_1)_{k}}) \ar[r] & \Gr_{n+1}(T_1) \ar[r] & \Gr_n(T_1) \ar[r] & 0. 
}
$$
Par récurrence sur $n\geq 1$, $\Gr_n(i)$ est donc affine, en particulier quasi-compact ; son conoyau est donc représentable.
En utilisant le lemme du serpent, on obtient alors que $\Coker(\Gr_{n+1}(i))$ est extension de $\Coker(\Gr_n(i))$ par un quotient de $V(\omega_{(T_1)_{k}})$. Un tel quotient est encore un groupe vectoriel. 
\end{proof}

\begin{Rem*}\label{pasqctores}
Ainsi donc, le conoyau de $\Gr_n(f):\Gr_n(T_0)\ra \Gr_n(T_1)$ est représentable bien que ce morphisme ne soit pas quasi-compact en général.
\end{Rem*}

\begin{Prop*} \label{Cokerrep}
Supposons $R$ complet et $k$ algébriquement clos. Soit $f_K:L_{K,0}\ra L_{K,1}$ un morphisme surjectif de $K$-schémas en groupes lisses connexes. Soient 
$$
L_0=(L_{0,i})_{i\in\bbN}:=\N(L_{0,K})\quad\textrm{et}\quad L_1=(L_{1,j})_{j\in\bbN}:=\N(L_{1,K})
$$
les modèles de Néron \ref{foncNer} de $L_{K,0}$ et $L_{K,1}$. Soit 
$$
\xymatrix{
f_{i,j}:L_{0,i} \ar[r] & L_{1,j}
}
$$
prolongeant $f_K$. Alors pour tout $n\geq 1$, le conoyau de $\Gr_n(f_{i,j}):\Gr_n(L_{0,i})\ra\Gr_n(L_{1,j})$ dans la catégorie des faisceaux fppf sur $k$ est représentable. Le $k$-schéma en groupes $\Coker((f_{i,j})_k)$ est lisse de type fini, et pour tout $n\geq 1$, $\Coker(\Gr_{n+1}(f_{i,j}))$ est extension de $\Coker(\Gr_n(f_{i,j}))$ par un $k$-groupe vectoriel.

Si de plus $f_K$ est fini, alors pour tout $n\geq 1$, $\Gr_n(f_{i,j}):\Gr_n(L_{0,i})\ra\Gr_n(L_{1,j})$ est quasi-compact.
\end{Prop*}

\begin{proof}
Reprenant les notations de \ref{structuremn} et \ref{structuremncompl}, le morphisme $f_K$ induit un diagramme commutatif exact 
$$
\xymatrix{
0 \ar[r] & T_{K,0} \ar[r] \ar[d] & L_{K,0} \ar[r] \ar[d]^{f_K} & H_{K,0} \ar[r] \ar[d] & 0 \\
0 \ar[r] & T_{K,1} \ar[r] & L_{K,1} \ar[r] & H_{K,1} \ar[r] & 0, 
}
$$
qui se prolonge en un diagramme commutatif exact
$$
\xymatrix{
0 \ar[r] & T_0 \ar[r] \ar[d] & L_{0,i} \ar[r] \ar[d]^{f_{i,j}} & H_{0,i} \ar[r] \ar[d] & 0 \\
0 \ar[r] & T_1 \ar[r] & L_{1,j} \ar[r] & H_{1,j} \ar[r] & 0. 
}
$$
Pour tout $n\geq 1$, le diagramme commutatif 
\begin{equation}\label{GrnTLH}
\xymatrix{
0 \ar[r] & \Gr_n(T_0) \ar[r] \ar[d] & \Gr_n(L_{0,i}) \ar[r] \ar[d]^{\Gr_n(f_{i,j})} & \Gr_n(H_{0,i}) \ar[r] \ar[d] & 0 \\
0 \ar[r] & \Gr_n(T_1) \ar[r] & \Gr_n(L_{1,j}) \ar[r] & \Gr_n(H_{1,j}) \ar[r] & 0 
}
\end{equation}
est alors exact d'après \ref{Grsurjlisse}. Maintenant, d'après \ref{surjtoresdepmax}, le morphisme $T_{K,0}\ra T_{K,1}$ induit par $f_K$ sur les plus grands sous-tores déployés est surjectif. D'après \ref{CokerGrreptoresdep}, le conoyau de $\Gr_n(T_0)\ra \Gr_n(T_1)$ est donc représentable par un $k$-schéma en groupes affine de type fini. Comme par ailleurs $\Gr_n(H_{0,i})\ra \Gr_n(H_{1,j})$ est un morphisme de $k$-schémas en groupes de type fini, donc à noyau de type fini et à conoyau représentable de type fini, on en déduit en utilisant le lemme du serpent que le conoyau de $\Gr_n(f_{i,j})$ est représentable de type fini. Il est de plus lisse puisque $\Gr_n(L_{1,j})$ l'est (\ref{troncationlisse}). Enfin, pour tout $n\geq 1$, on a d'après \ref{troncationlisse} un diagramme commutatif exact
$$
\xymatrix{
0\ar[r] & V(\omega_{(L_{0,i})_k}) \ar[r] \ar[d]^{V((f_{i,j})_{k}^*)} & \Gr_{n+1}(L_{0,i}) \ar[r] \ar[d]^{\Gr_{n+1}(f_{i,j})} & \Gr_n(L_{0,i}) \ar[r] \ar[d]^{\Gr_{n}(f_{i,j})} & 0 \\
0\ar[r] & V(\omega_{(L_{1,j})_k}) \ar[r] & \Gr_{n+1}(L_{1,j}) \ar[r] & \Gr_n(L_{1,j}) \ar[r] & 0, 
}
$$
qui montre que $\Coker(\Gr_{n+1}(f_{i,j}))\ra\Coker(\Gr_n(f_{i,j}))$ est un épimorphisme dont le noyau est quotient du groupe vectoriel $V(\omega_{(L_{1,j})_k})$. 

Si de plus $f_K$ est fini, alors le morphisme $T_{K,0}\ra T_{K,1}$ l'est aussi, et par conséquent pour tout $n\geq 1$, 
$\Gr_n(T_0)\ra \Gr_n(T_1)$ est affine d'après \ref{CokerGrreptoresdep}. En appliquant le lemme du serpent au diagramme (\ref{GrnTLH}), on obtient donc que le noyau de $\Gr_n(f_{i,j}):\Gr_n(L_{0,i})\ra\Gr_n(L_{1,j})$ est quasi-compact. De plus, comme $\Coker(\Gr_n(f_{i,j}))$ est représentable, l'image de $\Gr_n(f_{i,j})$ est représentable par un sous-$k$-schéma en groupes fermé de $\Gr_n(L_{1,j})$. Il en résulte que $\Gr_n(f_{i,j})$ est quasi-compact.
\end{proof}

\chapter{La cohomologie des $K$-schémas en groupes finis}

Soient $K$ un corps discrètement valué complet de caractéristique $p>0$, $R$ son anneau d'entiers et $k$ son corps résiduel, que l'on suppose algébriquement clos. 

\section{Nullité du $H^i$ pour $i>1$}

\begin{Prop}\label{dimcohf}
Soit $G_K$ un $K$-schéma en groupes fini. Alors
$$\forall i>1,\quad H^i(K,G_K)=0.$$
\end{Prop}

\begin{proof} 

\emph{a) Le cas où $G_K$ est de type multiplicatif.} Il existe $K'/K$ finie \emph{galoisienne} telle que $G_{K'}$ soit diagonalisable (\cite{SGA 3} X 4.6 (i)). Le groupe $G_{K'}$ se réalise alors comme sous-groupe d'un $K'$-tore déployé $T'$.  Le groupe $G_K$ apparaît donc comme sous-groupe de la restriction de Weil $\Pi_{K'/K}T'$, qui est un tore. Or le quotient d'un tore est un tore, et les $K$-tores sont $\Gamma(K,\cdot)$-acycliques (cf. \cite{Beg} 3.2.1). D'où la proposition dans le cas multiplicatif.

\medskip

\emph{b) Le cas où $G_K$ est unipotent radiciel.} Alors $G$ possède une suite de composition dont les quotients successifs sont isomorphes à $\alpha_{p,K}$ (\cite{SGA 3} XVII 4.2.1). On peut donc supposer $G_K=\alpha_{p,K}$, et la proposition résulte alors de la $\Gamma(K,\cdot)$-acyclicité de $\bbG_{a,K}$.

\medskip

\emph{c) Le cas où $G_K$ est unipotent étale.} Alors il existe $K'/K$ finie galoisienne telle que $G_{K'}$ soit un $p$-groupe ordinaire. Celui-ci se réalise comme sous-groupe d'un produit $W'$ de $K'$-schémas en groupes de Witt de longueur finie, via des $K'$-isogénies d'Artin-Schreier. Par suite, le $K$-groupe $G_K$ se réalise comme sous-groupe de la restriction de Weil $\Pi_{K'/K}W'$. Maintenant, $W'$ est un $K'$-groupe lisse unipotent connexe déployé. Il en résulte que $\Pi_{K'/K}W'$ est un $K$-groupe lisse unipotent connexe déployé : en effet  $\Pi_{K'/K}\bbG_{a,K'}$ est produit de copies de $\bbG_{a,K}$ et la restriction de Weil le long de $K'/K$ est exacte. \emph{A fortiori}, le quotient de $\Pi_{K'/K}W'$ par $G_K$ est déployé (cf. \ref{quotresol} ci-dessous). Comme les $K$-groupes lisses unipotents connexes déployés sont $\Gamma(K,\cdot)$-acycliques, on en déduit la proposition dans le cas unipotent étale.

\medskip

\emph{d) Le cas où $G_K$ est unipotent.} Alors $G_K$ est extension d'un $K$-groupe étale par un $K$-groupe radiciel. On conclut à l'aide des deux cas précédents.

\medskip

\emph{e) Le cas général.} Le $K$-groupe $G_K$ est extension d'un $K$-groupe unipotent par un $K$-groupe de type multiplicatif 
(\cite{SGA 3} XVII 7.2.1 a)). D'où la proposition.
\end{proof}

\begin{Lem}\label{quotresol}
Soit $G_K$ un $K$-groupe  lisse unipotent connexe déployé. Alors tout quotient de $G_K$ est déployé.
\end{Lem}

\begin{proof} 
Soit $Q_K$ un quotient de $G_K$. Soit 
$$0= G_{K,0}\subset G_{K,1}\subset\cdots\subset G_{K,n}=G_K$$
une suite de composition dont les quotients successifs sont isomorphes à $\bbG_{a,K}$. Notant $Q_{K,i}$ l'image de $G_{K,i}$ dans $Q_K$, le groupe $Q_{K,i+1}/Q_{K,i}$ est quotient du groupe $G_{K,i+1}/G_{K,i}$. D'après \cite{SGA 3} XVII 2.3, $Q_{K,i+1}/Q_{K,i}$ est donc nul ou isomorphe à $\bbG_{a,K}$, ce qui achève la démonstration.
\end{proof}

Bien qu'inutile pour la suite, notons le 

\begin{Cor}\label{dimcoha}
Soit $G_K$ un $K$-schéma en groupes de type fini affine. Alors
$$\forall i > 1,\quad H^i(K,G_K)=0.$$
\end{Cor}

\begin{proof} \emph{Le cas où $G_K$ est de type multiplicatif.} Même preuve que celle du cas multiplicatif de la proposition \ref{dimcohf}.

\medskip

\emph{Le cas où $G_K$ est unipotent connexe.} Il existe une suite exacte
$$
\xymatrix{
0\ar[r] & G_{K}'\ar[r] & G_K\ar[r] & G_{K}'' \ar[r] & 0
}
$$
avec $G_{K}'$ radiciel et $G_{K}''$ lisse déployé (\cite{SGA 3} XVII 4.3.1). La proposition \ref{dimcohf} conclut dans ce cas.

\medskip

\emph{Le cas où $G_K$ est unipotent.} Il existe une suite exacte
$$
\xymatrix{
0\ar[r] & G_{K}'\ar[r] & G_K\ar[r] & G_{K}''\ar[r] & 0
}
$$
avec $G_{K}'$ connexe et $G_{K}''$ étale, d'où à nouveau le résultat souhaité en utilisant le cas précédent et \ref{dimcohf}.

\medskip

\emph{Le cas général.} On écrit $G_K$ comme extension d'un $K$-groupe unipotent par un $K$-groupe de type multiplicatif.
\end{proof}

\section{Structure quasi-algébrique sur le $H^0$}

\begin{Def}\label{DefcH0}
Soit $G_K$ un $K$-schéma en groupes fini. La \emph{structure de groupe quasi-algébrique sur $H^0(K,G_K)$} est l'unique groupe quasi-algébrique fini sur $k$ dont le groupe ordinaire sous-jacent est $H^0(K,G_K)$. 
\end{Def}

\begin{Not}\label{NotcH0}
On note $\cH^0(G_K)$ la structure de groupe quasi-algébrique sur $H^0(K,G_K)$.
\end{Not}

\begin{Pt} \label{cH0fonc}
On a ainsi un \emph{foncteur}
$$
\xymatrix{
\cH^0:\textrm{($K$-schémas en groupes finis)} \ar[r] & \cQ.
}
$$
\end{Pt}

\begin{Lem}\label{schfinipt}
Soient $X_K$ un $K$-schéma fini et $k'/k$ une extension de corps. Alors l'application canonique
$$
\xymatrix{
X_K(K)\ar[r] & X_{K_{k'}}(K_{k'})
}
$$
(notation \ref{Kk'}) est bijective.
\end{Lem}

\begin{proof}
L'application considérée est évidemment injective, et il s'agit donc de montrer qu'elle est surjective. Soit $f:\Spec(K_{k'})\ra X_K$ un morphisme de $K$-schémas. Comme $\Spec(K_{k'})$ est connexe et réduit, $f$ se factorise par le réduit d'une composante connexe de $X_K$, qui est de la forme $\Spec(L)$ pour un corps $L$ fini sur $K$. On obtient ainsi des inclusions $K\subset L \subset K_{k'}$. La valuation discrète normalisée de $K_{k'}$ induit sur $L$ une valuation qui prolonge la valuation discrète normalisée de $K$. Par conséquent, l'extension $L/K$ est non ramifiée. Comme $k$ est algébriquement clos, on a donc nécessairement $L=K$.
\end{proof}

\begin{Cor}\label{cH0cb}
Si $k'/k$ est une extension de corps parfaite de $k$, alors on a un isomorphisme canonique
$$
\xymatrix{
\cH^0(G_K)_{k'}\ar[r] & \cH^0(G_{K_{k'}}).
}
$$
\end{Cor}

\section{Résolutions lisses}

\begin{Prop}\label{resolf}
Soient $G_K$ un $K$-schéma en groupes fini, et 
$$
\xymatrix{
0\ar[r] & G_K\ar[r] & L_K\ar[r] & L_{K}'\ar[r] & 0
}
$$
la résolution lisse standard de $G_K$ (\cite{Beg} 2.2.1). Alors $L_K$ et $L_{K}'$ sont connexes, et pour toute extension de corps algébriquement close $\kappa/k$, $L_{K_{\kappa}}$ et $L_{K_{\kappa}}'$ sont $\Gamma(K_{\kappa},\cdot)$-acycliques (\ref{Kk'}).
\end{Prop}

\begin{proof}
Par définition, $L_K$ est la restriction de Weil 
$$\Pi_{G_{K}^*/K}\bbG_{m,G_{K}^*}$$
où $G_{K}^*$ désigne le dual de Cartier de $G_K$. C'est donc un $K$-schéma en groupes lisse connexe (cf. par exemple \cite{CGP} Appendice A, A.5.2), et il en va donc de même de $L_{K}'$. Passons à l'assertion d'acyclicité. Comme la résolution lisse standard commute aux changements de base, en particulier à l'extension $K_{\kappa}/K$, il suffit de montrer que $L_K$ et $L_{K}'$ sont $\Gamma(K,\cdot)$-acycliques. Comme la cohomologie fppf d'un groupe lisse coïncide avec sa cohomologie étale, et $G_{K}^*$ est fini sur $K$, on a 
$$\forall i\geq 0,\quad H^i(K,L_K)=H^i(G_{K}^*,\bbG_{m,G_{K}^*}).$$
Montrons que 
\begin{equation}\label{Gmred}
\forall i\geq 1,\quad H^i(G_{K}^*,\bbG_{m,G_{K}^*})=H^i((G_{K}^*)_{\red},\bbG_{m,(G_{K}^*)_{\red}}).
\end{equation}
Quitte à décomposer l'immersion fermée $(G_{K}^*)_{\red}\hra G_{K}^*$ en épaississements d'ordre 1, on peut supposer que l'idéal quasi-cohérent $\cI$ de $\cO_{G_{K}^*}$ définissant cette immersion est de carré nul. Dans ce cas, on peut considérer la suite exacte exponentielle
$$
\xymatrix{
0\ar[r] & \cI \ar[r] & \bbG_{m,G_{K}^*} \ar[r] & \bbG_{m,(G_{K}^*)_{\red}} \ar[r] & 0
}
$$
sur le petit site étale de $G_{K}^*$. Comme $G_{K}^*$ est un schéma affine, $H^i(G_{K}^*,\cI)=0$ pour tout $i\geq 1$, d'où (\ref{Gmred}). Maintenant, le schéma $(G_{K}^*)_{\red}$ étant somme finie de spectres d'extensions finies du corps $K$, les hypothèses entraînent que 
$$H^i((G_{K}^*)_{\red},\bbG_{m,(G_{K}^*)_{\red}})=0$$
pour tout $i\geq 1$ (\cite{Beg} 3.2.1). D'où l'acyclicité de $L_K$. Celle de $L_{K}'$ est ensuite conséquence de \ref{dimcohf}.
\end{proof}

Les $K$-schémas en groupes $L_{K}$ et $L_{K}'$ de la résolution standard du $K$-schéma en groupes fini $G_K$ sont lisses, affines et connexes. Plus généralement, on a le fait suivant (néanmoins inutile pour la suite) :

\begin{Prop} \label{resola}
Soit $G_K$ un $K$-schéma en groupes de type fini affine. Il existe une suite exacte 
$$
\xymatrix{
0\ar[r] & G_K\ar[r] & L_{K} \ar[r] & L_{K}' \ar[r] & 0
}
$$
avec $L_K$ et $L_{K}'$ lisses affines connexes $\Gamma(K,\cdot)$-acycliques.
\end{Prop}

\begin{proof} 
D'après \ref{dimcoha}, il suffit de montrer qu'il existe un monomorphisme de $G_K$ dans un $K$-schéma en groupes lisse affine connexe $\Gamma(K,\cdot)$-acyclique $L_K$. 

\medskip

\emph{$G_K$ est de type multiplicatif.}
Le résultat a été vu dans la démonstration du cas multiplicatif de la proposition \ref{dimcoha}. 

\medskip

\emph{$G_K$ est unipotent connexe.}
D'après \cite{SGA 3} XVII 4.3.1, $G_K$ est extension d'un lisse déployé $G_{K}''$ par un radiciel $G_{K}'$. D'après la proposition \ref{resolf}, le $K$-schéma en groupes fini $G_{K}'$ se plonge dans un $L_{0,K}$ lisse affine connexe $\Gamma(K,\cdot)$-acyclique. En poussant par l'inclusion $G_{K}'\subseteq L_{0,K}$, on obtient une inclusion $G_K\subseteq L_K$ convenable :
\begin{displaymath}
\xymatrix{
0\ar[r]& L_{0,K} \ar[r]  & L_K \ar[r]  & G_{K}'' \ar[r] & 0 \\
0\ar[r] & G_{K}'\ar[u] \ar[r] & G_K \ar[r] \ar[u] & G_{K}''\ar[r] \ar@{=}[u] & 0 \\
&    0 \ar[u]     & 0. \ar[u] &  & 
}
\end{displaymath}

\medskip

\emph{$G_K$ est unipotent.} D'après le cas précédent, la composante neutre $G_{K}^0$ de $G_K$ se plonge dans un $L_{0,K}$ lisse affine connexe $\Gamma(K,\cdot)$-acyclique. En poussant par l'inclusion $G_{K}^0\subseteq L_{0,K}$, on obtient un diagramme :
\begin{displaymath}
\xymatrix{
0\ar[r]& L_{0,K} \ar[r]  & H_K \ar[r]  & G_K/G_{K}^0 \ar[r] & 0 \\
0\ar[r] & G_{K}^0\ar[u] \ar[r] & G_K \ar[r] \ar[u] & G_K/G_{K}^0\ar[r] \ar@{=}[u] & 0 \\
&    0 \ar[u]     & 0. \ar[u] &  & 
}
\end{displaymath}
Il suffit alors de prouver le résultat pour $H_K$, ce qui nous ramène au cas où $G_{K}^0$ est lisse $\Gamma(K,\cdot)$-acyclique.

Par ailleurs, on peut trouver un sous-groupe fini $F_K$ de $G_K$ qui rencontre chaque composante connexe de $G_K$. Pour le voir, on peut considérer une extension finie galoisienne $K'/K$ telle que chaque composante connexe de $G_{K'}$ possède un $K'$-point. On choisit une section ensembliste de $G_{K'}\ra G_{K'}/G_{K'}^0$, et on considère le sous-groupe $F_{K'}$ de $G_{K'}$ engendré par les points de cette section et par leurs conjugués sur $K$. Alors $F_{K'}$ est fini puisque $G_{K'}$ est  unipotent. Et $F_{K'}$ descend en un $K$-schéma en groupes fini $F_K$ convenable. 

On dispose alors  d'un diagramme commutatif exact
\begin{displaymath}
\xymatrix{
         &   0                          & 0           & & \\
         &  F_K/(G_{K}^0\cap F_K)  \ar[u] \ar[r]^>>>>>>{\sim} & G_{K}/G_{K}^0 \ar[u] & & \\
0\ar[r]& F_K \ar[r] \ar[u] & G_K \ar[r]  \ar[u] & G_K/F_K \ar[r] & 0 \\
0\ar[r] & G_{K}^0\cap F_K\ar[u] \ar[r] & G_{K}^0 \ar[r] \ar[u] & G_{K}^0/(G_{K}^0\cap F_K)\ar[r] \ar[u]_{\rotatebox{90}{$\sim$}} & 0 \\
&    0 \ar[u]     & 0. \ar[u] &  & 
}
\end{displaymath}
En particulier, $G_K/F_K$ est quotient de $G_{K}^0$. Comme ce dernier, c'est donc un $K$-groupe lisse affine connexe $\Gamma(K,\cdot)$-acyclique (en utilisant \ref{dimcohf}). D'autre part, d'après \ref{resolf}, il existe un plongement du $K$-schéma en groupes fini $F_K$ dans un $K$-groupe lisse affine connexe $\Gamma(K,\cdot)$-acyclique $L_{0,K}$. En poussant par l'inclusion $F_K\subseteq L_{0,K}$, on obtient
\begin{displaymath}
\xymatrix{
0\ar[r]& L_{0,K} \ar[r]  & L_K \ar[r]  & G_K/F_K \ar[r] & 0 \\
0\ar[r] & F_K\ar[u] \ar[r] & G_K \ar[r] \ar[u] & G_K/F_K\ar[r] \ar@{=}[u] & 0 \\
&    0 \ar[u]     & 0. \ar[u] &  & 
}
\end{displaymath}
Le $K$-groupe $L_K$ ainsi construit est  lisse affine connexe $\Gamma(K,\cdot)$-acyclique.

\medskip

\emph{Cas général.} $G_K$ est extension d'un groupe unipotent par un groupe de type multiplicatif, et cette extension devient triviale après une extension finie $K'/K$ (\cite{SGA 3} XVII 7.2.1). En particulier, d'après les cas multiplicatif et unipotent, il existe un plongement de $G_{K'}$ dans un $K'$-groupe lisse affine connexe $\Gamma(K',\cdot)$-acyclique $L_{K'}$. Il en résulte un plongement de $G_K$ dans la restriction de Weil de $L_{K'}$ à $K$. Celle-ci est un $K$-groupe  lisse affine connexe $\Gamma(K,\cdot)$-acyclique. 
\end{proof}

\begin{Rem}
Le $\delta$-foncteur $(H^i(K,\cdot))_{i\geq 0}$, de la catégorie abélienne des $K$-schémas en groupes de type fini affines dans celle des groupes ordinaires, est donc universel : cela résulte de \cite{G} II 2.2.1, puisque les $H^i(K,\cdot)$, $i>0$, sont effaçables par \ref{resola}.
\end{Rem}

\section{Projection des groupes algébriques lisses dans les groupes admissibles}

\begin{Pt}\label{foncNerIndGrpf}
On a défini en \ref{foncNer} un foncteur de Néron
$$
\xymatrix{
(\textrm{$K$-schémas en groupes lisses de type fini}) \ar[d]^{\N} \\
\Ind(\textrm{$R$-schémas en groupes lisses séparés}).
}
$$
Par ailleurs, le foncteur de Greenberg $\Gr$ \ref{green} se prolonge de façon unique en un foncteur
$$
\xymatrix{
\Ind(\textrm{$R$-schémas en groupes localement de type fini}) \ar[d]^{\Ind(\Gr)} \\
\Ind(\Pro(\textrm{$k$-schémas en groupes localement de type fini)}).
}
$$
Enfin, le foncteur de perfectisation
$$
\xymatrix{
\textrm{($k$-schémas en groupes localement de type fini)} \ar[d]^{(\cdot)^{\pf}} \\
\textrm{($k$-schémas en groupes parfaits localement algébriques)}
}
$$
se prolonge de manière unique en un foncteur
$$
\xymatrix{
\Ind(\Pro\textrm{($k$-schémas en groupes localement de type fini)}) \ar[d]^{\Ind(\Pro((\cdot)^{\pf}))} \\
\Ind(\Pro\textrm{($k$-schémas en groupes parfaits localement algébriques)}).
}
$$
En sous-entendant les plongements de catégories évidents, on peut donc considérer le foncteur composé
$$
\xymatrix{
(\textrm{$K$-schémas en groupes lisses de type fini}) \ar[d]^{\Ind(\Pro((\cdot)^{\pf}))\circ\Ind(\Pro(\Gr))\circ\N} \\
\Ind(\Pro\textrm{($k$-schémas en groupes parfaits localement algébriques)}),
}
$$
que l'on notera $\pi_*$.
\end{Pt}

\begin{Pt}\label{pi*kk'}
D'après \ref{Neroncb} et \ref{greencb}, le foncteur $\pi_*$ commute à toute extension de corps parfaite $k'/k$.
\end{Pt}

\begin{Prop}\label{proppi*}
\begin{enumerate}

\medskip

\item Soit $L$ un $R$-schéma en groupes lisse. Alors 
$$
\big(\Gr_n(L)^{\pf}(k)\big)_{n\in\bbN}=\big(\Gr_n(L)(k)\big)_{n\in\bbN}\in\Pro(\textrm{\emph{Groupes}})
$$ 
est à transitions surjectives. Pour tout $m\geq n\in\bbN$, le noyau de la transition 
$$\Gr_m(L)^{\pf}\ra \Gr_n(L)^{\pf}$$
est un $k$-schéma en groupes parfait connexe unipotent, et l'application 
$$\pi_0(\Gr_m(L))\ra \pi_0(\Gr_n(L))$$
qui s'en déduit est bijective. En particulier, 
\begin{itemize}
\medskip

\item le groupe proalgébrique associé à
$$
(\Gr(L)^{\pf})^0:=\big((\Gr_{n\in\bbN}(L)^0)^{\pf}\big)_{n\in\bbN}\in\Pro(\textrm{\emph{Groupes parfaits algébriques}})\cong\cP
$$
via l'équivalence \ref{Pab} est strictement connexe, i.e. appartient à $\cP^{00}$, et possède un ensemble de définition dont les éléments sont unipotents connexes ;

\medskip

\item le pro-groupe 
$$
\big(\pi_0(\Gr_n(L))\big)_{n\in\bbN}\in\Pro(\textrm{\emph{Groupes}})
$$
est isomorphe au groupe $\pi_0(L_k)$.
\end{itemize}

\medskip

\item Soit $L_K$ un $K$-schéma en groupes lisse de type fini. Soit $N(L_K)=(L_i)_{i\in I}$ le modèle de Néron \ref{foncNer} de $L_K$. Alors
$$
(\Gr(L_i)^{\pf}(k))_{i\in I}=(\Gr(L_i)(k))_{i\in I}=(L_i(R))_{i\in I}\in\Ind(\textrm{\emph{Groupes}})
$$
est à transitions injectives. Pour tout $i\geq j\in I$, le conoyau de la transition
$$
\Gr(L_i)^{\pf}\ra\Gr(L_j)^{\pf}
$$
dans la catégorie des faisceaux abéliens sur $(\parf/k)_{\et}$ est représentable par un $k$-schéma en groupes parfait localement algébrique. En particulier, 
\begin{itemize}

\medskip

\item le ind-groupe proalgébrique
$$
(\pi_*L_K)^0:=\big((\Gr(L_i)^{\pf})^0\big)_{i\in I}\in\Ind(\cP)
$$
est dans l'image de $\cI(\cP^{00})^{\ad}$ via le plongement \ref{IPexacte}, i.e. définit un groupe indproalgébrique admissible strictement connexe. 

\medskip

\item le ind-groupe 
$$
\big(\pi_0((L_i)_k)\big)_{i\in I}\in\Ind(\textrm{\emph{Groupes}})
$$
est isomorphe à un groupe de type fini.
\end{itemize}

\medskip

De plus, on a $(\pi_*L_K)(k)=L_K(K)$.
\end{enumerate}
\end{Prop}

\begin{proof}
Le point \emph{1.} résulte de \ref{troncationlisse} et \ref{pfexact}. La première assertion du point \emph{2.} résulte du fait que les $L_i$ sont séparés et les transitions $L_i\ra L_j$ pour $i\leq j\in I$ birationnelles. Comme elles sont également quasi-compactes, la seconde assertion du point \emph{2.} résulte de \ref{imrelisse} et \ref{pfexact}. En particulier (cf. \ref{troncationlisse} \emph{3.}), le conoyau de la transition 
$(\Gr(L_i)^{\pf})^0\ra(\Gr(L_j)^{\pf})^0$ est représentable par un $k$-schéma en groupes parfait algébrique. L'assertion relative au groupe des composantes connexes de la fibre spéciale du modèle de Néron a été vue en \ref{gctf}. Enfin, on a :
$$
(\pi_*L_K)(k)=\li_{i\in I}\lp_{n\in\bbN} \Gr_n(L)^{\pf}(k)=\li_{i\in I}L_i(R)=L_K(K).
$$
\end{proof}

\section{Lemmes d'image réciproque}

\begin{Lem}\label{lemimreproj}
Soit $f_{K}:X_K\ra Y_K$ un morphisme projectif de $K$-schémas de type fini. Si $F\subset Y_K(K)$ est un sous-ensemble borné dans $Y_K$ au sens de \cite{BLR} 1.1/2, alors $f_{K}(K)^{-1}(F)\subset X_K(K)$ est un sous-ensemble borné dans $X_K$.
\end{Lem}

\begin{proof}
Par définition, il existe un recouvrement fini $(V_j)_{j\in J}$ de $Y_K$ par des ouverts affines et une décomposition $F=\cup_{j\in J}F_j$ en sous-ensembles $F_j\subset V_j(K)$ telle que pour tout $j\in J$, $F_j$ soit borné dans $V_j$. Alors $(f_{K}^{-1}(V_j))_{j\in J}$ est un recouvrement fini de $X_K$ et 
$$
f_{K}(K)^{-1}(F)=\bigcup_{j\in J}f_{K}(K)^{-1}(F_j)
$$
est une décomposition telle que $f_{K}(K)^{-1}(F_j)\subset (f_{K}^{-1}(V_j))(K)$ pour tout $j\in J$. Pour montrer que $f_{K}(K)^{-1}(F)$ est borné dans $X_K$, on peut donc supposer $Y_K$ affine. Il existe alors un entier $n$ et une immersion fermée de $Y_K$-schémas 
$$
\xymatrix{
X_K \ar@{^{(}->}[r] \ar[dr]_{f_K} & \bbP_{K}^n\times_KY_K \ar[d] \\
& Y_K.
}
$$
Comme $\bbP_{K}^n(K)\times E\subset \bbP_{K}^n(K)\times Y_K(K)$ est borné dans $\bbP_{K}^n\times_KY_K$, on en déduit que $f_{K}(K)^{-1}(F)\subset X_K(K)$ est borné dans $X_K$.
\end{proof}

\begin{Cor}  \label{imreisog}
Soit $f_K:L_{K,0}\ra L_{K,1}$ une isogénie de $K$-schémas en groupes lisses connexes. Soient 
$$
L_0=(L_{0,i})_{i\in\bbN}:=\N(L_{0,K})\quad\textrm{et}\quad L_1=(L_{1,j})_{j\in\bbN}:=\N(L_{1,K})
$$
les modèles de Néron \ref{foncNer} de $L_{K,0}$ et $L_{K,1}$. Alors, pour tout $j\in\bbN$, il existe $i\in\bbN$ tel que
$$ 
f_{K}(K)^{-1}(L_{1,j}(R))\subset L_{0,i}(R).
$$
\end{Cor}

\begin{proof}
Supposons dans un premier temps que $L_{K,1}$ ne possède pas de sous-groupe isomorphe à $\bbG_{m,K}$ ; alors $L_{1,j}$ est de type fini sur $R$ pour tout $j\in\bbN$. Si donc $j\in\bbN$, l'ensemble $L_{1,j}(R)\subset L_{K,1}(K)$ est borné dans $L_{K,1}$. Comme $f_K:L_{K,0}\ra L_{K,1}$ est fini, donc projectif, l'image réciproque de $L_{1,j}(R)$ par $f_K(K):L_{K,0}(K)\ra L_{K,1}(K)$ est donc bornée de $L_{K,0}$ d'après \ref{lemimreproj}. Par conséquent, il existe un $R$-schéma de type fini $X$, de fibre générique isomorphe à $L_{K,0}$, tel que l'image de $X(R)\ra L_{K,0}(K)$ contienne $f_{K}(K)^{-1}(L_{1,j}(R))$. Quitte à lissifier le $R$-schéma $X$, on peut de plus supposer que $X$ est lisse sur $R$. Par définition du modèle de Néron $L_0$, il existe alors $i\in\bbN$ et un morphisme $X\ra L_{0,i}$ prolongeant l'isomorphismme $X_K\xrightarrow{\sim} L_K$. En particulier, 
$$ f_{K}(K)^{-1}(L_{1,j}(R))\subset L_{0,i}(R),$$
ce qui achève la démonstration dans le cas où $L_{K,1}$ ne contient pas $\bbG_{m,K}$.

Dans le cas général, on reprend les notations de \ref{structuremn}. Le morphisme $f_K:L_{K,0}\ra L_{K,1}$ induit un diagramme  commutatif exact
$$
\xymatrix{
0\ar[r] & T_{0,K} \ar[r] \ar[d] & L_{0,K} \ar[r] \ar[d]^{f_K} & H_{0,K} \ar[r] \ar[d] & 0 \\
0\ar[r] & T_{1,K} \ar[r] & L_{1,K} \ar[r] & H_{1,K} \ar[r] & 0.
}
$$ 
Comme $f_K$ est surjectif par hypothèse, $T_{0,K}\ra T_{1,K}$ est surjectif d'après \ref{surjtoresdepmax}. Le noyau de $H_{0,K}\ra H_{1,K}$ est donc quotient de celui de $f_K$, qui est fini par hypothèse ; par conséquent $H_{0,K}\ra H_{1,K}$ est fini. Soit $j\in\bbN$. D'après le cas traité précédemment, il existe $i\in\bbN$ tel que l'image réciproque de $H_{1,j}(R)\subset H_{K,1}(K)$ par $H_{K,0}(K)\ra H_{K,1}(K)$ soit contenue dans $H_{0,i}(R)$. Il en résulte que l'image réciproque de 
$$
L_{1,j}(R)=H_{1,j}(R)\times_{B_1(R)}\overline{L_1}(R)\subset H_{K,1}(K)\times_{B_{K,1}(K)}\overline{L_{K,1}}(K)=L_{K,1}(K)
$$
par $f_K(K):L_{K,0}(K)\ra L_{K,1}(K)$ est contenue dans 
$$
L_{0,i}(R)=H_{0,i}(R)\times_{B_0(R)}\overline{L_0}(R)\subset H_{K,0}(K)\times_{B_{K,0}(K)}\overline{L_{K,0}}(K)=L_{K,0}(K).
$$
\end{proof}

\section{Structure admissible sur le $H^1$}

\begin{Th} \label{cH1ad}
Soit $G_K$ un $K$-schéma en groupes fini. Alors il existe un unique objet 
$$\cH^1(G_K)\in\cI(\cP^{\pizf})$$ de fibre générique (\ref{deffibgenIP})
\begin{eqnarray*}
\cT^{\circ} & \ra & (\textrm{\emph{Groupes}}) \\
\kappa & \ra & H^1(K_{\kappa},G_{K_{\kappa}}) \ ;
\end{eqnarray*}
sa composante neutre est strictement connexe, son groupe des composantes connexes est un groupe fini, et c'est un objet admissible de $\cI\cP$. 

Si $f_K:L_{K,0}\ra L_{K,1}$ est une isogénie de $K$-schémas en groupes lisses connexes de noyau $G_K$, telle que
 $$H^1(K_{\kappa},L_{K_{\kappa},0})=0$$
pour toute extension algébriquement close $\kappa/k$, alors l'image de $\cH^1(G_K)$ dans 
$$\Ind\big(\Pro\textrm{\emph{\big(Faisceaux abéliens sur $(\parf/k)_{\et}$\big)}}\big)$$
réalise le conoyau de
$$
\pi_*f_K:\pi_*L_{K,0}\ra \pi_*L_{K,1}.
$$
\end{Th}

\begin{proof}
L'unicité résulte de \ref{ratIP0}. Pour l'existence, choisissons une isogénie de $K$-schémas en groupes lisses connexes $f_{K}:L_{K,0}\ra L_{K,1}$ de noyau $G_K$, telle que $H^1(K_{\kappa},L_{K_{\kappa},0})=0$ pour toute extension de corps algébriquement close $\kappa/k$ (par exemple, la résolution lisse standard de $G_K$, cf. \ref{resolf}). Soient 
$$
L_0=(L_{0,i})_{i\in\bbN}:=\N(L_{0,K})\quad\textrm{et}\quad L_1=(L_{1,j})_{j\in\bbN}:=\N(L_{1,K})
$$
les modèles de Néron \ref{foncNer} de $L_{K,0}$ et $L_{K,1}$. Représentons $\N(f_K)$ par un système compatible de morphismes de $R$-schémas en groupes
$$
\xymatrix{
f_i:L_{0,i} \ar[r] & L_{1,j_i}, & i\in\bbN ;
}
$$
en construisant la suite $(j_i)_{i\in\bbN}$ par récurrence, on peut la choisir croissante. D'après \ref{Cokerrep} et \ref{pfexact}, pour tous $i\in\bbN$ et $n\geq 1$, $\Coker(\Gr_n(f_i)^{\pf})$ est représentable par un $k$-schéma en groupes parfait algébrique, définissant donc un groupe quasi-algébrique $H_{n,i}$. Pour $i$ fixé, les transitions $H_{m,i}\ra H_{n,i}$, $m \geq n\geq 1$, sont surjectives à noyau connexe, de sorte que $H_i:=\lp_{n\geq 1}H_{n,i}$
est un groupe proalgébrique à composante neutre strictement connexe, et à groupes des composantes connexes fini égal à 
$$
\pi_0(H_{1,i}) = \Coker\big(\pi_0((f_i)_k):\pi_0((L_{0,i})_k)\ra\pi_0((L_{1,j_i})_k)\big), 
$$
qui est indépendant de $i$ d'après \ref{gctf}. De plus, pour tous $i\leq i'\in\bbN$, on a un diagramme commutatif de pro-objets de la catégorie des $k$-schémas en groupes parfaits localement algébriques 
$$
\xymatrix{
\Gr(L_{1,j_i})^{\pf} \ar[rr]^{\Gr(\tau_{j_i,j_{i'}})^{\pf}} \ar[d] && \Gr(L_{1,j_{i'}})^{\pf} \ar[rr] \ar[d] && \Coker(\Gr(\tau_{j_i,j_{i'}})^{\pf}) \ar[d] \\
\Coker(\Gr(f_i)^{\pf}) \ar[rr]^{\overline{\Gr(\tau_{j_i,j_{i'}})^{\pf}}} && \Coker(\Gr(f_{i'})^{\pf})  \ar[rr] && \Coker(\overline{\Gr(\tau_{j_i,j_{i'}})^{\pf}}). 
}
$$
En appliquant \ref{imrelisse} à l'isomorphisme $\tau_{j_i,j_{i'}}\otimes K:L_{1,j_i}\otimes K\ra L_{1,j_{i'}}\otimes K$ avec les modèles de Néron $\N(L_{1,j_i}\otimes K)=(L_{1,j})_{j\geq j_i}$ et $\N(L_{1,j_{i'}}\otimes K)=(L_{1,j'})_{j'\geq j_{i'}}$, on obtient que $\Coker(\Gr(\tau_{j_i,j_{i'}})^{\pf})$ est isomorphe à un groupe parfait algébrique, définissant en particulier un groupe quasi-algébrique. Comme il résulte du diagramme ci-dessus que le groupe proalgébrique $\Coker(H_i\ra H_{i'})$ est quotient de ce groupe quasi-algébrique, $\Coker(H_i\ra H_{i'})$ est lui aussi quasi-algébrique. 

Remarquons maintenant que le système projectif $(L_{0,i}(R/\fm^n))_{n\geq 1}$ étant à transitions surjectives, on a $H_i(k)=L_{1,j_i}(R)/\Ima f_i(R)$, et que d'après \ref{imreisog}, il existe $i'\geq i$ tel que 
$$
f_K(K)^{-1}(L_{1,j_i}(R))\subset L_{0,i'}(R),
$$ 
de sorte que 
$$
\Ker \big(L_{1,j_i}(R)/\Ima f_i(R)\ra\li_{i\in\bbN} L_{1,j_i}(R)/\Ima f_i(R)\big) 
$$
$$
= \\ \Ker \big(L_{1,j_i}(R)/\Ima f_i(R)\ra L_{1,j_{i'}}(R)/\Ima f_{i'}(R)\big).
$$
Autrement dit, le ind-objet $H:=(H_i)_{i\in\bbN}$ de $\cP$ vérifie la condition $(*)$ de \ref{condisoinjIP}. D'après \ref{IPexacte}, il existe donc $\widetilde{H}\in\cI\cP$ dont l'image dans $\Ind(\cP)$ est isomorphe à $H$. Par construction, le groupe indproalgébrique $\widetilde{H}$ possède un ensemble de définition formé de quotients des $H_i$, $i\in\bbN$. Il est donc à composante neutre strictement connexe, son groupe des composantes connexes est fini, et il vérifie la condition $(\ad)$ \ref{defIPf}. De plus, pour toute extension de corps algébriquement close $\kappa/k$, on a, avec \ref{pi*kk'} :
$$
\widetilde{H}(\kappa)=H(\kappa)=\li_{i\in\bbN}H_i(\kappa)=\li_{i\in\bbN}\Coker(f_i(R_{\kappa}))=\Coker(f_K(K_{\kappa}))=H^1(K_{\kappa},G_{K_{\kappa}})
$$
puisque $H^1(K_{\kappa},L_{0,K_{\kappa}})=0$ par hypothèse.
\end{proof}

\begin{Def}\label{DefcH1}
Soit $G_K$ un $K$-schéma en groupes fini. La \emph{structure de groupe admissible sur $H^1(K,G_K)$} est le groupe 
$\fQf$-admissible sur $k$ défini par \ref{cH1ad} et \ref{fQfIPfPIf}.
\end{Def}

\begin{Not}\label{NotcH1}
On note encore $\cH^1(G_K)$ la structure de groupe admissible sur $H^1(K,G_K)$. Le groupe ordinaire sous-jacent à 
$\cH^1(G_K)$ est $H^1(K,G_K)$, et plus généralement $\cH^1(G_K)(\kappa)=H^1(K_{\kappa},G_{K_{\kappa}})$ pour toute extension de corps algébriquement close $\kappa/k$.
\end{Not}

\begin{Pt}\label{cH1cb}
D'après \ref{cH1ad} et \ref{RemratIP0} (ou \ref{pi*kk'}), pour toute extension de corps parfaite $k'/k$, on a un isomorphisme canonique
$$
\xymatrix{
\cH^1(G_K)_{k'}\ar[r]^{\sim} & \cH^1(G_{K_{k'}}).
}
$$
\end{Pt}

\begin{Pt}\label{cH1foncQ}
D'après \ref{cH1ad} et \ref{ratIP0}, la formation de $\cH^1(G_K)$ est fonctorielle en $G_K$. On a ainsi défini un \emph{foncteur}
$$
\xymatrix{
\cH^1:\textrm{($K$-schémas en groupes finis)} \ar[r] & \fQf,
}
$$
induisant sur les points à valeurs dans une extension algébriquement close $\kappa/k$ arbitraire le foncteur
$$
\xymatrix{
H^1(K_{\kappa},\cdot):\textrm{($K$-schémas en groupes finis)} \ar[r] & \textrm{(Groupes)}.
}
$$
\end{Pt}

\section{Compatibilité aux restrictions des scalaires étales}

\begin{Pt}\label{K'k'}
Soit $K'/K$ une extension de corps finie. Alors $K'$ est un corps discrètement valué complet de caractéristique $p$, de même corps résiduel $k$ que $K$, puisque $k$ est algébriquement clos.  Si $k'/k$ est une extension de corps, on peut donc former l'extension $K_{k'}'/K'$ comme en \ref{Kk'}. On a alors  
$$R_{k'}':=R'\widehat{\otimes}_k k'=R'\widehat{\otimes}_R R\widehat{\otimes}_kk'=R'\widehat{\otimes}_R R_{k'}=R'\otimes_R R_{k'}$$
puisque $R\ra R'$ est fini, puis donc $$K_{k'}'=K'\otimes_K K_{k'}.$$
\end{Pt}

\begin{Prop}\label{cH1Res}
Soit $K'/K$ une extension de corps finie séparable. Soit $G_{K'}$ un $K'$-schéma en groupes fini. Alors il existe un isomorphisme 
$$
\xymatrix{
\cH^1(\Pi_{K'/K}G_{K'}) \ar[r]^<<<<{\sim} & \cH^1(G_{K'}) & \in \fQf,
}
$$
fonctoriel en $G_{K'}$, tel que pour toute extension algébriquement close $\kappa/k$, le morphisme induit sur les $\kappa$-points s'identifie à l'isomorphisme canonique
$$
\xymatrix{
H^1(K_{\kappa},\Pi_{K_{\kappa}'/K_{\kappa}}G_{K_{\kappa}'}) \ar[r]^<<<<{\sim} & H^1(K_{\kappa}',G_{K_{\kappa}'}).
}
$$
\end{Prop}

\begin{proof}
La collection des isomorphismes $H^1(K_{\kappa},\Pi_{K_{\kappa}'/K_{\kappa}}G_{K_{\kappa}'}) \xrightarrow{\sim} H^1(K_{\kappa}',G_{K_{\kappa}'})$, pour $\kappa/k$ variable, définit un isomorphisme 
$$
\xymatrix{
f_{\eta}:\cH^1(\Pi_{K'/K}G_{K'})_{\eta} \ar[r]^<<<<{\sim} & \cH^1(G_{K'})_{\eta}
}
$$
entre les fibres génériques des groupes indproalgébriques $\cH^1(\Pi_{K'/K}G_{K'})$ et $\cH^1(G_{K'})$. Comme ceux-ci ont un $\pi_0$ fini, il résulte de  \ref{ratIP0} que $f_{\eta}$ s'étend en un isomorphisme 
\begin{equation}
\xymatrix{
f:\cH^1(\Pi_{K'/K}G_{K'}) \ar[r]^<<<<{\sim} & \cH^1(G_{K'}).
}
\end{equation}
Celui-ci est fonctoriel en $G_{K'}$ par construction.
\end{proof}

\section{La suite exacte longue de cohomologie}

\begin{Def}\label{DefseIP}
Une suite de $\cI\cP$
$$
\xymatrix{
(G',S') \ar[r]^f & (G,S) \ar[r]^g & (G'',S'')
}
$$
est \emph{exacte (en $(G,S)$)} si $\Ima f=\Ker g$.
\end{Def}

\begin{Pt}
Par définition, le groupe indproalgébrique $\Ima f$ est un sous-objet de $(G,S)$ et $f$ induit un épimorphisme strict $(G',S')\ra \Ima f$ ; le monomorphisme $\Ima f\ra (G,S)$ n'est pas strict en général. Le groupe indproalgébrique $\Ker g$ est un sous-objet strict de $(G,S)$, et $g$ induit un monomorphisme $(G,S)/\Ker g\ra (G'',S'')$, qui n'est pas strict en général. La suite 
$$
\xymatrix{
(G',S') \ar[r]^f & (G,S) \ar[r]^g & (G'',S'')
}
$$
est exacte si et seulement si la suite
$$
\xymatrix{
0 \ar[r] & \Ima f \ar[r] & (G,S) \ar[r] & \Ima g \ar[r] & 0
}
$$ 
est une suite exacte courte de $\cI\cP$.
\end{Pt}

\begin{Def}\label{DefsefQf}\label{compsefQfIP}
Une suite de $\fQf$
$$
\xymatrix{
(G',S') \ar[r]^f & (G,S) \ar[r]^g & (G'',S'')
}
$$
est \emph{exacte (en $(G,S)$)} si la suite de $\cI\cP^{\ad}\subset\cI\cP$
$$
\xymatrix{
\indpro(G',S') \ar[r]^f & \indpro(G,S) \ar[r]^g & \indpro(G'',S'')
}
$$
est exacte (en $(G,S)$).
\end{Def}

\begin{Pt}\label{partial}
Soit
$$
\xymatrix{
0\ar[r] & G_{K}' \ar[r] & G_K \ar[r] & G_{K}'' \ar[r] & 0
}
$$
une suite exacte de $K$-schémas en groupes finis. Le bord 
$$
\xymatrix{
\partial:H^0(K,G_{K}'')\ar[r] & H^1(K,G_{K}')
}
$$
induit un morphisme de groupes admissibles $\partial:\cH^0(G_{K}'')\ra\cH^1(G_{K}')$ puisque $\cH^0(G_{K}'')$ est fini.
\end{Pt}

\begin{Prop}\label{cHsel}
Soit
$$
\xymatrix{
0\ar[r] & G_{K}' \ar[r] & G_K \ar[r] & G_{K}'' \ar[r] & 0
}
$$
une suite exacte de $K$-schémas en groupes finis. Alors la suite de $\fQf$
\begin{equation}\label{selenonce}
0\ra \cH^0(G_{K}') \ra  \cH^0(G_K) \ra \cH^0(G_{K}'') 
\xrightarrow{\partial} \cH^1(G_{K}') \ra \cH^1(G_K) \ra \cH^1(G_{K}'') \ra 0
\end{equation}
déduite de \ref{cH0fonc}, \ref{partial} et \ref{cH1foncQ}, est exacte. Si $\kappa/k$ est une extension de corps algébriquement close, alors la suite induite sur les $\kappa$-points est la suite exacte longue de cohomologie déduite de  
$$
\xymatrix{
0\ar[r] & G_{K_{\kappa}}' \ar[r] & G_{K_{\kappa}} \ar[r] & G_{K_{\kappa}}'' \ar[r] & 0.
}
$$
\end{Prop}

\begin{Lem} \label{cH1aff}
Soit $G_K$ un $K$-schéma en groupes fini. Alors les sous-groupes de définition de $\cH^1(G_K)\in\cI\cP$ sont à quotients affines.
\end{Lem}

\begin{proof}
Pour former $\cH^1(G_K)$, on peut utiliser la résolution lisse standard $f:L_{K,0}\ra L_{K,1}$ de $G_K$, où $L_{K,0}$ et 
$L_{K,1}$ sont affines. Posons $\N(L_{K,0})=(L_{0,i})_{i\in\bbN}$ et $\N(L_{K,1})=(L_{1,j})_{j\in\bbN}$. Pour tout $j\in\bbN$, $L_{1,j}^0$ est alors un $R$-schéma en groupes de type fini plat séparé, à fibre générique affine. D'après \cite{A} 2.3.1, il est donc affine, en particulier $(L_{1,j})_{k}^0$ est affine. Si $f_{i,j}:L_{0,i}\ra L_{1,j}$ prolonge $f_K$, la composante neutre du conoyau de $(f_{i,j})_k:(L_{0,i})_k \ra (L_{1,j})_k$ est quotient de $(L_{1,j})_{k}^0$, donc affine. Comme par ailleurs ce conoyau est de type fini, il est affine. Le lemme résulte alors de \ref{Cokerrep}.
\end{proof}

\begin{proof}[Démonstration de la proposition \ref{cHsel}.]
Si $\kappa/k$ est une extension de corps algébriquement close, la suite induite par (\ref{selenonce}) sur les $\kappa$-points est par construction la suite exacte longue de cohomologie 
\begin{equation*}
0\ra H^0(K_{\kappa},G_{K_{\kappa}}') \ra  H^0(K_{\kappa},G_{K_{\kappa}}) \ra H^0(K_{\kappa},G_{K_{\kappa}}'') 
\end{equation*}
\begin{equation*}
\xrightarrow{\partial} H^1(K_{\kappa},G_{K_{\kappa}}') \ra H^1(K_{\kappa},G_{K_{\kappa}}) \ra H^1(K_{\kappa},G_{K_{\kappa}}'') \ra 0 ;
\end{equation*}
l'exactitude en $H^1(K_{\kappa},G_{K_{\kappa}}'')$ résulte de \ref{dimcohf}. Ceci étant, l'exactitude de (\ref{selenonce}) en les quatre premiers termes (non triviaux) est claire. Choisissons ensuite des immersions fermées
$$
\xymatrix{
G_K \ar[r]^{\iota'} & L_{K,0}' & \textrm{et} & G_{K}'' \ar[r]^{\iota''} & L_{K,0}'', 
}
$$
avec $L_{K,0}'$ et $L_{K,0}''$ lisses connexes et vérifiant $H^1(K_{\kappa},L_{K_{\kappa},0}')=H^1(K_{\kappa},L_{K_{\kappa},0}'')=0$ pour toute extension algébriquement close $\kappa/k$. Posant $L_{K,0}:=L_{K,0}'\times L_{K,0}''$, on a alors un diagramme commutatif exact
\begin{equation*}
\xymatrix{
             & 0  \ar[d]                  & 0  \ar[d]                    & 0 \ar[d] \\
0  \ar[r] & G_{K}'   \ar[r]^{a} \ar[d]^{\iota'a} & G_{K}  \ar[r]^b \ar[d]^{(\iota',\iota''b)}  & G_{K}'' \ar[r] \ar[d]^{\iota''} & 0 \\
0  \ar[r] & L_{K,0}' \ar[r] \ar[d]^{f'} & L_{K,0}  \ar[r] \ar[d]^f  & L_{K,0}'' \ar[r] \ar[d]^{f''} & 0  \\
0 \ar[r]  & L_{K,1}' \ar[r] \ar[d] & L_{K,1}  \ar[r] \ar[d]  & L_{K,1}'' \ar[r] \ar[d] & 0\\
            & 0                   &  0                      & 0, 
}
\end{equation*}
dont la première ligne est scindée. Les deux dernières lignes se prolongent aux modèles de Néron \ref{foncNer} ; on peut représenter ce prolongement par une suite de diagrammes commutatifs
$$
\xymatrix{
0  \ar[r] & L_{0,i}' \ar[r] \ar[d]^{f_{i}'} & L_{0,i}  \ar[r] \ar[d]^{f_i}  & L_{0,i}'' \ar[r] \ar[d]^{f_{i}''} & 0  \\
            & L_{1,i}' \ar[r]                     & L_{1,i}  \ar[r]                   & L_{1,i}'', &    & i\in \bbN,
}
$$
dont la première ligne est exacte scindée et la seconde est un complexe. Pour tous $i\in\bbN$ et $n\geq 1$, on a alors un diagramme commutatif exact (\ref{pfexact})
$$
\xymatrix{
0  \ar[r] & \Gr_n(L_{0,i}')^{\pf} \ar[r] \ar[d]^{\Gr_n(f_{i}')^{\pf}} & \Gr_n(L_{0,i})^{\pf}  \ar[r] \ar[d]^{\Gr_n(f_i)^{\pf}}  & \Gr_n(L_{0,i}'')^{\pf} \ar[r] \ar[d]^{\Gr_n(f_{i}'')^{\pf}} & 0  \\
            & \Gr_n(L_{1,i}')^{\pf} \ar[r]  \ar[d]             & \Gr_n(L_{1,i})^{\pf}  \ar[r] \ar[d]                & \Gr_n(L_{1,i}'')^{\pf} \ar[d] & \\
            & \Coker(\Gr_n(f_{i}'')^{\pf}) \ar[r]  \ar[d] & \Coker(\Gr_n(f_i)^{\pf}) \ar[r] \ar[d] & \Coker(\Gr_n(f_{i}'')^{\pf}) \ar[d] & \\
            & 0                                                  &  0                                            & 0.
}
$$
La troisième ligne est un complexe de $k$-schémas en groupes parfaits algébriques, définissant donc un complexe de groupes quasi-algébriques
\begin{equation}\label{seHin}
\xymatrix{
H_{i,n}' \ar[r] & H_{i,n} \ar[r] & H_{i,n}''.
} 
\end{equation}
Posons $K_{i,n}:=\Ker (H_{i,n} \ra H_{i,n}'')\in\cQ$, $K_i:=\lp_{n\geq 1} K_{i,n}\in\cP$, et montrons que pour tout $i\in\bbN$, $\pi_0(K_i)$ est fini.

Notant $N_i$ le noyau de $L_{1,i}\ra L_{1,i}''$, on a un diagramme commutatif canonique
$$
\xymatrix{
&L_{0,i}' \ar[d]^{f_{i}'} \ar@/_1pc/[dd]_{\varphi_{i}':=} \\
&L_{1,i}' \ar[d] \ar[dr] \\
0 \ar[r] & N_i \ar[r] & L_{1,i}  \ar[r]  & L_{1,i}'',
}
$$
puis un système projectif de diagrammes commutatifs exacts
\begin{equation}\label{diagi}
\xymatrix{
0  \ar[r] & \Gr_n(L_{0,i}') \ar[r] \ar[d]^{\Gr_n(\varphi_{i}')} & \Gr_n(L_{0,i})  \ar[r] \ar[d]^{\Gr_n(f_i)}  & \Gr_n(L_{0,i}'') \ar[r] \ar[d]^{\Gr_n(f_{i}'')} & 0  \\
0 \ar[r]  & \Gr_n(N_i) \ar[r]                     & \Gr_n(L_{1,i})  \ar[r]                  & \Gr_n(L_{1,i}''), &    & 
n \geq 1.
}
\end{equation}
Le morphisme $\Gr_n(f_i)$ étant quasi-compact d'après \ref{Cokerrep}, $\Gr_n(\varphi_{i}')$ l'est aussi, et par conséquent $\Coker(\Gr_n(\varphi_{i}'))$ est représentable par un $k$-schéma en groupes localement de type fini. Comme $\Gr_n(f_{i}'')$ est quasi-compact, le lemme du serpent montre que $\Coker(\Gr_n(\varphi_{i}'))$ est de type fini sur $k$, et définit un système compatible d'épimorphismes de $k$-schémas en groupes de type fini
$$
\Coker(\Gr_n(\varphi_{i}')) \lra \Ker\big(\Coker(\Gr_n(f_i))\ra \Coker(\Gr_n(f_{i}''))\big), \quad n\geq 1.
$$
Notant $C_{i,n}$ le groupe quasi-algébrique défini par le groupe parfait algébrique $\Coker(\Gr_n(\varphi_{i}'))^{\pf}$, on a donc un système compatible d'épimorphismes $C_{i,n}\ra K_{i,n}\in\cQ$, $n\geq 1$ (\ref{pfexact}), puis un épimorphisme $C_i:=\lp_{n\geq 1}C_{i,n}\ra K_i\in\cP$. Pour montrer que $\pi_0(K_i)$ est fini, il suffit donc de montrer que $\pi_0(C_i)$ l'est. Introduisons pour cela la lissification $\ell:\widetilde{N_i}\ra N_i$ du $R$-schéma en groupes $N_i$, et notons 
$\widetilde{\varphi_{i}'}:L_{0,i}'\ra \widetilde{N_i}$ le relèvement de $\varphi_{i}':L_{0,i}'\ra N_i$. Alors 
$\Gr_n(\widetilde{\varphi_{i}'})$ est quasi-compact puisque $\Gr_n(\varphi_{i}')$ l'est, et on a un système projectif de diagrammes commutatifs exacts de $k$-schémas en groupes localement de type fini
$$
\xymatrix{
\Gr_n(L_{0,i}')\ar[rr]^{\Gr_n(\widetilde{\varphi_{i}'})} \ar@{=}[d] && \Gr_n(\widetilde{N_i}) \ar[d]^{\Gr_n(\ell)} \ar[r] & \Coker(\Gr_n(\widetilde{\varphi_{i}'})) \ar[r] \ar[d]^{\overline{\Gr_n(\ell)}} & 0 \\
\Gr_n(L_{0,i}')\ar[rr]^{\Gr_n(\varphi_{i}')} && \Gr_n(N_i) \ar[r] & \Coker(\Gr_n(\varphi_{i}')) \ar[r] & 0, & n\geq 1.
}
$$
Le conoyau de $\overline{\Gr_n(\ell)}$ est isomorphe à celui de $\Gr_n(\ell)$, qui est représentable puisque $\Gr_n(\ell)$ est quasi-compact ; $\Coker(\Gr_n(\widetilde{\varphi_{i}'}))$ est donc extension d'un sous-$k$-schéma en groupes fermé de 
$\Coker(\Gr_n(\varphi_{i}'))$ par un quotient de $\Ker(\Gr_n(\ell))$, et est donc de type fini sur $k$. Notant 
$\widetilde{C_{i,n}}$ le groupe quasi-algébrique défini par $\Coker(\Gr_n(\widetilde{\varphi_{i}'}))^{\pf}$, le système compatible de morphismes de $k$-schémas en groupes parfaits algébriques
$$
\Coker(\Gr_n(\widetilde{\varphi_{i}'}))^{\pf}\lra \Coker(\Gr_n(\varphi_{i}'))^{\pf},\quad n\geq 1,
$$
définit un morphisme de groupes proalgébriques $\widetilde{C_i}:=\lp_{n\geq 1}\widetilde{C_{i,n}}\ra C_i$. Comme $L_{0,i}'$ est lisse sur $R$, le diagramme commutatif 
$$
\xymatrix{
\lp_n\Gr_n(L_{0,i}')^{\pf}(k)\ar[rr]^>>>>>>>>>>{\lp_n\Gr_n(\widetilde{\varphi_{i}'})^{\pf}(k)} \ar@{=}[d] && \lp_n\Gr_n(\widetilde{N_i})^{\pf}(k) \ar[d]\ar[r] & \lp_n\Coker(\Gr_n(\widetilde{\varphi_{i}'}))^{\pf}(k)\ar[r] \ar[d] & 0 \\
\lp_n\Gr_n(L_{0,i}')^{\pf}(k)\ar[rr]^>>>>>>>>>>{\lp_n\Gr_n(\varphi_{i}')^{\pf}(k)} && \lp_n\Gr_n(N_i)^{\pf}(k) \ar[r] & 
\lp_n\Coker(\Gr_n(\varphi_{i}'))^{\pf}(k) \ar[r] & 0 \\
}
$$
est exact, et la deuxième flèche verticale est bijective puisqu'elle s'identifie à $\widetilde{N_i}(R)\ra N_{i}(R)$. Il en résulte que le morphisme de groupes proalgébriques $\widetilde{C_i}\ra C_i$ est un isomorphisme. Mais 
$L_{0,i}'$ et $\widetilde{N_i}$ étant lisses sur $R$, la projection $\widetilde{C_i}\ra \widetilde{C_{i,0}}$ induit un bijection sur les groupes de composantes connexes (cf. \ref{troncationlisse}). Ainsi $\pi_0(\widetilde{C_i})$, et donc $\pi_0(C_i)$ puis $\pi_0(K_i)$, sont finis.

Maintenant, les complexes (\ref{seHin}) de $\cQ$ définissent un complexe 
$$
\xymatrix{
(H_{i}')_{i\in\bbN} \ar[r] & (H_i)_{i\in\bbN} \ar[r] & (H_{i}'')_{i\in\bbN}
}
$$
de $\Ind(\cP)$, isomorphe à celui déduit du complexe
\begin{equation} \label{selen5}
\xymatrix{
\cH^1(G_{K}') \ar[r] & \cH^1(G_{K}) \ar[r] &  \cH^1(G_{K}'')
}
\end{equation}
de $\cI\cP$. Le ind-objet $(K_i)_{i\in\bbN}$ de $\cP$ est donc isomorphe à celui déduit du noyau de $\cH^1(G_{K}) \ra \cH^1(G_{K}'')$. En particulier, ce noyau possède un ensemble de définition formé de groupes proalgébriques à $\pi_0$ fini. De plus, ses sous-groupes de définition sont proaffines d'après \ref{cH1aff}. Comme pour toute extension de corps algébriquement close $\kappa/k$, la suite déduite de (\ref{selen5}) en prenant les 
$\kappa$-points est exacte, il résulte alors de \ref{foncfibgenIPcons00} appliqué au morphisme
$$
\Ima\big(\cH^1(G_{K}') \ra \cH^1(G_{K})\big) \lra \Ker\big(\cH^1(G_{K}) \ra \cH^1(G_{K}'')\big)
$$
que la suite (\ref{selenonce}) est exacte en $\cH^1(G_{K})$.

Enfin, $\cH^1(G_{K}'')$ possède un ensemble de définition formé de groupes proalgébriques à $\pi_0$ fini d'après \ref{cH1ad}, et ses sous-groupes de définition sont proaffines d'après \ref{cH1aff}. Comme $\cH^1(G_K)(\kappa)\ra\cH^1(G_{K}'')(\kappa)$ est surjective pour toute extension de corps algébriquement close $\kappa/k$, il résulte de \ref{foncfibgenIPcons00} appliqué au morphisme 
$$
\Ima\big(\cH^1(G_{K}) \ra \cH^1(G_{K}'')\big) \lra \cH^1(G_{K}'')
$$
que la suite (\ref{selenonce}) est exacte en $\cH^1(G_{K}'')$.
\end{proof}

\section{La partie continue de la structure admissible sur le $H^1$ est unipotente}

\begin{Prop}\label{cH1type}
Soit $G_K$ un $K$-schéma en groupes fini. Alors $\cH^1(G_K)^0\in\fUf$, et en particulier $\cH^1(G_K)\in\fUf$ si $G_K$ est $p$-primaire. 
De plus :
\begin{enumerate}

\medskip

\item si $G_K$ est de type multiplicatif, alors $\cH^1(G_K)$ est proalgébrique ;

\medskip

\item si $G_K$ est étale, alors $\cH^1(G_K)$ est indalgébrique ;

\medskip

\item si $G_K$ est biradiciel, alors $\cH^1(G_K)$ n'est ni proalgébrique ni indalgébrique ;

\medskip 

\item si $G_K$ est unipotent, alors $\cH^1(G_K)$ est connexe (donc strictement connexe).

\end{enumerate}
\end{Prop}

\begin{proof}
Soit $G_K$ un $K$-schéma en groupes fini. Alors $G_K$ est produit d'un $K$-groupe fini $G_{K}^{(p')}$ d'ordre premier à $p$ par un $K$-groupe fini $p$-primaire. Après une extension finie galoisienne de $K$, le groupe $G_{K}^{(p')}$ devient diagonalisable ; comme pour $n$ premier à 
$p$ $H^1(K,\mu_{n,K})=\bbZ/n$ d'après les hypothèses sur $R$ et $k$, la suite spectrale de Hochschild-Serre entraîne que $H^1(K,G_{K}^{(p')})$ est fini. Comme le foncteur $\cH^1$ est additif, il suffit donc pour démontrer la proposition de considérer le cas où $G_K$ est $p$-primaire.

\medskip 

\emph{a) Le cas où $G_K=\mu_p$.} On peut choisir pour résolution de $\mu_p$ la suite exacte de Kummer 
$$
\xymatrix{
0\ar[r] & \mu_p \ar[r] & \bbG_{m,K} \ar[r]^p & \bbG_{m,K} \ar[r] & 0.
}
$$
Le prolongement de $p:\bbG_{m,K}\ra\bbG_{m,K}$ au modèle de Néron de $\bbG_m$ est $p:\cG_m\ra\cG_m$, dont l'image par le foncteur de Greenberg est $p:\Gr(\cG_m)\ra\Gr(\cG_m)$. Par conséquent 
$$
\cH^1(\mu_p)=\Coker\big(p:\Gr(\cG_m)^{\pf}\ra\Gr(\cG_m)^{\pf}\big)
$$
\medskip
est extension de $\bbZ/p$ par un $k$-groupe provectoriel (cf. \ref{CokerGrreptoresdep}), en particulier est proalgébrique unipotent.

\medskip

\emph{b) Le cas où $G_K$ $p$-primaire de type multiplicatif.} Soit $K'/K$ une extension finee galoisienne telle que $G_{K'}$ soit diagonalisable. On a alors $\cH^1(G_{K'})\in\cP(\cU)$ d'après \emph{a)}, \ref{cHsel} puis $\cH^1(\Pi_{K'/K}G_{K'})\in\cP(\cU)$ d'après \ref{cH1Res}. Appliquant $\cH^1$ à l'immersion fermée canonique $G_K\ra \Pi_{K'/K}G_{K'}$, on en déduit avec \ref{cHsel} que $\cH^1(G_K)$ est proalgébrique unipotent.

\medskip

\emph{c) Le cas où $G_K=\bbZ/p$.} On peut choisir pour résolution de $\bbZ/p$ la suite exacte d'Artin-Schreier 
$$
\xymatrix{
0\ar[r] & \bbZ/p \ar[r] & \bbG_{a,K} \ar[r]^{1-F} & \bbG_{a,K} \ar[r] & 0.
}
$$
Le prolongement de $1-F:\bbG_{a,K}\ra\bbG_{a,K}$ au modèle de Néron $\cG_a$ \ref{defGa} de $\bbG_{a,K}$ est ind-représenté par le système 
 $$
\xymatrix{
\cG_{a,i} \ar[rr]^{t^{ip-i}-F} && \cG_{a,ip}, & i\in\bbN,
} 
$$
dont l'image par l'indisé du foncteur de Greenberg est ind-représentée par le système de morphismes de pro-schémas en groupes sur $k$ 
$$
\xymatrix{
\Gr(\cG_{a,i}) \ar[rr]^{\Gr(t^{ip-i}-F)} && \Gr(\cG_{a,ip}), & i\in\bbN.
} 
$$
Comme $t^{ip-i}-F:\bbG_{a,R}\ra\bbG_{a,R}$ est fini fidèlement plat, il résulte de \ref{fibgenetQ} que $\Coker(\Gr(t^{ip-i}-F)^{\pf})$ est un groupe parfait algébrique unipotent connexe (de dimension $ip-i$). D'où 
$$
\cH^1(\bbZ/p)=\bigg(\Coker\big(\Gr(t^{ip-i}-F)^{\pf}:\Gr(\cG_{a,i})^{\pf}\ra\Gr(\cG_{a,ip})^{\pf}\big)\bigg)_{i\in\bbN} \in\cI^0(\cU).
$$

\medskip

\emph{d) Le cas où $G_K$ est $p$-primaire étale.} Soit $K'/K$ une extension finie galoisienne telle que $G_{K'}$ soit constant. On a alors $\cH^1(G_{K'})\in\cI^0(\cU)$ d'après \emph{c)}, \ref{cHsel} et \ref{I0exacte}, puis $\cH^1(\Pi_{K'/K}G_{K'})\in\cI^0(\cU)$ d'après \ref{cH1Res}. Comme l'épimorphisme $\Pi_{K'/K}G_{K'}\ra G_K$ Cartier dual de l'immersion fermée canonique $G_{K}^*\ra \Pi_{K'/K}G_{K'}^*$ induit un épimorphisme strict $\cH^1(\Pi_{K'/K}G_{K'}) \ra \cH^1(G_K)$ d'après \ref{cHsel}, on en déduit $\cH^1(G_K)\in\cI^0(\cU)$.

\medskip

\emph{e) Le cas où $G_K=\alpha_p$.} On peut choisir pour résolution de $\alpha_p$ le Frobenius relatif de 
$\bbG_{a,K}$ 
$$
\xymatrix{
0\ar[r] & \alpha_p \ar[r] & \bbG_{a,K} \ar[r]^F & \bbG_{a,K} \ar[r] & 0.
}
$$
Le prolongement de $F:\bbG_{a,K}\ra\bbG_{a,K}$ au modèle de Néron $\cG_a$ \ref{defGa} de $\bbG_{a,K}$ est ind-représenté par le système 
$$
\xymatrix{
\cG_{a,i} \ar[r]^{F} & \cG_{a,ip}, & i\in\bbN,
} 
$$
dont l'image par l'indisé du foncteur de Greenberg est ind-représentée par le système de morphismes de pro-schémas en groupes sur $k$
$$
\xymatrix{
\Gr(\cG_{a,i}) \ar[r]^{\Gr(F)} & \Gr(\cG_{a,ip}), & i\in\bbN.
} 
$$
D'après \ref{ExsR} \emph{2.}, le conoyau de $\Gr(F)^{\pf}:\Gr(\bbG_{a,R})^{\pf}\ra\Gr(\bbG_{a,R})^{\pf}$ est le $k$-groupe provectoriel de dimension infinie $\cH^1(\alpha_{p,R})$. D'où
$$
\cH^1(\alpha_p)=\bigg(\Coker\big(\Gr(F)^{\pf}:\Gr(\cG_{a,i})^{\pf}\ra\Gr(\cG_{a,ip})^{\pf}\big)\bigg)_{i\in\bbN} \in\cI(\cP^{00})(\cU),
$$
et n'est ni proalgébrique ni indalgébrique.

\medskip

\emph{f) Le cas où $G_K$ est biradiciel.} Alors $G_K$ admet une suite de composition dont les quotients successifs sont isomorphes à $\alpha_{K,p}$, d'après \cite{SGA 3} XVII 4.2.1. D'après \emph{e)}, \ref{cHsel} et \ref{IP0exacte}, $\cH^1(G_K)\in\cI\cP$ est donc un groupe indproalgébrique unipotent strictement connexe qui n'est ni proalgébrique ni indaglégbrique.

\medskip

\emph{g) Le cas où $G_K$ est unipotent.} Alors $\cH^1(G_K)\in\cI(\cP^{00})(\cU)$ d'après \emph{d)}, \emph{f)}, \ref{cHsel} et \ref{IP0exacte}.

\medskip

\emph{h) Le cas où $G_K$ est $p$-primaire quelconque.} Alors $\cH^1(G_K)\in\cI\cP(\cU)$ \emph{b)}, \emph{g)}, et \ref{cHsel}.
\end{proof}

\begin{Rem}
Lorsque $G_K$ est $p$-primaire de type multiplicatif, on peut aussi montrer directement, c'est-à-dire sans dévisser au cas de $\mu_p$, que $\cH^1(G_K)\in\cP(\cU)$. Choisissons en effet une isogénie de $K$-tores $f_K:L_{K,0}\ra L_{K,1}$ de noyau $G_K$, nécessairement $\Gamma(K_{\kappa},\cdot)$-acycliques pour toute extension de corps algébriquement close $\kappa/k$ (cf.\cite{Beg} 3.2.1). Reprenant dans ce cas les notations générales de la démonstration de \ref{Cokerrep}, on a alors un diagramme commutatif exact de modèles de Néron
$$
\xymatrix{
0\ar[r] & T_0 \ar[r] \ar[d] & L_0 \ar[r] \ar[d]^f & H_0 \ar[r] \ar[d] & 0 \\
0\ar[r] & T_1 \ar[r] & L_1 \ar[r] & H_1 \ar[r] & 0.
}
$$
Comme $f_K$ est une isogénie $p$-primaire de $K$-tores, $T_{0,K}\ra T_{1,K}$ et $H_{0,K}\ra H_{1,K}$ aussi, et induit donc une isogénie $p$-primaire entre les tores maximaux de $(T_0)_k$ et $(T_1)_k$ d'une part, et de $(H_0)_k$ et $(H_1)_k$ d'autre part. Par conséquent $\Coker(f_k)$ est extension d'un $k$-schéma en groupes lisse unipotent par un $k$-schéma en groupes fini unipotent, et est donc unipotent. Il résulte alors de \ref{Cokerrep} que $\cH^1(G_K)\in\cP(\cU)$. 
\end{Rem}

\begin{Exs}\label{cH1Exs}
Les cas de $\mu_p$, $\bbZ/p$ et $\alpha_p$ explicités dans la démonstration de \ref{cH1type} admettent les généralisations suivantes. Soit $n\geq 1$. 

\medskip

\emph{1. $G_K=\mu_{p^n}$.} On peut choisir pour résolution de $\mu_{p^n}$ la suite exacte de Kummer 
$$
\xymatrix{
0\ar[r] & \mu_{p^n} \ar[r] & \bbG_{m,K} \ar[r]^{p^n} & \bbG_{m,K} \ar[r] & 0.
}
$$
Le prolongement de $p^n:\bbG_{m,K}\ra\bbG_{m,K}$ au modèle de Néron de $\bbG_m$ est $p^n:\cG_m\ra\cG_m$, dont l'image par le foncteur de Greenberg est $p^n:\Gr(\cG_m)\ra\Gr(\cG_m)$. D'où 
$$
\cH^1(\mu_{p^n})=\Coker\big(p^n:\Gr(\cG_m)^{\pf}\ra\Gr(\cG_m)^{\pf}\big) \in\cP(\cU).
$$
On a donc le diagramme commutatif exact de faisceaux abéliens représentables sur $(\parf/k)_{\et}$
$$
\xymatrix{
0\ar[r] & \Gr(\bbG_{m,R})^{\pf} \ar[r] \ar[d]^{p^n} & \Gr(\cG_m)^{\pf} \ar[r]^v \ar[d]^{p^n} & \bbZ \ar[d]^{p^n} \ar[r] & 0\\
0\ar[r] & \Gr(\bbG_{m,R})^{\pf} \ar[r] \ar[d] & \Gr(\cG_m)^{\pf} \ar[r]^v \ar[d] & \bbZ \ar[r] \ar[d] & 0\\
0\ar[r] & \cH^1(\mu_{p^n,R}) \ar[r] & \cH^1(\mu_{p^n}) \ar[r]^v & \bbZ/p^n \ar[r] & 0.
}
$$
Comme $\cH^1(\mu_{p^n,R})$ est connexe et $\bbZ/p^n$ fini, la dernière ligne s'identifie à la suite exacte canonique
$$
\xymatrix{
0\ar[r] & \cH^1(\mu_{p^n})^0  \ar[r] & \cH^1(\mu_{p^n})  \ar[r] & \pi_0(\cH^1(\mu_{p^n})) \ar[r] & 0.
}
$$ 
En particulier, $\cH^1(\mu_{p^n})^0(k)=R^*/(R^*)^{p^n}=H^1(R,\mu_{p^n,R})$.

\medskip

\emph{2. $G_K=\bbZ/p^n$.} On peut choisir pour résolution de $\bbZ/p^n$ la suite exacte d'Artin-Schreier-Witt 
$$
\xymatrix{
0\ar[r] & \bbZ/p^n \ar[r] & \bbW_{n,K} \ar[r]^{1-F} & \bbW_{n,K} \ar[r] & 0.
}
$$
Le prolongement de $1-F:\bbW_{n,K}\ra\bbW_{n,K}$ au modèle de Néron $\cW_n$ \ref{defWn} de $\bbW_{n,K}$ est ind-représenté par le système 
$$
\xymatrix{
\cW_{n,i} \ar[rr]^{[t^{ip-i}]-F} && \cW_{n,ip}, & i\in\bbN,
} 
$$
dont l'image par l'indisé du foncteur de Greenberg est ind-représentée par le système de morphismes de pro-schémas en groupes sur $k$ 
$$
\xymatrix{
\Gr(\cW_{n,i}) \ar[rr]^{\Gr([t^{ip-i}]-F)} && \Gr(\cW_{n,ip}), & i\in\bbN.
} 
$$
D'où 
$$
\cH^1(\bbZ/p^n)=\bigg(\Coker\big(\Gr([t^{ip-i}]-F)^{\pf}:\Gr(\cW_{n,i})^{\pf}\ra\Gr(\cW_{n,ip})^{\pf}\big)\bigg)_{i\in\bbN} \in\cI^0(\cU).
$$

\medskip

\emph{3. $G_K=\alpha_{p^n}$.} On peut choisir pour résolution de $\alpha_{p^n}$ le $n$-ième itéré du Frobenius relatif de 
$\bbG_{a,K}$ 
$$
\xymatrix{
0\ar[r] & \alpha_{p^n} \ar[r] & \bbG_{a,K} \ar[r]^{F^n} & \bbG_{a,K} \ar[r] & 0.
}
$$
Le prolongement de $F^n:\bbG_{a,K}\ra\bbG_{a,K}$ au modèle de Néron $\cG_a$ \ref{defGa} de $\bbG_{a,K}$ est ind-représenté par le système 
$$
\xymatrix{
\cG_{a,i} \ar[r]^{F^n} & \cG_{a,ip^n}, & i\in\bbN,
} 
$$
dont l'image par l'indisé du foncteur de Greenberg est ind-représentée par le système de morphismes de pro-schémas en groupes sur $k$
$$
\xymatrix{
\Gr(\cG_{a,i}) \ar[r]^{\Gr(F^n)} & \Gr(\cG_{a,ip^n}), & i\in\bbN.
} 
$$
D'où 
$$
\cH^1(\alpha_{p^n})=\bigg(\Coker\big(\Gr(F^n)^{\pf}:\Gr(\cG_{a,i})^{\pf}\ra\Gr(\cG_{a,ip^n})^{\pf}\big)\bigg)_{i\in\bbN} \in\cI(\cP^{00})(\cU).
$$
\end{Exs}

\begin{Pt}\label{cH1foncU}
En composant le foncteur $\cH^1$ \ref{cH1foncQ} avec le foncteur composante neutre, on obtient d'après \ref{cH1type} et \ref{cH1ad} un foncteur 
$$
\xymatrix{
(\cH^1)^0:\textrm{($K$-schémas en groupes finis)} \ar[r] & \fUf^{00},
}
$$
et en composant $\cH^1$ avec le foncteur groupe des composantes connexes, un foncteur
$$
\xymatrix{
\pi_0(\cH^1):\textrm{($K$-schémas en groupes finis)} \ar[r] & \textrm{(Groupes finis)}.
}
$$
De plus, les foncteurs $\cH^1$, $(\cH^1)^0$ et $\pi_0(\cH^1)$ induisent des foncteurs
$$
\xymatrix{
\cH^1:\textrm{($K$-schémas en groupes finis $p$-primaires)} \ar[r] & \fUf,
}
$$
$$
\xymatrix{
(\cH^1)^0:\textrm{($K$-schémas en groupes finis $p$-primaires)} \ar[r] & \fUf^{00},
}
$$
$$
\xymatrix{
\pi_0(\cH^1):\textrm{($K$-schémas en groupes finis $p$-primaires)} \ar[r] & \textrm{(Groupes finis $p$-primaires)}.
}
$$
\end{Pt}

\section{Lien avec la cohomologie des $R$-schémas en groupes finis plats}

\begin{Lem}\label{cH1RcH1K}
Soit $G$ un $R$-schéma en groupes fini plat. Alors le morphisme de groupes canonique $H^1(R,G)\ra H^1(K,G_K)$ est un morphisme
$$
\xymatrix{
\cH^1(G) \ar[r] & \cH^1(G_K)\in\cI\cP^{\ad},
}
$$
et induit un isomorphisme de $\cH^1(G)$ sur une origine de $\cH^1(G_K)$.
\end{Lem}

\begin{proof}
Soient $f:L_0\ra L_1$ la résolution lisse standard de $G$. Le morphisme $f_K:L_{K,0}\ra L_{K,1}$ induit sur les fibres génériques se prolonge aux modèles de Néron $\N(L_{K,0})$ et 
$\N(L_{K,1})$ \ref{foncNer} en un morphisme $\N(f_K)$. Par ailleurs, les $R$-groupes lisses $L_0$ et $L_1$ étant de type fini sur $R$, l'identité de $L_{K,0}$ (resp. $L_{K,1}$) se prolonge en un morphisme $\alpha_0:L_0\ra \N(L_{K,0})$ (resp. $\alpha_1:L_1\ra \N(L_{K,1})$). On obtient ainsi un diagramme commutatif de ind-schémas en groupes sur $R$
\begin{equation}\label{LN(L)}
\xymatrix{
L_0 \ar[r]^f \ar[d]^{\alpha_0} & L_1 \ar[d]^{\alpha_1} \\
\N(L_{K,0}) \ar[r]^{\N(f_K)} & \N(L_{K,1}). 
}
\end{equation}
Appliquant les foncteurs $\Ind(\Gr)$ puis $\Ind(\Pro((\cdot)^{\pf}))$, puis passant aux conoyaux, on obtient un morphisme de ind-pro-schémas en groupes parfaits algébriques
$$
\xymatrix{
\Coker(\pi_*f)\ar[r] & \Coker(\pi_*f_K)
}
$$
(\ref{cH1foncQRbis} et \ref{foncNerIndGrpf}), induisant le morphisme canonique $H^1(R,G)\ra H^1(K,G_K)$ sur les $k$-points. 
Celui-ci est donc bien indrpoalgébrique pour les structures proalgébrique $\cH^1(G)$ et indproalgébrique $\cH^1(G_K)$ d'après \ref{cH1foncQRbis} et \ref{cH1ad} (ce qui résulte aussi de \ref{ratIP0}). Il est de plus injectif (puisqu'un $G$-torseur étant propre, il est trivial s'il l'est génériquement), de sorte qu'il identifie $\cH^1(G)$ à un sous-groupe de définition de $\cH^1(G_K)\in\cI\cP$. Pour voir que $\cH^1(G)$ est une origine de $\cH^1(G_K)$, posons $\N(L_{K,0})=(L_{0,i})_{i\in\bbN}$, $\N(L_{K,0})=(L_{1,j})_{j\in\bbN}$, et représentons (\ref{LN(L)}) par un diagramme commutatif
$$
\xymatrix{
L_0 \ar[r]^f \ar[d]^{\alpha_{0,i}} & L_1 \ar[d]^{\alpha_{1,j_i}} \\
L_{0,i} \ar[r]^{f_i} & L_{1,j_i}. 
}
$$
Pour tout $n\geq 1$, on a alors le diagramme commutatif exact
$$
\xymatrix{
\Gr_n(L_0) \ar[r]^{\Gr_n(f)} \ar[d]^{\Gr_n(\alpha_{0,i})} & \Gr_n(L_1) \ar[d]^{\Gr_n(\alpha_{1,j_i})} \ar[r] & \Coker(\Gr_n(f)) \ar[r] \ar[d]^{\overline{\Gr_n(\alpha_{1,j_i})}} & 0 \\
\Gr_n(L_{0,i}) \ar[r]^{\Gr_n(f_i)} \ar[d] & \Gr_n(L_{1,j_i})  \ar[r] \ar[d] & \Coker(\Gr_n(f_i)) \ar[r] & 0 \\
\Coker(\Gr_n(\alpha_{0,i})) \ar[r]^{\overline{\Gr_n(f_i)}} \ar[d] & \Coker(\Gr_n(\alpha_{1,j_i})) \ar[d] &  & \\
0 & 0, & &
}
$$
identifiant $\Coker(\overline{\Gr_n(\alpha_{1,j_i})})$ à $\Coker(\overline{\Gr_n(f_i)})$. Les morphismes $\alpha_{0,i}$ et $\alpha_{1,j_i}$ sont birationnels par construction, et $\Coker(\Gr_n(\alpha_{0,i}))$ et $\Coker(\Gr_n(\alpha_{1,j_i}))$ sont représentables puisque $L_0$ et $L_1$ sont quasi-compacts. D'après \ref{imrelisse}, les transitions des systèmes projectifs $(\Coker(\Gr_n(\alpha_{0,i})))_{n\geq 1}$ et  
$(\Coker(\Gr_n(\alpha_{1,j_i})))_{n\geq 1}$ sont donc des isomorphismes pour tout $n$ assez grand ; il en va donc de même des transitions du système projectif $(\Coker(\overline{\Gr_n(\alpha_{1,j_i})}))_{n\geq 1}$. Comme pour $i$ variable, les pro-objets $(\Coker(\Gr_n(f_i)))_{n\geq 1}$ définissent un ensemble de définition de $\cH^1(G_K)\in\cI\cP$, on obtient bien que $\cH^1(G)$ est une origine de $\cH^1(G_K)$.
\end{proof}

\begin{Pt}
Le morphisme $\cH^1(G)\ra \cH^1(G_K)$ est donc un monomorphisme strict de $\cI\cP^{\ad}$, et le quotient 
$\cH^1(G_K)/\cH^1(G)$ est un objet de $\cI$. 
\end{Pt}

\chapter{L'isomorphisme de dualité}

Soient $K$ un corps discrètement valué complet de caractéristique $p>0$, $R$ son anneau d'entiers et $k$ son corps résiduel, que l'on suppose algébriquement clos. 

\section{La flèche sur les parties discrètes}

\begin{Nots}
\medskip
\begin{itemize}
\item Si $G_K$ est un $K$-schéma en groupes fini, annulé par $p^n$, $n\in\bbN$, on note
$$G_{K}^*:=\cH om_{\bbZ}(G_K,\bbG_m)=\cH om_{\bbZ/p^n}(G_K,\mu_{p^n})$$ 
le \emph{dual de Cartier} de $G_K$.
\medskip
\item Si $F$ un groupe ordinaire fini (resp. un groupe quasi-algébrique fini), annulé par $p^n$, $n\in\bbN$, on note 
$$F^{\vee_P}:=\Hom_{\bbZ}(F,\bbQ/\bbZ)=\Hom_{\bbZ/p^n}(F,\bbZ/p^n)$$
le \emph{dual de Pontryagin} de $F$. 
\end{itemize}
\end{Nots}

\begin{Pt} 
Soit $G_K$ un $K$-schéma en groupes fini, annulé par $p^n$, $n\in\bbN$. On définit une flèche 
\begin{equation}\label{flecheK1k}
\xymatrix{
H^{1}(K,G_K) \ar[r] & H^0(K,G_{K}^*)^{\vee_P}
}
\end{equation}
par la composition 
$$
\xymatrix{
H^{1}(K,G_K)=\Ext_{\bbZ/p^n}^1(\bbZ/p^n,G_K) \ar[r]^>>>>{(\cdot)^*} & \Ext_{\bbZ/p^n}^1(G_{K}^*,\mu_{p^n}) \ar[r]^<<<<{\partial} & \Hom_{\bbZ/p^n}(H^0(K,G_{K}^*),H^1(K,\mu_{p^n})) \ar[d]^{v}\\
&&\Hom_{\bbZ/p^n}(H^0(K,G_{K}^*),\bbZ/p^n)\ar@{=}[d]\\
&&H^0(K,G_{K}^*)^{\vee_P}.
}
$$
\end{Pt}

\begin{Pt}
Le morphisme 
$$
\xymatrix{
\Ext_{\bbZ/p^n}^1(\bbZ/p^n,G_K) \ar[r]^>>>>{(\cdot)^*} & \Ext_{\bbZ/p^n}^1(G_{K}^*,\mu_{p^n})
}
$$
associe à la classe d'isomorphisme d'une extension
$$
\xymatrix{
0\ar[r] & G_K \ar[r] & E_K \ar[r] & \bbZ/p^n \ar[r] & 0, 
}
$$
où $p^nE_K=0$, la classe d'isomorphisme de l'extension obtenue par dualité de Cartier
$$
\xymatrix{
0\ar[r] & \mu_{p^n} \ar[r] & E_{K}^* \ar[r] & G_{K}^* \ar[r] & 0.
}
$$
\end{Pt}

\begin{Pt}
Le morphisme 
$$
\xymatrix{
\Ext_{\bbZ/p^n}^1(G_{K}^*,\mu_{p^n}) \ar[r]^<<<<{\partial} & \Hom_{\bbZ/p^n}(H^0(K,G_{K}^*),H^1(K,\mu_{p^n})) 
}
$$
associe à la classe d'isomorphisme d'une extension
$$
\xymatrix{
0\ar[r] & \mu_{p^n} \ar[r] & E_{K}^* \ar[r] & G_{K}^* \ar[r] & 0,
}
$$ 
où $p^nE_{K}^*=0$, le bord 
$$
\xymatrix{
H^0(K,G_{K}^*)\ar[r]^{\partial} & H^1(K,\mu_{p^n})
}
$$
de la suite exacte longue de cohomologie qui s'en déduit.
\end{Pt}

\begin{Pt}
Le morphisme
$$
\xymatrix{ 
\Hom_{\bbZ/p^n}(H^0(K,G_{K}^*),H^1(K,\mu_{p^n})) \ar[r]^<<<<v & \Hom_{\bbZ/p^n}(H^0(K,G_{K}^*),\bbZ/p^n) 
}
$$
associe à un morphisme
$$
\xymatrix{
H^0(K,G_{K}^*)\ar[r]^{\partial} & H^1(K,\mu_{p^n})
}
$$
le morphisme composé
$$
\xymatrix{
H^0(K,G_{K}^*)\ar[r]^{\partial} & H^1(K,\mu_{p^n}) \ar[r]^v & \bbZ/p^n,
}
$$
où
$$
\xymatrix{
H^1(K,\mu_{p^n})=K^*/(K^*)^{p^n} \ar[r]^<<<<<v & \bbZ/p^n
}
$$
est induit par la valuation de $K$.
\end{Pt}

\begin{Pt} 
On définit une autre flèche 
\begin{equation}\label{oppflecheK1k}
\xymatrix{
H^{1}(K,G_K) \ar[r] & H^0(K,G_{K}^*)^{\vee_P}
}
\end{equation}
par la composition
$$
\xymatrix{
H^{1}(K,G_K) \ar[r]^>>>>>{\ev_{.}} &  \Hom_{\bbZ/p^n}(H^0(K,G_{K}^*),H^1(K,\mu_{p^n})) \ar[d]^{v}\\
&\Hom_{\bbZ/p^n}(H^0(K,G_{K}^*),\bbZ/p^n)\ar@{=}[d]\\
&H^0(K,G_{K}^*)^{\vee_P}.
}
$$
\end{Pt}

\begin{Pt}
Le morphisme 
$$
\xymatrix{
H^{1}(K,G_K) \ar[r]^>>>>>{\ev_{.}} &  \Hom_{\bbZ/p^n}(H^0(K,G_{K}^*),H^1(K,\mu_{p^n})) 
}
$$
associe à un élément $a\in H^1(K,G_K)$ le morphisme
\begin{eqnarray*}
H^0(K,G_{K}^*)=\Hom_{\bbZ/p^n}(G_K,\mu_{p^n}) & \lra & H^1(K,\mu_{p^n}) \\
f & \lmapsto & \ev_a(f):=f_*(a),
\end{eqnarray*}
où $f_*:=H^1(K,f):H^1(K,G_K)\ra H^1(K,\mu_{p^n})$.  
\end{Pt}

\begin{Prop}\label{compflecheK1k}
Les flèches (\ref{flecheK1k}) et (\ref{oppflecheK1k}) sont opposées. Plus précisément,  $\partial\circ (\cdot)^*=-\ev_{.}$.
\end{Prop}

\begin{Lem}\label{calculEKamupn}
Soit $a\in K^*\setminus (K^*)^{p^n}$. L'identification 
$$K^*/(K^*)^{p^n}=H^1(K,\mu_{p^n})=\Ext_{\bbZ/p^n}^1(\bbZ/p^n,\mu_{p^n})$$
fait correspondre à la classe de $a$ la classe de l'extension
\begin{equation} \label{extEKamupn}
\xymatrix{
0\ar[r] & \mu_{p^n} \ar[r] & E_{K,a} \ar[r] & \bbZ/p^n \ar[r] & 0,
}
\end{equation}
où $E_{K,a}$ est le $K$-schéma en groupes dont l'algèbre de Hopf est la suivante :
\medskip
\begin{itemize}
\item 
la $K$-algèbre de $E_{K,a}$ est :
$$
A_{K,a}:=K[X]\Big/\prod_{0\leq i\leq p^n-1}(X^{p^n}-a^i)\ ; 
$$
\item 
posant
$$
\delta_{a^i}:=\prod_{\substack{0\leq i'\leq p^n-1,\\ i'\neq i}}\bigg(\frac{X^{p^n}-a^{i'}}{a^i-a^{i'}}\bigg)\ \in A_{K,a}
$$
pour tout $0\leq i \leq p^n-1$, la co-multiplication de $A_{K,a}$ est donnée par :
\begin{eqnarray*}
A_{K,a} & \lra & A_{K,a}\otimes_K A_{K,a} \\
X & \lmapsto & \sum_{0\leq i,j \leq p^n-1}a^{-\left\lfloor\frac{i+j}{p^n}\right\rfloor}X\delta_{a^i}\otimes X\delta_{a^j}\ ;
\end{eqnarray*}
\item
la co-unité de $A_{K,a}$ est donnée par :
\begin{eqnarray*}
A_{K,a} & \lra & K \\
X & \lmapsto & 1\ ;
\end{eqnarray*}
\item
la co-inversion de $A_{K,a}$ est donnée par :
\begin{eqnarray*}
A_{K,a} & \lra & A_{K,a} \\
X & \lmapsto & X^{p^n-1}(\delta_1+\sum_{i=1}^{p^n-1}a^{1-i}\delta_{a^i}).
\end{eqnarray*}
\end{itemize}
De plus, le dual de Cartier de $E_{K,a}$ est $E_{K,a^{-1}}$, avec pour élément diagonal 
$$
\sum_{0\leq i,j\leq p^n-1}X^j\delta_{a^i}\otimes X^i\delta_{(a^{-1})^j}\ \in A_{K,a}\otimes_K A_{K,a^{-1}}. 
$$
L'extension obtenue à partir de (\ref{extEKamupn}) par dualité de Cartier est l'extension
\begin{equation} \label{extEK*amupn}
\xymatrix{
0\ar[r] & \mu_{p^n} \ar[r] & E_{K,a^{-1}} \ar[r] & \bbZ/p^n \ar[r] & 0.
}
\end{equation}
\end{Lem}

\begin{proof}
Pour faire les calculs, on utilise la $K$-base de $A_{K,a}$ donnée par 
$$
X^j\delta_{a^i},\quad 0\leq i,j\leq p^n-1.
$$ 
Dans cette base,
\medskip
\begin{itemize}
\item  la co-multiplication de $A_{K,a}$ est donnée par :
\begin{eqnarray*}
A_{K,a} & \lra & A_{K,a}\otimes_K A_{K,a} \\
X^j\delta_{a^i} & \lmapsto &\sum_{0\leq \ell,\ell'\leq p^n-1}a^{-j\left\lfloor\frac{\ell+\ell'}{p^n}\right\rfloor}\delta_{i,\ell+\ell'\modd p^n} X^j\delta_{a^{\ell}}\otimes X^j\delta_{a^{\ell'}},
\end{eqnarray*}
où $\delta_{\alpha,\beta}$ est le symbole de Kronecker ;
\medskip
\item
la co-unité de $A_{K,a}$ est donnée par :
\begin{eqnarray*}
A_{K,a} & \lra & K \\
X^j\delta_{a^i} & \lmapsto & \delta_{i,0}\ ;
\end{eqnarray*}
\item
la co-inversion de $A_{K,a}$ est donnée par :
\begin{eqnarray*}
A_{K,a} & \lra & A_{K,a} \\
X^j\delta_{a^i} & \lmapsto & X^{p^n-j}(\delta_{i,0}\delta_1+\sum_{\ell=1}^{p^n-1}a^{j-\ell}\delta_{i,p^n-\ell}\delta_{a^{\ell}}).
\end{eqnarray*}
\end{itemize}
\end{proof}

\begin{proof}[Démonstration de la proposition \ref{compflecheK1k}]
Rappelons que l'identification 
$$
\xymatrix{
\Ext_{\bbZ/p^n}^1(\bbZ/p^n,G_K) \ar[r]^<<<<{\sim} &  H^1(K,G_K)
}
$$
fait correspondre à la classe d'isomorphisme d'une extension
$$
\xymatrix{
0\ar[r] & G_K \ar[r] & E_{K} \ar[r] & \bbZ/p^n \ar[r] & 0 
}
$$
(où $p^nE_{K}=0$) l'image de $1$ par le bord 
$$
\xymatrix{
\partial:\bbZ/p^n\ar[r] & H^1(K,G_K)
}
$$
de la suite exacte longue de cohomologie qui s'en déduit. Géométriquement, $\partial(1)$ est la classe d'isomorphisme du $G_{K}$-torseur $(E_{K})_1$, fibre en $1$ de $E_{K} \ra \bbZ/p^n$. 

Ceci étant, soit $a\in\Ext_{\bbZ/p^n}^1(\bbZ/p^n,G_K)$ la classe d'isomorphisme d'une extension
$$
\xymatrix{
0\ar[r] & G_K \ar[r] & E_{K,a} \ar[r] & \bbZ/p^n \ar[r] & 0, 
}
$$
et soit $f\in H^0(K,G_{K}^*)=\Hom_{\bbZ/p^n}(G_K,\mu_{p^n})$. Considérons l'image directe par $f$ de l'extension $E_{K,a}$ :
\begin{equation}\label{diagf}
\xymatrix{
0\ar[r] & G_K \ar[r] \ar[d]^f & E_{K,a} \ar[d] \ar[r] & \bbZ/p^n \ar[r] \ar@{=}[d] & 0 \\
0\ar[r] & \mu_{p^n} \ar[r] & f_*E_{K,a} \ar[r] & \bbZ/p^n \ar[r] & 0. 
}
\end{equation}
L'image par $f_*:H^1(K,G_K)\ra H^1(K,\mu_{p^n})$  de la classe d'isomorphisme du $G_K$-torseur $(E_{K,a})_1$ est la classe d'isomorphisme du $\mu_{p^n}$-torseur $(f_*E_{K,a})_1$. Autrement dit, $\ev_a(f)=\ \textrm{classe de}\ (f_*E_{K,a})_1$. Par ailleurs, appliquant le foncteur de dualité de Cartier au diagramme (\ref{diagf}), on obtient un diagramme
commutatif exact
\begin{equation}
\xymatrix{
0\ar[r] & \mu_{p^n} \ar[r] & E_{K,a}^* \ar[r] & G_{K}^* \ar[r] & 0 \\
0\ar[r] & \mu_{p^n} \ar[r] \ar@{=}[u] & (f_*E_{K,a})^* \ar[r] \ar[u] & \bbZ/p^n \ar[r] \ar[u] & 0, 
}
\end{equation}
où $\bbZ/p^n\ra G_{K}^*$ envoie $1$ sur $f\in G_{K}^*(K)$. Passant aux suites exactes de cohomologie, on en déduit un diagramme commutatif
$$
\xymatrix{
G_{K}^*(K)\ar[r]^<<<<{\partial} & H^1(K,\mu_{p^n})  \\
\bbZ/{p^n} \ar[r]^<<<<<<{\partial} \ar[u] & H^1(K,\mu_{p^n}). \ar@{=}[u]
}
$$
Autrement dit, $(\partial\circ (\cdot)^*)_a(f)$ est la classe d'isomorphisme du $\mu_{p^n}$-torseur $((f_*E_{K,a})^*)_1$. D'après \ref{calculEKamupn}, il s'agit de l'opposé de la classe d'isomorphisme du $\mu_{p^n}$-torseur $(f_*E_{K,a})_1$. On a donc finalement $(\partial\circ (\cdot)^*)_a(f)=-\ev_a(f)$.
\end{proof}

\begin{Cor} \label{annulneutre}
La flèche (\ref{flecheK1k}) s'annule sur le sous-groupe 
$$\cH^1(G_K)^0(k)\subseteq \cH^1(G_K)(k)=H^1(K,G_K).$$
\end{Cor}

\begin{proof}
Soit $f\in H^0(K,G_{K}^*)=\Hom_{\bbZ/p^n}(G_K,\mu_{p^n})$. Appliquant le sous-foncteur 
$(\cH^1)^0$ de $\cH^1$ à $f$ (\ref{cH1foncU}), on obtient un diagramme commutatif
$$
\xymatrix{
\cH^1(G_K)^0\ar[r]^{\cH^1(f)^0} \ar[d] & \cH^1(\mu_{p^n})^0 \ar[d] \\
\cH^1(G_K)\ar[r]^{\cH^1(f)} & \cH^1(\mu_{p^n}).
}
$$
En prenant les $k$-points, on obtient le diagramme commutatif 
$$
\xymatrix{
\cH^1(G_K)^0(k)\ar[r] \ar[d] & \cH^1(\mu_{p^n})^0(k) \ar[d] \\
H^1(K,G_K)\ar[r]^{f_*} & H^1(K,\mu_{p^n}).
}
$$
Par conséquent, si $a\in\cH^1(G_K)^0(k)$, alors $\ev_a(f):=f_*(a)\in \cH^1(\mu_{p^n})^0(k)$. Comme 
$$\cH^1(\mu_{p^n})^0(k)=R^*/(R^*)^{p^n}$$ 
(\ref{cH1Exs} \emph{1.}), on a alors $(v\circ\ev_a)(f)=0$. Autrement dit, la flèche (\ref{oppflecheK1k}) s'annule sur $\cH^1(G_K)^0(k)$. D'où le corollaire.
\end{proof}

\begin{Pt}
La flèche  (\ref{flecheK1k}) induit donc une flèche 
\begin{equation}\label{flecheK1kq}
\xymatrix{
\cH^1(G_K)(k)/\cH^1(G_K)^0(k) \ar[r] & H^0(K,G_{K}^*)^{\vee_P}.
}
\end{equation}
La source d'identifie à $\pi_0(\cH^1(G_K))(k)$ d'après \ref{foncfibgenIPexact}, et $\pi_0(\cH^1(G_K))$ est un groupe quasi-algébrique fini d'après \ref{cH1ad}. La donnée de la flèche ci-dessus est donc équivalente à la donnée d'une flèche dans $\cQ$
\begin{equation} \label{fdisc}
\xymatrix{
f_{\disc}(G_K):\pi_0(\cH^1(G_K)) \ar[r] & \cH^0(G_{K}^*)^{\vee_P}.
}
\end{equation}
La flèche $f_{\disc}(G_K)$ ainsi construite est fonctorielle en $G_K$.
\end{Pt}

\section{La flèche sur les parties continues} 

\begin{Pt}\label{Ptprelimfl2}
Rappelons que si $G_K$ est un $K$-schéma en groupes fini, alors $\cH^1(G_K)^0\in\fUf^{00}$ (\ref{cH1ad}). En particulier, on peut considérer son dual de Serre $(\cH^1(G_K)^0)^{\vee_S}\in\fUf^{00}$ (\ref{dualSfUfzz}). Si $n\in\bbN$ est tel que $p^nG_K=0$, alors pour toute extension de corps algébriquement close $\kappa/k$, on a 
$$
(\cH^1(G_K)^0)^{\vee_S}(\kappa)=\Ext_{\fUnf_{\kappa}}^1(\cH^1(G_{K_{\kappa}})^0,\bbZ/p^n)
$$
d'après \ref{dualSexplicit4} et \ref{cH1cb}.
\end{Pt}

\begin{Pt}
Soit $G_K$ un $K$-schéma en groupes fini, annulé par $p^n$, $n\in\bbN$. On définit une flèche
\begin{equation} \label{flecheK2k}
\xymatrix{
\cH^1(G_{K})^0(k) \ar[r] & (\cH^1(G_{K}^*)^0)^{\vee_S}(k)
}
\end{equation}
par la composition :
$$
\xymatrix{
\cH^1(G_{K})^0(k) \ar[d]^{\rotatebox{90}{$\sim$}} & \\
\Ima\bigg(\cH^1(G_{K})^0(k) \ar@{^{(}->}[r] & H^1(K,G_K)=\Ext_{\bbZ/p^n}^1(\bbZ/p^n,G_K) \ar[r]^>>>>>>>>{(\cdot)^*} & \Ext_{\bbZ/p^n}^1(G_{K}^*,\mu_{p^n})\bigg)  \ar[d]^>>>>>{\iota^*\circ \ov_*\circ\cH^1} \\ 
& & \Ext_{\fUnf}^1\big(\cH^1(G_{K}^*)^0,\bbZ/p^n\big) \ar@{=}[d] \\
& & (\cH^1(G_{K}^*)^0)^{\vee_S}(k).
}
$$
\end{Pt}

\begin{Pt}
Notant $\Ext_{\bbZ/p^n}^1(G_{K}^*,\mu_{p^n})^0$ l'image de $\cH^1(G_{K})^0(k)$ dans $\Ext_{\bbZ/p^n}^1(G_{K}^*,\mu_{p^n})$, le morphisme  
$$
\xymatrix{
\Ext_{\bbZ/p^n}^1(G_{K}^*,\mu_{p^n})^0 \ar[rr]^>>>>>>>>>>{\iota^*\circ \ov_*\circ\cH^1} & & \Ext_{\fUnf}^1\big(\cH^1(G_{K}^*)^0,\bbZ/p^n)\big)
}
$$
associe à la classe d'isomorphisme d'une extension
\begin{equation}\label{ext0}
\xymatrix{
0\ar[r] & \mu_{p^n} \ar[r] & E_{K}^* \ar[r] & G_{K}^* \ar[r] & 0,
}
\end{equation} 
où $p^nE_{K}^*=0$, la troisième ligne du diagramme commutatif exact de $\fUnf$
\begin{equation}\label{diagdeffleche}
\xymatrix{
& \cH^1(\mu_{p^n}) \ar[r] \ar[d]_{v} & \cH^1(E_{K}^*) \ar[r] \ar[d] & \cH^1(G_{K}^*) \ar[r] \ar@{=}[d] & 0 \\
0 \ar[r] & \bbZ/p^n \ar[r]  & \ov_*\cH^1(E_{K}^*) \ar[r] & \cH^1(G_{K}^*) \ar[r] & 0 \\
0 \ar[r] & \bbZ/p^n \ar[r] \ar@{=}[u] & \iota^*\ov_*\cH^1(E_{K}^*) \ar[r] \ar[u] & \cH^1(G_{K}^*)^0 \ar[r] \ar[u]_{\iota} & 0,
}
\end{equation}
où $\iota:\cH^1(G_{K}^*)^0\ra\cH^1(G_{K}^*)$ est le morphisme canonique. Montrons que ce diagramme est bien défini et exact. D'après \ref{cHsel}, on a une suite exacte de $\fUnf$
$$
\xymatrix{
\cH^0(G_{K}^*)  \ar[r]^<<<<<{\partial} & \cH^1(\mu_{p^n}) \ar[r] &  \cH^1(E_{K}^*) \ar[r] &\cH^1(G_{K}^*) \ar[r] & 0.
}
$$
D'après \ref{annulneutre}, le fait que l'extension (\ref{ext0}) provienne d'un élément de $\cH^1(G_{K})^0(k)\subset \cH^1(G_{K})(k)$ entraîne que l'image de $\partial$ est contenue dans le noyau de $v:\cH^1(\mu_{p^n})\ra\bbZ/p^n$. En particulier, $v$ se factorise par $\cH^1(\mu_{p^n})/\Ima(\partial)$, et on a donc un diagramme commutatif exact
$$
\xymatrix{
0\ar[r] & \cH^1(\mu_{p^n})/\Ima(\partial) \ar[r] \ar[d]^{\ov}& \cH^1(E_{K}^*) \ar[r] \ar[d] &\cH^1(G_{K}^*) \ar[r] \ar@{=}[d] & 0 \\
0 \ar[r] & \bbZ/p^n \ar[r]  & \ov_*\cH^1(E_{K}^*) \ar[r] & \cH^1(G_{K}^*) \ar[r] & 0,
}
$$
définissant la deuxième ligne exacte du diagramme (\ref{diagdeffleche}). La troisième ligne est obtenue par image réciproque par $\iota$.
\end{Pt}

\begin{Pt}\label{fcontconstr}
Si $\kappa/k$ est une extension de corps algébriquement close, alors la flèche (\ref{flecheK2k}) pour $G_{K_{\kappa}}$ fournit une flèche
\begin{equation} \label{flecheK2kappa}
\xymatrix{
\cH^1(G_K)^0(\kappa) \ar[r] & (\cH^1(G_{K}^*)^0)^{\vee_S}(\kappa),
}
\end{equation}
d'après \ref{cH1cb} et \ref{Ptprelimfl2}. La collection de flèches obtenue en faisant varier $\kappa/k$ définit une transformation naturelle
\begin{equation}\label{flecheK2eta}
\xymatrix{
(\cH^1(G_K)^0)_{\eta} \ar[r] & ((\cH^1(G_{K}^*)^0)^{\vee_S})_{\eta}
}
\end{equation}
entre les fibres génériques des groupes indproalgébriques $\cH^1(G_K)^0$ et $(\cH^1(G_{K}^*)^0)^{\vee_S}$ (\ref{deffibgenIP}). Alors, d'après \ref{ratIP0}, il existe un unique morphisme dans $\fUnf$
\begin{equation}\label{fcont}
\xymatrix{
f_{\cont}(G_K):\cH^1(G_K)^0 \ar[r] & (\cH^1(G_{K}^*)^0)^{\vee_S}
}
\end{equation}
induisant (\ref{flecheK2eta}) (et un tel morphisme est uniquement déterminé par ses $k$-points). La flèche $f_{\cont}(G_K)$ ainsi construite est fonctorielle en $G_K$.
\end{Pt}

\section{Lemmes d'effacement dans la catégorie des $K$-schémas en groupes finis}

Les lemmes de coeffacement \ref{coefface} et d'effacement \ref{efface} qui suivent sont les analogues en égale caractéristique $p>0$ des lemmes \cite{Beg} 6.2.1 et 6.2.2 de Bégueri.

\begin{Notss}
On note $\cC(n)$ la catégorie des $K$-schémas en groupes finis annulés par $p^n$, et $\cC(n)_0$ la sous-catégorie pleine de $\cC(n)$ dont les objets sont les $G_K\in\cC(n)$ tels que $G_K(K)=0$. 
\end{Notss}

\begin{Pt}
La catégorie $\cC(n)$ est abélienne et $\cC(n)_0$ est épaisse dans $\cC(n)$, en particulier exacte.
\end{Pt}

\begin{Lem}\label{coefface}
Tout objet de $\cC(n)$ est quotient d'un objet de $\cC(n)_0$.
\end{Lem}

\begin{proof}
Soit $G_K\in\cC(n)$.

\medskip

\emph{a) Le cas $G_K=\bbZ/p^m$ avec $0\leq m\leq n$.} Soit $t$ une uniformisante de $K$. La classe de $t$ dans 
$$
K^*/(K^*)^{p^m}=H^1(K,\mu_{p^m})=\Ext_{\bbZ/p^m}^1(\bbZ/p^m,\mu_{p^m})
$$
correspond à la classe d'isomorphisme d'une extension
$$
\xymatrix{
0\ar[r] & \mu_{p^m} \ar[r] & C_0 \ar[r] & \bbZ/p^m \ar[r] & 0,
}
$$
avec $C_0\in\cC(m)\subset\cC(n)$. Le bord $\delta:\bbZ/p^m\ra H^1(K,\mu_{p^m})$ déduit de cette extension s'identifie alors au morphisme de groupes
\begin{eqnarray*}
\bbZ/p^m & \ra & K^*/(K^*)^{p^m} \\
1 & \mapsto & t\ \modd\ (K^*)^{p^m},
\end{eqnarray*}
qui est injectif. Par conséquent $0=\mu_{p^m}(K)\xrightarrow{\sim} C_0(K)$, et donc $C_0\in\cC(n)_0$.

\medskip

\emph{b) Le cas où $G_K\in\cC(n)$ est étale.}  Soit $K'/K$ une extension finie séparable telle que le $K'$-schéma en groupes $G_{K'}$ soit constant. L'immersion fermée canonique $G_{K}^*\ra\Pi_{K'/K}G_{K'}^*$ induit alors, par bidualité de Cartier, un épimorphisme $\Pi_{K'/K}G_{K'}\ra G_{K}$. Or d'après \emph{a)}, il existe $C_{0}\in(\cC(n)_0)_{K'}$ et un épimorphisme $C_{0}\ra G_{K'}$, et donc un épimorphisme $\Pi_{K'/K}C_{0}\ra \Pi_{K'/K}G_{K'}$. D'où par composition un épimorphisme $\Pi_{K'/K}C_{0}\ra G_{K}$, et $\Pi_{K'/K}C_{0}\in\cC(n)_0$ par construction.

\medskip

\emph{c) Le cas $G_K\in\cC(n)$ quelconque.} Notant $G_{K}^0$ la composante neutre de $G_K$, on a une extension de $\cC(n)$
$$
\xymatrix{
0\ar[r] & G_{K}^0 \ar[r] & G_K \ar[r] & G_{K}^{\et} \ar[r] & 0
}
$$
avec $G_{K}^{\et}$ étale. D'après \emph{b)}, il existe $C_0\in\cC(n)_0$ et un épimorphisme $C_0\ra G_{K}^{\et}$. En prenant l'image réciproque de l'extension précédente par cet épimorphisme, on obtient un diagramme commutatif exact
$$
\xymatrix{
0\ar[r] & G_{K}^0 \ar[r] \ar@{=}[d] & \widetilde{G_K} \ar[r] \ar[d] & C_0 \ar[r] \ar[d] & 0 \\
0\ar[r] & G_{K}^0 \ar[r] & G_K \ar[r] & G_{K}^{\et} \ar[r] & 0,
}
$$
où $\widetilde{G_K}\in\cC(n)_0$ et $\widetilde{G_K}\ra G_K$ est un épimorphisme.
\end{proof}

\begin{Not}
On note $\cC(n)_{00}$ la sous-catégorie pleine de $\cC(n)_0$ dont les objets sont les $G_K\in\cC(n)$ tels que $G_K(K)=G_{K}^*(K)=0$.
\end{Not}

\begin{Pt}
La catégorie $\cC(n)_{00}$ est épaisse dans $\cC(n)$, en particulier exacte.
\end{Pt}

\begin{Rem}\label{etquot00}
Il résulte de la démonstration de \ref{coefface} que tout $G_{K}\in\cC(n)$ étale est quotient d'un objet de $\cC(n)_{00}$. En effet :

\medskip

\emph{a)} Montrons d'abord que si $t$ est une uniformisante de $K$, alors le $K$-schéma en groupes $C_0$ construit au point \emph{a)} de la démonstration de \ref{coefface} est un objet de $\cC(m)_{00}$. L'image de la classe d'isomorphisme de l'extension 
$$
\xymatrix{
0\ar[r] & \mu_{p^m} \ar[r] & C_0 \ar[r] & \bbZ/p^m \ar[r] & 0
}
$$
par la dualité de Cartier $(\cdot)^*:\Ext_{\bbZ/p^m}^1(\bbZ/p^m,\mu_{p^m})\xrightarrow{\sim}\Ext_{\bbZ/p^m}^1(\bbZ/p^m,\mu_{p^m})$ est un élément de même ordre, i.e. $p^m$, du groupe $\Ext_{\bbZ/p^m}^1(\bbZ/p^m,\mu_{p^m})=H^1(K,\mu_{p^m})$ (on sait même d'après \ref{calculEKamupn} que cet élément est la classe de $t^{-1}$). Le bord $\bbZ/p^m \ra H^1(K,\mu_{p^m})$ de 
l'extension duale est donc lui aussi injectif, i.e. $C_{0}^*(K)=0$, d'où $C_0\in\cC(m)_{00}$.

\medskip

\emph{b)} Soit maintenant $G_{K}\in\cC(n)$ étale. D'après \emph{a)}, il existe une extension finie séparable $K'/K$ et un épimorphisme $\Pi_{K'/K}C_{0}\ra G_K$ avec $C_{0}\in(\cC(n)_{00})_{K'}$. Mais alors $\Pi_{K'/K}C_{0}\in\cC(n)_{00}$, d'où le résultat.
\end{Rem}

\begin{Lem}\label{efface}
Pour tout $G_K\in\cC(n)_0$, il existe une suite exacte
$$
\xymatrix{
0\ar[r] & G_K \ar[r] & H_{K,0} \ar[r] & H_{K,1} \ar[r] & 0
}
$$
avec $H_{K,0},H_{K,1}\in\cC(n)_{00}$.
\end{Lem}

\begin{proof}
\emph{a) Le cas où $G_K$ est de type multiplicatif.} 
Soit $K'/K$ une extension finie séparable telle que le $K'$-schéma en groupes $G_{K'}$ soit diagonalisable. Il existe alors $r\in\bbN$, $n_1,\ldots,n_r\in\bbN$ inférieurs à $n$ et $\alpha_1,\ldots,\alpha_r\in\bbN$ tels que
$$
G_{K'}^*\simeq (\bbZ/p^{n_1})^{\alpha_1}\times\cdots\times(\bbZ/p^{n_r})^{\alpha_r}.
$$
Fixons un entier $i\in [1,r]$. D'après \ref{etquot00} \emph{a)}, il existe un épimorphisme $F_i\ra\bbZ/p^{n_i}$ avec $F_i\in(\cC(n_i)_{00})_{K'}$. D'après \ref{dimcohf}, cet épimorphisme induit une surjection $H^1(K',F_i)\ra H^1(K',\bbZ/p^{n_i})$. Comme $H^1(K',\bbZ/p^{n_i})$ n'est pas annulé par $p^{n_i-1}$, $H^1(K',F_i)$ non plus, et il existe donc une extension de $\cC(n_i)_{K'}$
$$
\xymatrix{
0\ar[r] & F_i \ar[r] & G_i \ar[r] & \bbZ/p^{n_i} \ar[r] & 0
}
$$
dont la classe dans $\Ext_{\bbZ/p^{n_i}}^1(\bbZ/p^{n_i},F_i)=H^1(K',F_i)$ est d'ordre $p^{n_i}$. Le bord $\bbZ/p^{n_i}\ra H^1(K',F_i)$ de cette extension est alors injectif, i.e. $0=F_i(K')\xrightarrow{\sim} G_i(K')$. Ceci étant, la suite exacte produit 
dans $\cC(n)_{K'}$
$$
\xymatrix{
0\ar[r] & F_{1}^{\alpha_1}\times\cdots\times F_{r}^{\alpha_r}  \ar[r] & G_{1}^{\alpha_1}\times\cdots\times G_{r}^{\alpha_r} \ar[r] & G_{K'}^*\ar[r] & 0
}
$$
donne lieu par dualité de Cartier à une suite exacte de $\cC(n)_{K'}$
$$
\xymatrix{
0\ar[r] & G_{K'}  \ar[r] & H \ar[r] & I \ar[r] & 0
}
$$
avec $I(K')=I^*(K')=H^*(K')=0$. Par restriction de Weil le long de $K'/K$, on en déduit une suite exacte de $\cC(n)$
$$
\xymatrix{
0\ar[r] & \Pi_{K'/K}G_{K'}  \ar[r] & \Pi_{K'/K}H \ar[r] & \Pi_{K'/K}I \ar[r] & 0.
}
$$
En la poussant par l'épimorphisme $\Pi_{K'/K}G_{K'}\ra G_K$ dual de l'immersion fermée canonique $G_{K}^*\ra \Pi_{K'/K}G_{K'}^*$, on obtient un diagramme commutatif exact de $\cC(n)$
$$
\xymatrix{
0\ar[r] & \Pi_{K'/K}G_{K'}  \ar[r] \ar[d] & \Pi_{K'/K}H \ar[r] \ar[d] & \Pi_{K'/K}I \ar[r] \ar@{=}[d] & 0 \\
0\ar[r] & G_K  \ar[r] & H_{K,0} \ar[r] & H_{K,1} \ar[r] & 0,
}
$$
avec $H_{K,0}^*(K)\subset (\Pi_{K'/K}H^*)(K)=H^*(K')=0$,
$$H_{K,1}(K)=I(K')=0\textrm{ et } H_{K,1}^*(K)=(\Pi_{K'/K}I^*)(K)=I^*(K')=0.$$
Donc $H_{K,1}\in\cC(n)_{00}$, puis $H_{K,0}\in\cC(n)_{00}$ puisque $G_K(K)=0$.

\medskip

\emph{b) Le cas où $G_K$ est radiciel.}
Notant $G_{K,m}$ le plus grand sous-groupe de type multiplicatif de $G_K$, on a une extension de $\cC(n)$
$$
\xymatrix{
0\ar[r] & G_{K,m} \ar[r] & G_K \ar[r] & G_{K,u} \ar[r] & 0
}
$$
avec $G_{K,u}$ unipotent radiciel, donc biradiciel, en particulier dans $\cC(n)_{00}$. D'après \emph{a)}, on a un diagramme commutatif exact de $\cC(n)$
$$
\xymatrix{
           & 0 \ar[d]                      & 0 \ar[d]             & & \\
0\ar[r] & G_{K,m} \ar[r] \ar[d] & G_K \ar[r] \ar[d] & G_{K,u} \ar[r] \ar@{=}[d] & 0 \\
0\ar[r] & H_{K,0} \ar[r] \ar[d] & H_{K,2} \ar[r] \ar[d] & G_{K,u} \ar[r] & 0 \\
           & H_{K,1} \ar@{=}[r] \ar[d] & H_{K,1}\ar[d] & & \\
           & 0                                & 0 & &
}
$$
avec $H_{K,0},H_{K,1}\in\cC(n)_{00}$. On a aussi $H_{K,2}\in\cC(n)_{00}$ puisque cette catégorie est épaisse dans $\cC(n)$.

\medskip

\emph{c) Le cas $G_K\in\cC(n)_0$ quelconque.} Notant $G_{K}^0$ la composante neutre de $G_K$, on a une extension de $\cC(n)$
$$
\xymatrix{
0\ar[r] & G_{K}^0 \ar[r] & G_K \ar[r] & G_{K}^{\et} \ar[r] & 0
}
$$
avec $G_{K}^{\et}$ étale. D'après \ref{etquot00}, il existe un épimorphisme $F\ra G_{K}^{\et}$ avec $F\in\cC(n)_{00}$. En tirant l'extension précédente par cet épimorphisme, on obtient un diagramme commutatif exact de $\cC(n)$
$$
\xymatrix{
0\ar[r] & G_{K}^0 \ar[r] & G_K \ar[r] & G_{K}^{\et} \ar[r] & 0 \\
0\ar[r] & G_{K}^0 \ar[r] \ar@{=}[u] & G_{K}' \ar[r] \ar[u] & F \ar[r] \ar[u] & 0.
}
$$
Appliquant alors \emph{b)} à $G_{K}^0$, on trouve un diagramme commutatif exact de $\cC(n)$
$$
\xymatrix{
           & 0 \ar[d]                    & 0 \ar[d]                  &                              & \\
0\ar[r] & G_{K}^0 \ar[r] \ar[d] & G_{K}' \ar[r] \ar[d] & F \ar[r] \ar@{=}[d] & 0 \\
0\ar[r] & H_{K,0} \ar[r] \ar[d] &  H_{K,2} \ar[r] \ar[d] & F \ar[r] & 0 \\
          & H_{K,1} \ar@{=}[r] \ar[d]  & H_{K,1}  \ar[d]                 &             & \\
          & 0                              & 0
}
$$
avec $H_{K,0},H_{K,1}\in\cC(n)_{00}$, puis $H_{K,2}\in\cC(n)_{00}$ par épaisseur. En poussant la deuxième extension verticale par l'épimorphisme $G_{K}'\ra G_K$, on obtient finalement un diagramme commutatif exact de $\cC(n)$
$$
\xymatrix{
0\ar[r] & G_{K}' \ar[r] \ar[d] & H_{K,2} \ar[r] \ar[d] & H_{K,1} \ar[r] \ar@{=}[d] & 0 \\
0\ar[r] & G_K \ar[r]  &  H_{K,3} \ar[r] & H_{K,1} \ar[r] & 0 
}
$$
avec $H_{K,3}(K)=0$ puisque $G_K\in\cC(n)_0$ par hypothèse et $H_{K,1}\in\cC(n)_0$, et $H_{K,3}^*(K)\subset H_{K,2}^*(K)=0$. D'où $H_{K,3}\in\cC(n)_{00}$.
\end{proof}

\begin{Pt}\label{effprovisoire}
La démonstration de \ref{coefface} montre que dans  l'énoncé, on peut remplacer $\cC(n)_0$ par la sous-catégorie pleine  
$\cC(n)_{0}'$ dont les objets sont les $G_K\in\cC(n)$ tels que $G_K(K)=0$ et $\cH^1(G_{K}^*)$ est connexe, i.e. les $G_K\in\cC(n)_0$ tels que $f_{\disc}(G_K)$ est un isomorphisme. De même, la démonstration de \ref{efface} montre que dans l'énoncé, si l'on remplace $\cC(n)_0$ par $\cC(n)_{0}'$, alors on peut remplacer $\cC(n)_{00}$ par la sous-catégorie pleine $\cC(n)_{00}'$ dont les objets sont les $G_K\in\cC(n)$ tels que $G_K(K)=G_{K}^*(K)=0$ et $\cH^1(G_{K}^*)$, $\cH^1(G_K)$ connexes, i.e. les $G_K\in\cC(n)_{00}$ tels que $f_{\disc}(G_K)$ et $f_{\disc}(G_{K}^*)$ sont des isomorphismes\footnote{Une conséquence du théorème de dualité sera que $\cC(n)_{0}'=\cC(n)_0$ et $\cC(n)_{00}'=\cC(n)_{00}$.}.
\end{Pt}
 
\section{Lemmes d'injectivité dans les catégories de groupes unipotents}

\begin{Lem}\label{annulExt1}
Soit $W_n$ le groupe quasi-algébrique unipotent des vecteurs de Witt de longueur $n$. 
\begin{enumerate}
\medskip

\item $W_n$ est un objet injectif de $\cU(n)$, et tout objet de $\cU(n)$ admet un monomorphisme dans un produit fini de copies de $W_n$ ; en particulier, $\cU(n)$ possède assez d'injectifs.

\medskip

\item $W_n$ est un objet injectif de $\Pro(\cU(n))$, et plus généralement, tout produit dans $\Pro(\cU(n))$ d'objets injectifs de $\cU(n)$ est injectif dans $\Pro(\cU(n))$. Tout objet de $\Pro(\cU(n))$ admet un monomorphisme dans un produit d'injectifs de $\cU(n)$. En particulier, $\Pro(\cU(n))$ possède assez d'injectifs.

\medskip

\item Si 
$$
\xymatrix{
0 \ar[r] & G' \ar[r] & G \ar[r] & G'' \ar[r] & 0
}
$$
est une suite exacte de $\Ind(\cU(n))$ avec $G''$ dans l'image essentielle de $\cI(\cU(n))$, alors la suite de groupes ordinaires
$$
\xymatrix{
0\ar[r] & \Hom_{\Ind(\cU(n))}(G'',W_n) \ar[r] & \Hom_{\Ind(\cU(n))}(G,W_n)  
}
$$
$$
\xymatrix{
\ar[r] & \Hom_{\Ind(\cU(n))}(G',W_n) \ar[r] & 0
}
$$
est exacte.

\medskip

\item Plus généralement, si 
$$
\xymatrix{
0 \ar[r] & G' \ar[r] & G \ar[r] & G'' \ar[r] & 0
}
$$
est une suite exacte de $\Ind(\Pro(\cU(n)))$ avec $G''$ dans l'image essentielle de $\cI\cP(\cU(n))$, alors la suite de groupes
ordinaires
$$
\xymatrix{
0\ar[r] & \Hom_{\Ind(\Pro(\cU(n)))}(G'',W_n) \ar[r] & \Hom_{\Ind(\Pro(\cU(n)))}(G,W_n)  
}
$$
$$
\xymatrix{
\ar[r] & \Hom_{\Ind(\Pro(\cU(n)))}(G',W_n) \ar[r] & 0
}
$$
est exacte.
\end{enumerate}
\end{Lem}

\begin{proof}
\emph{1.} 
Ce résultat est contenu dans \cite{DG} V \S 1 n${}^{\circ}$ 2. Considérons en effet un diagramme dans $\cU(n)$ de la forme
$$
\xymatrix{
0 \ar[r] & G' \ar[r] \ar[d]^f & G \ar[r] & G'' \ar[r] & 0 \\
           & W_n, & &
}
$$
dont la première ligne est exacte. D'après \ref{descrseQ} et \ref{reppf}, on peut supposer que ce diagramme se trouve dans la catégorie des groupes algébriques lisses unipotents annulés par $p^n$. Alors si $G''=\bbG_{a,k}$, $f$ se relève à $G$ d'après \emph{loc. cit.} 2.3 et 2.4, en notant que $p^nG=V_{G/k}^nF_{G/k}^nG=0$ implique $V_{G/k}^n=0:G^{(p^n)}\ra G$ puisque $G$ est lisse. Raisonnant par récurrence sur la longueur d'une suite de composition de $G''$ dont les quotients successifs sont isomorphes à $\bbG_{a,k}$, on peut encore relever $f$ à $G$ si $G''$ est connexe. Remarquons par ailleurs que $\Ext_{\bbZ/p^n}^1(\bbZ/p^n,W_n)=0$, puis que $\Ext_{\bbZ/p^n}^1(\bbZ/p^m,W_n)=0$ pour tout $1\leq m \leq n$ puisque le morphisme 
$$
\xymatrix{
\Hom_{\bbZ/p^n}(\bbZ/p^n,W_n)\ar[r] & \Hom_{\bbZ/p^n}(\bbZ/p^{n-m},W_n)
}
$$
induit par l'injection canonique $\bbZ/p^{n-m}\ra\bbZ/p^n$ est surjectif par $p$-divisibilité de $W_n(k)$ ; ainsi $\Ext_{\bbZ/p^n}^1(G'',W_n)=0$ pour $G''$ unipotent étale annulé par $p^n$, et donc $f$ se relève à $G$ si $G''$ est de ce type. On obtient finalement que $f$ se relève à $G$ quel que soit $G''$ en écrivant celui-ci comme extension d'un étale par un connexe, ce qui montre que $W_n$ est injectif dans $\cU(n)$. Enfin, d'après \emph{loc. cit.} 2.5, tout groupe algébrique unipotent $G$ sur $k$ est plongeable dans une puissance de $k$-schémas en groupes de Witt de longueur finie, et d'après \emph{loc. cit.} 2.4, cette longueur peut être prise égale à $n$ si $G$ est lisse et annulé par $p^n$.
 
\medskip

\emph{2.} Soit $\{I_e\}_{e\in E}$ une famille d'objets injectifs de $\cU(n)$, et montrons que $I:=\prod_{e\in E}I_e\in\Pro(\cU(n))$ est injectif. Soit 
$$
\xymatrix{
0 \ar[r] & G' \ar[r] & G \ar[r] & G'' \ar[r] & 0
}
$$
une suite exacte de $\Pro(\cU(n))$. D'après \cite{KS} 8.6.6 (a), une telle suite est pro-représentable par un système de suites exactes de $\cU(n)$
$$
\xymatrix{
0 \ar[r] & G_{m}' \ar[r] & G_m \ar[r] & G_{m}'' \ar[r] & 0, & m\in M.
}
$$
Appliquant pour chaque $e\in E$ le foncteur exact $\Hom_{\cU(n)}(\cdot,I_e)$, on obtient un ensemble indexé par $E$ de système inductifs indexés par $M$ de suites exactes de groupes ordinaires
$$
\xymatrix{
0 \ar[r] & \Hom_{\cU(n)}(G_{m}'',I_e) \ar[r] & \Hom_{\cU(n)}(G_m,I_e) \ar[r] & \Hom_{\cU(n)}(G_{m}',I_e) \ar[r] & 0.
}
$$
Passant à la limite inductive sur $M$, on obtient pour tout $e\in E$ la suite exacte
$$
\xymatrix{
0 \ar[r] & \Hom_{\Pro(\cU(n))}(G'',I_e) \ar[r] & \Hom_{\Pro(\cU(n))}(G,I_e) \ar[r] & \Hom_{\Pro(\cU(n))}(G',I_e) \ar[r] & 0.  
}
$$
Prenant enfin le produit sur $E$, on obtient la suite 
$$
\xymatrix{
0 \ar[r] & \Hom_{\Pro(\cU(n))}(G'',I) \ar[r] & \Hom_{\Pro(\cU(n))}(G,I) \ar[r] & \Hom_{\Pro(\cU(n))}(G',I) \ar[r] & 0,  
}
$$
qui est donc exacte, et $I$ est donc bien un injectif dans $\Pro(\cU(n))$. Par ailleurs, tout objet de $\Pro(\cU(n))$ est évidemment plongeable dans un produit d'objets de $\cU(n)$, qui à son tour est plongeable dans un produit d'injectifs de 
$\cU(n)$ d'après \emph{1.} ; un tel produit étant injectif d'après ce qui précède, $\Pro(\cU(n))$ possède donc assez d'injectifs.

\medskip

\emph{3.} Soit 
$$
\xymatrix{
0 \ar[r] & G' \ar[r] & G \ar[r] & G'' \ar[r] & 0
}
$$
une suite exacte de $\Ind(\cU(n))$ avec $G''$ dans l'image essentielle de $\cI(\cU(n))$. D'après \emph{loc. cit.}, une telle suite est ind-représentable par un système de suite exactes de $\cU(n)$
$$
\xymatrix{
0 \ar[r] & G_{i}' \ar[r] & G_i \ar[r] & G_{i}'' \ar[r] & 0, & i\in I.
}
$$
Appliquant le foncteur $\Hom_{\cU(n)}(\cdot,W_n)$, exact d'après \emph{1.}, on obtient un système inductif de suites exactes de groupes ordinaires
$$
\xymatrix{
0 \ar[r] & \Hom_{\cU(n)}(G_{i}'',W_n) \ar[r] & \Hom_{\cU(n)}(G_i,W_n) \ar[r] & \Hom_{\cU(n)}(G_{i}',W_n) \ar[r] & 0,\ i\in I.
}
$$
Comme $G''$ est dans l'image essentielle de $\cI(\cU(n))$ et $\Hom_{\cU(n)}(\cdot,W_n)$ exact, le système projectif 
$(\Hom_{\cU(n)}(G_{i}'',W_n))_{i\in I}$ vérifie la condition de Mittag-Leffler. Passant à la limite projective sur $I$ (essentiellement dénombrable par hypothèse), on obtient la suite
$$
\xymatrix{
0 \ar[r] & \Hom_{\Ind(\cU(n))}(G'',W_n) \ar[r] & \Hom_{\Ind(\cU(n))}(G,W_n) \ar[r] & \Hom_{\Ind(\cU(n))}(G',W_n) \ar[r] & 0,   
}
$$
qui est donc exacte.

\medskip

\emph{4.} Si 
$$
\xymatrix{
0 \ar[r] & G' \ar[r] & G \ar[r] & G'' \ar[r] & 0
}
$$
est une suite exacte de $\Ind(\Pro(\cU(n)))$ avec $G''$ dans l'image essentielle de $\cI\cP(\cU(n))$, alors comme $\Hom_{\Pro(\cU(n))}(\cdot,W_n)$ est exact d'après \emph{2.}, on obtient comme en \emph{3.} que la suite de groupes ordinaires
$$
\xymatrix{
0\ar[r] & \Hom_{\Ind(\Pro(\cU(n)))}(G'',W_n) \ar[r] & \Hom_{\Ind(\Pro(\cU(n)))}(G,W_n)  
}
$$
$$
\xymatrix{
\ar[r] & \Hom_{\Ind(\Pro(\cU(n)))}(G',W_n) \ar[r] & 0
}
$$
est exacte. 
\end{proof}

\begin{Lem}\label{annulExt2}
Soit $G$ un objet de la catégorie abélienne $\Ind(\Pro(\cU(n)))$ contenu dans l'image essentielle de $\cI\cP(\cU(n))$. Alors pour tout $i\neq 0$,
$$
\Hom_{D(\Ind\Pro(\cU(n)))}(G,W_n[i])=0,
$$
et pour tout $i\neq 0,1$,
$$
\Hom_{D(\Ind\Pro(\cU(n)))}(G,\bbZ/p^n[i])=0.
$$ 
\end{Lem}

\begin{proof}
Posons $P:=\Pro(\cU(n))$ pour alléger. Nous allons utiliser que le bifoncteur  
$$
\Hom_{\Ind(P)}:\Ind(P)^{\circ}\times \Ind(P) \lra (\textrm{Groupes})
$$
est dérivable à droite en un foncteur 
$$
R^+\Hom_{\Ind(P)}:D^-(\Ind(P))^{\circ}\times D^+(\Ind(P)) \lra D^+(\textrm{Groupes}),
$$
tel que 
$$
H^iR^+\Hom_{\Ind(P)}(G,H)= \Hom_{D(\Ind(P))}(G,H[i])
$$
pour tous $G\in D^-(\Ind(P))$, $H\in D^+(\Ind(P))$ et $i\in\bbZ$. Un tel foncteur dérivé existe d'après \cite{KS} 15.3.8, puisque  $P$ possédant assez d'injectifs d'après \ref{annulExt1} \emph{2.}, $\Ind(P)$ possède assez de quasi-injectifs (au sens de \emph{loc. cit.} 15.2.1) d'après \emph{loc. cit.} 15.2.7. De plus, comme $W_n$ est injectif dans $P$ d'après \ref{annulExt1} \emph{2.}, donc quasi-injectif dans $\Ind(P)$, le complexe $R^+\Hom_{\Ind(P)}(G,W_n)$ pour $G\in \Ind(P)$ peut être explicité  de la manière suivante (\cite{KS} 15.3.8) : on choisit un quasi-isomorphisme 
$$
\xymatrix{
\widetilde{G}^{\cdot}\ar[r] & G,
}
$$
où $\widetilde{G}^{\cdot}$ est un complexe de $\Ind(P)$ concentré en degrés $\leq 0$ et dont les composantes sont des sommes directes dans $\Ind(P)$ d'objets de $P$ ; alors
$$
R^+\Hom_{\Ind(P)}(G,W_n) \simeq \Hom_{\Ind(P)}^{\cdot}(\widetilde{G}^{\cdot},W_n).
$$
Par construction des sommes directes dans $\Ind(P)$, les transitions du ind-objet $\widetilde{G}^i\in\Ind(P)$ de degré $i\in\bbZ$ de $\widetilde{G}$ sont des monomorphismes, de sorte que $\widetilde{G}^i$ se trouve dans l'image essentielle de $\cI\cP(\cU(n))$. Si de plus $G$ se trouve dans cette image essentielle, alors le complexe exact de $\Ind(P)$
$$
\xymatrix{
\cdots \ar[r] & \widetilde{G}^{-2} \ar[r] & \widetilde{G}^{-1} \ar[r] & \widetilde{G}^0 \ar[r] & G \ar[r] & 0
}
$$
se découpe en suites exactes courtes, toutes dans l'image essentielle de $\cI\cP(\cU(n))$. Comme le foncteur $\Hom_{\Ind(P)}(\cdot,W_n)$ est exact sur cette image essentielle d'après \ref{annulExt1} \emph{4.}, il en résulte que le complexe $\Hom_{\Ind(P)}^{\cdot}(\widetilde{G}^{\cdot},W_n)$ est exact dans ce cas, i.e.
$$
\Hom_{D(\Ind\Pro(\cU(n)))}(G,W_n[i])=H^iR^+\Hom_{\Ind(P)}(G,W_n)=0\quad \textrm{pour tout $i\neq 0$}. 
$$
En utilisant ensuite la suite exacte d'Artin-Schreier-Witt de $\cU(n)$
$$
\xymatrix{
0 \ar[r] & \bbZ/p^n \ar[r] & W_n \ar[r]^{F-\Id} & W_n \ar[r] & 0,
}
$$
on obtient 
$$
\Hom_{D(\Ind\Pro(\cU(n)))}(G,\bbZ/p^n[i])=0\quad \textrm{pour tout $i\neq 0,1$}. 
$$
\end{proof}

\section{La flèche dans la catégorie dérivée} 

\begin{Pt}\label{H1tri}
Soit $n\in\bbN$. D'après \ref{cH1foncU}, on a un foncteur 
$$
\xymatrix{
\cH^1:\cC(n) \ar[r] & \fUnf,  
}
$$
qui est additif. Il se prolonge donc de manière unique aux catégories homotopiques en un foncteur triangulé
$$
\xymatrix{
\cH^1:K(\cC(n)) \ar[r] & K(\fUnf).
}
$$
On note $\cN(\cC(n))$ le nul-système (au sens de \cite{KS} 10.2.1) de $K(\cC(n))$ dont les objets sont les complexes exacts de la catégorie abélienne $\cC(n)$, et $D(\cC(n)):=K(\cC(n))/\cN(\cC(n))$ la catégorie dérivée de $\cC(n)$. On note $\cN(\fUnf)$ le nul-système de $K(\fUnf)$ dont les objets sont les complexes exacts au sens de \ref{DefsefQf} de la catégorie exacte 
$\fUnf$, c'est-à-dire ceux dont l'image par le foncteur triangulé composé
$$
\xymatrix{
K(\fUnf) \ar[r] & K(\Ind(\cP)) \ar[r]^{H} & \Ind(\cP) 
}
$$
est nulle. Puis on note $D(\fUnf):=K(\fUnf)/\cN(\fUnf)$ la \emph{catégorie dérivée de $\fUnf$}.
\end{Pt}

\begin{Lem}\label{LH1}
Le foncteur triangulé $\cH^1:K^-(\cC(n)) \ra  K^-(\fUnf)$ \ref{H1tri} admet un dérivé gauche
$$
\xymatrix{
L\cH^1:D^-(\cC(n)) \ar[r] & D^-(\fUnf),
}
$$
c'est-à-dire est localisable à gauche en un foncteur triangulé, tel que pour tout $G_{K}^{\cdot}\in K(\cC(n)_0)$, $L\cH^1(G_{K}^{\cdot})\simeq \cH^1(G_{K})^{\cdot}$. On a alors pour tout $G_K\in\cC(n)$ :
$$
\left\{ \begin{array}{ll}
H^i(L\cH^1(G_K))=0 & \textrm{si $i\neq -1,0$}, \\
H^0(L\cH^1(G_K))=\cH^1(G_K), &  \\
H^{-1}(L\cH^1(G_K))=\cH^0(G_K) ; & 
\end{array} \right.
$$
en particulier $L\cH^1(G_K)\in D_{\fUnf}^{[-1,0]}(\fUnf)$, et on a un triangle distingué dans $D(\fUnf)$ :
$$
\xymatrix{
\cH^0(G_K)[1] \ar[r] & L\cH^1(G_K) \ar[r] & \cH^1(G_K) \ar[r]^<<<{+1} &.
}
$$
De plus, si $\kappa/k$ est une extension de corps algébriquement close, le carré
$$
\xymatrix{
D^-(\cC(n)) \ar[r]^{L\cH^1} \ar[d]^{(\cdot)_{K_{\kappa}}} & D^-(\fUnf) \ar[d]^{(\cdot)_{\kappa}}\\
D^-(\cC(n)_{K_{\kappa}}) \ar[r]^{L\cH^1}  & D^-(\fUnf_{\kappa}) 
}
$$
est commutatif.
\end{Lem}

\begin{proof}
D'après \ref{coefface}, la sous-catégorie additive pleine $\cC(n)_0$ de $\cC(n)$ est génératrice. De plus, si 
$$
\xymatrix{
0 \ar[r] & G_K \ar[r] & G_K \ar[r] & G_{K}'' \ar[r] & 0
}
$$
est une suite exacte de $\cC(n)$ avec $G_K, G_{K}''\in \cC(n)_0$, alors $G_K\in \cC(n)_0$ et la suite 
$$
\xymatrix{
0 \ar[r] & \cH^1(G_K) \ar[r] & \cH^1(G_K) \ar[r] & \cH^1(G_{K}'') \ar[r] & 0
}
$$
de $\fUnf$ est exacte d'après \ref{cHsel}. D'après \cite{KS} 13.3.8 appliqué au foncteur $\cH^1:\cC(n)\ra\Ind(\cP)$, la catégorie $K^-(\cC(n)_0)$ est projective relativement au foncteur $\cH^1:K^-(\cC(n))\ra K^-(\Ind(\cP))$, ce qui par définition signifie que $K^-(\cC(n)_0)$ est projective relativement au foncteur triangulé $\cH^1:K^-(\cC(n))\ra K^-(\fUnf)$ et aux sous-catégories $\cN^-(\cC(n))$ et $\cN^-(\fUnf)$. D'où le foncteur dérivé $L\cH^1$ cherché d'après \cite{KS} 10.3.3 (b) et (10.3.2). De plus si $G_K\in\cC(n)$, il existe d'après \ref{coefface} une suite exacte
$$
\xymatrix{
0\ar[r] & H_{K,-1} \ar[r] & H_{K,0} \ar[r] & G_K \ar[r] & 0
}
$$
avec $H_{K,0}\in\cC(n)_0$ et donc aussi ($H_{K,-1}\in\cC(n)_0$). Le calcul de la cohomologie de $L\cH^1(G_K)$ en résulte avec \ref{cHsel}. La compatibilité de $L\cH^1$ aux extensions de corps algébriquement closes $\kappa/k$ résulte quant à elle de \ref{cH0cb} et \ref{cH1cb}.
\end{proof}

\begin{Pt}
Composant $L\cH^1:D^-(\cC(n))\ra D^-(\fUnf)$ par les foncteurs canoniques 
$$ 
\xymatrix{
D^-(\fUnf) \ar[r] & D^-(\Pro\Ind(\fUnf)) & \textrm{ou} & D^-(\fUnf) \ar[r] & D^-(\Ind\Pro(\fUnf)) 
}
$$
induits par les foncteurs exacts $\proind$ (\ref{fQfPIf}) et $\indpro$ (\ref{fQfIPf}) (cf. \ref{fQfexacte}), on obtient des foncteurs triangulés
$$ 
\xymatrix{
D^-(\cC(n)) \ar[r] & D^-(\Pro\Ind(\fUnf)) & \textrm{et} & D^-(\cC(n)) \ar[r] & D^-(\Ind\Pro(\fUnf)), 
}
$$
encore notés $L\cH^1$.
\end{Pt}

\begin{Pt}\label{Cderiv}
Le foncteur de dualité de Cartier
$$
\xymatrix{
(\cdot)^*:\cC(n) \ar[r] & \cC(n) ^{\circ} 
}
$$
étant exact, il se dérive trivialement en un foncteur triangulé 
$$
\xymatrix{
(\cdot)^*:D(\cC(n)) \ar[r] & D(\cC(n))^{\circ}.  
}
$$
\end{Pt}

\begin{Pt}
D'après \ref{PIS}, on a un foncteur de dualité de Serre
$$
\xymatrix{
(\cdot)^{\vee_S}:D^+(\Ind\Pro(\cU(n)))^{\circ} \ar[r] & D^-(\Pro\Ind(\cU(n))).
}
$$
\end{Pt}

\begin{Pt}\label{flecheH100'}
Pour tout $G_K\in\cC(n)_{00}'$, $\cH^1(G_K)$ et $\cH^1(G_{K}^*)$ sont connexes par définition, d'où par \ref{fcont} une flèche dans $\fUnf$
$$
\xymatrix{
f_{\cont}(G_K):\cH^1(G_K)\ar[r] & \cH^1(G_{K}^*)^{\vee_S}.
}
$$
\end{Pt}

\begin{Lem}\label{constrfderiv}
Il existe une unique transformation naturelle de foncteurs de source $D^b(\cC(n))$ et de but  $D^b(\Pro\Ind(\cU(n)))$
\begin{equation}\label{fderiv}
\xymatrix{
f_{\deriv}:L\cH^1\ar[r] & (\cdot)^{\vee_S}\circ L\cH^1 \circ (\cdot)^*,
}
\end{equation}
telle que pour tout $G_{K}^{\cdot}\in K^b(\cC(n)_{00}')$, le morphisme
$$
\xymatrix{
\cH^1(G_{K})^{\cdot} \ar[r] & (\cH^1(G_{K}^*)^{\vee_S})^{\cdot}
}
$$
qui s'en déduit par \ref{LH1} soit $f_{\cont}(G_{K})^{\cdot}$. De plus, si $\kappa/k$ est une extension de corps algébriquement close, le carré
$$
\xymatrix{
(\cdot)_{\kappa}\circ L\cH^1\ar[rr]^<<<<<<<<<<{(f_{\deriv})_{\kappa}} \ar[d]^{\rotatebox{90}{$\sim$}} && (\cdot)_{\kappa}\circ (\cdot)^{\vee_S}\circ L\cH^1 \circ (\cdot)^* \ar[d]^{\rotatebox{90}{$\sim$}} \\
L\cH^1\ar[rr]^<<<<<<<<{f_{\deriv}}\circ (\cdot)_{K_{\kappa}}  && (\cdot)^{\vee_S}\circ L\cH^1 \circ (\cdot)^* \circ (\cdot)_{K_{\kappa}} 
}
$$
est commutatif.
\end{Lem}

\begin{proof}
D'après \ref{coefface}, \ref{effprovisoire} et \cite{KS} 10.2.7 (ii), le foncteur canonique 
$$
\xymatrix{
K^b(\cC(n)_{0}')/\cN^b(\cC(n)_{0}')\ar[r] & D^b(\cC(n))
}
$$ 
est une équivalence de catégories. D'après \ref{efface}, \ref{effprovisoire} et \cite{KS} 10.2.7 (i), le foncteur canonique 
$$
\xymatrix{
K^b(\cC(n)_{00}')/\cN^b(\cC(n)_{00}')\ar[r] & K^b(\cC(n)_{0}')/\cN^b(\cC(n)_{0}')
}
$$ 
est une équivalence de catégories. D'où l'équivalence de catégories composée
$$
\xymatrix{
\iota:K^b(\cC(n)_{00}')/\cN^b(\cC(n)_{00}')\ar[r]^<<<<<{\sim} & D^b(\cC(n)).
}
$$

Maintenant, les flèches $f_{\cont}(G_K)$ pour $G_K\in\cC(n)_{00}'$ variable définissent une transformation naturelle de foncteurs de source $K^b(\cC(n)_{00}')$ et de but $D^b(\Pro\Ind(\cU(n)))$
$$
\xymatrix{
f_{\cont}:\cH^1\ar[r] & (\cdot)^{\vee_S}\circ \cH^1 \circ (\cdot)^*,
}
$$
\emph{a fortiori} une transformation naturelle de foncteurs de source $K^b(\cC(n)_{00}')/\cN^b(\cC(n)_{00}')$ et de but $D^b(\Pro\Ind(\cU(n)))$
$$
\xymatrix{
f_{\cont}:\cH^1\ar[r] & (\cdot)^{\vee_S}\circ \cH^1 \circ (\cdot)^*.
}
$$
Comme $L\cH^1\circ\iota\simeq \cH^1$ d'après \ref{LH1}, et 
$$
((\cdot)^{\vee_S}\circ L\cH^1\circ (\cdot)^*)\circ\iota\simeq (\cdot)^{\vee_S}\circ \cH^1\circ (\cdot)^*
$$
d'après \ref{Cderiv}, \ref{LH1} et le fait que $\cH^1(G_{K}^*)\in\fUnf^{00}\subset\Ind\Pro(\cU^0)$ pour tout $G_K\in\cC(n)_{00}'$, $f_{\cont}$ définit la transformation naturelle $f_{\deriv}$ souhaitée grâce au foncteur $\iota^{-1}$. Comme de plus $f_{\cont}$ commute aux extensions de corps algébriquement closes $\kappa/k$ d'après \ref{fcont}, $f_{\deriv}$ également.
\end{proof}

\begin{Prop}\label{HiLcH}
Soit $G_K\in\cC(n)$. 
\begin{enumerate}

\medskip

\item On a un triangle distingué dans $D(\Pro\Ind(\fUnf))$
$$
\xymatrix{
\cH^1(G_{K}^*)^{\vee_S} \ar[r] & (L\cH^1(G_{K}^*))^{\vee_S} \ar[r] & \cH^0(G_{K}^*)^{\vee_P} \ar[r]^<<<<{+1} & \ ;
}
$$
en particulier $(L\cH^1(G_{K}^*))^{\vee_S}\in D_{\fUnf}^{[-1,0]}(\Pro\Ind\cU(n))$ et
$$
\left\{ \begin{array}{ll}
(H^0((L\cH^1(G_{K}^*))^{\vee_S}))^0=(\cH^1(G_{K}^*)^0)^{\vee_S}\textrm{ et }\pi_0(H^0((L\cH^1(G_{K}^*))^{\vee_S}))=\cH^0(G_{K}^*)^{\vee_P}, &  \\
H^{-1}((L\cH^1(G_{K}^*))^{\vee_S})=(\pi_0(\cH^1(G_{K}^*)))^{\vee_P}. & 
\end{array} \right.
$$

\medskip

\item Le morphisme de groupes ordinaires  $H^{-1}(f_{\deriv}(G_K))(k)$ s'identifie au composé
$$
\xymatrix{
\cH^0(G_K)(k)=\Hom_{\cC(n)}(\bbZ/p^n,G_K) \ar[r]^<<<<<<<{(\cdot)^*} & \Hom_{\cC(n)}(G_{K}^*,\mu_{p^n}) \ar[d]^{\cH^1} \\ 
& \Hom_{\fUnf}(\cH^1(G_{K}^*),\cH^1(\mu_{p^n})) \ar[d]^{v\circ}\\
& \Hom_{\fUnf}(\cH^1(G_{K}^*),\bbZ/p^n).
}
$$

\medskip

\item Le morphisme de groupes ordinaires $H^0(f_{\deriv}(G_K))(k)$ s'identifie au composé
$$
\xymatrix{
\cH^1(G_K)(k)=\Hom_{D(\cC(n))}(\bbZ/p^n,G_K[1]) \ar[r]^<<<<<<<<<<{(\cdot)^*} & \Hom_{D(\cC(n))}(G_{K}^*,\mu_{p^n}[1]) \ar[d]^{L\cH^1} \\
& \Hom_{D(\Ind\Pro\cU(n))}(L\cH^1(G_{K}^*),L\cH^1(\mu_{p^n})[1]) \ar[d]^{v[1]\circ} \\ 
& \Hom_{D(\Ind\Pro\cU(n))}(L\cH^1(G_{K}^*),\bbZ/p^n[1]).
}
$$
\end{enumerate}
\end{Prop}

\begin{proof} Pour alléger, on pose $D:=D(\Ind\Pro(\cU(n))$.

\medskip

\emph{1.} Le triangle distingué 
$$
\xymatrix{
\cH^1(G_{K}^*)^{\vee_S} \ar[r] & (L\cH^1(G_{K}^*))^{\vee_S} \ar[r] & \cH^0(G_{K}^*)^{\vee_P} \ar[r]^<<<<{+1} & \ ;
}
$$
s'obtient en appliquant la dualité de Serre \ref{PIS} à l'image dans $D$ du triangle distingué \ref{LH1} 
$$
\xymatrix{
\cH^0(G_{K}^*)[1]\ar[r] & L\cH^1(G_{K}^*)) \ar[r] & \cH^1(G_{K}^*) \ar[r]^<<<<{+1} & 
}
$$
et en utilisant que $(\cH^0(G_{K}^*)[1])^{\vee_S}=(\cH^0(G_{K}^*)^{\vee_P}[1])[-1]=\cH^0(G_{K}^*)^{\vee_P}$.

\medskip

\emph{2. et 3. a) Le cas où $G_K\in\cC(n)_{00}'$.} On a 
$$H^{-1}(L\cH^1(G_K))(k)=\cH^0(G_K)(k)=G_K(K)=\Hom_{\cC(n)}(\bbZ/p^n,G_K)=0$$
et
$$H^{-1}((L\cH^1(G_{K}^*))^{\vee_S})(k)=(\pi_0(\cH^1(G_{K}^*)))^{\vee_P}=\Hom_{\fUnf}(\cH^1(G_{K}^*),\bbZ/p^n)=0,$$
d'où le point \emph{2.} dans ce cas. De plus, 
$$H^0(L\cH^1(G_K))(k)=\cH^1(G_K)(k)=H^1(K,G_K)=\Ext_{\cC(n)}^1(\bbZ/p^n,G_K)=\Hom_{D(\cC(n))}(\bbZ/p^n,G_K[1])$$
et d'après \ref{dualSexplicit3},
$$H^0((L\cH^1(G_{K}^*))^{\vee_S})(k)=H^0(\cH^1(G_{K}^*)^{\vee_S})(k)=\cH^1(G_{K}^*)^{\vee_S}(k)$$
$$=\Ext_{\Ind\Pro(\cU(n))}^1(\cH^1(G_{K}^*),\bbZ/p^n)=\Hom_{D}(\cH^1(G_{K}^*),\bbZ/p^n[1])=\Hom_{D}(L\cH^1(G_{K}^*),\bbZ/p^n[1]).$$
Comme $f_{\deriv}(G_K)=f_{\cont}(G_K)$ par définition \ref{constrfderiv}, $H^0(f_{\deriv}(G_K))(k)$ s'identifie bien dans ce cas à la composée du point \emph{3.} de l'énoncé.

\medskip

\emph{2. et 3. b) Le cas où $G_K\in\cC(n)_{0}'$.} On a 
$$H^{-1}(L\cH^1(G_K))(k)=\cH^0(G_K)(k)=G_K(K)=\Hom_{\cC(n)}(\bbZ/p^n,G_K)=0$$
et
$$H^{-1}((L\cH^1(G_{K}^*))^{\vee_S})(k)=(\pi_0(\cH^1(G_{K}^*)))^{\vee_P}=\Hom_{\fUnf}(\cH^1(G_{K}^*),\bbZ/p^n)=0,$$
d'où le point \emph{2.} dans ce cas. De plus, d'après \ref{efface} et \ref{effprovisoire}, il existe une suite exacte de 
$\cC(n)$
$$
\xymatrix{
0\ar[r] & G_K \ar[r] & H_{K,0} \ar[r] & H_{K,1} \ar[r] & 0
}
$$
avec $H_{K,0},H_{K,1}\in\cC(n)_{00}'$. Il en résulte un diagramme commutatif
$$
\xymatrix{
0\ar[r] & \cH^1(G_K) \ar[r] \ar[d]^{H^0(f_{\deriv}(G_K))} & \cH^1(H_{K,0}) \ar[r] \ar[d]^{H^0(f_{\deriv}(H_{K,0}))} & \cH^1(H_{K,1}) \ar[r] \ar[d]^{H^0(f_{\deriv}(H_{K,1}))} & 0 \\
0 \ar[r] & H^0((L\cH^1(G_{K}^*))^{\vee_S}) \ar[r] & H^0((L\cH^1(H_{K,0}^*))^{\vee_S}) \ar[r] & H^0((L\cH^1(H_{K,1}^*))^{\vee_S}) \ar[r] & 0 ;
}
$$
la première ligne est exacte d'après \ref{cHsel} puisque $\cH^0(H_{K,1})=0$, et la seconde ligne est exacte d'après \emph{1.} puisque $H^{-1}((L\cH^1(H_{K,1}^*))^{\vee_S})=H^{-1}((\cH^1(H_{K,1}^*))^{\vee_S})=0$. Sur les $k$-points, le carré de droite s'identifie d'après \emph{3. a)} au carré de droite du diagramme commutatif suivant dont les flèches verticales sont celles du point \emph{3.} de l'énoncé :
$$
\xymatrix{
0\ar[r] & \Hom_{\cC(n)}(\bbZ/p^n,G_K[1]) \ar[r] \ar[d] & \Hom_{\cC(n)}(\bbZ/p^n,H_{K,0}[1]) \ar[d]  \\
0 \ar[r] & \Hom_{D}(L\cH^1(G_{K}^*),\bbZ/p^n[1]) \ar[r] &  \Hom_{D}(L\cH^1(H_{K,0}^*),\bbZ/p^n[1])   \\
& \Hom_{\cC(n)}(\bbZ/p^n,H_{K,0}[1]) \ar[r] \ar[d] & \Hom_{\cC(n)}(\bbZ/p^n,H_{K,1}[1]) \ar[r] \ar[d] & 0 \\
& \Hom_{D}(L\cH^1(H_{K,0}^*),\bbZ/p^n[1])  \ar[r] &  \Hom_{D}(L\cH^1(H_{K,1}^*),\bbZ/p^n[1])  \ar[r] & 0.
}
$$
La première ligne de ce diagramme s'identifie aux $k$-points de la première ligne du diagramme précédent, en particulier elle est exacte. La seconde ligne est exacte puisque $\Hom_{D}(L\cH^1(H_{K,1}^*),\bbZ/p^n)=\Hom_{\cI\cP(\cU(n))}(\cH^1(H_{K,1}^*),\bbZ/p^n)=0$ ($\cH^1(H_{K,1}^*)$ étant connexe) d'une part, et puisque 
$$
\xymatrix{
H^0((L\cH^1(H_{K,0}^*))^{\vee_S})(k) \ar[r] & H^0((L\cH^1(H_{K,1}^*))^{\vee_S})(k)
}
$$ 
est surjective d'autre part. Par conséquent, la première flèche verticale s'identifie à $H^0(f_{\deriv}(G_K))(k)$, d'où le point \emph{3.} de l'énoncé dans ce cas. Notons au passage que : 
\begin{itemize}

\medskip

\item $\Hom_{D}(L\cH^1(G_{K}^*),\bbZ/p^n[2])$ s'injecte dans 
$$
\Hom_{D}(L\cH^1(H_{K,0}^*),\bbZ/p^n[2])=\Hom_{D}(\cH^1(H_{K,0}^*),\bbZ/p^n[2]),
$$
et par conséquent est nul d'après \ref{annulExt2} ;

\medskip

\item $H^{-1}((L\cH^1(G_{K}^*))^{\vee_S})$ s'injecte d'après \emph{1.} dans 
$$H^{-1}((L\cH^1(H_{K,0}^*))^{\vee_S})=H^{-1}((\cH^1(H_{K,0}^*))^{\vee_S})=0$$
et par conséquent est nul.
\end{itemize}

\medskip

\emph{2. et 3. c) Le cas $G_K\in\cC(n)$ quelconque.} D'après \ref{coefface} et \ref{effprovisoire}, il existe une suite exacte de $\cC(n)$
$$
\xymatrix{
0 \ar[r] & H_{K,-1} \ar[r] & H_{K,0} \ar[r] & G_K \ar[r] & 0
}
$$
avec $H_{K,-1}, H_{K,0} \in\cC(n)_{0}'$. Il en résulte un diagramme commutatif
$$
\xymatrix{
0\ar[r] & \cH^0(G_K) \ar[r] \ar[d]^{H^{-1}(f_{\deriv}(G_K))} & \cH^1(H_{K,-1}) \ar[r] \ar[d]^{H^0(f_{\deriv}(H_{K,-1}))} & \cH^1(H_{K,0}) \ar[d]^{H^0(f_{\deriv}(H_{K,0}))} \\
0 \ar[r] & H^{-1}((L\cH^1(G_{K}^*))^{\vee_S}) \ar[r] & H^0((L\cH^1(H_{K,-1}^*))^{\vee_S}) \ar[r] & H^0((L\cH^1(H_{K,0}^*))^{\vee_S}) \\
& \cH^1(H_{K,0}) \ar[r] \ar[d]^{H^0(f_{\deriv}(H_{K,0}))} & \cH^1(G_K) \ar[r] \ar[d]^{H^0(f_{\deriv}(G_K))} & 0 \\
& H^0((L\cH^1(H_{K,0}^*))^{\vee_S}) \ar[r] & H^0((L\cH^1(G_{K}^*))^{\vee_S}) \ar[r] & 0 \ ;
}
$$
la première ligne est exacte d'après \ref{cHsel} puisque $\cH^0(H_{K,0})=0$, et la seconde ligne est exacte d'après \emph{1.} puisque $H^{-1}((L\cH^1(H_{K,0}^*))^{\vee_S})=0$ comme noté à la fin de \emph{2. et 3. b)}. Sur les $k$-points, le carré central s'identifie d'après \emph{3. b)} au carré central du diagramme commutatif suivant dont les flèches verticales sont celles des points \emph{2.} et \emph{3.} de l'énoncé :
$$
\xymatrix{
0\lra \Hom_{\cC(n)}(\bbZ/p^n, G_K) \ar[r] \ar[d] & \Hom_{\cC(n)}(\bbZ/p^n,H_{K,-1}[1]) \ar[r] \ar[d] & \Hom_{\cC(n)}(\bbZ/p^n,H_{K,0}[1]) \ar[d] \\
0 \lra \Hom_{\fUnf}(\cH^1(G_{K}^*),\bbZ/p^n) \ar[r] & \Hom_{D}(L\cH^1(H_{K,-1}^*),\bbZ/p^n[1]) \ar[r] &  \Hom_{D}(L\cH^1(H_{K,0}^*),\bbZ/p^n[1])  \\
 \Hom_{\cC(n)}(\bbZ/p^n,H_{K,0}[1]) \ar[r] \ar[d] & \Hom_{\cC(n)}(\bbZ/p^n,G_{K}[1]) \ar[r] \ar[d]  & 0\\
 \Hom_{D}(L\cH^1(H_{K,0}^*),\bbZ/p^n[1])  \ar[r] &  \Hom_{D}(L\cH^1(G_{K}^*),\bbZ/p^n[1])  \ar[r] & 0. 
}
$$
La première ligne de ce diagramme s'identifie aux $k$-points de la première ligne du diagramme précédent, en particulier elle est exacte. La seconde ligne est exacte puisque $\cH^1(H_{K,0}^*)$ étant connexe, $\Hom_{\fUnf}(\cH^1(H_{K,0}^*),\bbZ/p^n)=0$, et puisque $\Hom_{D}(L\cH^1(H_{K,-1}^*),\bbZ/p^n[2])=0 $ comme noté à la fin de \emph{2. et 3. b)}. Par conséquent, la première flèche verticale s'identifie à $H^{-1}(f_{\deriv}(G_K))(k)$ et la dernière à $H^0(f_{\deriv}(G_K))(k)$. D'où le point les points \emph{2.} et \emph{3.}
\end{proof}

\begin{Cor} \label{disccontderiv}
Soit $G_K\in\cC(n)$. 
\begin{enumerate}

\medskip

\item On a 
$$
f_{\disc(G_K)}=-H^{-1}(f_{\deriv}(G_{K}^*))^{\vee_P}:\pi_0(\cH^1(G_K))\ra \cH^0(G_{K}^*)^{\vee_P}.
$$

\medskip

\item On a un diagramme commutatif exact de groupes ordinaires
$$
\xymatrix{
0 \ar[r] & \cH^1(G_K)^0(k) \ar[r] \ar[d]^{f_{\cont}(G_K)(k)} & \cH^1(G_K)(k)  \ar[d]^{H^0(f_{\deriv}(G_K))(k)}  \\
0\ar[r] & \Hom_{D(\Ind\Pro(\fUnf))}(\cH^1(G_{K}^*),\bbZ/p^n[1]) \ar[r] & \Hom_{D(\Ind\Pro(\fUnf))}(L\cH^1(G_{K}^*),\bbZ/p^n[1])  \\
&\cH^1(G_K)(k) \ar[r] \ar[d]^{H^0(f_{\deriv}(G_K))(k)} & \pi_0(\cH^1(G_K))(k) \ar[r] \ar[d]^{f_{\disc}(G_K)(k)} & 0 \\
&\Hom_{D(\Ind\Pro(\fUnf))}(L\cH^1(G_{K}^*),\bbZ/p^n[1]) \ar[r] & \Hom_{\fUnf}(\cH^0(G_{K}^*),\bbZ/p^n) \ar[r] & 0.
}
$$
\end{enumerate}
\end{Cor}

\begin{proof} On pose encore $D:=D(\Ind\Pro(\cU(n))$ pour alléger.

\medskip

\emph{1.} Résulte de \ref{HiLcH} \emph{2.} et \ref{compflecheK1k}.

\medskip

\emph{2.} 
La seconde ligne du diagramme est déduite du triangle distingué de $D$
$$
\xymatrix{
\cH^0(G_{K}^*)[1] \ar[r] & L\cH^1(G_{K}^*) \ar[r] & \cH^1(G_{K}^*) \ar[r]^<<<{+1} &
}
$$
par application de $\Hom_{D}(\cdot,\bbZ/p^n[1])$. Elle est exacte puisque $\cH^0(G_{K}^*)$ et $\bbZ/p^n$ étant concentrés en degré 0, $\Hom_D(\cH^0(G_{K}^*)[2],\bbZ/p^n[1])=0$, et puisque 
$$\Hom_D(\cH^1(G_{K}^*)[-1],\bbZ/p^n[1])=\Hom_D(\cH^1(G_{K}^*),\bbZ/p^n[2])=0$$
d'après \ref{annulExt2}. La commutativité du diagramme résulte de \ref{HiLcH} \emph{3.}
\end{proof}

\section{L'isomorphisme de dualité pour les $K$-schémas en groupes finis}

\subsection{\'Enoncé}

\begin{Th*}\label{Thp}
Soit $G_{K}\in\cC(n)$. Alors la flèche (\ref{fderiv}) 
$$
\xymatrix{
f_{\deriv}(G_{K}):L\cH^1(G_{K}) \ar[r] & (L\cH^1(G_{K}^*))^{\vee_S} & \in D(\Pro\Ind(\cU(n))),
}
$$
la flèche (\ref{fdisc})
$$
\xymatrix{
f_{\disc}(G_{K}):\pi_0(\cH^1(G_{K})) \ar[r] & \cH^0(G_{K}^*)^{\vee_P} & \in \cU(n),
}
$$
et la flèche (\ref{fcont})
$$
\xymatrix{
f_{\cont}(G_{K}):\cH^1(G_{K})^0 \ar[r] & (\cH^1(G_{K}^*)^0)^{\vee_S} & \in \fUnf
}
$$
sont des isomorphismes.
\end{Th*}

\begin{Pt*}\label{fdisccontgen}
Si $G_K$ est un $K$-schéma en groupes fini non nécessairement $p$-primaire, alors $G_K$ est produit d'un $K$-schéma en groupes fini $G_{K}^{(p')}$ d'ordre premier à $p$ par un $K$-schéma en groupes fini $p$-primaire $G_{K}^{(p)}$. Le groupe $H^1(K,G_{K}^{(p')})$ est fini (cf. \ref{cH1type}), en particulier 
$$
\cH^1(G_K)^0=\cH^1(G_{K}^{(p)})^0,
$$
et
$$
\pi_0(\cH^1(G_K))=\pi_0(\cH^1(G_{K}^{(p)}))\times\pi_0(\cH^1(G_{K}^{(p')}))=\pi_0(\cH^1(G_{K}^{(p)}))\times\cH^1(G_{K}^{(p')}). 
$$
En remplaçant $p^n$ par un entier annulant $G_K$ dans (\ref{flecheK1k}), on définit une flèche 
\begin{equation}\label{fdiscgen}
\xymatrix{
f_{\disc}(G_K):\pi_0(\cH^1(G_K)) \ar[r] & \cH^0(G_{K}^*)^{\vee_P},
}
\end{equation}
fonctorielle en $G_K$, et coïncidant avec $f_{\disc}(G_{K}^{(p')})$ (\ref{fdisc}) sur $G_{K}^{(p')}$. On pose par ailleurs
\begin{equation}\label{fcontgen}
\xymatrix{
f_{\cont}(G_K):=f_{\cont}(G_{K}^{(p)}):\cH^1(G_K)^0=\cH^1(G_{K}^{(p)})^0\ar[r] & (\cH^1(G_{K}^*)^0)^{\vee_S}=(\cH^1((G_{K}^{(p)})^*)^0)^{\vee_S}.
}
\end{equation}
On a alors la version suivante de l'énoncé $p$-primaire \ref{Thp} :
\end{Pt*}

\begin{Th*}\label{Th}
Soit $G_K$ un $K$-schéma en groupes fini. Le morphisme de groupes finis (\ref{fdiscgen})
$$
\xymatrix{
f_{\disc}(G_K):\pi_0(\cH^1(G_K)) \ar[r] & \cH^0(G_{K}^*)^{\vee_P} 
}
$$
et le morphisme de groupes admissibles unipotents strictement connexes (\ref{fcontgen})
$$
\xymatrix{
f_{\cont}(G_K):\cH^1(G_K)^0\ar[r] & (\cH^1(G_{K}^*)^0)^{\vee_S}
}
$$
sont des isomorphismes.
\end{Th*}

\subsection{Démonstration par dévissage}

\begin{Pt*}
Le théorème \ref{Th} est conséquence de la version $p$-primaire \ref{Thp} et de \ref{fdiscmultunip}.
\end{Pt*}

\begin{Pt*} \label{deviss1}
Si 
$$
\xymatrix{
0 \ar[r] & G_{K}' \ar[r] & G_K \ar[r] & G_{K}'' \ar[r] & 0
}
$$
est une suite exacte de $\cC(n)$, alors pour que $f_{\deriv}(G_K)$ soit un isomorphisme, il suffit que $f_{\deriv}(G_{K}')$ et $f_{\deriv}(G_{K}'')$ soient des isomorphismes.
\end{Pt*}

\begin{Pt*} \label{deviss2}
Pour que $f_{\deriv}(G_K)$ soit un isomorphisme, il faut et il suffit que $H^{-1}(f_{\deriv}(G_K))$ et $H^0(f_{\deriv}(G_K))$ soient des isomorphismes de $\fUnf^{\pizf}$ d'après \ref{LH1} et \ref{HiLcH} \emph{1.} Si $\kappa/k$ est une extension de corps algébriquement close, alors d'après \ref{constrfderiv}, on a pour $i=-1,0$,
$$
H^{i}(f_{\deriv}(G_K))(\kappa)=H^{i}(f_{\deriv}(G_K))_{\kappa}(\kappa)=H^{i}(f_{\deriv}(G_K)_{\kappa})(\kappa)=H^{i}(f_{\deriv}(G_{K_{\kappa}}))(\kappa).
$$
Pour que $H^{-1}(f_{\deriv}(G_K))$ et $H^0(f_{\deriv}(G_K))$ soient des isomorphismes, il faut et il suffit donc d'après \ref{foncfibgenIPcons00} qu'elles induisent des bijections au niveau des $k$-points. D'après \ref{disccontderiv}, il suffit pour cela que $f_{\disc}(G_{K}^*)$, $f_{\disc}(G_K)$ et $f_{\cont}(G_K)$ soient des isomorphismes.
\end{Pt*}

\begin{Pt*} \label{deviss4}
D'après \ref{fdiscmultunip} ci-dessous, $f_{\disc}(G_K)$ est un isomorphisme si $G_K\in\cC(n)$ est de type multiplicatif ou unipotent.
\end{Pt*}

\begin{Pt*} \label{deviss5}
D'après \ref{deviss2}, \ref{deviss4} et \ref{fcontalphap} ci-dessous, $f_{\deriv}(G_K)$ est un isomorphisme si $G_K=\alpha_p$.
\end{Pt*}

\begin{Pt*} \label{deviss6}
D'après \ref{deviss1}, \ref{deviss2}, \ref{deviss4} et \ref{fcontzsurp} ci-dessous, $f_{\deriv}(G_K)$ est un isomorphisme si $G_K\in\cC(n)$ est unipotent étale constant. D'après \ref{deviss4} et \ref{disccontderiv}, $f_{\disc}(G_K)$ et $f_{\cont}(G_K)$ sont aussi des isomorphismes dans ce cas.
\end{Pt*}

\begin{Pt*} \label{deviss8}
Soit $G_K\in\cC(n)$ unipotent étale. D'après \ref{deviss4}, $f_{\disc}(G_K)$ et $f_{\disc}(G_{K}^*)$ sont des isomorphismes. Montrons que $f_{\cont}(G_K)$ est injectif. Soit 
$$a\in \cH^1(G_K)(k)=H^1(K,G_K)=\Ext_{\bbZ/p^n}^1(\bbZ/p^n,G_K)$$
($\cH^1(G_K)$ est connexe d'après \ref{cH1type} \emph{4.}), représenté par une extension 
\begin{equation}\label{Exta}
\xymatrix{
0\ar[r] & G_K \ar[r] & E_{K,a} \ar[r] & \bbZ/p^n \ar[r] & 0\ ;
}
\end{equation}
soit 
$$
\xymatrix{
0\ar[r] & \mu_{p^n} \ar[r] & E_{K,a}^* \ar[r] & G_{K}^* \ar[r] & 0
}
$$
l'extension Cartier duale, puis 
$$
\xymatrix{
0\ar[r] & \cH^1(\mu_{p^n}) \ar[r] & \cH^1(E_{K,a}^*) \ar[r] & \cH^1(G_{K}^*) \ar[r] & 0
} 
$$
la suite de $\fUnf$ qui s'en déduit. Celle-ci est exacte d'après \ref{cHsel} puisque $G_{K}^*$ étant de type multiplicatif $p$-primaire, $G_{K}^*(K)=0$. Notant $\iota:\cH^1(G_{K}^*)^0\ra\cH^1(G_{K}^*)$ l'inclusion canonique, $f_{\cont}(G_K)(a)$ est par définition l'image directe de l'extension 
$$
\xymatrix{
0\ar[r] & \cH^1(\mu_{p^n}) \ar[r] & \iota^*\cH^1(E_{K,a}^*) \ar[r] & \cH^1(G_{K}^*)^0 \ar[r] & 0
}
$$
par le morphisme $v:\cH^1(\mu_{p^n})\ra\bbZ/p^n$, et celui-ci s'identifie à la projection canonique $\cH^1(\mu_{p^n})\ra\pi_0(\cH^1(\mu_{p^n}))$ d'après \ref{cH1Exs} \emph{1}. Si $f_{\cont}(G_K)(a)=0$, le morphisme de groupes finis
$$
\xymatrix{
\pi_0(\cH^1(\mu_{p^n}))\ar[r] & \pi_0(\iota^*\cH^1(E_{K,a}^*))
}
$$
est donc bijectif, puis donc 
$$
\xymatrix{
\pi_0(\cH^1(\mu_{p^n}))\ar[r] & \pi_0(\cH^1(E_{K,a}^*))
}
$$
injectif. D'après \ref{deviss4}, la flèche 
$$
\xymatrix{
(\bbZ/p^n)^{\vee_P} \ar[r] & (E_{K,a}(K))^{\vee_P}
}
$$
est donc injective, puis 
$$
\xymatrix{
E_{K,a}(K) \ar[r] &  \bbZ/p^n
}
$$
surjective. Comme $a\in H^1(K,G_K)$ est l'image de $1\in\bbZ/p^n$ par le bord de l'extension (\ref{Exta}), on a donc $a=0$, et $f_{\cont}(G_K)$ est bien injectif. Choisissons maintenant une extension fine séparable $K'/K$ telle que $G_{K'}$ soit constant. Notant $N_{K}^*$ le conoyau de l'immersion fermée canonique de $G_{K}^*$ dans la restriction de Weil $\Pi_{K'/K}G_{K'}^*$, la suite exacte de $\cC(n)$ 
$$
\xymatrix{
0 \ar[r] & N_K \ar[r] & \Pi_{K'/K}G_{K'} \ar[r] & G_K \ar[r] & 0 
}
$$
induit un carré commutatif
$$
\xymatrix{
\cH^1(\Pi_{K'/K}G_{K'}) \ar[r] \ar[d]^{f_{\cont}(\Pi_{K'/K}G_{K'})} & \cH^1(G_K) \ar[d]^{f_{\cont}(G_K)}  \\
(\cH^1(\Pi_{K'/K}G_{K'}^*)^0)^{\vee_S} \ar[r] & (\cH^1(G_{K}^*)^0)^{\vee_S}. 
}
$$
La seconde ligne est un épimorphisme strict : appliquant $(\cdot)^{\vee_S}$ au triangle distingué défini par la suite exacte 
$$
\xymatrix{
0 \ar[r] & \cH^1(G_{K}^*) \ar[r] & \cH^1(\Pi_{K'/K}G_{K'}^*) \ar[r] & \cH^1(N_{K}^*) \ar[r] & 0
}
$$
(cf. \ref{cHsel}, $N_{K}^*(K)=0$ puisque $N_{K}^*$ est de type multiplicatif $p$-primaire), on obtient la suite exacte
$$
\xymatrix{
(\cH^1(N_{K}^*)^0)^{\vee_S} \ar[r] & (\cH^1(\Pi_{K'/K}G_{K'}^*)^0)^{\vee_S} \ar[r] & (\cH^1(G_{K}^*)^0)^{\vee_S} \ar[r] & 0,
}
$$
en particulier $(\cH^1(\Pi_{K'/K}G_{K'}^*)^0)^{\vee_S} \ra (\cH^1(G_{K}^*)^0)^{\vee_S}$ est bien un épimorphisme strict. Comme la flèche $f_{\cont}(\Pi_{K'/K}G_{K'})$ est un isomorphisme d'après \ref{cH1Res} et \ref{deviss6}, il en résulte que $f_{\cont}(G_K)$ est un épimorphisme strict, et finalement un isomorphisme. Puis $f_{\deriv}(G_K)$ est un isomorphisme d'après \ref{deviss2}.
\end{Pt*}

\begin{Pt*} \label{deviss7}
D'après \ref{deviss1}, \ref{deviss2}, \ref{deviss4} et \ref{fcontmup} ci-dessous, $f_{\deriv}(G_K)$ est un isomorphisme si $G_K\in\cC(n)$ est diagonalisable $p$-primaire. D'après \ref{deviss4} et \ref{disccontderiv}, $f_{\disc}(G_K)$ et $f_{\cont}(G_K)$ sont aussi des isomorphismes dans ce cas.
\end{Pt*}

\begin{Pt*} \label{deviss9}
Soit $G_K\in\cC(n)$ de type multiplicatif $p$-primaire. D'après \ref{deviss4}, $f_{\disc}(G_K)$ et $f_{\disc}(G_{K}^*)$ sont des isomorphismes. Montrons que $f_{\cont}(G_K)$ est injectif. Choisissons une extension finie séparable $K'/K$ telle que 
$G_{K'}$ soit diagonalisable. Notant $C_K$ le conoyau de l'immersion fermée canonique de $G_K$ dans la restriction de Weil $\Pi_{K'/K}G_{K'}$, la suite exacte de $\cC(n)$ 
$$
\xymatrix{
0 \ar[r] & G_K \ar[r] & \Pi_{K'/K}G_{K'} \ar[r] & C_K \ar[r] & 0 
}
$$
induit un carré commutatif 
$$
\xymatrix{
\cH^1(G_K)^0 \ar[r] \ar[d]^{f_{\cont}(G_K)} & \cH^1(\Pi_{K'/K}G_{K'})^0 \ar[d]^{f_{\cont}(\Pi_{K'/K}G_{K'})}  \\
\cH^1(G_{K}^*)^{\vee_S} \ar[r] & \cH^1(\Pi_{K'/K}G_{K'}^*)^{\vee_S}
}
$$
($\cH^1(G_{K}^*)$ et $\cH^1(\Pi_{K'/K}G_{K'}^*)$ sont connexes d'après \ref{cH1type} \emph{4.}), dont la première ligne est un monomorphisme d'après \ref{cHsel} puisque $C_K(K)=0$, $C_K$ étant de type multiplicatif $p$-primaire. Comme $f_{\cont}(\Pi_{K'/K}G_{K'})$ est un isomorphisme d'après \ref{cH1Res} et \ref{deviss7}, $f_{\cont}(G_K)$ est donc bien injectif. En notant maintenant $N_{K}^*$ le conoyau de l'immersion fermée canonique de $G_{K}^*$ dans $\Pi_{K'/K}G_{K'}^*$, la suite exacte de $\cC(n)$ 
$$
\xymatrix{
0 \ar[r] & N_K \ar[r] & \Pi_{K'/K}G_{K'} \ar[r] & G_K \ar[r] & 0 
}
$$
induit un carré commutatif
$$
\xymatrix{
\cH^1(\Pi_{K'/K}G_{K'})^0 \ar[r] \ar[d]^{f_{\cont}(\Pi_{K'/K}G_{K'})} & \cH^1(G_K)^0 \ar[d]^{f_{\cont}(G_K)}  \\
\cH^1(\Pi_{K'/K}G_{K'}^*)^{\vee_S} \ar[r] & \cH^1(G_{K}^*)^{\vee_S}.
}
$$
La seconde ligne est un épimorphisme strict : appliquant $(\cdot)^{\vee_S}$ aux triangles distingués définis par les suites exactes (cf. \ref{cHsel})
$$
\xymatrix{
0 \ar[r] & \cH^1(G_{K}^*)/\Ima\partial \ar[r] & \cH^1(\Pi_{K'/K}G_{K'}^*) \ar[r] & \cH^1(N_{K}^*) \ar[r] & 0
}
$$
$$
\xymatrix{
0 \ar[r] & \Ima\partial \ar[r] & \cH^1(G_{K}^*) \ar[r] & \cH^1(G_{K}^*)/\Ima\partial \ar[r] & 0,
}
$$
on obtient des suites exactes
$$
\xymatrix{
0 \ar[r] & \cH^1(N_{K}^*)^{\vee_S} \ar[r] & \cH^1(\Pi_{K'/K}G_{K'}^*)^{\vee_S} \ar[r] & (\cH^1(G_{K}^*)/\Ima\partial)^{\vee_S}  \ar[r] & 0
}
$$
$$
\xymatrix{
0 \ar[r] & (\Ima\partial)^{\vee_P} \ar[r] & (\cH^1(G_{K}^*)/\Ima\partial)^{\vee_S} \ar[r] & \cH^1(G_{K}^*)^{\vee_S}  \ar[r] & 0,
}
$$
qui montrent que $\cH^1(\Pi_{K'/K}G_{K'}^*)^{\vee_S}\ra \cH^1(G_{K}^*)^{\vee_S}$ est composé d'épimorphismes stricts. Comme la flèche $f_{\cont}(\Pi_{K'/K}G_{K'})$ est un isomorphisme d'après \ref{cH1Res} et \ref{deviss7}, il en résulte que $f_{\cont}(G_K)$ est un épimorphisme strict, et finalement un isomorphisme. Puis $f_{\deriv}(G_K)$ est un isomorphisme d'après \ref{deviss2}.
\end{Pt*}

\begin{Pt*}
En conclusion, pour tout $G_K\in\cC(n)$, $f_{\deriv}(G_K)$ est un isomorphisme d'après \ref{deviss1},  \ref{deviss5}, \ref{deviss8} et \ref{deviss9}. D'après \ref{disccontderiv}, $f_{\disc}(G_K)$ et $f_{\cont}(G_K)$ sont donc aussi des isomorphismes pour tout $G_{K}\in\cC(n)$, d'où finalement le théorème \ref{Thp}.
\end{Pt*}

\subsection{Les cas multiplicatif et unipotent discrets}

\begin{Prop*} \label{fdiscmultunip}
Soit $G_K$ un $K$-schéma en groupes fini de type multiplicatif ou unipotent. Alors la flèche (\ref{fdiscgen})
$$
\xymatrix{
f_{\disc}(G_K):\pi_0(\cH^1(G_K))\ar[r] & \cH^0(G_{K}^*)^{\vee_P} & \in \cQ 
}
$$
est un isomorphisme.
\end{Prop*}

\begin{proof}
S'agissant d'une flèche dans $\cQ$, il est équivalent de montrer que la flèche induite sur les $k$-points, i.e. la flèche (\ref{flecheK1kq}), est un isomorphisme de groupes ordinaires. 

\medskip

\emph{a) Le cas $G_{K}=\mu_{p^n}$, $n\in\bbN$.} La flèche (\ref{oppflecheK1k}) associe alors à $a\in H^1(K,\mu_{p^n})$ le morphisme
\begin{eqnarray*}
\bbZ/p^n & \lra & \bbZ/p^n \\
1 & \lmapsto & v(\Id_*(a))=v(a).
\end{eqnarray*}
Autrement dit, la flèche (\ref{oppflecheK1k}) coïncide dans ce cas avec le morphisme $v:H^1(K,\mu_{p^n})\ra\bbZ/p^n$
induit par la valuation de $K$, et (\ref{flecheK1k}) est son opposé d'après \ref{compflecheK1k}. La flèche 
(\ref{flecheK1kq}) est donc bien une bijection pour $G_K=\mu_{p^n}$ d'après \ref{cH1Exs} \emph{1}.

\medskip

\emph{b) Le cas $G_{K}=\mu_m$, $m\in\bbN$ premier à $p$.} On vérifie encore dans ce cas que la flèche (\ref{fdiscgen}) s'identifie à la multiplication par $-1$ sur $\bbZ/m$, en particulier est un isomorphisme.

\medskip

\emph{c) Le cas où $G_K$ est diagonalisable.} Résulte de \emph{a)} et \emph{b)}.

\medskip

\emph{d) Le cas où $G_K$ est de type multiplicatif.} Soit $K_1/K$ une extension finie séparable telle que le $K_1$-schéma en groupes $G_{K_1}$ soit diagonalisable. Si $n$ est un entier tel que $nG_K=0$, on peut alors trouver une immersion fermée $G_{K_1}^*\ra(\bbZ/n)^{\alpha_1}$ pour un certain $\alpha_1\in\bbN$. Prenant sa restriction de Weil $\Pi_{K_1/K}G_{K_1}^*\ra\Pi_{K_1/K}(\bbZ/n)^{\alpha_1}$ et composant avec l'inclusion canonique $G_{K}^*\subset \Pi_{K_1/K}G_{K_1}^*$, on obtient une immersion fermée $G_{K}^*\ra\Pi_{K_1/K}(\bbZ/n)^{\alpha_1}$ qui, par dualité de Cartier, fournit un épimorphisme $\Pi_{K_1/K}(\mu_{n})^{\alpha_1}\ra G_{K}$. Appliquant la même construction au noyau de cet épimorphisme, on obtient une suite exacte
$$
\xymatrix{
\Pi_{K_2/K}(\mu_{n})^{r_2}\ar[r] & \Pi_{K_1/K}(\mu_{n})^{r_1} \ar[r] & G_K \ar[r] & 0.
}
$$
Appliquant le foncteur $\cH^1$, on obtient une suite de $\cI\cP$
$$
\xymatrix{
\cH^1(\Pi_{K_2/K}(\mu_{n})^{r_2}) \ar[r] & \cH^1(\Pi_{K_1/K}(\mu_{n})^{r_1}) \ar[r] & \cH^1(G_K) \ar[r] & 0,
}
$$
qui est exacte d'après \ref{cHsel}. Le foncteur $\pi_0:\cI(\cP^{\pizf})\ra\cI$ étant exact à droite d'après \ref{pi0exactd}, on en déduit une suite exacte dans $\cI$,
$$
\xymatrix{
\pi_0(\cH^1(\Pi_{K_2/K}(\mu_{n})^{r_2})) \ar[r] & \pi_0(\cH^1(\Pi_{K_1/K}(\mu_{n})^{r_1})) \ar[r] & \pi_0(\cH^1(G_K)) \ar[r] & 0,
}
$$
qui est ici une suite exacte de groupes quasi-algébriques finis. Par ailleurs, les dualités de Cartier et de Pontryagin étant exactes, et le foncteur 
$\Gamma(K,\cdot)$ exact à gauche,
on obtient finalement un diagramme commutatif de $\cQ$ à lignes exactes 
$$
\xymatrix{
\pi_0(\cH^1(\Pi_{K_2/K}(\mu_{n})^{r_2})) \ar[r] \ar[d] & \pi_0(\cH^1(\Pi_{K_1/K}(\mu_{n})^{r_1})) \ar[r] \ar[d] & \pi_0(\cH^1(G_K)) \ar[r] \ar[d] & 0\\
\cH^0((\Pi_{K_2/K}(\mu_{n})^{r_2})^*)^{\vee_P} \ar[r] & \cH^0((\Pi_{K_1/K}(\mu_{n})^{r_1})^*)^{\vee_P} \ar[r] & \cH^0(G_{K}^*)^{\vee_P} \ar[r] & 0,
}
$$
dans lequel les flèches verticales sont données par (\ref{fdisc}). Les deux premières sont des isomorphismes d'après \ref{cH1Res} et le cas \emph{c)}. La troisième flèche verticale est donc un isomorphisme, ce qui prouve le théorème dans le cas où $G_K$ est de type multiplicatif.

\medskip

\emph{e) Le cas où $G_K$ est unipotent.} Alors $\cH^1(G_K)$ est connexe d'après \ref{cH1type} \emph{4.}, i.e. 
$\pi_0(\cH^1(G_K))=0$, et $G_{K}^*$ est radiciel, donc $\cH^0(G_{K}^*)=0$.
\end{proof}

\subsection{L'isomorphisme différentiel du corps de classes}

\begin{Not*}
On note $\fm$ l'idéal maximal de $R$.
\end{Not*}

\begin{Def*}
Soit $M$ un $R$-module de type fini. La \emph{structure provectorielle $\fm$-adique} sur le $k$-vectoriel $M$ est la structure d'objet de $\cP(\cV)$ engendrée par l'ensemble de sous-$k$-vectoriels 
$$
\{\fm^nM \ |\ n\in\bbN\}.
$$ 
\end{Def*}

\begin{Pt*}\label{Kvectad}
Soit $V$ un $K$-vectoriel de dimension finie. Si $M$ et $M'$ sont deux réseaux de $V$, alors il existe $i\in\bbN$ tel que 
$M\subset \fm^{-i}M'$. De plus, si $M\subset M'$, alors le morphisme d'inclusion $M\ra M'$ est provectoriel pour les structures 
$\fm$-adiques. Il en résulte qu'il existe au plus une structure indprovectorielle sur $V$ engendrée par un ensemble $S$ de sous-$k$-vectoriels vérifiant :

\medskip

\begin{itemize}
\item les éléments de $S$ sont des réseaux de $V$ munis de leur structure provectorielle $\fm$-adique ;
\item l'ensemble $S$ est stable par multiplication par $\fm^{-i}$ pour tout $i\in\bbN$.
\end{itemize}

\medskip

Par ailleurs, une telle structure existe : si l'on choisit un réseau $M$ de $V$, l'ensemble de sous-$k$-vectoriels de $M$
$$
\{\fm^{-i}M\ |\ i\in\bbN\}
$$
engendre une structure convenable sur $V$. D'après ce qui précède, elle ne dépend pas du choix du réseau $M$, et commute au produits. Elle est de plus admissible, i.e. $V\in\cI\cP(\cV)^{\ad}$, puisque le quotient d'un réseau par un sous-réseau est de dimension finie sur $k$. La structure sur $V$ d'objet de $\fVf$ qui s'en déduit par l'équivalence $\ad:\cI\cP(\cV)^{\ad}\xrightarrow{\sim}\fVf$ est engendrée par l'ensemble de sous-$k$-vectoriels de $V$
$$
\{\fm^{\alpha}M\ |\ \alpha\in\bbZ\}.
$$
L'image de $V\in\fVf$ par le foncteur naturel $\fVf\ra\fUuf$ est l'image par le foncteur $\pi_*$ \ref{foncNerIndGrpf} du $K$-groupe vectoriel défini par $V$, que l'on notera simplement $\pi_*V$.
\end{Pt*}

\begin{Def*}
Soit $V$ un $K$-vectoriel de dimension finie. La \emph{structure admissible $\fm$-adique} sur le $k$-vectoriel $V$ est la structure d'objet de $\fVf$ définie en \ref{Kvectad}.
\end{Def*}

\begin{Pt*}
L'action de $K$ sur $\Omega_{K}^1$
$$
\xymatrix{
K\times\Omega_{K}^1 \ar[r] & \Omega_{K}^1,
}
$$
composée avec le résidu $\Res:\Omega_{K}^1\ra k$, définit un accouplement $k$-bilinéaire
\begin{equation}\label{accres}
\xymatrix{
K\times \Omega_{K}^1  \ar[r] & k.
}
\end{equation}
Celui-ci induit un isomorphisme dans $\fVf$
$$ 
\xymatrix{
\Omega_{K}^1\ar[r]^>>>>>{\sim} & \Hom_{\fVf}(K,k)=K^{\vee_L}
}
$$
(\ref{dualLexplicit3}). D'après \ref{dualLdualSad}, il en résulte un isomorphisme dans $\fUuf$
\begin{equation}\label{recipK}
\xymatrix{
\pi_*\Omega_{K}^1\ar[r]^>>>>>{\sim} & (\pi_*K)^{\vee_S}.
}
\end{equation}
Pour toute extension de corps algébriquement close $\kappa/k$, celui-ci donne lieu au diagramme commutatif
\begin{equation}\label{compFC}
\xymatrix{
\Ext_{\fUuf_{\kappa}}^1((\pi_*K)_{\kappa},\bbZ/p)=(\pi_*K)^{\vee_S}(\kappa) \ar[rr]^{(\pi_*F)^{\vee_S}(\kappa)} && (\pi_*K)^{\vee_S}(\kappa)=\Ext_{\fUuf_{\kappa}}^1((\pi_*K)_{\kappa},\bbZ/p)  \\
\Omega_{K_{\kappa}}^1=(\pi_*\Omega_{K}^1)(\kappa)\ar[rr]^{C} \ar[u]^{\rotatebox{90}{$\sim$}} && (\pi_*\Omega_{K}^1)(\kappa)=\Omega_{K_{\kappa}}^1,\ar[u]^{\rotatebox{90}{$\sim$}} 
}
\end{equation}
où $F:\bbG_{a,K}\ra\bbG_{a,K}$ est le Frobenius relatif de $\bbG_{a,K}$ et $C:\Omega_{K_{\kappa}}^1\ra \Omega_{K_{\kappa}}^1/dK_{\kappa}\xrightarrow{\sim}\Omega_{K_{\kappa}}^1$ est l'opérateur de Cartier.
\end{Pt*}

\begin{Lem*} \label{sefmupexacte}
On a la suite exacte suivante dans $\fUuf$ :
$$
\xymatrix{
0 \ar[r] & \cH^1(\mu_p) \ar[r]^{\dlog}  & \pi_*\Omega_{K}^1 \ar[r]^{C-1} & \pi_*\Omega_{K}^1 \ar[r] & 0 & =:F(\mu_p).
}
$$
\end{Lem*}

\begin{proof}
C'est un cas particulier de \ref{sef2AMRK} \emph{2.} ci-dessous.
\end{proof}

\begin{Pt*}\label{vRes}
On a le diagramme commutatif exact de $\fUuf$ suivant :
$$
\xymatrix{
0 \ar[r] & \cH^1(\mu_p) \ar[r]^{\dlog} \ar[d]^{v} & \pi_*\Omega_{K}^1 \ar[r]^{C-1} \ar[d]^{\Res} & \pi_*\Omega_{K}^1 \ar[r] \ar[d]^{\Res} & 0 & =:F(\mu_p)\\
0 \ar[r] & \bbZ/p \ar[r] & \bbG_{a,k}^{\pf} \ar[r]^{1-F} & \bbG_{a,k}^{\pf}  \ar[r] & 0 & =:\epsilon.
}
$$
Par conséquent,
$$
v_*F(\mu_p)\cong \Res^*\epsilon\ \in\Ext_{\fUuf}^1(\pi_*\Omega_K,\bbZ/p).
$$
Si $\kappa/k$ est une extension algébriquement close de $k$, l'extension de groupes obtenue en prenant les $\kappa$-points de cette extension de $\fUuf$ est la deuxième ligne du diagramme 
$$
\xymatrix{
0 \ar[r] & K_{\kappa}^*/(K_{\kappa}^*)^p \ar[r]^{\dlog} \ar[d]^{v} & \Omega_{K_{\kappa}}^1 \ar[r]^{C-1} \ar[d] & \Omega_{K_{\kappa}}^1 \ar[r] \ar@{=}[d] & 0 & =F(\mu_p)(\kappa)\\
0 \ar[r] & \bbZ/p \ar[r] \ar@{=}[d] & \Omega_{K_{\kappa}}^1/\Ker(C-1)\cap \Omega_{R_{\kappa}}^1 \ar[r]^>>>>>{C-1} \ar[d]^{\Res} & \Omega_{K_{\kappa}}^1 \ar[r] \ar[d]^{\Res} & 0 \\
0 \ar[r] & \bbZ/p \ar[r] & \kappa \ar[r]^{1-F} & \kappa \ar[r] & 0 & =\epsilon(\kappa).
}
$$
\end{Pt*}

\begin{Pt*}
On a les variantes suivantes de (\ref{recipK})
\begin{equation}\label{variantesrecipK1}
\xymatrix{
\pi_*R \ar[r]^>>>>>>{\sim} & (\pi_*\Omega_{K}^1/\pi_*\Omega_{R}^1)^{\vee_S} & & (\textrm{resp. } \pi_*\Omega_{K}^1/\pi_*\Omega_{R}^1 \ar[r]^>>>>>>{\sim} & (\pi_*R)^{\vee_S}),
}
\end{equation}
\begin{equation}\label{variantesrecipK2}
\xymatrix{
\pi_*\Omega_{R}^1 \ar[r]^>>>>>>{\sim} & (\pi_*K/\pi_*R)^{\vee_S} & & (\textrm{resp. } \pi_*K/\pi_*R \ar[r]^>>>>>>{\sim} & (\pi_*\Omega_{R}^1)^{\vee_S}),
}
\end{equation}
donnant lieu à des diagrammes commutatifs analogues à (\ref{compFC}).
\end{Pt*}

\subsection{Le cas de $\alpha_p$ continu}

\begin{Prop*} \label{fcontalphap}
La flèche (\ref{fcont})
$$
\xymatrix{
f_{\cont}(\alpha_p):\cH^1(\alpha_p) \ar[r] & \cH^1(\alpha_p)^{\vee_S} & \in \fUuf
}
$$
est un isomorphisme.
\end{Prop*}

\begin{Lem*}\label{calculEKaalphap}
Soit $a\in K$. L'identification 
$$K/K^p=H^1(K,\alpha_p)=\Ext_{\bbZ/p}^1(\bbZ/p,\alpha_p)$$
fait correspondre à la classe de $a$ la classe de l'extension
\begin{equation}\label{extEKaalphap}
\xymatrix{
0\ar[r] & \alpha_{p} \ar[r] & E_{K,a} \ar[r] & \bbZ/p \ar[r] & 0,
}
\end{equation}
où $E_{K,a}$ est le $K$-schéma en groupes dont l'algèbre de Hopf est la suivante :
\medskip
\begin{itemize}
\item 
la $K$-algèbre de $E_{K,a}$ est :
$$
A_{K,a}:=K[X,Y]\Big/(X^p-aY,Y^p-Y)\ ; 
$$
\item 
posant 
$$
\delta_{\lambda}:=\prod_{\substack{\mu\in\bbF_p\setminus\{\lambda\}}}\bigg(\frac{Y-\mu}{\lambda-\mu}\bigg)\ \in A_{K,a}
$$
pour tout $\lambda\in\bbF_p$, la co-multiplication de $A_{K,a}$ est donnée par :
\begin{eqnarray*}
A_{K,a} & \lra & A_{K,a}\otimes_K A_{K,a} \\
X & \lmapsto &\sum_{\lambda,\mu\in\bbF_p}\delta_{\mu}\otimes X\delta_{\lambda-\mu}+X\delta_{\mu}\otimes\delta_{\lambda-\mu} \\
Y & \mapsto & \sum_{\lambda,\mu\in\bbF_p}\lambda\delta_{\mu}\otimes\delta_{\lambda-\mu} \ ;
\end{eqnarray*}
\item
la co-unité de $A_{K,a}$ est donnée par :
\begin{eqnarray*}
A_{K,a} & \lra & K \\
X & \lmapsto & 0 \\
Y & \lmapsto & 0 \ ;
\end{eqnarray*}
\item
la co-inversion de $A_{K,a}$ est donnée par :
\begin{eqnarray*}
A_{K,a} & \lra & A_{K,a} \\
X & \lmapsto & -X \\
Y & \lmapsto & \sum_{\lambda\in\bbF_p}\lambda\delta_{-\lambda}.
\end{eqnarray*}
\end{itemize}
De plus, le dual de Cartier $E_{K,a}^*$ de $E_{K,a}$ est le $K$-schéma en groupes dont l'algèbre de Hopf est la suivante : 
\medskip
\begin{itemize}
\item 
la $K$-algèbre de $E_{K,a}^*$ est :
$$
A_{K,a}^*:=K[U,V]\Big/(U^p,V^p-1),
$$
avec pour élément diagonal
$$
\sum_{0\leq i\leq p-1,\ \lambda\in\bbF_p}X^i\delta_{\lambda}\otimes\frac{U^i}{i!}V^{\lambda} \in A_{K,a}\otimes_K A_{K,a}^*\ ;
$$
\item 
la co-multiplication de $A_{K,a}^*$ est donnée par :
\begin{eqnarray*}
A_{K,a}^* & \lra & A_{K,a}^*\otimes_K A_{K,a}^* \\
U & \lmapsto & U\otimes 1 + 1\otimes U \\
V & \lmapsto & V\otimes V\bigg(1+ a\sum_{i=1}^{p-1}\frac{U^i}{i!}\otimes\frac{U^{p-i}}{(p-i)!}\bigg)\ ;
\end{eqnarray*}
\item
la co-unité de $A_{K,a}^*$ est donnée par :
\begin{eqnarray*}
A_{K,a}^* & \lra & K \\
U & \lmapsto & 0\\
V & \lmapsto & 1\ ;
\end{eqnarray*}
\item
la co-inversion de $A_{K,a}^*$ est donnée par :
\begin{eqnarray*}
A_{K,a}^* & \lra & A_{K,a}^* \\
U & \lmapsto & -U\\
V & \lmapsto & V^{p-1}. 
\end{eqnarray*}
\end{itemize}
L'extension obtenue à partir de (\ref{extEKaalphap}) par dualité de Cartier est l'extension
\begin{equation}\label{extEK*aalphap}
\xymatrix{
0\ar[r] & \mu_{p} \ar[r]^<<<<<{\alpha} & E_{K,a}^* \ar[r]^{\beta} & \alpha_p \ar[r] & 0,
}
\end{equation}
où $\alpha(v)=(0,v)$ et $\beta(u,v)=u$.

En particulier, le $K$-schéma en groupes fini $E_{K,a}^*$ est de hauteur 1, i.e. annulé par Frobenius. Les dérivations 
\begin{eqnarray*}
\partial_1:A_{K,a}^* & \lra & A_{K,a}^* \\
U & \lmapsto & 1 \\
V & \lmapsto & a\frac{U^{p-1}}{(p-1)!}V, 
\end{eqnarray*}
\begin{eqnarray*}
\partial_2:A_{K,a}^* & \lra & A_{K,a}^* \\
U & \lmapsto & 0 \\
V & \lmapsto & V 
\end{eqnarray*}
sont invariantes, et forment une $K$-base de la $p$-algèbre de Lie de $E_{K,a}^*$. La matrice de la puissance $p$-ième symbolique dans la base $\{\partial_1,\partial_2\}$ est
$$
\left( \begin{array}{cc}
0 & 0 \\
a & 1 
\end{array} \right).
$$
\end{Lem*}

\begin{proof}
Pour faire les calculs, on utilise la $K$-base de $A_{K,a}$ donnée par 
$$
X^i\delta_{\lambda},\quad 0\leq i\leq p-1,\ \lambda\in\bbF_p.
$$ 
Dans cette base,
\medskip
\begin{itemize}
\item  la co-multiplication de $A_{K,a}$ est donnée par :
\begin{eqnarray*}
A_{K,a} & \lra & A_{K,a}\otimes_K A_{K,a} \\
X^i\delta_{\lambda} & \lmapsto &\sum_{0\leq j\leq i,\ \mu\in\bbF_p}\binom{i}{j}X^j\delta_{\mu}\otimes X^{i-j}\delta_{\lambda-\mu}\ ;
\end{eqnarray*}
\item
la co-unité de $A_{K,a}$ est donnée par :
\begin{eqnarray*}
A_{K,a} & \lra & K \\
X^i\delta_{\lambda} & \lmapsto & \delta_{i,0}\delta_{\lambda,0},
\end{eqnarray*}
où $\delta_{\alpha,\beta}$ est le symbole de Kronecker ;
\medskip
\item
la co-inversion de $A_{K,a}$ est donnée par :
\begin{eqnarray*}
A_{K,a} & \lra & A_{K,a} \\
X^i\delta_{\lambda} & \lmapsto & (-1)^iX^i\delta_{-\lambda}.
\end{eqnarray*}
\end{itemize}
\end{proof}

Ceci étant, rappelons le résultat suivant de la théorie de Artin-Milne \cite{AM}.

\begin{Th*} \label{se2} \emph{(Artin-Milne \cite{AM} 2.4)} 
Soit $X$ un $k$-schéma réduit possédant localement une $p$-base finie\footnote{Dans \emph{loc. cit.}, $X$ est supposé lisse sur $k$, mais les arguments donnés restent valides sous cette hypothèse plus faible.}. Soit $\alpha: X_{\FPPF}\ra X_{\et}$ le morphisme canonique du gros topos fppf de $X$ dans le petit topos étale de $X$. Soit $G$ un schéma en groupes fini localement libre sur $X$ de hauteur $1$. On a une suite exacte de faisceaux abéliens sur $X_{\et}$
\begin{equation}\label{sef2AM}
\xymatrix{
0\ar[r] & R^1\alpha_*G \ar[r]^<<<<<{=:\dlog_G} & \fgg\otimes_{\cO_X}Z^1F_{X^{(p^{-1})}/k*}\Omega_{X^{(p^{-1})}}^{\cdot}\ar[r]^<<<<<{\psi_0-\psi_1} & \fgg\otimes_{\cO_X}\Omega_{X}^1\ar[r] & 0,
}
\end{equation}
où $\fgg$ est la $p$-algèbre de Lie de $G$, et
\begin{itemize}

\medskip

\item le morphisme $\psi_0$ égale $1\otimes C$, où
$$C:Z^1F_{X^{(p^{-1})}/k*}\Omega_{X^{(p^{-1})}}^{\cdot}\ra \Omega_{X}^1$$
est l'opération de Cartier ;

\medskip

\item le morphisme $\psi_1$ égale $(\cdot)^{(p)}\otimes W^*$, où
$$(\cdot)^{(p)}: \fgg\ra \fgg$$
est l'opération de puissance $p$-ième symbolique, et
$$W^*:Z^1F_{X^{(p^{-1})}/k*}\Omega_{X^{(p^{-1})}}^{\cdot}\ra \Omega_{X}^1$$
est induit par $W$.
\end{itemize}
\end{Th*}

\begin{Prop*}\label{sef2AMRK}
\begin{enumerate}

\medskip

\item Soit $G$ un $R$-schéma en groupes fini plat de hauteur 1. Alors la suite des sections globales de (\ref{sef2AM}) est une suite exacte de $\cP$
\begin{equation} \label{sef2AMR}
\xymatrix{
0\ar[r] & \cH^1(G)  \ar[r]^<<<<<{\dlog_G} & \pi_*(\fgg\otimes_{R}F_{R^{(p^{-1})}/k*}\Omega_{R^{(p^{-1})}}^{1})\ar[r]^<<<<<{\psi_0-\psi_1} & \pi_*(\fgg\otimes_R\Omega_R^1)\ar[r] & 0. 
}
\end{equation}

\medskip

\item Soit $G_K$ un $K$-schéma en groupes fini de hauteur 1. Alors la suite des sections globales de (\ref{sef2AM}) est une suite exacte de $\fUuf$
\begin{equation} \label{sef2AMK}
\xymatrix{
0\ar[r] & \cH^1(G_K)  \ar[r]^<<<<<{\dlog_{G_K}} & \pi_*(\fgg_K\otimes_{K}F_{K^{(p^{-1})}/k*}\Omega_{K^{(p^{-1})}}^{1})\ar[r]^<<<<<{\psi_0-\psi_1} & \pi_*(\fgg_K\otimes_K\Omega_K^1)\ar[r] & 0. 
}
\end{equation}

\end{enumerate}
\end{Prop*}

\begin{proof}
En prenant les sections globales de (\ref{sef2AM}) pour $X=\Spec(R)$ (resp. $X=\Spec(K)$), on obtient la suite de groupes ordinaires sous-jacente à (\ref{sef2AMR}) (resp. (\ref{sef2AMK})). Comme 
$$H^1(R,R^1\alpha_*G)=H^2(R,G)=0$$
puisque $G$ est radiciel et d'après (\ref{dimcohfR}) (resp. $H^1(K,R^1\alpha_*G_K)=H^2(K,G_K)=0$ puisque $G_K$ est radiciel et d'après (\ref{dimcohf})), la suite de groupes ordinaires ainsi obtenue est exacte. De plus, pour toute extension de corps algébriquement close $\kappa/k$, on a une suite de groupes ordinaires analogue sur $R_{\kappa}$ (resp. $K_{\kappa}$), et ceci de façon naturelle en 
$\kappa/k$. D'après \ref{ratIP0}, les deux flèches de la suite (\ref{sef2AMR}) (resp. (\ref{sef2AMK})) sont donc proalgébriques (resp. indproalgébriques). Il en résulte que la suite (\ref{sef2AMR}) est exacte dans $\cP$, et il reste à montrer que la suite (\ref{sef2AMK}) est exacte dans $\cI\cP$.

Pour cela, on utilise que $G_K$ possède un modèle $G$ fini plat sur $R$, nécessairement de hauteur 1 (par exemple l'adhérence schématique de $G_K$ dans un $\bbG\bbL_{n,R}$ convenable, qui est finie sur $R$ parce que $G_K$ est radiciel). D'après la suite exacte (\ref{sef2AMR}) pour ce modèle, on obtient qu'une origine de $\pi_*(\fgg_K\otimes_K\Omega_K^1)$ est de définition dans  l'image indproalgébrique de $\psi_0-\psi_1:\pi_*(\fgg_K\otimes_{K}F_{K^{(p^{-1})}/k*}\Omega_{K^{(p^{-1})}}^1)\ra \pi_*(\fgg_K\otimes_K\Omega_K^1)$, et donc que ce dernier est un épimorphisme strict d'après \ref{critfoncfibgenIPcons}. De même, (\ref{sef2AMR}) montre avec \ref{origexact} que $\cH^1(G)$ est une origine du noyau indproalgébrique de $\psi_0-\psi_1:\pi_*(\fgg_K\otimes_{K}F_{K^{(p^{-1})}/k*}\Omega_{K^{(p^{-1})}}^1)\ra \pi_*(\fgg_K\otimes_K\Omega_K^1)$, et comme 
$\cH^1(G)$ est de définition dans $\cH^1(G_K)$ d'après \ref{cH1RcH1K},
il résulte de \ref{critfoncfibgenIPcons} que $\cH^1(G_K)$ est bien le noyau de $\psi_0-\psi_1$.
\end{proof}

En particulier :

\begin{Cor*}\label{calculH1EK*aalphap}
Dans la situation du lemme \ref{calculEKaalphap}, on a le diagramme commutatif exact de $\fUuf$ suivant :
$$
\xymatrix{
            &  0  \ar[d]                                                      & 0 \ar[d]                                                                & 0 \ar[d] \\                               
0 \ar[r] & \cH^1(\mu_p) \ar[r] \ar[dd]^{\dlog_{\mu_p}=\dlog} & \cH^1(E_{K,a}^*) \ar[r] \ar[dd]^{\dlog_{E_{K,a}^*}} & \cH^1(\alpha_p) \ar[r] \ar[dd]^{\dlog_{\alpha_p}=d} & 0 \\
\\
0 \ar[r] & \pi_*\Omega_{K}^1 \ar[r]^{i_2} \ar[dd]^{C-1} & \pi_*(\Omega_{K}^1)^{\oplus 2} \ar[r]^{p_1} \ar[dd]^{\binom{C\ 0}{-a\ C-1}} & \pi_*\Omega_{K}^1 \ar[r]\ar[dd]^{C} & 0 \\
\\
0 \ar[r] & \pi_*\Omega_{K}^1 \ar[r]^{i_2} \ar[d] & \pi_*(\Omega_{K}^1)^{\oplus 2} \ar[r]^{p_1}  \ar[d] & \pi_*\Omega_{K}^1 \ar[r] \ar[d] & 0 \\
           & 0                                                 &  0                                                                     & 0,
}
$$
où la première ligne est déduite de (\ref{extEK*aalphap}) par application de $\cH^1$, et les trois colonnes sont les suites exactes (\ref{sef2AMK}). 
 \end{Cor*}

\begin{proof}
Le diagramme est commutatif par fonctorialité en $G_K$ de la suite exacte (\ref{sef2AMK}). La première ligne est exacte d'après \ref{cHsel} puisque $\alpha_p(K)=0$. Les deux autres lignes sont scindées par définition. 
\end{proof}

\begin{proof}[Démonstration de la proposition \ref{fcontalphap}.]
Pour que $f_{\cont}(\alpha_p)$ soit un isomorphisme, il suffit d'après \ref{foncfibgenIPcons00} que sa fibre générique $f_{\cont}(\alpha_p)_{\eta}$ (\ref{deffibgenIP}) soit un isomorphisme. Comme pour toute extension de corps algébriquement close $\kappa/k$, on a $f_{\cont}(\alpha_p)(\kappa)=f_{\cont}(\alpha_{p,K_{\kappa}})(\kappa)$ par construction \ref{fcontconstr}, il suffit donc de montrer que le morphisme de groupes ordinaires $f_{\cont}(\alpha_p)(k)$ est bijectif. On commence par l'expliciter. Soit $a\in K$. D'après \ref{calculEKaalphap} et \ref{calculH1EK*aalphap}, l'image de $a$ par la composée
$$
\xymatrix{
K\ar[r] & K/K^p = H^1(K,\alpha_p)=\Ext_{\bbZ/p}^1(\bbZ/p,\alpha_p) \ar[r]^>>>>>>>>>>{(\cdot)^*} & \Ext_{\bbZ/p}^1(\alpha_p,\mu_p) \ar[d]^{\cH^1(\cdot)} \\
 & &  \Ext_{\fUuf}^1(\cH^1(\alpha_p),\cH^1(\mu_p))
}
$$
est une extension de $\fUuf$
\begin{equation*}\label{H1extEK*aalphap}
\xymatrix{                             
0 \ar[r] & \cH^1(\mu_p) \ar[r] & \cH^1(E_{K,a}^*) \ar[r] & \cH^1(\alpha_p) \ar[r] & 0, &=:h_a
}
\end{equation*}
qui admet la description différentielle suivante :
$$
\xymatrix{                            
0 \ar[r] & \cH^1(\mu_p) \ar[r] \ar@{=}[d] & \cH^1(E_{K,a}^*) \ar[r] \ar[d]^{\dlog_{E_{K,a}^*}}_{\rotatebox{90}{$\sim$}} & \cH^1(\alpha_p) \ar[r] \ar[d]^{d}_{\rotatebox{90}{$\sim$}} & 0 \\
0 \ar[r] & \cH^1(\mu_p) \ar[r]^<<<<<{\dlog} \ar@{=}[d] & \pi_*\Omega_{K}^1{}_{\searrow C-1}\times_{\pi_*\Omega_{K}^1}{}_{a\swarrow}\Ker(C) \ar[r] \ar[d] & \Ker (C) \ar[r] \ar[d]^{a} & 0 \\
0 \ar[r] & \cH^1(\mu_p) \ar[r]^{\dlog} & \pi_*\Omega_{K}^1 \ar[r]^{C-1} & \pi_*\Omega_{K}^1 \ar[r] & 0. & =:F(\mu_p)
}
$$
Autrement dit, 
$$h_a\cong(ad)^*F(\mu_p).$$ 
En particulier, $d^*F(\mu_p)\cong h_1$ est triviale par construction. Comme $ad=d(a\cdot)-\cdot da$, on a donc
$$h_a\cong-(\cdot da)^*F(\mu_p).$$
En combinant avec (\ref{vRes}), on obtient
$$
v_*h_a\cong -(\Res(\cdot da))^*\epsilon.
$$
Ainsi donc, l'image de $a\in K$ par la composée 
$$
\xymatrix{
K\ar[r] & K/K^p = H^1(K,\alpha_p)=\Ext_{\bbZ/p}^1(\bbZ/p,\alpha_p) \ar[r]^>>>>>>>>>>{(\cdot)^*} & \Ext_{\bbZ/p}^1(\alpha_p,\mu_p) \ar[d]^{\cH^1} \\
 & &  \Ext_{\fUuf}^1(\cH^1(\alpha_p),\cH^1(\mu_p)) \ar[d]^{v_*}\\
 & &  \Ext_{\fUuf}^1(\cH^1(\alpha_p),\bbZ/p) 
}
$$
est la classe d'isomorphisme de l'extension $-(\Res(\cdot da))^*\epsilon$, i.e. 
$$
f_{\cont}(\alpha_p)(a\modd K^p)=-\textrm{ classe de }(\Res(\cdot da))^*\epsilon.
$$

Maintenant, d'après (\ref{compFC}), on a un diagramme commutatif exact
$$
\xymatrix{
0\ar[r] & \Ext_{\fUuf}^1(\cH^1(\alpha_p),\bbZ/p)\ar[r] & \Ext_{\fUuf}^1(\pi_*K,\bbZ/p) \ar[rr]^{\Ext_{\fUuf}^1(\pi_*F,\bbZ/p)} & & \Ext_{\fUuf}^1(\pi_*K,\bbZ/p) \ar[r] & 0 \\
0 \ar[r] & K/K^p \ar[r]^d \ar[u]^{\rotatebox{90}{$\sim$}} & \Omega_{K}^1\ar[rr]^{C} \ar[u]^{\rotatebox{90}{$\sim$}} && \Omega_{K}^1\ar[u]^{\rotatebox{90}{$\sim$}}  \ar[r] & 0,
}
$$
où la première flèche verticale est donnée par
\begin{eqnarray*}\label{areconnalphap}
K/K^p = H^1(K,\alpha_p) & \lra &  \Ext_{\fUuf}^1(\cH^1(\alpha_p),\bbZ/p) \\
a \modd K^p & \lmapsto & \textrm{classe de }(\Res(\cdot da))^*\epsilon.
\end{eqnarray*}
Comme cette flèche est bijective, on en conclut que $f_{\cont}(\alpha_p)(k)$ l'est également.
\end{proof}

\subsection{Le cas de $\bbZ/p$ continu}

\begin{Prop*} \label{fcontzsurp}
La flèche (\ref{fcont})
$$
\xymatrix{
f_{\cont}(\bbZ/p):\cH^1(\bbZ/p) \ar[r] & (\cH^1(\mu_p)^0)^{\vee_S} & \in \fUuf
}
$$
est un isomorphisme.
\end{Prop*}

\begin{Lem*}\label{calculEKaZ/p}
Soit $a\in K$. L'identification 
$$K/(1-F)(K)=H^1(K,\bbZ/p)=\Ext_{\bbZ/p}^1(\bbZ/p,\bbZ/p)$$
fait correspondre à la classe de $a$ la classe de l'extension
\begin{equation}\label{extEKaZ/p}
\xymatrix{
0\ar[r] & \bbZ/p \ar[r] & E_{K,a} \ar[r] & \bbZ/p \ar[r] & 0,
}
\end{equation}
où $E_{K,a}$ est le $K$-schéma en groupes dont l'algèbre de Hopf est la suivante :
\medskip
\begin{itemize}
\item 
la $K$-algèbre de $E_{K,a}$ est :
$$
A_{K,a}:=K[X,Y]\Big/(X-X^p-aY,Y^p-Y)\ ; 
$$
\item 
posant 
$$
\delta_{\lambda}:=\prod_{\substack{\mu\in\bbF_p\setminus\{\lambda\}}}\bigg(\frac{Y-\mu}{\lambda-\mu}\bigg)\ \in A_{K,a}
$$
pour tout $\lambda\in\bbF_p$, la co-multiplication de $A_{K,a}$ est donnée par :
\begin{eqnarray*}
A_{K,a} & \lra & A_{K,a}\otimes_K A_{K,a} \\
X & \lmapsto &\sum_{\lambda,\mu\in\bbF_p}\delta_{\mu}\otimes X\delta_{\lambda-\mu}+X\delta_{\mu}\otimes\delta_{\lambda-\mu} \\
Y & \mapsto & \sum_{\lambda,\mu\in\bbF_p}\lambda\delta_{\mu}\otimes\delta_{\lambda-\mu} \ ;
\end{eqnarray*}
\item
la co-unité de $A_{K,a}$ est donnée par :
\begin{eqnarray*}
A_{K,a} & \lra & K \\
X & \lmapsto & 0 \\
Y & \lmapsto & 0 \ ;
\end{eqnarray*}
\item
la co-inversion de $A_{K,a}$ est donnée par :
\begin{eqnarray*}
A_{K,a} & \lra & A_{K,a} \\
X & \lmapsto & -X \\
Y & \lmapsto & \sum_{\lambda\in\bbF_p}\lambda\delta_{-\lambda}.
\end{eqnarray*}
\end{itemize}
De plus, le dual de Cartier $E_{K,a}^*$ de $E_{K,a}$ est le $K$-schéma en groupes dont l'algèbre de Hopf est la suivante : 
\medskip
\begin{itemize}
\item 
la $K$-algèbre de $E_{K,a}^*$ est :
$$
A_{K,a}^*:=K[U,V]\Big/(U^p,V^p-1),
$$
avec pour élément diagonal
$$
\sum_{0\leq i\leq p-1,\ \lambda\in\bbF_p}X^i\delta_{\lambda}\otimes\frac{U^i}{i!}V^{\lambda} \in A_{K,a}\otimes_K A_{K,a}^*\ ;
$$
\item 
la co-multiplication de $A_{K,a}^*$ est donnée par :
\begin{eqnarray*}
A_{K,a}^* & \lra & A_{K,a}^*\otimes_K A_{K,a}^* \\
U & \lmapsto & U\otimes 1 + 1\otimes U + \sum_{i=1}^{p-1}\frac{U^i}{i!}\otimes\frac{U^{p-i}}{(p-i)!}\\
V & \lmapsto & V\otimes V\bigg(1- a\sum_{i=1}^{p-1}\frac{U^i}{i!}\otimes\frac{U^{p-i}}{(p-i)!}\bigg)\ ;
\end{eqnarray*}
\item
la co-unité de $A_{K,a}^*$ est donnée par :
\begin{eqnarray*}
A_{K,a}^* & \lra & K \\
U & \lmapsto & 0\\
V & \lmapsto & 1\ ;
\end{eqnarray*}
\item
la co-inversion de $A_{K,a}^*$ est donnée par :
\begin{eqnarray*}
A_{K,a}^* & \lra & A_{K,a}^* \\
U & \lmapsto & -U\\
V & \lmapsto & V^{p-1}. 
\end{eqnarray*}
\end{itemize}
L'extension obtenue à partir de (\ref{extEKaZ/p}) par dualité de Cartier est l'extension
\begin{equation}\label{extEK*aZ/p}
\xymatrix{
0\ar[r] & \mu_{p} \ar[r]^<<<<<{\alpha} & E_{K,a}^* \ar[r]^{\beta} & \mu_p \ar[r] & 0,
}
\end{equation}
où $\alpha(v)=(0,v)$ et $\beta(u,v)=\sum_{i=0}^{p-1}\frac{u^i}{i!}$.

En particulier, le $K$-schéma en groupes fini $E_{K,a}^*$ est de hauteur 1, i.e. annulé par Frobenius. Les dérivations 
\begin{eqnarray*}
\partial_1:A_{K,a}^* & \lra & A_{K,a}^* \\
U & \lmapsto & 1 + \frac{U^{p-1}}{(p-1)!}\\
V & \lmapsto & -a\frac{U^{p-1}}{(p-1)!}V, 
\end{eqnarray*}
\begin{eqnarray*}
\partial_2:A_{K,a}^* & \lra & A_{K,a}^* \\
U & \lmapsto & 0 \\
V & \lmapsto & V 
\end{eqnarray*}
sont invariantes, et forment une $K$-base de la $p$-algèbre de Lie de $E_{K,a}^*$. La matrice de la puissance $p$-ième symbolique dans la base $\{\partial_1,\partial_2\}$ est
$$
\left( \begin{array}{cc}
1 & 0 \\
-a & 1 
\end{array} \right).
$$
\end{Lem*}

\begin{proof}
Pour faire les calculs, on utilise la $K$-base de $A_{K,a}$ donnée par 
$$
X^i\delta_{\lambda},\quad 0\leq i\leq p-1,\ \lambda\in\bbF_p.
$$ 
Dans cette base,
\medskip
\begin{itemize}
\item  la co-multiplication de $A_{K,a}$ est donnée par :
\begin{eqnarray*}
A_{K,a} & \lra & A_{K,a}\otimes_K A_{K,a} \\
X^i\delta_{\lambda} & \lmapsto &\sum_{0\leq j\leq i,\ \mu\in\bbF_p}\binom{i}{j}X^j\delta_{\mu}\otimes X^{i-j}\delta_{\lambda-\mu}\ ;
\end{eqnarray*}
\item
la co-unité de $A_{K,a}$ est donnée par :
\begin{eqnarray*}
A_{K,a} & \lra & K \\
X^i\delta_{\lambda} & \lmapsto & \delta_{i,0}\delta_{\lambda,0},
\end{eqnarray*}
où $\delta_{\alpha,\beta}$ est le symbole de Kronecker ;
\medskip
\item
la co-inversion de $A_{K,a}$ est donnée par :
\begin{eqnarray*}
A_{K,a} & \lra & A_{K,a} \\
X^i\delta_{\lambda} & \lmapsto & (-1)^iX^i\delta_{-\lambda}.
\end{eqnarray*}
\end{itemize}
\end{proof}

Par ailleurs, le lemme \ref{sef2AMRK} admet le 

\begin{Cor*}\label{calculH1EK*aZ/p}
Dans la situation du lemme \ref{calculEKaZ/p}, on a le diagramme commutatif exact de $\fUuf$ suivant :
$$
\xymatrix{
            &  0  \ar[d]                                                      & 0 \ar[d]                                                                & 0 \ar[d] \\                               
0 \ar[r] & \cH^1(\mu_p) \ar[r] \ar[dd]^{\dlog_{\mu_p}=\dlog} & \cH^1(E_{K,a}^*) \ar[r] \ar[dd]^{\dlog_{E_{K,a}^*}} & \cH^1(\mu_p) \ar[r] \ar[dd]^{\dlog_{\mu_p}=\dlog} & 0 \\
\\
0 \ar[r] & \pi_*\Omega_{K}^1 \ar[r]^{i_2} \ar[dd]^{C-1} & \pi_*(\Omega_{K}^1)^{\oplus 2} \ar[r]^{p_1} \ar[dd]^{\binom{C-1\ 0}{a\ C-1}} & \pi_*\Omega_{K}^1 \ar[r]\ar[dd]^{C-1} & 0 \\
\\
0 \ar[r] & \pi_*\Omega_{K}^1 \ar[r]^{i_2} \ar[d] & \pi_*(\Omega_{K}^1)^{\oplus 2} \ar[r]^{p_1}  \ar[d] & \pi_*\Omega_{K}^1 \ar[r] \ar[d] & 0 \\
           & 0                                                 &  0                                                                     & 0,
}
$$
où la première ligne est déduite de (\ref{extEK*aZ/p}) par application de $\cH^1$, et les trois colonnes sont les suites exactes  (\ref{sef2AMK}). 
\end{Cor*}

\begin{proof}
Le diagramme est commutatif par fonctorialité en $G_K$ de la suite exacte \ref{sef2AMK}. La première ligne est exacte d'après \ref{cHsel} puisque $\mu_p(K)=0$. Les deux autres lignes sont scindées par définition. 
\end{proof}

\begin{proof}[Démonstration de la proposition \ref{fcontzsurp}.]
Comme dans la démonstration de \ref{fcontalphap}, il suffit de montrer que l'application $f_{\cont}(\bbZ/p)(k)$ est bijective, et on commence par l'expliciter. Soit $a\in K$. D'après \ref{calculEKaZ/p} et \ref{calculH1EK*aZ/p}, l'image de $a$ par la composée
$$
\xymatrix{
K\ar[r] & K/(1-F)(K) = H^1(K,\bbZ/p)=\Ext_{\bbZ/p}^1(\bbZ/p,\bbZ/p) \ar[r]^>>>>>>>>>>{(\cdot)^*} & \Ext_{\bbZ/p}^1(\mu_p,\mu_p) \ar[d]^{\cH^1} \\
 & &  \Ext_{\fUuf}^1(\cH^1(\mu_p),\cH^1(\mu_p))
}
$$
est une extension de $\fUuf$
\begin{equation*}\label{H1extEK*Z/p}
\xymatrix{                             
0 \ar[r] & \cH^1(\mu_p) \ar[r] & \cH^1(E_{K,a}^*) \ar[r] & \cH^1(\mu_p) \ar[r] & 0, &=:h_a
}
\end{equation*}
qui admet la description différentielle suivante :
$$
\xymatrix{                            
0 \ar[r] & \cH^1(\mu_p) \ar[r] \ar@{=}[d] & \cH^1(E_{K,a}^*) \ar[r] \ar[d]^{\dlog_{E_{K,a}^*}}_{\rotatebox{90}{$\sim$}} & \cH^1(\mu_p) \ar[r] \ar[d]^{\dlog}_{\rotatebox{90}{$\sim$}} & 0 \\
0 \ar[r] & \cH^1(\mu_p) \ar[r]^<<<<<{\dlog} \ar@{=}[d] & \pi_*\Omega_{K}^1{}_{\searrow C-1}\times_{\pi_*\Omega_{K}^1}{}_{-a\swarrow}\Ker(C-1) \ar[r] \ar[d] & \Ker (C-1) \ar[r] \ar[d]^{-a} & 0 \\
0 \ar[r] & \cH^1(\mu_p) \ar[r]^{\dlog} & \pi_*\Omega_{K}^1 \ar[r]^{C-1} & \pi_*\Omega_{K}^1 \ar[r] & 0. & =:F(\mu_p)
}
$$
Autrement dit
$$
h_a\cong - (a\dlog)^*F(\mu_p).
$$
En combinant avec (\ref{vRes}), on obtient 
$$
v_*h_a\cong -(\Res(a\dlog))^*\epsilon.
$$
Puis 
$$
\iota^*v_*h_a\cong -(\Res(a\dlog\circ\iota))^*\epsilon.
$$
Ainsi donc, l'image de $a\in K$ par la composée 
$$
\xymatrix{
K\ar[r] & K/(1-F)(K) = H^1(K,\bbZ/p)=\Ext_{\bbZ/p}^1(\bbZ/p,\bbZ/p) \ar[r]^>>>>>>>>>>{(\cdot)^*} & \Ext_{\bbZ/p}^1(\mu_p,\mu_{p}) \ar[d]^{\cH^1} \\
 & &  \Ext_{\fUuf}^1(\cH^1(\mu_p),\cH^1(\mu_p)) \ar[d]^{v_*}\\
 & &  \Ext_{\fUuf}^1(\cH^1(\mu_p),\bbZ/p) \ar[d]^{\iota^*} \\
 & &  \Ext_{\fUuf}^1(\cH^1(\mu_p)^0,\bbZ/p)
}
$$
est la classe d'isomorphisme de l'extension $-(\Res(a\dlog\circ\iota))^*\epsilon$, i.e.
$$
f_{\cont}(\bbZ/p)(a \modd (1-F)(K))=-\textrm{ classe de }(\Res(a\dlog\circ\iota))^*\epsilon. 
$$

Maintenant, d'après (\ref{variantesrecipK2}), on a un diagramme commutatif exact
$$
\xymatrix{
            &                                                                                                         &&                                                                 &
\Ext_{\fUuf}^1(\cH^1(\mu_p)^0,\bbZ/p) \\           
0\ar[r] & \Ext_{\fUuf}^1(\pi_*\Omega_{R}^1,\bbZ/p)\ar[rr]^{\Ext_{\fUuf}^1(1-C,\bbZ/p)} && \Ext_{\fUuf}^1(\pi_*\Omega_{R}^1,\bbZ/p) \ar[r] & \Ext_{\fUuf}^1(\Ker(1-C),\bbZ/p) \ar[r] \ar[u]^{\rotatebox{90}{$\sim$}}_{(\dlog\circ\iota)^*} & 0 \\
0 \ar[r] & K/R \ar[rr]^{\overline{1-F}} \ar[u]^{\rotatebox{90}{$\sim$}} && K/R \ar[r] \ar[u]^{\rotatebox{90}{$\sim$}} & K/(1-F)(K)\ar[u]^{\rotatebox{90}{$\sim$}}  \ar[r] & 0,
}
$$
où la troisième flèche verticale est donnée par
\begin{eqnarray*}\label{areconnZ/p}
K/(1-F)(K) = H^1(K,\bbZ/p) & \lra &  \Ext_{\fUuf}^1(\cH^1(\mu_p)^0,\bbZ/p) \\
a \modd (1-F)(K) & \lmapsto & \textrm{classe de } (\dlog\circ\iota)^*\Res(a\cdot)^*\epsilon=\textrm{classe de }(\Res(a\dlog\circ\iota))^*\epsilon.
\end{eqnarray*}
Comme cette flèche est bijective, on en conclut que $f_{\cont}(\bbZ/p)(k)$ l'est également.
\end{proof}

\subsection{Le cas de $\mu_p$ continu}\label{LecasGKmup}

\begin{Prop*} \label{fcontmup}
La flèche (\ref{fcont})
$$
\xymatrix{
f_{\cont}(\mu_p):\cH^1(\mu_p)^0 \ar[r] & \cH^1(\bbZ/p)^{\vee_S} & \in \fUuf
}
$$
est un isomorphisme.
\end{Prop*}

\begin{Lem*}\label{calculH1EK*amup}
Soit $a\in K^*\setminus (K^*)^p$. En reprenant les notations de \ref{calculEKamupn} dans le cas $n=1$, on a le diagramme commutatif exact de $\fUuf$ suivant :
$$
\xymatrix{                               
0\ar[r] & \cH^0(\bbZ/p)  \ar[r] \ar@{=}[d] & \cH^1(\mu_p) \ar[r] \ar[d]_{\rotatebox{90}{$\sim$}}^{\dlog_{\mu_p}=\dlog} & \cH^1(E_{K,a}^*) \ar[r] \ar[d]_{\rotatebox{90}{$\sim$}}^{\dlog_{E_{K,a}^*}} & \cH^1(\bbZ/p) \ar[r] \ar[d]_{\rotatebox{90}{$\sim$}} & 0 \\
0\ar[r] & \cH^0(\bbZ/p)  \ar[r] & \Ker(C-1) \ar[r] & \pi_*\Omega_{K}^1/K^p\dlog a \ar[r]^<<<<<{\frac{C-1}{\dlog a}} & \pi_*K/\pi_*(1-F)(\pi_*K) \ar[r] & 0,  
}
$$
où la première ligne est la suite exacte longue \ref{cHsel} déduite de l'extension (\ref{extEK*amupn}). 
\end{Lem*}

\begin{proof}
On utilise le théorème de comparaison plat/cristallin de Trihan-Vauclair \cite{TV} 8.14, dans le cas particulier du groupe $p$-divisible $G$ défini par $(a^{-1}\ \modd\ (K^*)^{p^n})_{n\in\bbN^*}\in\lp K^*/(K^*)^{p^n}$ (de sorte que $G(1)=E_{K,a}^*$). Soit $\bbD$ le cristal de Dieudonné contravariant de $G$. Choisissant un $p$-anneau $K_2$ de longueur $2$ de corps résiduel $K$, on note $\cS_{1}(\bbD)$ le complexe de faisceaux sur le petit site étale de $\Spec(K)$, concentré en degrés $0$ et $1$, égal à
$$
\xymatrix{
(\Fil^1\bbD)_{K_2}/p \ar[r]^<<<<{\overline{\nabla}} & \bbD_K\otimes_K\Omega_{K}^{1}, 
}
$$
où 
$$
(\Fil^1\bbD)_{K_2}:=\Ker\big(\bbD_{K_2}\ra\bbD_K\xrightarrow{\sim}(F_{K}^*\bbD_K)^{\nabla}\ra (F_{K}^*\bbD_K/\Ima (V_{\bbD})_K)^{\nabla}\big).
$$
D'après le théorème de comparaison, on a une suite exacte longue  
$$
0 \ra \alpha_*E_{K,a}^* \ra \cH^0(\cS_{1}(\bbD)) \xrightarrow{1-\frac{F_{\bbD}\otimes F_{K_2}}{p}}  \cH^0(\bbD_K\otimes_K\Omega_{K}^{\cdot}) 
$$
\begin{equation}\label{selcomp}
\ra   R^1\alpha_*E_{K,a}^* \xrightarrow{=:\dlog_{E_{K,a}^*}} \cH^1(\cS_{1}(\bbD))  \xrightarrow{1-\frac{F_{\bbD}\otimes F_{K_2}}{p}}  \cH^1(\bbD_K\otimes_K\Omega_{K}^{\cdot}).
\end{equation}

Ceci étant, rappelons la description du cristal de Dieudonné de $E_{K,a}$ donnée par de Jong \cite{dJ} 4.3.5. Appliquant $\bbD(\cdot)$ à l'extension (\ref{extEKamupn}), on obtient une extension
$$
\xymatrix{
0 \ar[r] & \bbD(\bbZ/p)\ar[r] & \bbD(E_{K,a}) \ar[r] & \bbD(\mu_p) \ar[r] & 0.
}
$$
Evaluant sur $K$, on obtient la suite exacte
$$
\xymatrix{
0 \ar[r] & (K.1,d,1,0) \ar[r] & (K.1\oplus K.\widetilde{1},\nabla,(F_{\bbD})_K,(V_{\bbD})_K) \ar[r] & (K.\widetilde{1},d,0,1) \ar[r] & 0,
}
$$
où
\begin{tabular}{@{} l @{} l @{} l @{}}
$\nabla=\left( \begin{array}{cc}
0 & \dlog a \\
0 & 0 
\end{array} \right),$ & $\quad(F_{\bbD})_K=\left( \begin{array}{cc}
1 & 0 \\
0 & 0 
\end{array} \right),$ & $\quad(V_{\bbD})_K=\left( \begin{array}{cc}
0 & 1 \\
0 & 0 
\end{array} \right).$
\end{tabular}

En particulier,
$$
(\Fil^1\bbD)_{K_2}=\Ker(K_2.1\oplus K_2.\widetilde{1} \ra K.1)=pK_2.1\oplus K_2.\widetilde{1},
$$
et $\cS_{1}(\bbD)$ est donné par
$$
\xymatrix{
pK_2.1\oplus K.\widetilde{1} \ar[rr]^{\nabla=\binom{0\ \dlog a}{0\ 0}} && \Omega_{K}^1.1\oplus\Omega_{K}^1.\widetilde{1}.
}
$$
On obtient donc les diagrammes commutatifs 
$$
\xymatrix{
\cH^0(\cS_{1}(\bbD)) \ar[rr]^{1-\frac{F_{\bbD}\otimes F_{K_2}}{p}} \ar@{=}[d]&&\cH^0(\bbD_K\otimes_K\Omega_{K}^{\cdot})\\
pK_2.1 \ar[rr]^{1-\frac{F_{K_2}}{p}} && d^{-1}(K^p\dlog a) \ar[u]_{\rotatebox{90}{$\sim$}} \ar@{=}[d] \\
K \ar[u]_{\rotatebox{90}{$\sim$}}^{\underline{p}} \ar[rr]^{-F_K}_{\sim} && K^p 
}
$$
et 
$$
\xymatrix{
\cH^1(\cS_{1}(\bbD)) \ar[rr]^{1-\frac{F_{\bbD}\otimes F_{K_2}}{p}} \ar@{=}[d]&&\cH^1(\bbD_K\otimes_K\Omega_{K}^{\cdot}) \ar@{=}[d]\\
\Omega_{K}^1.1\oplus\Omega_{K}^1.\widetilde{1}/\nabla(K.\widetilde{1})  \ar[rr]^{\binom{1-\frac{F_{K_2}}{p}\ 0}{0\ 1}} && \Omega_{K}^1.1\oplus\Omega_{K}^1.\widetilde{1}/\Ima(\nabla).
}
$$
Par conséquent,
$$
\xymatrix{
R^1\alpha_*E_{K,a}^* \ar[r]_>>>>>{\sim}^<<<<{\dlog_{E_{K,a}^*}} & \Ker \bigg(\Omega_{K}^1.1\oplus\Omega_{K}^1.\widetilde{1}/\nabla(K.\widetilde{1})  \ar[rr]^{\binom{1-\frac{F_{K_2}}{p}\ 0}{0\ 1}} && \Omega_{K}^1.1\oplus\Omega_{K}^1.\widetilde{1}/\Ima(\nabla)\bigg).
}
$$

Maintenant, en utilisant la description de $\nabla$, et le fait que $\frac{F_{K_2}}{p}:\Omega_{K_2}^1\ra F_{K*}\Omega_{K}^1$ est un relèvement de l'opérateur de Cartier inverse $C^{-1}:\Omega_{K}^1\ra\Omega_{K}^1/dK$, on vérifie que la formule $\omega\ra(\omega,0)$ induit un isomorphisme de
$$
\Omega_{K}^1/K^p\dlog a
$$
sur le noyau ci-dessus. Ainsi donc, l'extension (\ref{extEK*amupn}) induit un diagramme commutatif exact de faisceaux sur le petit site étale de $\Spec(K)$
$$
\xymatrix{
0\ar[r] & \bbZ/p \ar[rr] \ar@{=}[d] && R^1\alpha_*\mu_p \ar[r] \ar[d]_{\rotatebox{90}{$\sim$}}^{\dlog_{\mu_p}=\dlog} & R^1\alpha_*E_{K,a}^* \ar[r] \ar[d]_{\rotatebox{90}{$\sim$}}^{\dlog_{E_{K,a}^*}} & 0 \\
0\ar[r] & \bbZ/p \ar[rr]^{1\mapsto\dlog a^{-1}}  && \Ker(C-1)\ar[r]  & \Omega_{K}^1/K^p\dlog a \ar[r] & 0.
}
$$
En prenant les sections globales, on obtient avec \ref{dimcohf} le diagramme commutatif exact de groupes ordinaires
$$
\xymatrix{                               
0\ar[r] & H^0(K,\bbZ/p)  \ar[r] \ar@{=}[d] & H^1(K,\mu_p) \ar[r] \ar[d]_{\rotatebox{90}{$\sim$}}^{\dlog_{\mu_p}=\dlog} & H^1(K,E_{K,a}^*) \ar[r] \ar[d]_{\rotatebox{90}{$\sim$}}^{\dlog_{E_{K,a}^*}} & H^1(K,\bbZ/p) \ar[r] \ar[d]_{\rotatebox{90}{$\sim$}} & 0 \\
0\ar[r] & H^0(K,\bbZ/p)  \ar[r] & \Ker(C-1) \ar[r] & \Omega_{K}^1/K^p\dlog a \ar[r]^<<<<<{\frac{C-1}{\dlog a}} & K/(1-F)(K) \ar[r] & 0,  
}
$$
dont les flèches verticales sont des isomorphismes. En remplaçant $k$ par une extension algébriquement close $\kappa/k$ arbitraire, et donc $K$ par $K_{\kappa}$, le lemme résulte alors de \ref{ratIP0} et \ref{cHsel}.
\end{proof}

\begin{proof}[Démonstration de la proposition \ref{fcontmup}.]
Comme dans la démonstration de \ref{fcontalphap}, il suffit de montrer que l'application $f_{\cont}(\mu_p)(k)$ est bijective, et on commence par l'expliciter. Soit $a\in R^*\setminus (R^*)^p$. D'après \ref{calculEKamupn} et \ref{calculH1EK*amup}, l'image de $a$ par la composée
$$
\xymatrix{
R^* \ar[d] & \\
R^*/(R^*)^p \ar@{=}[d] & \\
\cH^1(\mu_p)^0(k) \ar[d]^{\rotatebox{90}{$\sim$}} & \\
\Ima\bigg(\cH^1(\mu_p)^0(k) \ar@{^{(}->}[r] & H^1(K,\mu_p)=\Ext_{\bbZ/p}^1(\bbZ/p,\mu_p) \ar[r]^>>>>>>>>>{(\cdot)^*} & \Ext_{\bbZ/p}^1(\bbZ/p,\mu_p)\bigg)  \ar[d]^>>>>>{\ov\circ\cH^1} \\ 
& & \Ext_{\fUuf}^1\big(\cH^1(\bbZ/p),\bbZ/p\big), 
}
$$
i.e. $f_{\cont}(\mu_p)(a\modd (R^*)^p)$, est la classe d'isomorphisme d'une extension de $\fUuf$
\begin{equation*}\label{H1extEK*mup}
\xymatrix{                             
0 \ar[r] & \bbZ/p \ar[r] & \cH^1(E_{K,a}^*)/\cH^1(\mu_p)^0  \ar[r] & \cH^1(\bbZ/p) \ar[r] & 0, &=:h_a
}
\end{equation*}
qui admet la description différentielle suivante :
$$
\xymatrix{                               
0\ar[r] & \bbZ/p \ar[r] \ar@{=}[d] & \cH^1(E_{K,a}^*)/\cH^1(\mu_p)^0 \ar[r] \ar[d]_{\rotatebox{90}{$\sim$}}^{\dlog_{E_{K,a}^*}} & \cH^1(\bbZ/p) \ar[r] \ar[d]_{\rotatebox{90}{$\sim$}} & 0 \\
0\ar[r] & \bbZ/p  \ar[r] & \pi_*\Omega_{K}^1/(K^p\dlog a +\dlog(R^*/(R^*)^p)) \ar[r]^<<<<<{\frac{C-1}{\dlog a}} & \pi_*K/\pi_*(1-F)(\pi_*K) \ar[r] & 0,  
}
$$
où l'image de $1\in\bbZ/p$ par la première flèche horizontale de la seconde ligne est la classe de $\dlog t$, d'après le diagramme commutatif exact
$$
\xymatrix{
0\ar[r] & R^*/(R^*)^p \ar[r] \ar[d]_{\rotatebox{90}{$\sim$}}^{\dlog} & K^*/(K^*)^p \ar[r]^v \ar[d]_{\rotatebox{90}{$\sim$}}^{\dlog} & \bbZ/p \ar[r] \ar[d]^{\rotatebox{90}{$\sim$}} & 0 \\
0\ar[r] & \Ker(C-1)\cap\Omega_{R}^1 \ar[r] & \Ker(C-1) \ar[r] & \Ker(C-1)/\Ker(C-1)\cap\Omega_{R}^1 \ar[r] & 0.
}
$$

Par ailleurs, d'après (\ref{variantesrecipK2}), on a un diagramme commutatif exact
$$
\xymatrix{
            &R^*/(R^*)^p \ar[d]_{\rotatebox{90}{$\sim$}}^{\dlog}                     &                                           &&                                  & \\                                                                 
0\ar[r] & \Ker(C-1)\cap\Omega_{R}^1 \ar[r] \ar[d]_{\rotatebox{90}{$\sim$}} & \Omega_{R}^1 \ar[rr]^{C-1} \ar[rr] \ar[d]_{\rotatebox{90}{$\sim$}} && \Omega_{R}^1 \ar[r]  \ar[d]_{\rotatebox{90}{$\sim$}} & 0 \\
0 \ar[r] & \Ext_{\fUuf}^1(\cH^1(\bbZ/p),\bbZ/p)  \ar[r] & \Ext_{\fUuf}^1(\pi_*K/\pi_*R,\bbZ/p) \ar[rr]^{\Ext_{\fUuf}^1(\overline{\pi_*(1-F)},\bbZ/p)} \ar[rr] && \Ext_{\fUuf}^1(\pi_*K/\pi_*R,\bbZ/p)  \ar[r] & 0,
}
$$
où la première flèche verticale est donnée par
\begin{eqnarray*}\label{areconnmup}
R^*/(R^*)^p = \cH^1(\mu_p)^0(k) & \lra &  \Ext_{\fUuf}^1(\cH^1(\bbZ/p),\bbZ/p) \\
a \modd (R^*)^p & \lmapsto &\textrm{classe de } \Res(\cdot\dlog a)^*\epsilon.
\end{eqnarray*}
Nous allons montrer que l'extension $\Res(\cdot\dlog a)^*\epsilon$ est isomorphe à l'extension $h_a$ ; comme la flèche $a \modd (R^*)^p \mapsto \textrm{classe de }\Res(\cdot\dlog a)^*\epsilon$ est bijective par construction, cela achèvera de démontrer que $f_{\cont}(\mu_p)(k)$ est bijective.

L'image de $\dlog a$ par la flèche $\Omega_{R}^1\ra \Ext_{\fUuf}^1(\pi_*K/\pi_*R,\bbZ/p)$ est la première ligne du diagramme commutatif exact
$$
\xymatrix{
0\ar[r] & \bbZ/p \ar[r] \ar@{=}[d] & \pi_*\Omega_{K}^1/\pi_*\Omega_{R}^1{}_{\searrow C-1}\times_{\pi_*\Omega_{K}^1/\pi_*\Omega_{R}^1} {}_{\swarrow \cdot\dlog a}\pi_*K/\pi_*R \ar[r] \ar[d] & \pi_*K/\pi_*R \ar[r] \ar[d]^{\cdot\dlog a} & 0 \\
0\ar[r] & \bbZ/p \ar[r] \ar@{=}[d] & \pi_*\Omega_{K}^1/\pi_*\Omega_{R}^1 \ar[r]^{C-1} \ar[d]^{\Res}  & \pi_*\Omega_{K}^1/\pi_*\Omega_{R}^1 \ar[r] \ar[d]^{\Res} & 0 \\
0\ar[r] & \bbZ/p \ar[r] & \bbG_{a,k}^{\pf} \ar[r]^{1-F} & \bbG_{a,k}^{\pf} \ar[r] & 0 &=:\epsilon.
}
$$
En utilisant que $(C-1)(\Omega_{R}^1)=\Omega_{R}^1$ (puisque $H^1(R,R^1\alpha_*\mu_p)=H^2(R,\mu_p)=0$), et \ref{ratIP0}, on vérifie que le morphisme
\begin{eqnarray*}
\pi_*\Omega_{K}^1/(C-1)^{-1}(R\dlog a)\cap\Omega_{R}^1 & \lra &  \pi_*\Omega_{K}^1/\pi_*\Omega_{R}^1{}_{\searrow C-1}\times_{\pi_*\Omega_{K}^1/\pi_*\Omega_{R}^1} {}_{\swarrow \cdot\dlog a}\pi_*K/\pi_*R \\
\overline{\omega} & \lmapsto & (\overline{\omega},\frac{(C-1)(\omega)}{\dlog a})
\end{eqnarray*}
est un isomorphisme. En utilisant que $(1-F)(R)=R$ (puisque $H^1(R,\bbZ/p)=0$), on vérifie que 
$$
(C-1)^{-1}(R\dlog a)\cap\Omega_{R}^1=R^p\dlog a + \dlog(R^*/(R^*)^p).
$$
On a donc le diagramme commutatif exact
$$
\xymatrix{
0\ar[r] & \bbZ/p \ar[r] \ar@{=}[d] & \pi_*\Omega_{K}^1/\pi_*\Omega_{R}^1{}_{\searrow C-1}\times_{\pi_*\Omega_{K}^1/\pi_*\Omega_{R}^1} {}_{\swarrow \cdot\dlog a}\pi_*K/\pi_*R \ar[r] \ar[d]_{\rotatebox{90}{$\sim$}} & \pi_*K/\pi_*R \ar[r] \ar@{=}[d] & 0 \\
0\ar[r] & \bbZ/p \ar[r] \ar@{=}[d] & \pi_*\Omega_{K}^1/(R^p\dlog a + \dlog(R^*/(R^*)^p)) \ar[r]^<<<<<<<<{\frac{C-1}{\dlog a}} \ar[d] & \pi_*K/\pi_*R \ar[r] \ar[d] & 0 \\
0\ar[r] & \bbZ/p \ar[r] & \pi_*\Omega_{K}^1/(K^p\dlog a + \dlog(R^*/(R^*)^p)) \ar[r]^<<<<<<<<{\frac{C-1}{\dlog a}}  & \pi_*K/\pi_*(1-F)(\pi_*K) \ar[r] & 0.
} 
$$
D'après la description différentielle de $h_a$ vue en début de démonstration, on a donc obtenu
$$
\Res(\cdot\dlog a)^*\epsilon=h_a\quad\in\Ext_{\fUuf}^1(\cH^1(\bbZ/p),\bbZ/p).
$$
\end{proof}

\end{document}